%% file: long-ca-ms-24.tex
\documentclass[11pt,article]{memoir}
\setlrmargins{*}{*}{1}
\checkandfixthelayout

\counterwithout{section}{chapter}
\counterwithout{figure}{chapter}
\firmlists

\usepackage{amssymb}

\usepackage{float} 
\usepackage[backref,hyperindex,colorlinks,citecolor=blue,bookmarksnumbered]{hyperref}
\usepackage{memhfixc} 
\usepackage[numbered]{bookmark} 
\usepackage[all]{hypcap} 
\usepackage[algo2e,algosection,tworuled,noend,noline]{algorithm2e}
\usepackage{cheatpf2}
\usepackage{gacs-algo}
\usepackage[theorem-section]{gacs}


\newcommand{\glo}[1]{\gloproc#1.)))}
\def\gloproc#1@#2.))){\glossary(#1){#2}{}} 
\makeglossary
\makeindex

\numberwithin{equation}{section} 

\newcommand\qq{\qquad}
\newcommand\nn{\nonumber}

\renewcommand\evof[1]{\mathopen\{\verythinmathskip#1\verythinmathskip\mathclose\}}

\newcommand\leor[1] {\stackrel{\scriptstyle #1}{<}}
\newcommand\geor[1] {\stackrel{\scriptstyle #1}{>}}

\newcommand\prl {\hbox{$\mskip\medmuskip\parallel\mskip\medmuskip$}}

\renewcommand{\var}[1]{\ensuremath{\mathord{\textit{#1\/}}}}
 
 \newcommand{\cns}[1]{\ensuremath{\mathinner{\textsf{\upshape{#1}}}}}
 \newcommand{\rul}[1]{\ensuremath{\mathord{\texttt{\slshape{#1\/}}}}}
 \newcommand{\prg}[1]{\ensuremath{\mathord{\texttt{#1}}}}
 \newcommand{\fld}[1]{\ensuremath{\mathord{\textsf{\slshape{#1\/}}}}}

 \newcommand{\cc}{\sqcup} 

\newcommand{\prt}{\boldsymbol{\pi}_{\text{t}}}
\newcommand{\prs}{\boldsymbol{\pi}_{\text{s}}}

\newcommand\Amp {\mathbf{\Psi}}
\newcommand{\ba}{{\mathbf{a}}}
\newcommand\B {{B}} 
\newcommand\e {{e}} 
\newcommand\f {{f}} 
\newcommand\g {{g}} 
\newcommand\GF {\mathrm{GF}} 
\newcommand\grFac {\rho}
\newcommand\h {{h}} 
\newcommand\bb {\mathbf{b}}
\newcommand\bv {\delta} 
\newcommand\bw {w} 
\newcommand\D{D} 
\newcommand\hierCode{\bH}
\newcommand{\numDirAff}{\Dg}
\newcommand\p{{p}} 
\newcommand\nP{P}
\newcommand{\bq}{{\mathbf{q}}} 
\newcommand\Q {{Q}}
\newcommand{\br}{{\mathbf{r}}} 
\newcommand{\red}{\delta}
\newcommand{\rep}{R}
\newcommand{\nS}{S} 
\renewcommand{\bs}{{\mathbf{s}}} 
\newcommand{\sites}{{\Lambda}}
\newcommand{\states}{{\mathbb{S}}}
\newcommand\bu{{\mathbf{u}}} 
\newcommand\T {{T}} 
\newcommand\Tu {T^{\bullet}}   
\newcommand\Tus {T^{\bullet *}}
\newcommand\Tl {{T_{\bullet}}}  
\newcommand\Tls {T_{\bullet}^{*}}
\newcommand\U {{U}} 
\newcommand{\Vdam}{V}
\newcommand\x{{\mathbf{x}}} 

\newcommand{\bsize}{\beta}
\newcommand{\bdist}{\gamma}

\newcommand\Addr {\fld{Addr}}
\newcommand\Adj {\fld{Adj}}
\newcommand\AdjCol {\fld{AdjCol}}
\newcommand\Age {\fld{Age}}
\newcommand\All {\fld{All}}
\newcommand{\Color}{\fld{Color}}
\newcommand\PlCommands {\fld{Pl-commands}}
\newcommand\Cpt {\fld{Cpt}}
\newcommand\Creating{\fld{Creating}}
\newcommand\Cur {\fld{Cur}}
\newcommand{\Decaying}{\fld{Decaying}}
\newcommand\Depth{\fld{Depth}}
\newcommand\Det {\fld{Det}}
\newcommand\Dominant {\fld{Dominant}}
\newcommand{\Dying}{\fld{Dying}}
\newcommand\Doomed{\fld{Doomed}}
\newcommand\F {\fld{F}}
\newcommand{\Frozen}{\fld{Frozen}}
\newcommand\G {\fld{G}}
\newcommand\GermEdge {\fld{G-edge}}
\newcommand\GermSize {\fld{G-size}}
\newcommand\GrowDir {\fld{Grow-dir}}
\newcommand\Growing{\fld{Growing}}
\newcommand\fH {\fld{H}}
\newcommand{\Healing}{\fld{Healing}}
\newcommand\Hold {\fld{Hold}}
\newcommand\Info {\fld{Info}}
\newcommand\Input {\fld{Input}}
\newcommand\Kind {\fld{Kind}}
\newcommand\fL {\fld{L}}

\newcommand\Fromaddr {\fld{Fromaddr}}
\newcommand\slice{\fld{slice}}
\newcommand\Mail {\fld{Mail}}

\newcommand\MainBit{\fld{Main-bit}}
\newcommand\Memory {\fld{Memory}}

\newcommand\PayloadFromNb {\fld{Payload-from-nb}}
\newcommand\PayloadToNb {\fld{Payload-to-nb}}
\newcommand\New {\fld{New}}

\newcommand\Output {\fld{Output}}
\newcommand\Payload {\fld{Payload}}
\newcommand\PlRedun {\fld{Pl-redun}}
\newcommand\PlUpd {\fld{Pl-upd}}

\newcommand\Prev {\fld{Prev}}
\newcommand\Rand {\fld{Rand}}
\newcommand\Receiving{\fld{Receiving}}
\newcommand\Redun{\fld{Redun}}
\newcommand\Retrieved{\fld{Retrieved}}

\newcommand\StageEnd {\fld{Stage-end}}
\newcommand\StageStart {\fld{Stage-start}}

\newcommand\Toaddr {\fld{Toaddr}}
\newcommand\TempPayloadToNb {\fld{Temp-Payload-to-nb}}
\newcommand\ToRefresh{\fld{To-refresh}}

\newcommand\Util{\fld{Util}}
\newcommand\Vote{\fld{Vote}}
\newcommand\Wait {\fld{Wait}}
\newcommand\Work {\fld{Work}}

\newcommand{\lArg}{\fld{\_Arg}}
\newcommand{\lAdjArg}{\fld{\_Adj-arg}}
\newcommand{\lCpt}{\fld{\_Cpt}}
\newcommand{\lColor}{\fld{\_Color}}
\newcommand\lDecoded{\fld{\_Decoded}}
\newcommand\lEncoded{\fld{\_Encoded}}

\newcommand\lInfo{\fld{\_Info}}
\newcommand\lMail{\fld{\_Mail}}
\newcommand\lParam{\fld{\_Param}}
\newcommand\lPayload{\fld{\_Payload}}
\newcommand\lPlRedun{\fld{\_Pl-redun}}
\newcommand\lPlUpd {\fld{\_Pl-upd}}
\newcommand\lToEncode{\fld{\_To-encode}}
\newcommand\lToRefresh{\fld{\_To-refresh}}
\newcommand\lRetrieved{\fld{\_Retrieved}}
\newcommand\lSimOutput{\fld{\_Sim-output}}
\newcommand\lToUpdate{\fld{\_To-update}}
\newcommand\lUpdated{\fld{\_Updated}}
\newcommand\lVote{\fld{\_Vote}}

\newcommand{\vActive}{\var{active}}
\newcommand{\vPassive}{\var{passive}}
\newcommand\ACA {\var{ACA}}
\newcommand\CA {\var{CA}}
\newcommand\PCA {\var{PCA}}
\newcommand\CCA {\var{CCA}}

\newcommand{\adjNb}{\ensuremath{\mathord{\textrm{adj-nb}}}}

\newcommand{\AMed}{\var{AMed}}

\newcommand{\Bad}{\var{Bad}}
\newcommand\Channel {\var{Channel}}

\newcommand{\Configs}{\var{Configs}}
\newcommand\cp{\var{Cap}}
\newcommand{\Dam}{\var{Dam}}
\newcommand{\Damage}{\var{Damage}}
\newcommand\depth{\var{depth}}
\newcommand\Edge{\var{Edge}}
\newcommand\Fit{\var{Fit}}

\newcommand\Ext {\var{Ext}}
\newcommand{\Frame}{\var{Frame}}
\newcommand{\Germ}{\var{Germ}}
\newcommand{\GermGrowing}{\var{G-growing}}
\newcommand{\GermLess}{\var{G-less}}
\newcommand\Growth {\var{Growth}}
\newcommand{\Height}{\var{Height}}
\newcommand{\Histories}{\var{Histories}}
\newcommand\Init {\var{Init}}
\newcommand\Interpr {\var{Interpr}}
\newcommand{\konst}{K}
\newcommand{\legal}{\var{legal}}
\newcommand\Latent {\var{Latent}}
\newcommand{\loc}{\var{loc}}
\newcommand{\MailFree}{\var{Mail-free}}
\newcommand{\MailLegal}{\var{Mail-legal}}
\newcommand{\Med}{\var{Med}}
\newcommand\Member {\var{Member}}
\newcommand{\MyRules}{\var{My-rules}}
\newcommand{\nLocMaint}{\var{n}_{\text{l-m}}}
\newcommand{\NonGermLess}{\var{Non-germ-less}}
\newcommand\Param{\var{Param}}
\newcommand\peer{\var{peer}}
\newcommand{\PlTrans} {\var{Pl-trans}}
\newcommand\Prim{}
\def\Prim-var{\var{Prim-var}}
\newcommand{\Rectangles}{\var{Rectangles}}
\newcommand{\Rob}{\var{Rob}}
\newcommand\Safe {\var{Safe}}
\newcommand{\Space}{\var{space}}
\newcommand\Trajs {\var{Trajs}}
\newcommand\trans {\var{Tr}}
\newcommand{\Undef}{\var{Undef}}
\newcommand\univ {\var{Univ}}
\newcommand\Vac {\var{Vac}}
\newcommand\Vacant {\var{Vacant}}
\newcommand\VacantString {\var{Vacant-str}}
\newcommand{\Xposed}{\var{Xposed}}

\newcommand{\Adapt}{\rul{Adapt}}
\newcommand{\Birth}{\rul{Birth}}
\newcommand{\Broadcast}{\rul{Broadcast}} 
\newcommand{\rCheck}{\rul{Check}} 
\newcommand\Compute{\rul{Compute}}
\newcommand\Copy {\rul{Copy}}
\newcommand\Create{\rul{Create}}
\newcommand{\Decay}{\rul{Decay}}
\newcommand\Decode{\rul{Decode}}
\newcommand{\Die}{\rul{Die}}
\newcommand\Encode{\rul{Encode}}
\newcommand\Eval {\rul{Eval}}
\newcommand\Extend {\rul{Extend}}
\newcommand\Finish {\rul{Finish}}
\newcommand{\Freeze}{\rul{Freeze}}
\newcommand\GermGrow {\rul{Germ-grow}}
\newcommand\Grow {\rul{Grow}}
\newcommand\rActive {\rul{active}} 
\newcommand\passive {\rul{passive}}
\newcommand{\Heal}{\rul{Heal}}
\newcommand{\HealRevive}{\rul{Heal-revive}}
\newcommand{\HealSync}{\rul{Heal-sync}}
\newcommand{\Initialize}{\rul{Initialize}}
\newcommand{\Interpret}{\rul{Interpret}}
\newcommand{\Legalize}{\rul{Legalize}}
\newcommand{\Lift}{\rul{Lift}}
\newcommand{\LocMaintain}{\rul{Loc-maintain}}
\newcommand{\March}{\rul{March}}
\newcommand{\MoveMail}{\rul{Move-mail}}
\newcommand{\PostMail}{\rul{Post-mail}}
\newcommand{\PassDominant}{\rul{Pass-dominant}}
\newcommand{\PassGermSize}{\rul{Pass-germ-size}}
\newcommand{\ProcPayload}{\rul{Proc-payload}}
\newcommand{\Purge}{\rul{Purge}}
\newcommand\ReceiveMail {\rul{Receive-mail}}
\newcommand\RefreshPayload {\rul{Refresh-payload}}
\newcommand\Retrieve {\rul{Retrieve}}
\newcommand{\Send}{\rul{Send}} 

\newcommand\UpdateLocMaint{\rul{Update-loc-maint}}
\newcommand\UpdatePacket{\rul{Update-packet}}
\newcommand\UpdatePayload{\rul{Update-payload}}
\newcommand\Write {\rul{Write}}
\newcommand{\WriteParam}{\rul{Write-param}}

\newcommand{\attrT}{\cns{attrib-t}}
\newcommand{\computeStart}{\cns{compute-start}}
\newcommand{\destrS}{\cns{destr-s}}
\newcommand{\destrT}{\cns{destr-t}}
 \newcommand{\gapLb}{\cns{gap-lb}}
\newcommand{\growStart}{\cns{grow-start}}
\newcommand{\growEnd}{\cns{grow-end}}
\newcommand{\healSpan}{\cns{heal-span}}
\newcommand{\healT}{\cns{heal-t}}
\newcommand{\islSz}{\cns{island-size}}
\newcommand\maxDepth{\cns{max-depth}}
\newcommand{\payloadStart}{\cns{proc-payload-start}}
\newcommand{\purgeT}{\cns{purge-t}}
\newcommand{\splitT}{\cns{split-t}} 
\newcommand{\splitN}{\cns{split-n}}

\begin{document}
 \title
{Reliable Cellular Automata with Self-Organization}

\author{Peter G\'acs
\thanks{Partially supported by NSF grant CR-920484.  The author also
thanks the IBM Alma-den Research Center and the Center for Wiskunde and
Informatica (Amsterdam) for their support during the long gestation of
this project.}
\\Computer Science Department\\ Boston University
     \\ gacs@bu.edu}




\date{\today~(last printing)}

\maketitle

\begin{abstract}
  In a probabilistic cellular automaton in which all local transitions have
positive probability, the problem of keeping a bit of information
indefinitely is nontrivial, even in an infinite automaton.  Still, there is
a solution in 2 dimensions, and this solution can be used to construct a
simple 3-dimensional discrete-time universal fault-tolerant cellular
automaton.  This technique does not help much to solve the following
problems: remembering a bit of information in 1 dimension; computing in
dimensions lower than 3; computing in any dimension with non-synchronized
transitions.

Our more complex technique organizes the cells in blocks that perform a
reliable simulation of a second (generalized) cellular automaton.  The
cells of the latter automaton are also organized in blocks, simulating even
more reliably a third automaton, etc.  Since all this (a possibly infinite
hierarchy) is organized in ``software'', it must be under repair all the
time from damage caused by errors.  A large part of the problem is
essentially self-stabilization recovering from a mess of arbitrary size and
content.  The present paper constructs an asynchronous one-dimensional
fault-tolerant cellular automaton, with the further feature of
``self-organization''.  The latter means that unless a large amount of
input information must be given, the initial configuration can be chosen
homogeneous.

This is a corrected and strengthened version of the journal article~\cite{GacsSorg01}.
 \end{abstract}

 \newpage
 \tableofcontents

 \section{Introduction}

A \emph{cellular automaton} is a homogenous array of identical, locally
communicating finite-state automata.
Traditionally, the model is called \df{interacting particle system} 
when time is continuous.
By choosing the transition function appropriately, 
a cellular automaton can perform an arbitrary computation.
Indeed, for any one-tape Turing machine one can construct a one-dimensional
cellular automaton simulating it step-for-step.

Fault-tolerant information storage and computation in cellular automata
is a natural and challenging mathematical problem but there are also some
arguments indicating an eventual practical significance of the subject,
since there are advantages in uniform structure for parallel computers.

Fault-tolerant cellular automata belong to the larger
category of reliable computing devices built from unreliable
components, in which the error probability of the individual
components is not required to decrease as the size of the device
increases.
In such a model it is essential that the faults are assumed to be
transient: they change the local state but not the local transition
function.

A fault-tolerant computer of this kind must use massive parallelism.
Indeed, information stored anywhere during computation is subject to
decay, and therefore must be actively maintained.
It does not help to run two computers simultaneously, comparing
their results periodically, since if the computers are sufficiently large,
faults will occur in both of them
between comparisons with high probability.
The self-correction mechanism must be built into each part of the
computer.
In cellular automata, it must be a property of the transition function
of the cells.

Due to the homogeneity of cellular automata, since large groups of errors
can destroy large parts of any kind of structure,
``self-stabilization''\footnote{In distributed computing, ``self-stabilization''
  refers to techniques of restoring some structure from an arbitrarily corrupted state;
  however, no new faults are assumed to occur during restoration.}
techniques are needed in conjunction with traditional error-correction.

\subsection {Historical remarks}

The problem of reliable computation with unreliable components was
addressed in~\cite{VonNeum56} in the context of Boolean circuits.
Von Neumann's solution, as well as its improved versions 
in~\cite{DobrOrtyUp77} and \cite{PippRelBool85}, rely on high
connectivity and non-uniform constructs.
The best currently known result of this type is 
in~\cite{SpielmanFTPrl96} where redundancy has been substantially decreased
for the case of computations whose computing time is larger than the
storage requirement.

Of particular interest to us are those probabilistic cellular
automata in which all local transition probabilities are positive (let
us call such automata \df{noisy}), since such an automaton is
obtained by way of ``perturbation'' from a deterministic cellular
automaton.
The automaton may have, for example,
two distinguished initial configurations:
\begin{align*}
 \xi_{0},  \xi_{1}  
\end{align*}
in which all cells have state 0 and in which all have
state 1 (there may be other states besides 0 and 1).
Let \( p_{i}(x,t) \) be the probability that, starting from initial
configuration \( \xi_{i} \), the state of cell \( x \) at time \( t \) is \( i \).
If \( p_{i}(x,t) \) is bigger than, say, \( 2/3 \) for all \( x,t \) then we can say
that the automaton remembers the initial configuration forever.

Informally speaking, a probabilistic cellular automaton is called
\df{ergodic} if it eventually forgets all information about its
initial configuration.
Finite noisy cellular automata are always ergodic.
In the example above, one can define the ``relaxation time'' as the time
by which the probability decreases below \( 2/3 \).
If an infinite automaton is ergodic then the relaxation time of the
corresponding finite automaton is bounded independently of size.
A minimal requirement of fault-tolerance is therefore that the
infinite automaton be non-ergodic.

The difficulty in constructing non-ergodic noisy one-dimensional cellular
automata is that eventually large blocks of errors which we might call
``islands'' will randomly occur.
We can try to design a transition function that (except for a small
error probability) attempts to decrease these islands.
It is a natural idea that the function should replace the state of each
cell, at each transition time, with the majority of the cell states
in some neighborhood.
However, majority voting among the five nearest neighbors (including the
cell itself) seems to lead to an ergodic transition function, even in two
dimensions, if the ``failure'' probabilities are not symmetric with respect
to the interchange of 0's and 1's; it has not been proved to be non-ergodic
even in the symmetric case.
Perturbations of the one-dimensional majority voting function were
actually shown to be ergodic in~\cite{Gray82} and \cite{Gray87}.

Noisy cellular automata remembering a bit forever
(in the sense defined above) in dimensions 2 and higher were
constructed in~\cite{Toom80}.
The paper~\cite{GacsReif3dim88} recognizes that adding one more
dimension, Toom's idea can be used not just for
remembering a bit but for simulating reliably an arbitrary computation.
It designs a simple three-dimensional fault-tolerant cellular automaton that
simulates arbitrary one-dimensional arrays.
The theorem will be spelled out precisely in Section~\ref{sec:1dim.comp}.
Toom's original proof was simplified and adapted to strengthen these
results in~\cite{BermSim88} (see also~\cite{GacsToom95}).  

\begin{remark}
A three-dimensional fault-tolerant cellular automaton cannot be built to
arbitrary size in the physical space.
Indeed, there will be an (inherently irreversible) error-correcting
operation on the average in every constant number of steps in each cell.
This will produce a steady flow of heat from each cell that needs
therefore a separate escape route for each cell.
 \end{remark}

A simple one-dimensional deterministic cellular automaton
eliminating finite islands in the absence of failures was defined
in~\cite{GacsKurdLevIslands78} (see also~\cite{DeSaMaes92}). 
It is now known (see~\cite{ParkThes96}) that perturbation (at least, in a
strongly biased way) makes this automaton ergodic.

\subsection {Hierarchical constructions}\label{sec:hier-constr}

  The limited geometrical possibilities in one dimension suggest that
only some non-local organization can cope with the task of
eliminating finite islands.
  Indeed, imagine a large island of 1's in the 1-dimensional ocean of
0's.
  Without additional information, cells at the left end of this
island will not be able to decide locally whether to move the
boundary to the right or to the left.
  This information must come from some global organization that,
given the fixed size of the cells, is expected to be hierarchical.
  The ``cellular automaton'' in~\cite{Tsir77} gives such a
hierarchical organization.
  It indeed can hold a bit of information indefinitely.
  However, the transition function is not uniform either in space or time:
the hierarchy is ``hardwired'' into the way the transition function
changes.

  The paper~\cite{Gacs1dim86} constructs a non-ergodic one-dimensional
cellular automaton working in discrete time, using some ideas from
the very informal paper~\cite{Kurd78} of Georgii Kurdyumov.
  Surprisingly, it seems even today that in one dimension, the keeping
of a bit of information requires all the organization
needed for general fault-tolerant computation.
  The paper~\cite{Gacs2dim89} constructs a two-dimensional fault-tolerant
cellular automaton.
  In the two-dimensional work, the space requirement of the reliable
implementation of a computation is only a constant times greater than
that of the original version.
  (The time requirement increases by a logarithmic factor.)

  In both papers, the cells are organized in blocks that perform a
fault-tolerant simulation of a second, generalized cellular automaton.
  The cells of the latter automaton are also organized in blocks,
simulating even more reliably a third generalized automaton, etc.
  In all these papers (including the present one), since all this
organization is in ``software'', that is  it is encoded into the states of the
cells, it must be under repair all the time from breakdown caused by
errors.
  In the two-dimensional case, Toom's transition function simplifies the
repairs.

\subsection {New features}

\paragraph{Asynchrony}
In the three-dimensional fault-tolerant cellular automaton 
of~\cite{GacsReif3dim88}, the components must work in discrete time and switch
simultaneously to their next state.
This requirement is unrealistic for arbitrarily large arrays.
A more natural model for asynchronous probabilistic cellular automata
is that of a continuous-time Markov process.
This is a much stronger assumption than allowing an adversary
scheduler, but it still leaves a lot of technical problems to be
solved.
Informally, it allows cells to choose in each moment,
whether to update at the present
time, independently of the choice their neighbors make.

The paper~\cite{BermSim88} gives a simple method to implement
arbitrary computations on asynchronous machines with otherwise
perfectly reliable components.
A two-dimensional asynchronous fault-tolerant cellular automaton
was constructed in~\cite{Wang90}.
Experiments combining this technique with the error-correction
mechanism of~\cite{GacsReif3dim88} were made, among others,
in~\cite{BenGriHeJayMuk90}. 

The present paper constructs a one-dimensional asynchronous
fault-tolerant cellular automaton, thus completing the refutation of
the so-called Positive Rates Conjecture in~\cite{Liggett85}.

\paragraph{Self-organization}
  Most hierarchical constructions, including our earlier ones, start from a complex,
hierarchical initial configuration (in case of an infinite system, and
infinite hierarchy).
The present paper avoids this: its transition function increases the height of the
hierarchy as the computation length requires it.

\paragraph{Proof method simplification}
  Several methods have emerged that help managing the complexity of a large
construction but the following two are the most important.
 \begin{itemize}
  \item 
  A number of ``interface'' concepts is introduced (generalized
simulation, generalized cellular automaton) helping to separate the
levels of the infinite hierarchy, and making it possible to speak
meaningfully of a single pair of adjacent levels.

  \item
  Though the construction is large, its problems are presented one at
a time.
  For example, the messiest part of the self-stabilization is the so-called
Attribution Lemma, showing how after a while all cells can be attributed to
some large organized group (colony), and thus no ``debris'' is in the way
of the creation of new colonies.
  This lemma relies mainly on the Freeze and Decay rules, and will be
proved before introducing many other major rules.
  Other parts of the construction that are not possible to ignore are
used only through ``interface conditions'' (specifications).
 \end{itemize}
  
  We believe that the new result and the new method of presentation
will serve as a firm basis for other new results.
  An example of a problem likely to yield to the new framework is the
\emph{growth rate of the relaxation time} as a function of the size of a
finite cellular automaton.
  At present, the relaxation time of all known cellular automata either
seems to be bounded (ergodic case) or grows exponentially.
  We believe that our constructions will yield examples for other,
intermediate growth rates.

 \subsection{Overview of the paper}
 \begin{bullets}

  \item Section~\ref{sec:CA} defines probabilistic cellular automata.

  \item Section~\ref{sec:results}  and spells out the main theorems for
    discrete time and infinite space.

  \item Section~\ref{sec:codes} introduces block codes and block simulations
    using colonies, and also generalized
    cellular automata (called \df{abstract media}), allowing more general kinds
    of simulation.

 \item Section~\ref{sec:hier} defines hierarchical codes and a hierarchy
    of simulations (called \df{amplifier}).
   It also explains the main technical problems of the construction and the
ways to solve them:
  \begin{enumerate}[label = {\bfseries\textendash}]
   \item correction of structural damage by destruction followed by
rebuilding from the neighbors;
   \item a ``hard-wired'' program;
   \item ``legalization'' of all locally consistent structures;
  \end{enumerate}
 
  \item
    Section~\ref{sec:lat} extends the main discrete-time theorems to the case
    of finite space.

  \item Section~\ref{sec:med} defines \df{media}, a specialization of abstract
media with the needed stochastic structure.
Along with media, we will define \df{canonical simulations}, whose form
guarantees that they are simulations between media.
We will give the basic examples of media with the basic simulations
between them.

  The section also defines variable-period media and formulates the main
  theorems for continuous time.

\item Section~\ref{sec:sync} shows the method we will use to simulate
  a discrete-time cellular automaton by a variable-period one.

  \item 
  Section~\ref{sec:sim} develops some simple simulations, to be used either
directly or as a paradigm.
  The example transition function defined here will correct any set of
faults in which no two faults occur close to each other.

  We also develop the language used for defining our transition function in
the rest of the paper.

  \item 
  A class of media called \df{robust media} for which nontrivial fault-tolerant simulations
exist will be defined in Section~\ref{sec:rob}.
  In these, cells are not necessarily adjacent to each other.
  The transition function can erase as well as create cells.

  The set of space-time points with ``bad'' values is called the
``damage''.
The Restoration Property in Condition~\ref{cond:restor}
requires that at any point of a trajectory, damage
occurs (or persists) only with small probability (\( \eps \)) \glo{eq:greek@\( \eps \)}. 
  The Computation Property requires that the trajectory obey the
transition function in the absence of damage.

As a history \( \eta \) \glo{y.greek@\( \eta \)} of medium \( M_{1} \) simulates
a history \( \eta^{*} \)
of a medium \( M_{2} \)\glo{m.cap@\( M_{k} \)}, we define the damage
of \( \eta^{*} \) \glo{y.greek.star@\( \eta^{*} \)}
in terms of that of \( \eta \) essentially as follows.
Damage occurs at a certain point \( (x,t) \) of \( \eta^{*} \) if
within a certain space-time rectangle in the past of \( (x,t) \), the damage of
\( \eta \) cannot be covered by a small rectangle of a certain size.
  This is saying, essentially, that damage occurs at least ``twice'' in \( \eta \).
  The Restoration Property for \( \eta \) with \( \eps \) will then guarantee
that the damage in \( \eta^{*} \) also satisfies a restoration property 
with \( \approx \eps^{2} \).

\item Section~\ref{sec:amp} defines the kind of amplifiers to be built.
  The \emph{main lemma}, called the Amplifier Lemma, says that these
  amplifiers exist.
  The rest of the section applies the main lemma to the proof of one of the
main theorems.

  \item 
  Section~\ref{sec:plan} gives an overview of an amplifier.
  As indicated above, the restoration property will be satisfied
automatically.
  In order to satisfy the computation property, the general framework of the
program will be similar to the outline in Section~\ref{sec:sim}.
  However, besides the single-error fault-tolerance property achieved
there, it will also have a
\df{self-stabilization}\index{self!stabilization} property.
  This means that a short time after the occurrence of arbitrary damage,
the configuration enables us to interpret it in terms of colonies.
  (In practice, pieces of incomplete colonies will eliminate themselves.)
  In the absence of damage, therefore, the colony structure will recover
from the effects of earlier damage, that is predictability in the simulated
configuration is restored.

\item Section~\ref{sec:consistency} defines the kind of local consistency
  needed for a colony to function.

  \item 
  Section~\ref{sec:kill} gives the rules for killing, creation, 
  growth, including the growth of germs which are precursors of colonies.
  It proves the basic lemmas about space-time paths connecting live cells.
  It also defines the healing rule; due to the need to restore some
  local clock values consistently with the neighbors, this rule
  is somewhat elaborate.

  \item
  Section~\ref{sec:gaps} defines the decay rule and shows that a large
gap will eat up a whole colony.

  \item
  Section~\ref{sec:attrib} proves the Attribution Lemma that traces back
each non-germ cell to a full colony.
  This lemma expresses the ``self-stabilization'' property mentioned above.
  The proof relies on Section~\ref{sec:gap.bad} showing that if a gap
will not be healed promptly then it grows.

  \item 
  Section~\ref{sec:heal} proves the Healing Lemma, showing how the effect of
a small amount of damage will be corrected.

  \item
  Section~\ref{sec:commun} introduces and applies the communication
rules needed to prove the Computation Property in simulation.
  These are rather elaborate, due to the need to communicate with not
completely reliable neighbor colonies asynchronously.

  \item
    Section~\ref{sec:comp} introduces and uses the error-correcting
    computation rules.

  \item
    Section~\ref{sec:sim-robust} proves the main properties of the reliable
    simulation defined in the preceding sections.

 \item 
   Section~\ref{sec:sorg} shows how the germ-growth rules lead to self-organization,
   and proves the main theorems that use self-organization.

  \item 
  The concluding remarks in Section~\ref{sec:concl} hint at some applications and questions.

 \end{bullets}

  The above constructions will be carried out for the case when the cells
work asynchronously (with variable time between switchings).
  This does not introduce any insurmountable difficulty but makes life
harder at several steps: more care is needed in the updating and correction
of the counter field of a cell, and in the communication between neighbor
colonies.
  The analysis in the proof also becomes more involved.

\section{Cellular automata}\label{sec:CA}

In the introductory sections, we confine ourselves to one-dimensional
infinite cellular automata.
Let us introduce some notation to be used throughout.

\begin{notation}
Let \( \bbR \) be the set of real numbers, and
 \( \bbZ_{m} \)\glo{z.cap.z@\( \bbZ_{m} \)}  
the set of remainders modulo \( m \).
For \( m=\infty \), this is the set \( \bbZ \) of integers.

We will use the notation
\begin{align}\label{eq:lem}
   f(n)\lem g(n)
\end{align}
for what usually is written as \( f(n)=O(g(n)) \): that is for the fact that \( f(n)\le c g(n) \)
for some constant \( c \).
We write \( f(n)\eqm g(n) \) if \( f(n)\lem g(n) \) and \( g(n)\le f(n) \).

The standard mathematical notation for open
intervals and for pairs is the same; we hope that the context will make the
meaning always clear.
We will use the same notation for intervals of integers as for those of real
numbers: the context 
will make it clear, whether \( \clint{a}{b} \) or \( \clint{a}{b}\cap \bbZ \) is
understood.
Given a set \( A \) of space or space-time and a real number \( c \), we write
 \[
  c A = \setof{cv : v\in A},
 \]
and \( \ol A \) for its closure: for example, \( \ol{\rint{a}{b}}=\clint{a}{b} \).
Given two space-time sets \( A,B \), we denote
 \begin{align}\label{eq:set-sum}
 A+B = \setOf{a+b}{a\in A,b\in B}.
 \end{align}
 
Some lists of assertions are denoted by
(a), (b), \( \dots \) and some by (1), (2), \( \dots \).
The attempted convention is that for a 
list properties that all hold or are required (conjunction) the
items are labeled with (a),(b), \( \dots \) while if the list is a list of
several possible cases (disjunction) then the items are labeled with (1),
(2), \( \dots \).

Maxima and minima will sometimes be denoted by \( \vee \) and
\( \wedge \)\glo{"+@\( \vee,\wedge \)}.
  We will write \( \log \)\glo{log@\( \log \)} for \( \log_{2} \).
  For two strings \( u,v \), we will denote by%
\glo{"+@\( \cc \)}
 \begin{equation}\label{eq:concat}
   u\cc v
 \end{equation}
  their concatenation.
\end{notation}

\subsection{Deterministic cellular automata}\label{sec:det-ca}

Let us give here the most frequently used definition of cellular
automata.
Later, we will use a certain generalization.
First, the notions associated without considering time.

 \begin{definition}[Space and states]\label{def:space}
  The set \( \sites \)\glo{c.cap@\( \sites \)} of
 \df{sites}\index{site} has the form \( \bbZ_{m} \) for finite or infinite \( m \).
This will mean that in the finite case, we take \df{periodic boundary
conditions}.
In a space-time vector \( (x,t) \), we will always write the space coordinate
first.
For a space-time set \( E \), we will denote its space- and time projections
by \glo{p.greek@\( \prs, \prt \)}
 \begin{equation}\label{eq:proj}
 \prs E, \prt E
 \end{equation}
respectively.
We will have a finite set \( \states \)\glo{s.cap.s@\( \states \)} of states, the potential
states of each site.
A \df{configuration}\index{configuration} is a function
 \glo{x.greek@\( \xi(x) \)}
 \[ 
   \xi(x)
 \]
  for \( x\in\sites \).
Here, \( \xi(x) \) is the state of site \( x \).
 \end{definition}

Now, the notations needed for reasoning about the evolution of configurations.
For uniform treatment, it is useful to view even discrete-time cellular automata
as working in continuous time, but of course they will change their states only
at certain predetermined instants.

 \begin{definition}[Space-time]\label{def:space-time}
The time of work of our cellular automata will be the interval
\( \lint{0}{\infty} \).
  Our \df{space-time} is given by 
 \[
 \sites\times\lint{0}{\infty}.
 \]
A \df{history}\index{history} is a
space-time function \( \eta(x,t) \)\glo{y.greek.func@\( \eta(x,t) \)} which for each \( t \)
defines a configuration.
If in a history \( \eta \) we have \( \eta(x,v)=s_{2} \) and
\( \eta(x,t)=s_{1}\ne s_{2} \) for all \( t<v \) sufficiently close to \( v \) then we
can say that there was a
\df{switch}\index{switch} from state \( s_{1} \) to state \( s_{2} \) at time
\( v \).
For ordinary discrete-time cellular automata, we allow only histories
in which all switching times are natural numbers
\( 0,1,2,\dots \).
The time 0 is considered a switching time.
If there is an \( \eps \) such that \( \eta(c,t) \) is constant for
\( a-\eps<t<a \) then this constant value will be denoted by
 \glo{-@\( \eta(c, a-) \)}
 \begin{equation} \label{eq:a.minus} 
  \eta(c, a-).
 \end{equation}
The \df{subconfiguration} \( \xi(D') \) of a configuration \( \xi \)
defined on \( D\spsq D' \) is the restriction of \( \xi \) to \( D' \).
Sometimes, we write
 \[
  \eta(V)
 \]
 for the sub-configuration over the space-time set \( V \).
 \end{definition}

Cellular automata describe a kind of dynamic for evolutions.

  \begin{definition}[Deterministic cellular automata]
  A
 \df{deterministic cellular automaton}\index{cellular automaton!deterministic}%
\glo{ca@\( \CA(\trans, \sites) \)}%
 \[
   \CA(\trans, \sites).
 \]
  is determined by a   transition function\index{transition function}
\( \trans:\states^{3}\to\states \)\glo{tr@\( \trans, \ol\trans \)} and the set \( \sites \) of
sites.
We will omit \( \sites \) from the notation when it is obvious from the
context.
A history \( \eta \) is a
\df{trajectory}\index{trajectory!deterministic cellular automaton} of this
automaton if 
 \[
   \eta(x,t)=\trans(\eta(x-1,t-1), \eta(x,t-1), \eta(x+1,t-1))
 \]
holds for all \( x,t \) with \( t>0 \).
For a history \( \eta \) let us write
 \begin{equation}\label{eq:trans}
   \ol\trans(\eta, x,t)=\trans(\eta(x-1,t-1), \eta(x,t-1), \eta(x+1,t-1))
 \end{equation}
for the ``intended'' value of \( \eta(x,t) \).
  \end{definition}

Given a configuration \( \xi \) over the space \( \sites \) and a transition
function, there is a unique trajectory \( \eta \) with the given transition
function and the initial configuration \( \eta(\cdot,0)=\xi \).

\subsection{Fields of a local state}\label{sec:fields}

The history of a deterministic cellular automaton can be
viewed as a ``computation''.
Moreover, every imaginable computation can be performed by an
appropriately chosen cellular automaton function.
This is not the place to explain the meaning of this statement if
it is not clear to the reader. 
But it becomes maybe clearer if we point out that a better known model of
computation, the Turing machine\index{Turing machine}, can be considered a
special cellular automaton.

\begin{definition}[Capacity]
We will deal, from now on, only with cellular automata in which the set
\( \states \) of local states consists of binary strings of some fixed length
\( \cp=\nm{\states} \)\glo{"|@\( \nm{\states} \)} called the
\df{capacity}\index{capacity} of the sites.
Thus, if the automaton has 16 possible states then its states can be
considered binary strings of length 4.
\end{definition}

  If \( \nm{\states}>1 \) then the information represented by the state can be
broken up naturally into parts.
  It will greatly help reasoning about a transition rule if it assigns
different functions to some of these parts; a typical ``computation'' would
indeed do so.

 \begin{definition}[Fields]\label{def:fields}
Nonempty subsets of the set \( \set{0,\dots,\nm{\states}-1} \) will be called
 \df{fields}\index{field}.
Some of these subsets will have special names.
Let \glo{all@\( \All \)}
 \[
   \All = \set{0,\dots,\nm{\states}-1}.
 \]
If \( s=\tupof{s(i) : i\in\All} \) is a bit string and \( \F=\set{i_{1},\dots,i_{k}} \)
is a field with \( i_{j}<i_{j+1} \) then we will write \glo{.@\( s.\F \)}
 \[
   s.\F=\tup{s(i_{1}),\dots,s(i_{k})}
 \]
 for the bit string that is called \df{field} \( \F \) of the state.
 The field consisting of the bits \( s(i),s(i+1),\dots,s(j-1) \) of the state will be
 denoted
 \begin{align}\label{eq:F[i:j]}
   \F[i\frTo j].
 \end{align}
The array obtained from joining 
the same fields (say the Mail field) of all the different cells is called 
a \df{track}\index{track}, (say the Mail track).
(The terminology recalls tracks of a magnetic tape.)
We will generally define fields that 
are either disjoint or contained in each other (but there may be some
exceptions.)
We will call the number of bits \( |\F| \) in a field its \( \df{width} \), and
and the proportion \( |F|/\cp \) its \df{relative width}.
The width and relative width of a \df{track} is defined to be the same as that of the
corresponding field.
 \end{definition}

  \begin{example}
  If the capacity is 12 we could subdivide the interval \( \clint{0}{11} \) into
subintervals of lengths 2,2,1,1,2,4 respectively and call these fields
the input, output, mail coming from left, mail coming from right,
memory and workspace.
  We can denote these as
 \( \Input \)\glo{input@\( \Input \)},
 \( \Output \)\glo{output@\( \Output \)},
 \( \Mail_{j} \)\glo{mail@\( \Mail \)} (\( j=-1,1 \)),
 \( \fld{Work} \)\glo{work@\( \fld{Work} \)} and
 \( \Memory \)\glo{memory@\( \Memory \)}.
  If \( s \) is a state then \( s.\Input \) denotes the first two bits of \( s \),
\( s.\Mail_{1} \) means the sixth bit of \( s \), etc.
Treating these fields differently means we may impose some useful
restrictions on the transition function.
We might require the following, calling \( \Mail_{-1} \) the ``right-directed
mail field'':
 \begin{quote}
  The information in \( \Mail_{-1} \) moves always to the right.
  More precisely, in a trajectory \( \eta \), the only part of the state
\( \eta(x,t) \) that depends on the state \( \eta(x-1,t-1) \) of the left
neighbor is the right-directed mail field \( \eta(x,t).\Mail_{-1} \).
  This field, on the other hand, depends only on the right-directed
mail field of the left neighbor and the workspace field
\( \eta(x,t-1).\fld{Work} \).
  The memory depends only on the workspace.
 \end{quote}
  Confining ourselves to computations that are structured in a similar
way make reasoning about them in the presence of faults much easier.
  Indeed, in such a scheme, the effects of a fault can propagate only
through the mail fields and can affect the memory field only if the
workspace field's state allows it.
 \end{example}

The following concepts related to the transition function will play a role.
 \glo{legal@\( \legal \)}%
\begin{definition}[Legality]\label{def:legal}
  Given a transition function \( \trans \), we say that state \( s' \) is a
  \df{legal successor} of state \( s \) if there are states \( r,t \) with
  \( s'=\trans(r,s,t) \).
  We define
  \begin{align*}
 \legal(s,s') = 1
\end{align*}
if this holds and 0 if it does not.  
\end{definition}

  \subsection{Probabilistic cellular automata}\label{sec:PCA}

In probabilistic cellular automata, the transitions allow randomness.

 \begin{definition}
  A \df{random history}\index{history!random} is a pair
\( (\mu,\eta) \)\glo{m.greek@\( \mu \)} where \( \mu \) is a probability measure over some
measurable space\index{measurable space}
 \( (\Og,\cA) \)\glo{o.greek@\( \Og,\og \)}\glo{a.cap@\( \cA,\cA(W),\cA_{t} \)}
 together with a measurable function \( \eta(x,t,\og) \)
 which is a history for every value of \( \og\in\Og \).
We will generally omit \( \og \) from the arguments of \( \eta \).
When we omit the mention of \( \mu \) we will use \( \Prob \)\glo{prob@\( \Prob \)}
to denote it.
If it does not lead to confusion, for some property of the form
\( \setof{\eta\in R} \), the quantity
 \( \mu\setOf{\og}{\eta(\cdot,\cdot,\og)\in R} \) will be written as usual, as
 \[
   \mu\evof{\eta\in R} .
 \]
 \end{definition}

\begin{notation}
  We will denote the expected value of \( f \) with respect to \( \mu \) by
 \glo{eq:cap @\( \Expv_{\mu} f \)}
 \[
   \Expv_{\mu} f
 \]
  where we will omit \( \mu \) when it is clear from the context. 
\end{notation}

 \begin{definition}[Events]\label{def:events}
 A function \( f(\eta) \) with values 0, 1 (that is an indicator function) and
measurable in the \( \sigma \)-algebra \( \cA \) will be called an
 \df{event function}\index{event function}  over \( \cA \).
Let \( W \) be any subset of space-time that is the union of some rectangles.
  Then 
 \[
  \cA(W)
 \]
  denotes the \( \sigma \)-algebra generated by events of the form
 \[
  \evof{\eta(x,t)=s \txt{ for } t_{1}\le t < t_{2}}
 \]
 for \( s\in\states \), \( \pair{x}{t_{i}}\in W \).
 \end{definition}

Sometimes, we may want to refer to events not necessarily expressible by
the history \( \eta \); still, we need a sense in which they become knowable
over time.
The notion of filtration captures this.

\begin{definition}[Filtration]\label{def:filtration}
For all times \( t \), we assume the existence of a \( \sigma \)-algebra
 \[
  \cA_{t} \supseteq \cA(\sites\times\clint{0}{t}),\quad \cA_{<t}=\bigcup_{u<t} \cA_{u}.
 \]
The system \( \{\cA_{t}\} \) is required to be increasing:
 \[
  t < u \imp \cA_{t} \subseteq \cA_{u}.
 \]
The function \( t\mapsto\cA_{t} \) is our \df{filtration}.
\end{definition}

\begin{example}
  The system \( t\mapsto\cA(\sites\times\clint{0}{t}) \) is a filtration.
\end{example}

Transition matrices will play the role of transition functions.

 \begin{definition}[Probabilistic cellular automaton]
A \df{transition matrix}\index{transition matrix}
\( \bP(s, (r_{-1},r_{0},r_{1})) \)\glo{p.cap@\( \bP(s, (r_{-1},r_{0},r_{1})) \)}
is a function \( \bP:\states^{4}\to\clint{0}{1} \) with the property
 \( \sum_{s}\bP(s,\br)=1 \).
For an arbitrary history \( \eta \), and space-time point 
\( (x,t) \), denote
 \begin{equation}\label{eq:Trans-prob}
 \ol\bP(\eta, s, x, t)
    = \bP(s, (\eta(x-1,t),\eta(x,t),\eta(x+1,t))).
 \end{equation}
  We will omit the parameter \( \eta \) when it is clear from the context.

A \df{probabilistic cellular automaton}\index{cellular
automaton!probabilistic}%
 \glo{PCA@\( \PCA(\bP, \sites) \)}%
 \[
  \PCA(\bP, \sites)
 \]
  is characterized by saying which random histories are
considered \df{trajectories}\index{trajectory}.
Now a trajectory is not a single history (sample path)
but a distribution over histories that satisfies the
following condition, saying
that the random history \( \eta \) is a
trajectory if and only if the following holds.
Let \( s_{1},\dots, s_{n}\in\states \).
 \begin{equation*}
 \Prob\Bigparen{\bigcap_{i=1}^{n}\{\eta(x_{i},t)=s_{i}\} \Bigm\vert \cA_{t-1}}
 = \prod_{i=1}^{n} \ol\bP(\eta, s_{i}, x_{i}, t-1).
 \end{equation*}
This expression uses the notion of conditional probability over sigma-algebras.
For a more elementary statement of the same property, 
let \( x_{0},\dots,x_{n+1} \) be given with \( x_{i+1}=x_{i}+1 \).
Let us fix an arbitrary history \( \zeta \) and an arbitrary
event \( \cH\ni\zeta \) in \( \cA_{t-1} \).
  Then we require
 \begin{multline*}
 \Prob\Bigparen{\bigcap_{i=1}^{n}\{\eta(x_{i},t)=\zeta(x_{i},t)\} \cap
  \cH\cap\bigcap_{i=0}^{n+1}\{\eta(x_{i},t-1)=\zeta(x_{i},t-1)\}}
 \\ =
 \Prob\Bigparen{\cH\cap\bigcap_{i=0}^{n+1}\{\eta(x_{i},t-1)=\zeta(x_{i},t-1)\}}
 \prod_{i=1}^{n} \ol\bP(\eta, \zeta(x_{i},t), x_{i}, t-1).
 \end{multline*}
A probabilistic cellular automaton is \df{noisy}\index{noisy} if
\( \bP(s,\br)>0 \) for all \( s,\br \).
Bandwidth can be defined for transition probabilities just as for
transition functions.
\end{definition}

 \begin{example}\label{xmp:rand-gen}
  As a simple example, consider a deterministic cellular automaton
with a ``random number generator''.
  Let the local state be a record with two fields, \( \Det \)\glo{det@\( \Det \)}
and \( \Rand \) where \( \Rand \)\glo{rand@\( \Rand \)} consists of a single bit.
  In a trajectory \( (\mu,\eta) \), the field \( \eta.\Det(x,t+1) \) is
computed by a deterministic transition function from \( \eta(x-1,t) \),
\( \eta(x,t) \), \( \eta(x+1,t) \), while \( \eta.\Rand(x,t+1) \) is obtained by
coin-tossing.
 \end{example}

A trajectory of a probabilistic cellular automaton is a discrete-time
Markov process.
If the set of sites consists of a single site then \( \bP(s, \br) \) is the
transition probability matrix of this \df{finite Markov
chain}\index{Markov!chain}.
We get a finite Markov chain as long as the number of sites is finite.

\subsection{Continuous-time probabilistic cellular automata}
 \label{sec:Markov}

For later reference, let us define here (1-dimensional)
probabilistic cellular automata in which the sites make a random
decision ``in each moment'' on whether to make a transition to another
state or not.
A systematic theory of such systems and an overview of many results
available in 1985 can be found in~\cite{Liggett85}.
Here, we indicate an elementary construction
similar to the one in~\cite{Gray82}.

 \begin{definition}[Transition rate matrix]
Let us call a \df{transition rate matrix} a function
\( \bR:\states^{4}\to\lint{0}{\infty} \),\glo{rem:cap@\( \bR(s, \br) \)}%
written as 
 \begin{align*}
 \bR(s, \br)\ge 0,
 \end{align*}
with the normalization property \( \bR(r_{0}, \tup{r_{-1},r_{0},r_{1}})=0 \) 
for all \( \br \).
Its elements are called the \df{transition rates}\index{transition rate}.
 \end{definition}

We will obtain the continous-time process as the limit of
certain discrete processes:

\begin{definition}\label{def:discr-approx}
  Consider a generalization 
 \begin{align*}
  \PCA(\bP, \B , \delta, \sites)
 \end{align*} 
of probabilistic
cellular automata in which the sites are at positions \( i\B  \) for some fixed
\( \B  \) called the \df{body size} and integers \( i \), and the switching times are
at \( 0,\delta,2\delta,3\delta,\dots \) for some small positive \( \delta \)\glo{def:greek@\( \delta \)}.
\end{definition}

The body size parameter \( \B \) will make more sense later, in the context of simulations.
The parameter \( \delta \) is of more interest now, as the limit \( \delta\to 0 \) will be considered.
Fixing \( \sites,\bP \), let
 \begin{align*}
 M_{\delta}=\PCA(\bP,1,\delta,\sites)
 \end{align*}
\glo{m.cap.delta@\( M_{\delta} \)}
with \( \bP(s,\br)=\delta\bR(s,\br) \) when
\( s\ne r_{0} \) and \( 1-\delta\sum_{s'\ne r_{0}}\bR(s',\br) \) otherwise.
(The definition is sound when \( \delta \) is small enough to make the last
expression nonnegative.)
It can be shown
that with any fixed initial configuration \( \eta(\cdot,0) \), the distributions
of trajectories \( \eta_{\delta} \) of \( M_{\delta} \) will converge to a certain random
process \( \eta \) which is the continuous-time probabilistic cellular
automaton with these rates.
For the sense of convergence and other constructions of the same process,
see~\cite{Liggett85}, \cite{Gray82} and the works quoted there.

 \begin{definition}[Interacting particle system]
The random process defined by the above procedure 
 will be denoted by \glo{cca@\( \CCA(\bR, \sites) \)}
 \[
   \CCA(\bR, \sites),
 \]
and called the \df{(continuous-time) interacting particle
system}\index{interacting particle system}
with the given rate matrix.
We call this system \df{noisy} if \( \bR(s,\br)>0 \) for all \( s\ne r_{0} \).
 \end{definition}

It can be shown 
that the process defined this way is a Markov process\index{Markov!process},
that is if we fix the past before some time \( t_{0} \) then the conditional
distribution of the process after \( t_{0} \) will only depend on the
configuration at time \( t_{0} \).
(For a more general definition allowing simultaneous change in a
finite number of sites, see~\cite{Liggett85}.)

\subsection{Perturbation}\label{sec:perturb}

  Intuitively, a deterministic cellular automaton is fault-tolerant if
even after it is ``perturbed'' into a probabilistic cellular
automaton, its trajectories can keep the most important properties of
a trajectory of the original deterministic cellular automaton.
Formally, we introduce perturbed versions of our cellular automata.

 \paragraph{Discrete time}
The definition is somewhat more straightforward in the discrete-time case, since
then only the transition function (or probability matrix) will be perturbed.

 \begin{definition}[Random perturbation of a deterministic system]
  We will say that a random history \( (\mu,\eta) \) is a
\df{trajectory} of the \( \eps \)-\df{perturbation}\index{perturbation}
 \glo{cae@\( \CA_{\eps}(\trans, \sites) \)}
 \[
  \CA_{\eps}(\trans, \sites)
 \] 
  of the transition function \( \trans \) if the following holds.
For all \( x_{0},\dots x_{n+1},t \) with \( x_{i+1}=x_{i}+1 \) 
for all \( 0<i_{1}<\cdots<i_{k}<n \),
 \begin{equation*}
 \mu\Bigparen{\bigcap_{j=1}^{k}\{\eta(x_{i_{j}},t)\ne\ol\trans(\eta, x_{i_{j}},t-1)\}
 \Bigm\vert \cA_{t-1}} \le\eps^{k}.
 \end{equation*}
Here is the same property without using conditional probability over sigma-algebras:
for events \( \cH \) in \( \cA_{<(t-1)} \) with \( \mu(\cH)>0 \), 
 \[
 \mu\Bigparen{\bigcap_{j=1}^{k}
    \{\eta(x_{i_{j}},t)\ne\ol\trans(\eta, x_{i_{j}},t-1)\} \Bigm|
  \cH\cap\bigcap_{i=0}^{n+1}\{\eta(x_{i},t-1)=s_{i}\}}\le \eps^{k}.
 \]
 \end{definition}

Note that the model \( \CA_{\eps}(\trans, \sites) \) defined this way
is not a probabilistic cellular automaton, since even if the initial
configuration is fixed there are many random processes that are
accepted as its trajectories.
Given any probabilistic cellular automaton \( \PCA(\bP, \sites) \) such
that \( \bP(s,\br, \sites)\ge 1-\eps \) whenever \( s=\trans(\br) \), the
trajectories of this are accepted as trajectories of \( \CA_{\eps}(\trans, \sites) \);
however, these do not exhaust all the possibilities.
We may think of the trajectory of a perturbation as a process created by
an ``adversary'' who is trying to defeat whatever conclusions we want to
make about the trajectory, and is only restricted by the inequalities that
the distribution of the trajectory must satisfy.

This new freedom becomes important in our later construction.

\paragraph{Continuous time}
There are several choices of how to generalize perturbation to  continous time.

 \begin{definition}[Random perturbation of a continous-time system]
  By the \df{\( \eps \)-perturbation} of a continuous-time interacting particle
system with transition rates given by \( \bR(s,\br) \), we understand the
following: in the above construction of a process, perturb the matrix
elements \( \bR(s,\br) \) by some arbitrary amounts smaller than \( \eps \).
 \end{definition}

Note that this is a more modest kind of perturbation since the perturbed
process is again a continuous-time interacting particle system, just with changed
parameters.
Other perturbations are imaginable, but not worth the trouble 
defining formally now.

\section{Some results}\label{sec:results}

We state some of the main results in this section; but a more detailed and formalized
description of them will be given in Sections~\ref{sec:lat} and~\ref{sec:med}.
For simplicity the theorems here are only for infinite space.
A noisy cellular automaton on a finite space is ergodic, so eventually it forgets
all about its initial state.
However, as will be shown in the later sections, even in finite spaces
the ``relaxation time'', that is the time
before which information as well as computation results of the reliable cellular automaton
can be trusted, grows (almost) exponentially as a function of the space size.

\subsection{Information storage}

 \paragraph{Remembering a few bits}
The simplest task we may want to assign a cellular automaton is to store some
information.

 \begin{definition}[Remembering some bits]
  Suppose that the bit string that is a local state has some field \( \F \)
(it can for example be the first two bits of the state).
We will say that \( \trans \) \df{remembers}\index{remembering} field \( \F \)
over an infinite set of sites \( \sites \) if
there is an \( \eps>0 \) such that for each string \( s\in\{0,1\}^{|\F|} \) there
is a configuration \( \xi_{s} \) such that for all
trajectories \( (\mu,\eta) \) of the \( \eps \)-perturbation
 \( \CA_{\eps}(\trans, \sites) \) with \( \eta(\cdot,0)=\xi_{s} \), for all \( x,t \) we
 have
 \begin{equation}\label{eq:remember-field}
  \mu\evof{\eta(x,t).\F\ne s}< f(t)+O(\eps)\text{ where } \lim_{t\to\infty}f(t)=0.
\end{equation}
  We define similarly the notions of remembering a field for a
probabilistic transition matrix \( \bP \) and a probabilistic
transition rate matrix \( \bR \).

Let us call a configuration \( \xi \) \df{homogeneous} if there is a state \( q\in\states \) such
that \( \xi(x)=q \) for all \( x \).
We say that \( \trans \) remembers \( \F \) in a \df{self-organizing way}
if the initial configurations \( \xi_{s} \) can be made homogeneous.
\end{definition}
For some \( c>1 \) let
\begin{equation}\label{eq:h_0}
\begin{aligned}
   h_{0}(N,c)&=c^{(\log N)^{1/2}}.
\end{aligned}
\end{equation}
\begin{theorem}[Remembering a field in discrete time]\label{thm:1dim.nonerg}
  For any constant \( c_{1}>1 \) 
 there is a one-dimensional probabilistic transition matrix that remembers a field,
 in a self-organizing way;  the function \( f(t) \) in~\eqref{eq:remember-field}
can be chosen as \( t^{-c_{1}}\).
If the space is finite with size \( N \) then in any cell, at time \( t \),
the probability of forgetting it at time \( t \) is bounded by \( \eps^{h_{0}(N,c_{2})} \)
for an appropriate constant \( c_{2} \).
 \end{theorem}
 Thus for a finite space, the memory lasts, even if not for exponential time, but
 for a time exponential in \( h_{0}(N,c_{2}) \).
 \begin{remark}\label{rem:sorg-poor}
   The error term in~\eqref{eq:remember-field} improves with time, but 
   why does it not decrease with \( \eps \)?
   Technically, in our proofs it is the price of the need for symmetry-breaking random
   choices in self-organization, but it still seems improvable.
 \end{remark}

In the proof of this theorem 
in~\cite{Gacs1dim86} the initial configuration \( \xi_{s} \) has an
infinite hierarchical structure, so the remembering was not self-organizing.
In Theorem~\ref{thm:1dim.nonerg} this is not necessary anymore.
We achieve the simplification of the result
by a technique we call \emph{self-organization}: the
hierarchy will be built up by the probabilistic transition during computation.    
The following continuous-time version of the same result is also new:

 \begin{theorem}[Remembering a field in continuous time]
\label{thm:1dim.nonerg.cont}
  There is a one-dimensional transition-rate matrix that remembers a field in a self-organizing way.
\end{theorem}

\paragraph{Remembering a long (possibly infinite) string}
In order to store information reliably in a computing device that deals with its
bits individually, it is necessary to add redundancy; otherwise
some bits can be lost in the very first step.
So in order to store a string the initial configuration of the cellular automaton
will contain it encoded by some error-correcting code \( (\fg_{*},\fg^{*}) \)
(the notation will be motivated in Section~\ref{sec:codes}).
A finite or infinite string \( \rho \)
in some alphabet \( \Sigma \) will be encoded into a string \( \fg_{*}(\rho) \)
in another alphabet \( \states \) of symbols, cell states of some
cellular automaton.
These states have fields, and one of them can be called \( \Info \).
The code we will use is such that if \( \xi=\fg_{*}(\rho) \) then for position \( n \)
in \( \rho \) we have \( \xi(x).\Info=\rho(x) \): the original word
is in the code explicitly in the \( \Info \) field, while the other fields serve for
error checks and other structure.
Though it will take some time to describe our codes (in the rest of the paper)
as they have a hierarchical
structure, in fact they are easy to compute: their computational complexity is low.
Here is a theorem about remembering forever an infinite string:

\begin{theorem}[Remembering a long string]\label{thm:remember-inf}
  Given some alphabet \( \Sigma \) there is a one-dimensional deterministic
  cellular automaton \( \CA(\trans) \) with a field \( \F \),
  a code \( \psi_{*}:\Sigma^{\bbZ_{N}}\to\states^{\bbZ_{N}} \) for finite
  or infinite \( N \) and a constant \( c>1 \)
  such that for sufficiently small \( \eps \),
  if \( \eta \) is the trajectory of an \( \eps \)-perturbation \( \CA_{\eps}(\trans) \)
  with \( \eta(\cdot,0)=\psi_{*}(\rho) \) then for all \( x,t \) we have
  \begin{align*}
    \Pbof{\eta(x,t).\F\ne\rho(x)} = O(\eps)+\eps^{h_{0}(N,c)}t.
  \end{align*}
\end{theorem}
As we see this code turns our cellular automaton an information-transmission device
from the present to the future, with \emph{constant capacity}: storing one symbol per cell.
(Some non-\( \Info \) fields of the cells will be used for error-checks.)

\subsection{Computation}

In order to accommodate computations that last forever,
and to formalize the idea of ``more and more output'',
let us introduce a special alphabet:

 \begin{definition}[Standard alphabet]\label{def:std-alphabet}
  Let us call the set
 \glo{s.greek.cap@\( \Sg_{0} \)}
 \begin{equation}\label{eq:std-alph}
 \Sg_{0}=\{0,1,\#,*\}
 \end{equation}
  the \df{standard alphabet}\index{standard!alphabet}.
 \end{definition}

Symbol \( \# \)\glo{,@\( *,\# \)} 
will be used to delimit binary strings, and \( * \) will serve as a
 ``don't-care'' symbol\index{don't-care symbol}.
The following definition formalizes this idea.

 \begin{definition}[Specification relation between strings]\label{def:preceq}
Each field \( \F \) of a cell state such that the field size is even, can be
considered not only a binary string but a string of (half as many) symbols
in the standard alphabet.
  If \( r,s \) are strings in \( (\Sg_{0})^{n} \) then%
 \glo{"<@\( \preceq \)}%
 \[
   r\preceq s
 \]
  will mean that \( s(i)=r(i) \) for all \( 0\le i<n \) such that \( r(i)\ne * \).
  For functions \( f, g \), with values in \( \Sigma_{0}^{*} \) we will write
\( f\preceq g \)  if for all \( s \) we have \( f(s)\preceq g(s) \).
 \end{definition}

Thus, a don't-care symbol \( r(i) \) imposes no restriction on \( s \) in
the relation \( r\preceq s \).




\begin{definition}\label{def:std-comp}
Let \( \trans \) be any deterministic
transition function with states in some alphabet \( \bbA \)
with distinguished fields \( \Input \), \( \Output \) in the standard alphabet.
We say that \( \trans \) has \df{monotonic output}\index{monotonic output} if
for all trajectories \( \eta \) of \( \CA(\trans) \) we have
 \[
  \eta(x,t).\Output \preceq \eta(x,t+1).\Output.
\]
  We call the transition function \( \trans \), all of whose fields have even
size (that is they are in the standard alphabet) together with some
distinguished fields \( \Input \), \( \Output \), a
\df{standard computing transition function} if it
 \begin{cjenum} 


  \item never changes \( \Input \);

  \item has monotonic output;

  \item does not change anything if
    the middle argument has all \( \# \)'s in its input field, or if all three
    arguments have \( * \) in all their fields.
    
 \end{cjenum}
We will call a cellular automaton with such a transition function a
\df{standard computing medium}\index{medium!standard computing}.
Our construction uses for this a special field called 
\begin{align*}
 \Color.
\end{align*}
It can take a constant number of values; in one case, only \( -1,0,1 \).
\end{definition}

The computation will occur in the segment of the infinite space
marked with color \( 0 \); the space to the left
and right of it will be marked with color \( -1 \) and \( 1 \) respectively.
The reason for this will be explained in Section~\ref{sec:plan}.
 
\begin{definition}[Init]\label{def:Init}
Let \( \Sigma_{0} \) be the standard input/output alphabet, and
\( \states \)  a set of states with \( \states.\Input=\Sigma_{0} \),
having another field \( \Color\in\{-1,0,1\} \), and
a distinguished ``latent'' state \( s_{0} \).
Let \( \psi_{*}:\Sigma_{0}^{*}\to\states^{*} \) be an encoding.
For any finite string \( \rho\in\Sigma_{0}^{n} \), let \( n'=|\psi_{*}(\rho)| \).
For any finite or infinite \( N>|\rho'| \), we construct the initial configuration
  \glo{init@\( \Init_{\psi_{*}}(\rho) \)}
	\begin{equation*}
      \xi'=\Init_{\psi_{*}}(\rho)\in\states^{N}
    \end{equation*}
    as follows.
    Define \( \xi \) by setting
    \( \xi(x)=\psi_{*}(\rho)(x) \) for each \( x\in\lint{0}{n'} \), and \( s_{0} \)
    for all other \( x \).
    In case \( N=\infty \),  obtain \( \xi' \) from \( \xi \) by setting
    \begin{align*}
    \xi'(x).\Color=
      \begin{cases}
    0  & \text{ for } 0\le x<n, 
\\    -1 & \text{ for } x<0, 
\\    1 & \text{ for } x\ge n.
      \end{cases}
    \end{align*}
    When the space is the finite one \( \bbZ_{N} \) then, with \( r=(N-n)/2 \),
 set \( \xi'(x).\Color=-1 \) for \( -r\le x<0 \),  and \( 1 \) for \( n\le x < n+r \).
\end{definition}

For some \( c>1 \) let
\begin{align}\label{eq:h_1}
  h_{1}(t,c)= c^{(\log t)^{1/2}\log\log t}.
\end{align}

\begin{sloppypar}
 \begin{theorem}[Reliable computation in dim 1] \label{thm:1dim.comp}
  Let \( \trans \) be a standard computing transition function over an alphabet \( \Sigma \),
  and \( \delta>0 \) a constant.
  There is
 \begin{itemize}
  \item a transition function \( \trans' \) with a state space \( \states \) having fields \( \Input \),
\( \Output \), \( \Color \), with \( \states.\Input=\states.\Output=\Sigma_{0} \).
\item a code \( \psi_{*}: \Sigma_{0}^{*}\to \states^{*} \), with  \( |\psi_{*}(\rho)|\le 2|\rho| \);
\item a constant \( c_{1}>1 \)   
 \end{itemize}
 such that the following holds for all \( \rho\in\Sigma_{0}^{*} \).
 Let \( \zeta(x,t) \) be any trajectory of \( \trans \) with
 \( \zeta(x,0).\Input=\rho(x) \) for \( x\in\lint{0}{|\rho|} \) and filled with \( * \)'s
 otherwise.
 Let \( \eta \) be any trajectory of an \( \eps \)-perturbation of \( \trans' \), with
 \( \eta(\cdot,0)=\Init_{\psi_{*}}(\rho) \).
 For all \( t \), all \( t'> t\cdot h_{1}(t,c_{1}) \), all \( x\in\bbZ \) we have
 for all sufficiently small \( \eps \):
 \[
 \Pbof{\zeta(x, t).\Output \not\preceq \eta(x, t').\Output} <\delta.
\]
\end{theorem} 
\end{sloppypar}
So the slowdown paid for reliability is somewhat worse than logarithmic:
by a factor \( h_{1}(t,c_{1}) \).
We will see that the code \( \psi_{*} \) is simple, hiding no complex computation.
The proof of this theorem relies on a move we can call \df{lifting}:
though the input configuration is a hierarchically encoded one, but only to the level
needed by the size of the input.
Higher levels will be built up as additional reliability is needed for longer computation time.

\begin{remarks}\label{rem:sorg-poorer}
  \begin{enumerate}
  \item 
  This note is similar to Remark~\ref{rem:sorg-poor}, but is even more important,
  as the error term does not even decrease with increasing \( t \).
\item The initialization via \( \Init \) avoids building up an infinite hierarchy in the
  initial configuration, leaving this to self-organization.
  But it introduces a huge asymmetry by coloring the left of the input with \( -1 \) and
  the right with \( 1 \).
  This has an effect similar to making the space infinite in only
  one direction, and placing the input into the beginning.
  \end{enumerate}
 \end{remarks}

  \section{Codes}\label{sec:codes}

For the moment, let us focus on the task of remembering a single constant-size
field of a cellular automaton, called
\begin{align*}
   \MainBit.
\end{align*}
\glo{mainbit@\( \MainBit \)}
The intention is that in the initial configuration, each cell's \( \MainBit \) is 
the same, and the transition function should try to keep this field constant.
We mentioned in Section~\ref{sec:hier-constr} that in one
dimension, even this simple task will require the construction and
maintenance of some non-local organization, as this seems the only way
to eliminate large islands.
The first idea is to store pieces of information redundantly in segments,
which we will call colonies.
But then colonies will be organized into supercolonies, and so on.
To simplify this idea, we will say that colonies encode cells of another 
cellular automaton (which can have its own colonies\( \dots \)).

 \subsection{Colonies}\label{sec:col}

A colony is a fixed-size segment of cells that are supposed to cooperate.

 \begin{definition}[Colony]
Let \( x \) be a site and \( \Q \) a positive integer.
The set of \( \Q \)\glo{q.cap@\( \Q \)} sites \( x+i \) for \( i\in\lint{0}{\Q} \) will be
called the \( \Q \)-\df{colony}\index{colony} with \df{base}\index{colony!base}
\( x \), and site \( x+i \) will be said to have \df{address}\index{address} \( i \) in
this colony.
 \end{definition}

Let us be given a configuration \( \xi \) of a cellular automaton \( M \)
with state set \( \states \).
The fact that \( \xi \) is ``organized into colonies''\index{organized into
colonies} will mean that one can break up the set of all sites into
non-overlapping colonies of size \( \Q \), using the information in the
configuration \( \xi \) in a translation-invariant way.
This will be achieved with the help of an address field.

 \begin{definition}[Address]
An \df{address field} is a field of our cells that we will denote
\( \Addr \)\glo{addr@\( \Addr \)}: we will always have it
when we speak about colonies.
The value \( \xi(x).\Addr \) is a binary string which will be interpreted
as an integer in \( \lint{0}{\Q} \).
 \end{definition}

Generally, we will assume that the \( \Addr \) field is large enough (its
size is at least \( \log \Q \)).
Then we could call a certain \( \Q \)-colony \( \cC \) a ``real'' colony of
\( \xi \) if for each element \( y \) of \( \cC \) with address \( i \) we have
\( \xi(y).\Addr=i \).
(But we do not introduce such a definition, since the address field of 
some cells could be faulty.)
Cellular automata working with colonies will not change the value of
the address field unless it seems to require correction.
In the absence of faults, if such a cellular automaton is started with a
configuration grouped into colonies then the sites can always use the
\( \Addr \) field to identify their colleagues within their colony.
Grouping into colonies seems to help preserve the \( \MainBit \)
field since each colony has this information in \( \Q \)-fold redundancy.
The transition function may somehow involve the colony members in a
coordinated activity of restoring this information from the
degradation caused by faults (for example with the help of some majority
operation).
This activity could be repeated periodically.

 \begin{definition}[Work period]
For a parameter \( \U  \) \glo{u.cap@\( \U  \)}
that we will fix in each case, let us call \( \U  \) steps of work of a colony a
\df{work period}\index{period!work}.
 \end{definition}

  The best we can expect from a transition function of the kind described
above is that unless too many faults happen during some colony work
period the \( \MainBit \) field of most sites in the colony will always
be the original one.
One can indeed write such transition functions; however,
they do not accomplish qualitatively much more than a local
majority vote for the \( \MainBit \) field among three neighbors.
Suppose that a group of failures changes the original content of the
\( \MainBit \) field in some colony, in so many sites that internal
correction is no more possible.
The information is not entirely lost since most probably, neighbor
colonies still have it.
But correcting the information in a whole colony with the help of other
colonies requires organization reaching wider than a single colony.
To arrange this broader activity also in the form of a cellular
automaton we use the notion of simulation with error-correction.

We will codify another convention:

\begin{definition}
Let us denote by \( M_{1} \) the fault-tolerant cellular
automaton to be built.  
\end{definition}

In the automaton \( M_{1} \), a colony \( \cC \)\glo{c.cap@\( \cC \)} with base \( x \) will be
involved in two kinds of activity during each of its work periods.
 \begin{description}

  \item[Simulation]\index{simulation}
   Manipulating the collective information of the colony in a way that can
be interpreted as the simulation of a single state transition of site \( x \)
of some cellular automaton \( M_{2} \).

  \item[Error-correction]\index{error!correction}
  Using the collective information (the state of \( x \) in \( M_{2} \)) to
correct each site within the colony as necessary.

 \end{description}
  Of course, even the sites of the simulated automaton \( M_{2} \) will not be
immune to errors.
  They must also be grouped into colonies simulating an automaton \( M_{3} \),
and so on; the organization must be a \emph{hierarchy of simulations}
(more precise definitions follow).

\subsection{Block codes}\label{sec:block-code}

  The notion of simulation relies on the notion of a \df{code}\index{code},
since the way the simulation works is that the simulated history
can be decoded from the simulating history.
We will develop a system of codes, starting from the simplest,
well-known example of a block code over strings.

 \begin{definition}[Code on strings]\label{def:code}
  A code \( \fg \)\glo{f.greek@\( \fg_{*},\fg_{*} \)} between two sets \( R,S \) is
a pair \( (\fg^{*},\fg_{*}) \) where \( \fg_{*}: R\to S \) is the
\df{encoding function} and \( \fg^{*}: S\to R \) is the \df{decoding function},
and the relation
 \[
  \fg^{*}(\fg_{*}(r)) = r
 \]
   holds.
We will be particularly interested in the example when 
for a positive integer \( \Q \) called the \df{block size} and some finite sets
\( \states_{1},\states_{2} \) we have \( R=\states_{2} \),
\( S=\states_{1}^{\Q} \).
Such a code is called a \df{block code}.
In a block code, strings of the form \( \fg_{*}(r) \) are called \df{codewords}.
The elements of a codeword \( s=\fg_{*}(r) \) are numbered as
\( s(0),\dots,s(\Q-1) \).

The codes \( \fg \) between sets \( R, S \) used in our simulations will have a
feature similar to the acceptance and rejection of 
Example~\ref{xmp:block-code}.
The set \( R \) will always have a subset of symbols called
\df{vacant}\index{vacant symbols}, 
and among them a distinguished special symbol \( \Vac \)\glo{vac@\( \Vac \)}.
An element \( s\in S \) will be called
\df{accepted}\index{accepted by decoding} by the decoding if
\( \fg^{*}(s)\ne \Vac \), otherwise it is called \df{rejected}.
   \end{definition}
Having more than one vacant symbol is just a matter of convenience.

A simple example code \( \gamma \) would be \( R=\{0,1\} \), \( S=R^{3} \),
\( \gamma_{*}(r)=(r,r,r) \) while \( \gamma^{*}((r,s,t)) \) is the majority of \( r,s,t \).

   \begin{remark}
     The notation \( (\fg^{*},\fg_{*}) \) to use for decoding and encoding is not in
common use: I find it suggestive, though, since \( \fg^{*} \) is something
like an inverse function without actually being one.
   \end{remark}

The following block code can be considered a significantly more complex,
paradigmatic example of the codes we will use later.

\begin{example}\label{xmp:block-code}
Suppose that \( \states_{1}=\states_{2}=\{0,1\}^{12} \) is the state set of both
cellular automata \( M_{1} \) and \( M_{2} \).
Let us introduce the fields \( s.\Addr \) and
\( s.\Info \)\glo{info@\( \Info \)} of a state \( r=(s_{0},\dots, s_{11}) \) in
\( \states_{1} \).
The \( \Addr \) field consists of the first 5 bits \( s_{0},\dots,s_{4} \),
while the \( \Info \) field is the last bit \( s_{11} \).
The other bits do not belong to any named field.
Let \( \Q=31 \).
Thus, we will use codewords of size \( 31 \), formed of the symbols (local
states) of \( M_{1} \), to encode local states of \( M_{2} \).
The encoding funcion \( \fg_{*} \) assigns a codeword
\( \fg_{*}(r)=(s(0),\dots, s(30)) \) of elements of \( \states_{1} \) to each
element \( r \) of \( \states_{2} \).
Let \( r=(r_{0},\dots, r_{11}) \).
We will set \( s(i).\Info=r_{i} \) for \( i=0,\dots,11 \).
The 5 bits in \( s(i).\Addr \) will denote the number \( i \) in binary
notation.
This did not determine all bits of the symbols \( s(0),\dots, s(30) \) in
the codeword.
In particular, the bits belonging to neither the \( \Addr \) nor the
\( \Info \) field are not determined, and the values of the \( \Info \) field
for the symbols \( s(i) \) with \( i\not\in\lint{0}{12} \) are not determined.
To determine \( \fg_{*}(r) \) completely, we could set these bits to 0.

The decoding function is simpler.
Given a word \( s=(s(0),\dots,s(30)) \) we first check whether it is a
``normal'' codeword, and as such, 
has \( s(0).\Addr=0 \) and \( s(i).\Addr\ne 0 \) for \( i\ne 0 \).
If yes then \( r=\fg^{*}(s) \) is defined by \( r_{i}=s(i).\Info \) for
\( i\in\lint{0}{12} \), and the word is accepted'.
Otherwise \( \fg^{*}(s)=0\cdots 0 \), and the word is rejected.

Informally, the symbols of the codeword use their first 5 bits to mark
their address within the codeword.
The last bit is used to remember their part of the information about
the encoded symbol.
 \end{example}

\begin{figure}
 \setlength{\unitlength}{0.2mm}
 \[
 \begin{picture}(600,150)
    \put(0, 0){\framebox(600, 150){}}
    \put(0, 20){\line(1, 0){600}}
    \put(0, 90){\line(1, 0){600}}
    \put(0, 110){\line(1, 0){600}}
    \put(0, 130){\line(1, 0){600}}

    \put(610, 0){\makebox(0,20)[l]{\( \var{Program} \)}}
    \put(610,20){\makebox(0,70)[l]{\( \var{Info} \)}}
    \put(610,90){\makebox(0,20)[l]{\( \var{Worksp} \)}}
    \put(610,110){\makebox(0,20)[l]{\( \var{Addr} \)}}
    \put(610,130){\makebox(0,20)[l]{\( \var{Age} \)}}

   \multiput(0,0)(20,0){20}{\line(0,1){3}}
   \put(200,0){\line(0, 1){150}}
   \multiput(400,0)(20,0){10}{\line(0,1){150}}
 \end{picture}
 \]
 \caption{Three neighbor colonies with their tracks}
 \label{fig:3-neighb}
\end{figure}

 \begin{example}\label{xmp:code.aggreg}
The trivial example here will not be really used as a code but rather as a
notational convenience.
For every symbol set \( \states_{1} \), blocksize \( \Q \) and \( \states_{2}=\states_{1}^{\Q} \),
there is a special block code \( \iota_{\Q} \)\glo{i:greek@\( \iota_{\Q} \)} called
\df{aggregation}\index{aggregation} defined by
 \[
  \iota_{\Q}^{*}((s(0),\dots,s(\Q-1)))=s(0)\cc\cdots\cc s(\Q-1),
 \]
and \( \iota_{\Q*} \) defined accordingly.
Thus, \( \iota_{\Q}^{*} \) is essentially the identity: it just aggregates \( \Q \)
symbols of \( \states_{1} \) into one symbol of \( \states_{2} \).
We use concatenation here since we identify all symbols with binary
strings.
\end{example}

\begin{definition}[Aggregating a field]\label{def:aggreg-field}
  Let \( \fg \) be a block code from \( \states_{2}\) to \( \states_{1}^{\Q} \),
  and \( \F_{1},\F_{2} \) be fields of \( \states_{i} \).
  We will say that \( \fg \) \df{aggregates} field \( \F_{1} \) into \( \F_{2} \)
  if 
\begin{align*}
 \fg^{*}((s(0),\dots,s(\Q-1))).\F_{2}=s(0).\F_{1}\cc\dots\cc s(\Q-1).\F_{1}.
\end{align*}
Hence for \( u,v\in\states_{2} \) the assumption \( u.\F_{2}=v.\F_{2} \)
implies that for all \( a\in\lint{0}{\Q} \) we have
\begin{align}\label{eq:aggreg-field-prop}
 \fg_{*}(u)(a).\F_{1}=\fg_{*}(v)(a).\F_{1} .
\end{align}
 \end{definition}

 We will use a strengthening of the aggregating property for a particular address \( a \).

 \begin{definition}[Controlling by a field]\label{def:control}
 Assume that code \( \fg \) aggregates field \( \F_{1} \) into \( \F_{2} \).
 We say that field \( \F_{1} \) \df{controls} address \( a \) for code \( \fg \)
 via a function \( \gamma:\states_{1}.\F_{1}\to\states_{1} \) if
 in all codewords \( w \) of the form
 \( w=\fg_{*}(u) \), we have \( w(a)=\gamma(w(a).\F_{1}) \).
\end{definition}

In general, a symbol \( w(a) \) of a codeword \( w=\fg_{*}(u) \)
may contain information in its fields other than \( \F_{1} \) about \( u \), and it is
typically important for some symbols of \( w \) to do it.
But there is no harm in freeing any one particular address \( a \), say \( a=1 \),
from this responsibility.

\subsection{Generalized cellular automata (media)}\label{sec:media} 

  A block code \( \fg \) could be used to define a code on configurations
between cellular automata \( M_{1} \) and \( M_{2} \).
Suppose that a configuration \( \xi \) of \( M_{2} \) is given.
Then we could define the configuration
 \( \xi_{*}=\fg_{*}(\xi) \)\glo{f.greek.conf@\( \fg_{*}(\xi),\fg^{*}(\xi) \)}
\glo{x.greek@\( \xi_{*},\xi^{*} \)}%
of \( M_{1} \) by setting for each cell \( x \) of \( \xi \) and \( 0\le i < \Q \),
  \[
   \xi_{*}(\Q x + i)=\fg_{*}(\xi(x))(i).
  \]
  The decoding function would be defined correspondingly.
  This definition of decoding is, however, unsatisfactory for our
purposes. 
  Suppose that \( \xi_{*} \) is obtained by encoding a configuration
\( \xi \) via \( \fg_{*} \) as before, and \( \zeta \) is obtained by shifting
\( \xi_{*} \): \( \zeta(x)=\xi_{*}(x-1) \).
Then the decoding of \( \zeta \) will return all vacant values since now the
strings \( (\zeta(\Q x),\cdots,\zeta(\Q x+\Q -1)) \) are not ``real'' colonies.
However, it will be essential for reasoning about error correction that
decoding should notice all parts of a configuration that form a colony, 
even if a shifted one.
With our current definition of cellular automata, the decoding function
could not be changed to do this.
Indeed, if \( \zeta^{*} \) is the configuration decoded from \( \zeta \) then
\( \zeta^{*}(0) \) corresponds to the value decoded from
\( \tup{\zeta(0),\dots,\zeta(\Q -1)} \), and \( \zeta^{*}(1) \) to the value decoded from
\( \tup{\zeta(\Q ),\dots,\zeta(2\Q -1)} \).
There is no site to correspond to the value decoded from
\( \tup{\zeta(1),\dots,\zeta(\Q )} \).

Our solution is to generalize the notion of cellular automata.
Let us give at once the most general definition which then we will
specialize later in varying degrees.
 
 An \df{medium}\index{medium}, or \df{generalized cellular automaton}
is given by the following ingredients:
 \[
 \states, \sites, \B, \Configs, \Histories, \Trajs, (\cA_{t})_{t\ge 0}.
\]
\( \states \) is the set of possible local states, and it has a subset of distinguished states
called \( \Vacant \), and among them one distinguished state
called \df{the vacant state}\index{state!vacant} \( \Vac\in\Vacant \).
(Definition~\ref{def:code} already hinted at the meaning of the vacant states.)
\( \sites \) is the set of sites as introduced in Definition~\ref{def:space},
 in our case the product of a few sets of the form \( \bbZ_{m} \) for finite or infinite \( m \).
The positive integer \( \B  \)\glo{b.cap@\( \B  \)} is called the \df{body size}.
In ordinary cellular automata, \( \B =1 \).
\( \Configs \)\glo{configs@\( \Configs \)}
is the set of functions \( \xi:\sites\to\states \) that are \df{configurations}.
\index{configuration}
If in a configuration \( \xi \) we have \( \xi(x)\not\in\Vacant \) then we will say
that there is a \df{cell}\index{cell} at site \( x \) in \( \xi \).
For a site \( x \), interval \( \lint{x}{x+\B } \) will be called the
\df{body}\index{body} of a possible cell with \df{base} \( x \).
In a \df{configuration}, cells must have non-intersecting bodies.
\( \Histories \)\glo{evols@\( \Histories \)}  is the set of functions
\( \eta:\sites\times\lint{0}{\infty}\to\states \) that are \df{histories}.
\index{history}
In a history \( \eta \) is is required that 
 \begin{cjenum}

  \item \( \eta(\cdot,t) \) is a space configuration for each \( t \);

  \item \( \eta(x,t) \) is a right-continuous function of \( t \);

  \item Each finite time interval contains only finitely many switching
times (see Definition~\ref{def:space-time}) for each site \( x \);

 \end{cjenum}
\( \Trajs \)\glo{trajs@\( \Trajs \)} is the set of random
histories \( (\mu,\eta) \) that are \df{trajectories}\index{trajectory}.
The increasing system \( (\cA_{t})_{t\ge 0} \) of \( \sigma \)-algebras was introduced in
Section~\ref{sec:PCA}.
The sets \( \states \), \( \sites \), \( \Configs \) and \( \Histories \) are 
defined by the set \( \Trajs \) implicitly---therefore
we may omit them from the notation, and write
 \[
   \Med(\B,\Trajs, \{\cA_{t}\}) =
   \Med(\states, \sites, \B, \Configs, \Histories, \Trajs, (\cA_{t})_{t\ge 0})
 \]
 for a medium defined by them.
We have a body size parameter \( \B \) because
some of our cellular automata will ``live in simulation'',
and it is very convenient for the simulated cell's body to coincide
with the interval of the ``base'' cells simulating it.

 \begin{definition}[Dwell period]
A \df{dwell period}\index{period!dwell} of \( \eta \)
is a tuple \( (x,s,t_{1},t_{2}) \)
such that \( x \) is a site, \( s \) is a nonvacant state, and \( 0\le t_{1}<t_{2} \)
are switching times with \( \eta(x,t_{1})=s \).
The rectangle \( \lint{x}{x+\B }\times \lint{t_{1}}{t_{2}} \) is the
\df{space-time body}\index{body!space-time} of the dwell period.
 \end{definition}
Ordinary cellular automata obey some further restrictions:

 \begin{definition}[Lattice configuration]
   A configuration is a \df{lattice
     configuration}\index{configuration!lattice} if all cells are at sites of
the form \( i\B  \) for integers \( i \).
We can also talk about \df{lattice histories}\index{history}: these have
space-time bodies of the form
 \[
   \lint{i\B }{(i+1)\B }\times\lint{j\T }{k\T }
 \]
  for integers \( i,j<k \).
 \end{definition}

A lattice history is a history of an ordinary cellular automaton, except 
that cells are given a common size \( \B  \) possibly different from 1, and the 
time step is given a size \( \T  \) possibly different from 1.

 \begin{definition}[Deterministic cellular automaton with size parameters]
\label{def:determ-ca-param}
A \df{deterministic cellular automaton}\index{cellular
automaton!deterministic}%
 \glo{ca@\( \CA(\trans, \B , \T , \sites) \)}%
 \[
   \CA(\trans, \B , \T , \sites)
 \]
 is determined by parameters 
\( \B ,\T >0 \) and a transition function \( \trans:\states^{3}\to\states \).
We may omit some obvious arguments from this notation.
A lattice history \( \eta \) with parameters \( \B ,\T  \) is a
\df{trajectory}\index{trajectory!deterministic cellular automaton} of this
automaton if
 \[
   \eta(x,t)=\trans(\eta(x-\B ,t-\T ), \eta(x,t-\T ), \eta(x+\B ,t-\T ))
 \]
 holds for all \( x,t \) with \( t\ge \T  \).
For a history \( \eta \) let us write 
 \begin{equation}\label{eq:trans-B}
   \ol\trans(\eta, x,t,\B )=\trans(\eta(x-\B ,t), \eta(x,t), \eta(x+\B ,t)).
 \end{equation}
 \end{definition}

Probabilistic cellular automata and perturbations are generalized
correspondingly as%
 \glo{pca@\( \PCA(\bP, \B , \T , \sites) \)}%
\begin{align}\label{eq:gen-pca}
   \PCA(\bP, \B , \T , \sites),\qq \CA_{\eps}(\trans, \B , \T , \sites).
  \end{align}
From now on, whenever we talk about a deterministic, probabilistic or
perturbed cellular automaton we understand one also having parameters
\( \B ,\T  \).

We will also consider cellular automata in which the dwell periods
do not have a fixed length:

 \begin{definition}[Variable-period cellular medium]
A medium is said to have a \df{variable period}\index{medium!variable-period}
if in its histories, not all dwell periods have necessarily the same size.
 \end{definition}

\paragraph{Block codes between cellular automata}
In a cellular abstract medium with body size \( \B  \), a
colony\index{colony!in a cellular medium} of size \( \Q  \) is going to be some
set of cells \( x+i\B  \) for \( i\in\lint{0}{\Q } \).
Thus, the union of its cell bodies occupies the interval \( \lint{x}{x+\Q \B } \).
In what follows we use decoding to recognize colonies (when reasoning about a 
configuration).

 \begin{definition}[Overlap-free]
A block code will be called
\df{overlap-free}\index{code!block!overlap-free} if for every string
\( \tup{s(0),\dots s(n-1)} \), and all \( i\le n-\Q  \), if both
strings \( \tup{s(0),\dots,s(\Q -1)} \)
and \( \tup{s(i+1),\dots,s(i+\Q -1)} \) are accepted then \( i\ge \Q  \).
 \end{definition}

In other words, a code is overlap-free if two accepted words cannot
overlap in a nontrivial way.

The simple tripling code \( \fg_{*}(x)=(x,x,x) \) with majority decoding
is not overlap-free, since for example if \( s(1)s(2)s(3)s(4)=0000 \) 
then both \( s(1)s(2)s(3) \) and \( s(2)s(3)s(4) \) are accepted.
On the other hand, the code \( \fg_{*}(x)=(0,x,1,x,2,x) \) 
with \( \fg^{*}(0,x,1,y,2,z)=\Maj(x,y,z) \) (and vacant in other cases)
is overlap-free.
The code in Example~\ref{xmp:block-code} is similarly overlap-free.
Overlap-free codes are used, among others, in~\cite{ItkisLevin94}.

A block code \( \fg \) of block size \( \Q  \) can be used to define a code on
configurations between cellular abstract media \( M_{1} \) and \( M_{2} \).

\begin{definition}[Block code on configuration]
Suppose that a configuration \( \xi \) of \( M_{2} \), which is a cellular medium
\( \AMed(\Q \B ) \), is given.
Then we define the encoded configuration \( \xi_{*}=\fg_{*}(\xi) \) of \( M_{1} \),
which is an \( \AMed(\B ) \), by setting for each cell \( x \) of \( \xi \) and
 \( 0\le i < \Q  \),
  \[
   \xi_{*}(x + i \B )=\fg_{*}(\xi(x))(i).
  \]
Suppose that a configuration \( \xi \) of \( M_{1} \) is given.
We define the decoded configuration
 \( \xi^{*}=\fg^{*}(\xi) \) of \( M_{2} \) as follows:
for site \( x \), set \( \xi'(x)=\fg^{*}(s) \) where
 \begin{equation}\label{eq:cand.colony}
   s = \tup{\xi(x), \xi(x + \B ), \dots, \xi(x + (\Q  - 1)\B )}.
 \end{equation}
We define \( \xi^{*}(x)=\xi'(x) \) if the latter is vacant or
there is no \( y \) closer than \( \Q \B  \) to \( x \)
with \( \xi'(y) \) non-vacant; otherwise \( \xi^{*} = \Vac \).
\end{definition}

\begin{proposition}
Given a block code \( \fg=(\fg_{*},\fg^{*}) \), its extension to configurations is
also a code in the sense that the equation
 \begin{align*}
   \fg^{*}(\fg_{*}(\xi))=\xi
 \end{align*}
still holds.
\end{proposition}

The proof of this statement is immediate.

As we know most configurations cannot be obtained by encoding.
Those that can merit a special name.

\begin{definition}[Code configuration]\label{def:code-config}
A configuration \( \xi \) is called a \df{code configuration} if
\(  \fg_{*}(\fg^{*}(\xi))=\xi \) and
the decoded configuration \( \fg^{*}(\xi) \) covers the space with adjacent cells.
\end{definition}

For example, if the code is the simple tripling \( \gamma_{*}(x)=(x,x,x) \) mentioned
above, then any code configuration would consist of triply repeated symbols.

Let \( \fg \) be an overlap-free code.
If \( \xi = \fg_{*}(\zeta) \) is a code configuration with \( \xi^{*}=\fg^{*}(\xi) \)
then \( \xi^{*}(x) \) is
nonvacant only at positions \( x \) where \( \zeta(x) \) is nonvacant.
If \( \xi \) is not a code configuration then it may happen
that in \( \xi^{*} \), the cells will 
not be exactly at a distance \( \Q  \B  \) apart.
Our construction still garantees that the distance of
cells in \( \fg^{*}(\xi) \) is at least \( \Q  \B  \).
This situation can be taken as one of the justifications for the notion of
cellular abstract media.

\subsection{Block simulations}\label{sec:block.sim}

A simulation is defined as a code allowing to decode the trajectory of
one cellular automaton from another.

 \begin{definition}[Block simulation between deterministic cellular automata]
Suppose that \( M_{1} \) and \( M_{2} \) are deterministic cellular automata
where \( M_{i}=\CA(\trans_{i},\B _{i},\T_{i}) \), and \( \fg \) is a block code with
 \[
  \B _{1}=\B ,\qq \B _{2}=\Q \B .
 \]
  The decoding function may be as simple as in Example~\ref{xmp:block-code}:
there is an \( \Info \) track and once the colony is 
accepted the decoding function depends only on this part of the
information in it.

For each history \( \eta \) of \( M_{1} \), we define the history
\( \eta^{*}=\fg^{*}(\eta) \) of \( M_{2} \) by setting
 \begin{equation}
  \eta^{*}(\cdot,t)=\fg^{*}(\eta(\cdot,t)).
 \end{equation}
  We will say that the code \( \fg \) is a
 \df{simulation}\index{simulation!block code} if for each configuration
\( \xi \) of \( M_{2} \), for the (unique) trajectory \( \eta \) 
of \( M_{1} \), such that \( \eta(\cdot,0)=\fg_{*}(\xi) \), the
history \( \eta^{*} \) is a trajectory of \( M_{2} \).
 \end{definition}  

  We can view \( \fg_{*} \) as an encoding of the initial configuration of
\( M_{2} \) into that of \( M_{1} \).
  Our requirements say that from every trajectory of \( M_{1} \) with
  a ``good'' initial configuration (that is a code configuration), the
 simulation-decoding results in a
trajectory of \( M_{2} \).

 \begin{example}
   Let us show one particular way in which the transition
function \( \trans_{1} \) can make the code \( \fg \) a simulation.
  Assume that
 \[
  \T_{1}=\T ,\ \T_{2}=\U \T 
 \]
  for some positive integer \( \U  \) called the
 \df{work period size}\index{period!work!size}.
Each cell of \( M_{1} \) will go through a period consisting of \( \U  \) steps
in such a way that the \( \Info \) field will be changed only in the last
step of this period.
The initial configuration \( \eta(\cdot,0)=\fg_{*}(\xi) \) is chosen in
such a way that each cell is at the beginning of its work period.
By the nature of the code, in the initial configuration, cells of
\( M_{1} \) are grouped into colonies.

Once started from such an initial configuration, during each work
period, each colony, in cooperation with its two neighbor colonies,
computes the new configuration.
With the block code in Example~\ref{xmp:block-code}, this may happen
as follows.
Let us denote by \( r_{-1},r_{0},r_{1} \) the value in the first 12 bits of
the \( \Info \) track in the left neighbor colony, in the colony itself
and in the right neighbor colony respectively.
First, \( r_{-1} \) and \( r_{1} \) are shipped into the middle colony.
Then, the middle colony computes \( s=\trans_{2}(r_{-1},r_{0},r_{1}) \)
where \( \trans_{2} \) is the transition function or \( M_{2} \) and stores it on
a memory track.
(It may help understanding how this happens if we think of the
possibilities of using some mail, memory and workspace tracks.)
Then, in the last step, \( s \) will be copied onto the \( \Info \) track.
 \end{example}

 \begin{example}[Aggregated transition]\label{xmp:sim.aggreg} 
Here is a trivial example of a block simulation which will be
applied, however, later in the paper.
Given a one-dimensional transition function \( \trans(x,y,z) \) with
state space \( \states \), we can define for all positive integers \( \Q  \) an
\df{aggregated} transition function\index{transition function!aggregated}
 \( \trans^{\Q }(u,v,w) \)\glo{tr@\( \trans^{\Q }(u,v,w) \)} as follows.
The state space of \( \trans^{\Q } \) is \( \states^{\Q } \).
Let \( \br_{j}=\tup{r_{j}(0),\dots, r_{j}(\Q -1)} \) for \( j=-1,0,1 \) be three
elements of \( \states^{\Q } \).
Concatenate these three strings to get a string of length \( 3\Q  \) and
apply the transition function \( \trans \) to each group of three
consecutive symbols to obtain a string of length \( 3\Q  - 2 \) (the end
symbols do not have both neighbors).
Repeat this \( \Q  \) times to get a string of \( \Q  \) symbols of \( \states \):
this is the value of \( \trans^{\Q }(\br_{-1},\br_{0},\br_{1}) \),
defining an \df{aggregated cellular automaton} \( \CA(\trans^{\Q},\Q,\Q) \) where
we used the notation in Definition~\ref{def:determ-ca-param}.

For \( M_{1}=\CA(\states,\trans,\B ,\T ) \) and
\( M_{2}=\CA(\states^{\Q },\trans^{\Q },\Q \B ,\Q \T ) \), the aggregation code \( \iota_{\Q } \)
defined in Example~\ref{xmp:code.aggreg} will be a block simulation of
\( M_{2} \) by \( M_{1} \) with a work period consisting of \( \U =\Q  \) steps.
If along with the transition function \( \trans \), there were some fields
\( \F,\G,\dots\subseteq\All \) also defined then we define, say, the field \( \F \) in
the aggregated cellular automaton as \( \bigcup_{i=0}^{\Q -1}(\F+i\nm{\states}) \).
Thus, if \( \br=r(0)\cc\dotsm\cc r(\Q -1) \) is a state of the aggregated
cellular automaton then
 \( \br.\F=r(0).\F\cc r(1).\F\cc\dotsm\cc r(\Q -1).\F \).

We may allow only less communication
with the neighbors, and then modify the above as follows to a function
\( \trans^{\Q,\bw } \) as follows, where we call \( \bv \) the \df{slowdown rate}, where
we assume that \emph{both \( \bv\Q \) and \( 1/\bv \) are integers}.
The state space is still \( \states^{\Q } \).
Let \( \br_{j}=\tup{r_{j}(0),\dots, r_{j}(\Q -1)} \) for \( j=-1,0,1 \) be three
elements of \( \states^{\Q } \).
Concatenate these three strings to get a string of length \( 3\Q  \) and
apply the transition function \( \trans \) to each group of three
consecutive symbols to obtain a string of length \( 3\Q  - 2 \) (the end
symbols do not have both neighbors).
Repeat this \emph{only \( \bw\Q  \) times}
to get a string of \( \Q  \) symbols of \( \states \):
this is the value of \( \trans^{\Q,\bv }(r_{-1},r_{0},r_{1}) \).
This defines the cellular automaton \( \CA(\trans^{\Q,\bv},\Q,\bv\Q) \).
The function \( \trans^{\Q} \) becomes the special case \( \trans^{\Q,1} \).
Denoting by \( s.\slice^{\bv,-1} \), \( s.\slice^{\bv,1} \)
the leftmost and rightmost  \( \bv|\Q| \) symbols of
a string \( s \) in \( \states^{\Q} \),
then in fact \( \trans^{\Q,\bv}(\br) \) depends only on
\( r_{-1}.\slice^{\bv,1},r_{0},r_{1}.\slice^{\bv,-1} \).
 \end{example}

To simulate \( n \) steps of a transition function \( \trans \)  over an interval \( \lint{a}{b} \)
start with a larger interval \( \lint{a-n}{b+n} \).
After step \( t \), we only need to work on the interval \( \lint{a-n+t}{b+n-t} \).

Cellular automata are ``computationally universal''
in the informal sense that Turing machines are: arbitrary computations 
can be implemented on appropriate cellular automata.
We will return to this point later, in Section~\ref{sec:func-prog}.
Also, there are cellular automata that are universal in a sense similar to
universal Turing machines.
The notion of block simulations give rise to a particularly simple
such universality property.

 \begin{definition}[Intrinsically universal cellular automata]\label{def:univ}
  A transition function \( \trans \) is
 \df{intrinsically universal}\index{transition function!intrinsically universal} if for every other
transition function \( \trans' \) there are \( \Q ,\U  \) and a block code \( \fg \) such
that \( \fg \) is a block simulation of \( \CA(\trans',\Q ,\U ) \) by
\( \CA(\trans,1,1) \).
\end{definition}

In the literature, intrinsic universality is defined in a slightly more general way:
see~\cite{OllingerRichard4states11}.
There are some extremely simple cellular automata that are computationally
universal but are probably not intrinsically so.
The most famous example is Cook's theorem for the computational
universality of Rule 110 in~\cite{Cook04}.

\begin{theorem}[Intrinsically universal cellular automata] \label{thm:univ}
 There is an intrinsically universal transition function.
 \end{theorem}

Here, we only give a proof sketch of existence.
A particularly simple intrinsically universal cellular automaton is defined 
in~\cite{OllingerRichard4states11}: it is one-dimensional, 
nearest-neighbor, with 4 states.

 \begin{proof}[Sketch of proof:]
This theorem is proved somewhat analogously to the theorem on the
existence of universal Turing machines.
If the intrinsically universal transition function is \( \trans \) then for simulating
another transition function \( \trans' \), the encoding demarcates colonies of
appropriate size with \( \Addr=0 \), and writes a string \( \var{Table} \) that is
the code of the transition table of \( \trans' \) onto a special track called
\( \fld{Prog} \)\glo{prog@\( \fld{Prog} \)} in each of these colonies.
The computation is just a table-look-up: the triple
\( (r_{-1},r_{0},r_{1}) \) mentioned in the above example must be looked up
in the transition table.
The transition function governing this activity does not depend on
the particular content of the \( \fld{Prog} \) track, and is therefore
independent of \( \trans' \).
For references to the first proofs of universality (in a technically
different but similar sense), see~\cite{BerlConwGuy82,TofMarg87}
 \end{proof}

Note that an intrinsically
universal cellular automaton cannot use codes similar to
Example~\ref{xmp:block-code}.
Indeed, in that example, the capacity of the cells of \( M_{1} \) is at
least the binary logarithm of the colony size, since each colony cell
contained its own address within the colony.
But if \( M_{1} \) is intrinsically
universal then the various simulations in which it
participates will have arbitrarily large colony sizes.

The size \( \Q  \) of the simulating colony in the above proof sketch
will generally be very large
also since the latter contains the whole table of the simulated
transition function.
There are many special cellular automata \( M_{2} \), however, whose
transition function can be described by a small computer program and
computed in relatively little space and time (linear in the size
\( \nm{\states_{2}} \)).
We will only deal with such automata, and will develop a simple
language to describe their programs.
(Any intrinsically universal transition function will simulate these with
correspondingly small \( \Q  \) and \( \U  \).)

\subsection{A single-fault-tolerant block simulation}\label{sec:ftol-block}

Here we outline a cellular automaton \( M_{1} \) that block-simulates a
cellular automaton \( M_{2} \) correctly as long as at most a single
error occurs in a colony work period of size \( \U  \).
The outline is very informal: it is only intended to give some
framework to refer to later: in particular, we add a few more fields
to the local states in addition to the ones introduced earlier.

The automaton \( M_{1} \) is not intrinsically universal, so the automaton
\( M_{2} \) cannot be chosen arbitrarily.
Among others, this is due to the fact that the address field of a
cell of \( M_{1} \) will hold its address within its colony
\footnote{We could avoid non-constant size fields, 
but since it is not necessary, for simplicity we keep them}.
But we will see later that intrinsic universality is not needed in this
context.

The cells of \( M_{1} \) will have, besides the \( \Addr \) field, also a
field \( \Age \)\glo{age@\( \Age \)}.
If no errors occur then in the \( i \)-th step of the colony work
period, each cell will have the number \( i \) in the field \( \Age \).
There are also fields called \( \Mail \), \( \Info \), \( \fld{Work} \), \( \Hold \),
\( \fld{Prog} \).

The \( \Info \) field holds the state of the represented cell of \( M_{2} \)
in three copies.
The \( \Hold \)\glo{hold@\( \Hold \)} field will hold parts of the final result
before it will be, in the last step of the work period, copied into
\( \Info \).
The role of the other fields is clear.

The program will be described from the point of view of a certain
colony \( \cC \).
Here is an informal description of the activities taking place in
the first third of the work period.

 \begin{enumerate}

  \item
From the three thirds of the \( \Info \) field, by majority vote, a
single string is computed.
Let us call it the \df{input string}.
This computation, as all others, takes place in the workspace field
\( \fld{Work} \); the \( \Info \) field is not affected.
The result is also stored in the workspace.

  \item
  The input strings computed in the two neighbor colonies are shipped
into \( \cC \) and stored in the workspace separately from each other and
the original input string.

  \item
  The workspace field behaves as a computationally
universal cellular automaton, and from the
three input strings and the \( \fld{Prog} \) field, computes the string that
would be obtained by the transition function of \( M_{2} \) from them.
This string will be copied to the first third of the \( \Hold \) track.

  \end{enumerate}

In the second part of the work period, the same activities will be
performed, except that the result will be stored in the second part
of the \( \Hold \) track.
Similarly with the third part of the work period.
In a final step, the \( \Hold \) field is copied into the \( \Info \) field.

The computation is coordinated with the help of the \( \Addr \) and
\( \Age \) fields.
It is therefore important that these are correct.
Fortunately, if a single fault changes such a field of a cell then
the cell can easily restore it using the \( \Addr \) and \( \Age \) fields of
its neighbors.
 
It is not hard to see that with such a program (transition function), if
the colony started with ``perfect'' information then a single fault will
not corrupt more than a third of the colony at the end of the work period.
On the other hand, if two thirds of the colony was correct at the
beginning of the colony work period and there is no fault during the
colony work period then the result will be ``perfect''.

 \subsection{General simulations}\label{sec:gen-sim}

The general notion of abstract media allows a very general definition of simulations:

 \begin{definition}[General simulation]\label{def:general-simulation}
A \df{simulation}\index{simulation} of medium \( M_{2} \) by 
medium \( M_{1} \) having the same system \( \{\cA_{t}\}_{t\ge 0} \) of \( \sigma \)-algebras 
is given by a pair
\glo{f.greek.cap@\( \Phi^{*} \)}
\begin{align}  \label{eq:Phi}
   \Phi=(\fg, \Phi^{*})  
\end{align}
where \( \Phi^{*} \) is a mapping of the set of histories of
\( M_{1} \) into those of \( M_{2} \) (the decoding), \( \fg_{*} \) is a mapping of
the set of configurations of \( M_{2} \) to the set of configurations of
\( M_{1} \) (the encoding for initialization).
We will assume to \( \fg_{*} \) always belongs also a decoding function \( \fg^{*} \),
so they form a code \( \fg=(\fg_{*},\fg^{*}) \).
Denote 
 \[
  \eta^{*}=\Phi^{*}(\eta).
 \]
 We require that for each configuration \( \xi \) of \( M_{2} \), each
 trajectory of \( \eta \) of \( M_{1} \) with
 \( \eta(\cdot,0) = \fg_{*}(\xi) \)
the decoded history\( \eta^{*} \) is a trajectory of \( M_{2} \)
with \( \eta^{*}(\cdot,0)=\xi \).
 \end{definition}


The actual simulations we will use will all have a certain locality property:
\begin{sloppypar}
\begin{definition}[Local simulation]
  A simulation will be called
\df{non-anticipating} \index{simulation!non-anticipating} if
for each \( t \) the random function
\( \eta^{*}(\cdot,t) \) is measurable in the \( \sigma \)-algebra \( \cA_{t} \)
(this says that it does not depend on the future).
Let \( \sigma \) be the shift on one-dimensional configurations and histories: that is
\( (\sigma\xi)(x)=\xi(x+1) \), \( (\sigma\eta)(x,t)=\eta(x+1,t) \).
A simulation \( (\fg,\Phi^{*}) \) is called \df{shift-invariant} if both \( \fg \) and
\( \Phi^{*} \) commute with the shift.
It is called \df{local}\index{simulation!local} 
if there is a finite space-time rectangle 
 \( V^{*}=I\times \rint{u}{0} \)\glo{v.cap@\( V^{*} \)} 
with \( u\ge 0 \) such that
\( \Phi^{*}(\eta)(w,t) \) depends only on \( \eta(\pair{w}{t}+V^{*}) \).
A local simulation is both non-anticipating and shift-invariant.

All our simulations will be local, unless stated otherwise.
It follows from the non-anticipating property that 
\( \eta^{*}(\cdot,0) \) depends only on
\( \eta(\cdot,0) \), and therefore the simulation always defines a decoding
function 
 \begin{align*}
  \fg^{*}
 \end{align*}
 on configurations (as already stipulated in Definition~\ref{def:general-simulation}). 
If \( u=0 \) and so the configuration \( \eta^{*}(\cdot,t) \) depends only on
the configuration \( \eta(\cdot,t) \), then the simulation will be called
\df{memoryless}\index{simulation!memoryless}.
\end{definition}  
\end{sloppypar}

Corollary~\ref{crl:var-sim} gives an example of non-local simulation.
Locality implies that a simulation is determined by a function defined on the 
set of configurations over \( V^{*} \).
Our eventual simulations will not be memoryless:
the decoding looks back on the history during
\( \rint{t-u}{t} \) but, being non-anticipating, still does not look ahead.
For a memoryless simulation, the simulation property is
identical to the one we gave at the beginning of
Section~\ref{sec:block.sim}.

 \paragraph{Simulation between perturbations}
  Our goal is to find nontrivial simulations between cellular automata
\( M_{1} \) and \( M_{2} \), especially when these are not deterministic.
If \( M_{1}, M_{2} \) are probabilistic cellular automata then the
simulation property would mean that whenever we have a trajectory
\( (\mu,\eta) \) of \( M_{1} \) the random history
\( \eta^{*} \) decoded from \( \eta \) would be a trajectory of \( M_{2} \).
There are hardly any nontrivial examples of this sort since in order
to be a trajectory of \( M_{2} \), the conditional probabilities of
\( \fg^{*}(\eta) \) must satisfy certain equations defined by
\( \bP_{2} \), while the conditional probabilities of \( \eta \) satisfy
equations defined by \( \bP_{1} \).

There is more chance of success in the case when \( M_{1} \) and \( M_{2} \) are
perturbations of some deterministic cellular automata since in this case,
only some inequalities must be satisfied.
The goal of improving reliability could be this.
For some universal transition function \( \trans_{2} \), and at least two
different initial configurations \( \xi_{i} \) (\( i=0,1 \)), find
\( \trans_{1},\Q ,\U ,c \) with
\( \B _{1}=\B  \), \( \B _{2}=\B \Q  \), \( \T_{1}=\T  \), \( \T_{2}=\T \U  \)
and a block simulation \( \Phi_{1} \) such that for all \( \eps>0 \), if \( \eps_{1}=\eps \),
\( \eps_{2}=c\eps^{2} \) and \( M_{k} \) is the perturbation
 \[
  \CA_{\eps_{k}}(\trans_{k},\B _{k},\T_{k}, \bZ)
 \]
 then \( \Phi_{1} \) is a simulation of \( M_{2} \) by \( M_{1} \).
The meaning of this is that even if we have to cope with the fault
probability \( \eps \), the simulation will compute \( \trans_{2} \) with a much
smaller fault probability \( c\eps^{2} \).
The hope is not unreasonable since in Section~\ref{sec:ftol-block},
we outlined a single-fault-tolerant block simulation while the
probability of several faults happening during one work period is only
of the order of \( (\Q \U \eps)^{2} \).
  However, it turns out that the only simply stated property of a
perturbation that survives noisy simulation is a certain initial stability
property (see below).

\paragraph{Trickle-down}
  Even if the above goal can be achieved, the reason for the existence of
the simulated more reliable abstract medium is to have feedback from it to
the simulating one.
Let us define the nature of this feedback via the notion of trickle-down.
We want the cell state simulated by a colony
to determine \emph{some properties} (fields) of the simulating cells in the colony;
however, we want to allow some other fields to vary.
For example there could be an \( \Age \) field that varies during any work period,
while the information of the big cell represented by the colony should
be reflected in its cells---via the code---all the time.

\begin{definition}[Trickle-down]\label{def:trickle-down-1}
Let \( \Phi=(\fg,\Phi^{*}) \) be a simulation as in~\eqref{eq:Phi} whose encoding
\( \fg_{*} \) is a block code with block size \( \Q  \), between abstract cellular
media \( M_{k} \) (\( k=1,2 \)).
Let \( \D_{k}>0 \) (\( k=1,2 \)) be some parameters, and \( \F_{1} \) 
a field of the medium \( M_{1} \).
We say that \( \Phi \) has the
\df{\( \eps \)-trickle-down property}\index{error!trickle-down} 
of field \( \F_{k+1} \) into field \( \F_{1} \) with respect
to \( \D_{1}, \D_{2} \) if the following holds for every configuration
\( \xi \) of \( M_{2} \) and every trajectory \( (\mu,\eta) \) of \( M_{1} \) with
\( \eta(\cdot, 0)=\fg_{*}(\xi) \).

Recall the notation \( \preceq \) in Definition~\ref{def:preceq}.
Let \( \eta^{1}=\eta \), \( \eta^{2}=\Phi^{*}(\eta) \).
Let \( x_{1},x_{2} \) be sites where \( x_{1} \)
has address \( a \) in the \( \Q  \)-colony with base \( x_{2} \), let \( t_{0} \) be some
time.
For \( k=1,2 \) let \( \cE_{k}(x) \) be the event that \( \eta^{k}(x,t) \)
is non-vacant and \( \eta^{k}(x,t).\F_{k} \)
does not change during \( \rint{t_{0}-\D_{k}}{t_{0}+\D_{k}} \).
Let \( \cE'_{1} \) be the event that \( \cE_{1}(x_{1}) \) holds and also
 \[
   \eta(x_{1}, t).\F_{1} \succeq \fg_{*}(\eta^{*}(x_{2}, t))(a).\F_{1}
 \]
 during this time.
Then \( \Prob\paren{\cE_{2}(x_{2})\setminus\cE'}<\eps \).
  \end{definition}

 The term ``trickle-down'' is used since information from a higher level cell
is used to determine something in a lower-level cell.
Informally, this means that for all \( x_{1},x_{2},a \) in the given
relation, the state \( \eta(x_{1},t) \) is with large probability what we
expect by encoding some \( \eta^{*}(x_{2}, t') \)
and taking the \( a \)-th symbol of the code-word.

There is yet another, simple error-resistance
property of a medium worth spelling out: namely, the time until the initial
state does not change, except with small probability.

\begin{definition}[Initial stability]
A cellular medium \( M \) 
is called \df{initially stable} for an initial configuration \( \xi \) 
with parameters \( (\eps,\T ) \)
if for each trajectory \( \eta \) of \( M \) with
\( \eta(\cdot,0)=\xi \), for each site \( x \), we have
 \[
  \Pbof{\exist{t < \T }\eta(x,t) \ne \eta(x,0)} < \eps.
 \]
\end{definition}

Initial stability is a trivial property for lattice cellular automata:

 \begin{example}
Consider the perturbation \( M=\CA_{\eps}(\trans,\B ,\T ) \)
of  a 1-dimensional deterministic cellular automaton.
Then for any configuration \( \xi \) and arbitrary \( \delta>0 \)
the medium \( M \) is initially stable for initial configuration \( \xi \) and parameters
\( (\delta,\T ) \).
Indeed, by definition no cell is changing it state before time \( \T  \).
 \end{example}

For variable-period cellular automata, and other media in which the length of
dwell periods is not fixed in advance, the property is less trivial.

\section{Hierarchy}\label{sec:hier}

The reliable cellular automaton we are building is working in a hierarchical
way.
Before understanding a dynamically working hierarchy, it is helpful to
understand thoroughly a static hierarchy, for example
a hierarchically constructed initial configuration.

 \subsection{Hierarchical codes} \label{sec:hier.codes}

In the present work, formally, a hiearchy will be defined via
a ``composite code''.
Let us introduce the basic operation of
the hierarchical structure arising in an amplifier.

\begin{definition}[Composition of codes]
If \( \fg,\psi \) are two codes then \( \fg\circ\psi \) is defined by
\( (\fg\circ\psi)_{*}(\xi)=\fg_{*}(\psi_{*}(\xi)) \) and
\( (\fg\circ\psi)^{*}(\zeta)=\psi^{*}(\fg^{*}(\zeta)) \).
It is assumed that \( \xi \) and \( \zeta \) are here configurations of the
appropriate cellular automata, that is the cell body sizes are in the
corresponding relation.
The code \( \fg\circ\psi \)\glo{"+@\( \circ \)} is called the
\df{composition}\index{code!composition} of \( \fg \) and \( \psi \).
\end{definition}

Here is a detailed example.

 \begin{example}
Let \( M_{1},M_{2},M_{3} \) have cell body sizes \( 1,31,31^{2} \)
respectively.
Let us use the code \( \fg \) from Example~\ref{xmp:block-code}.
The code \( \fg^{2}=\fg\circ\fg \) maps each cell \( c \) of \( M_{3} \) with body
size \( 31^{2} \) into a ``supercolony''\index{supercolony} of \( 31\cdot 31 \)
cells of body size 1 in \( M_{1} \).

Suppose that \( \zeta=\fg_{*}^{2}(\xi) \) is a configuration obtained by
encoding from a lattice configuration of body size \( 31^{2} \) in \( M_{3} \),
where the bases of the cells are at positions \( -480 + 31^{2} i \).
(We chose -480 only since \( 480=(31^{2}-1)/2 \) but we could have chosen any
other number.)
Then \( \zeta \) can be grouped into colonies of size \( 31 \) starting at
any of those bases.
Cell 55 of \( M_{1} \) belongs to the colony with base \( 47=-480+17\cdot 31 \)
and has address 8 in it.
Therefore the address field of \( \zeta(55) \) contains a binary
representation of 8.
The last bit of this cell encodes the 8-th bit the of cell (with base)
47 of \( M_{2} \) represented by this colony.
If we read together all 12 bits represented by the \( \Info \) fields of
the first 12 cells in this colony we get a state \( \zeta^{*}(47) \) (we count
from 0).
The cells with base \( -15+31j \) for \( j\in\bbZ \) with states
\( \zeta^{*}(-15+31j) \) obtained this way are also broken up into colonies.
In them, the first 5 bits of each state form the address and the
last bits of the first 12 cells, when put together, give back the
state of the cell represented by this colony.
Notice that these 12 bits were really drawn from \( 31^{2} \) cells of
\( M_{1} \).
Even the address bits in \( \zeta^{*}(47) \) come from different cells of
the colony with base 47.
Therefore the cell with state \( \zeta(55) \) does not contain information
allowing us to conclude that it is cell 55.
It only ``knows'' that it is the 8-th cell within its own colony
(with base 47) but does not know that its colony has address 17 within
its supercolony (with base \( -15\cdot 31 \)) since it has at most one
bit of that address.
 \end{example}

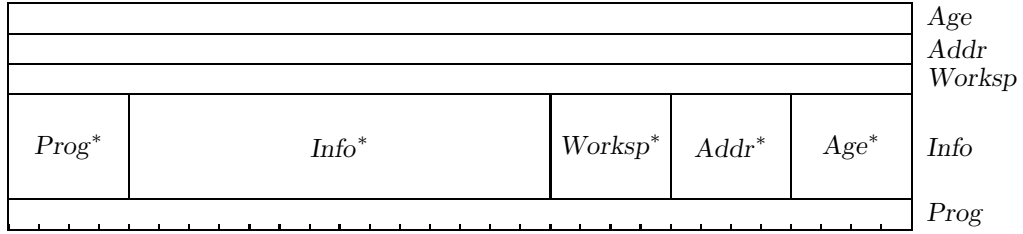
\begin{figure}
 \setlength{\unitlength}{0.2mm}
 \[
 \begin{picture}(600,150)
    \put(0, 0){\framebox(600, 150){}}

    \put(610, 0){\makebox(0,20)[l]{\( \var{Prog} \)}}

    \put(0,20){\line(1,0){600}}
    \put(610,20){\makebox(0,70)[l]{\( \var{Info} \)}}

    \put(0, 20){\makebox(80,70){\( \var{Prog}^{*} \)}}
    \put(80, 20){\line(0, 1){70}}
    \put(80, 20){\makebox(280,70){\( \var{Info}^{*} \)}}
    \put(360, 20){\line(0, 1){70}}
    \put(360, 20){\makebox(80,70){\( \var{Worksp}^{*} \)}}
    \put(440, 20){\line(0, 1){70}}
    \put(440, 20){\makebox(80,70){\( \var{Addr}^{*} \)}}
    \put(520, 20){\line(0, 1){70}}
    \put(520, 20){\makebox(80,70){\( \var{Age}^{*} \)}}

    \put(0,90){\line(1,0){600}}
    \put(610,90){\makebox(0,20)[l]{\( \var{Worksp} \)}}

    \put(0, 110){\line(1,0){600}}
    \put(610,110){\makebox(0,20)[l]{\( \var{Addr} \)}}

    \put(0, 130){\line(1,0){600}}
    \put(610,130){\makebox(0,20)[l]{\( \var{Age} \)}}

   \multiput(0,0)(20,0){30}{\line(0,1){3}}
 \end{picture}
 \]
 \caption{Fields of a cell simulated by a colony}
 \label{fig:fields}
\end{figure}

 \paragraph{Infinite composition}
  A code can form composition an arbitrary number of times with itself
or other codes.
In this way, a hierarchical, that is highly nonhomogenous, structure can
be defined using cells that have only a small number of states.

\begin{definition}[Hierarchical code]\label{def:xi^k}
A \df{hierarchical code}\index{code!hierarchical} is given by a sequence
 \glo{f.greek.hier@\( \fg_{k*} \)}%
\( \states_{k},\fg_{k*} \), (\( k\ge 1 \))
where \( \states_{k} \) is an alphabet and \( \fg_{k*} \) is a code on configurations.
In the present work we will also assume that \( \fg_{k*} \) always comes from a block
code 
 \begin{align*}
  \fg_{k*}:\states_{k+1}\to\states_{k}^{\Q_{k}}.
 \end{align*}
 \glo{q.cap@\( \Q_{k} \)}%
Since \( \states_{k} \) and \( \Q_{k} \) are implicitly defined by \( \fg_{k*} \) we can
refer to the code as just \( (\fg_{k*})_{k\ge 1} \).
  Consider a hierarchical code  \( (\fg_{k*})_{k\ge 1} \) with state sequence \( (\states_{k}) \)
  and block sizes \( (\Q_{k}) \).
  If for each \( \states_{k} \) a field \( \F_{k} \) is also given, and
  the code \( \fg_{k*} \) is aggregating from \( \F_{k} \) to \( \F_{k+1} \) as in
  Definition~\ref{def:aggreg-field}, then the hierarchical code \( (\fg_{k*})_{k\ge 1} \)
  will be called \df{aggregating} for the field sequence \( (\F_{k}) \).
  If also \( \F_{k} \) controls address 1 via a function \( \gamma_{k} \)
  (as by Definition~\ref{def:control}) then sequence of triples
  \begin{align}\label{eq:H}
   \hierCode = (\fg,\F_{k}, \gamma_{k})_{k\ge 1}
  \end{align}
  will be called a \df{centrally aggregating hierarchical code}, or shortly, a \df{code complex}.
\end{definition}

\begin{definition}
  Given a hierarchical code as above, for a configuration \( \xi \) we define the 
higher-order configurations \( \xi^{k} \) by the recursion
\begin{equation}\label{eq:xi^k}
  \begin{aligned}
       \xi^{1}    &=\xi,
\\   \xi^{k+1} &=\fg_{k}^{*}(\xi^{k}).      
  \end{aligned}
\end{equation}
The configuration \( \xi \) is called a \df{code configuration} for the given
hierarchical code if \( \xi^{k} \) is a code configuration for the code \( \fg_{k} \)
for each \( k \) in the sense of Definition~\ref{def:code-config}.
\end{definition}

\paragraph{Constructing code configurations}
The construction of a code configuration for a hierarchical code seems to
require a composition \( \fg_{1*}\circ\fg_{2*}\circ\cdots \).
What is the meaning of this?
We will want to compose the codes ``backwards'', that is in such a
way that from a configuration \( \xi^{1} \) of some
medium \( M_{1} \) with cell body size 1,
we can decode the configuration \( \xi^{2}=\fg_{1}^{*}(\xi^{1}) \) of \( M_{2} \) with
cell body size \( \B _{2}=\Q_{1} \), configuration \( \xi^3=\fg_{2}^{*}(\xi^{2}) \), of
\( M_{3} \) with body size \( \B _{3}=\Q_{1}\Q_{2} \), and so on.
Such constructions are not unknown, they were used for example  to define
``Morse sequences'' with applications in group theory as well in the theory
of quasicrystals (\cite{ItkisLevin94,Radin91}).
All known nonperiodic tiling sets are also based on a similar hierarchical
construction. 

 \begin{definition}\label{def:coord-hier}
Suppose that a code configuration \( \xi \) of a hierarchical code
\( (\fg_{k})_{k\ge 1} \) with block sizes \( \Q_{k} \) is given.
The sizes of the higher-order cells are defined by \glo{b.cap@\( B_{k} \)}%
\(  \B _{1}= 1 \), \( \B _{k+1} = Q_{k}\B_{k} \).
For a finite or infinite \( n \), let \glo{k.cap@\( K_{0}(n) \)}
\begin{align}\label{eq:K0}
 K_{0}(n) = \sup_{\B _{k}\le n} k.
\end{align}
If \( n<\infty \)  then this is the largest \( k \) for which a \( k \)-cell can fit into \( \bbZ_{n} \).
As long as \( k< K_{0}(n) \)
the origin of our space \( \sites \) will be contained in some cell \( x_{k} \) of level \( k \).
Since it is in our power to define the code configuration \( \xi \), for simplicity
we will always arrange that \( x_{k} \) is the second cell (the cell with address 1)
of the colony coding cell \( x_{k+1} \).
In general, we could use any addresses \( 0<a_{k}<\Q-1 \) instead of setting them always to 1,
but this is not important now.
Setting \( a_{k}=1 \) gives
\begin{align*}
  x_{1}=0, \quad  x_{k}  = -\B _{1} - \cdots - \B _{k-1} \text{ for } k>1.
\end{align*}
\end{definition}

In a hierarchical code as in Definition~\ref{def:coord-hier},
the sequence of states 
\begin{align*}
 s_{k}=\xi^{k}(x_{k})  
\end{align*}
obeys some restrictions:

 \begin{definition}[Fitted sequence]
A finite or infinite sequence
\( (s_{1},s_{2},\dots) \) will be called
\df{fitted}\index{fitted sequence} to a code complex
\( \hierCode \) as in~\eqref{eq:H} if
 \begin{align}\label{eq:fitted}
  \fg_{k*}(s_{k+1})(1)=s_{k}
 \end{align}
 for all \( k \).
 \end{definition}

Every fitted sequence allows the construction of a code configuration for the
hierarchical code.
(If we used addresses \( a_{k} \) then the condition above would be
\(   \fg_{k*}(s_{k+1})(a_{k})=s_{k} \), and the proposition below would still hold.)

  \begin{proposition}\label{propo:fitted}
    For a fitted sequence \( (s_{k})_{k\ge 1} \)
there are configurations \( \xi^{k} \) of \( M_{k} \) over \( \bbZ \) such that
for all \( k\ge 1 \) we have
 \begin{align*}
   \fg_{k*}(\xi^{k+1})	&=  \xi^{k},
\\ \xi^{k}(x_{k})  	        &= s_{k}.
 \end{align*}
 \end{proposition}
 \begin{proof}
   Let \( N \) (finite or infinite) be the size of our space,  \glo{x.greek@\( \xi^{i}_{k} \)}
   and let
 \begin{equation}\label{eq:xi^{k}_{k}}
  \xi^{k}_{k}
 \end{equation}
 be the configuration of \( M_{k} \) which has state \( s_{k} \) at site \( x_{k} \)
and arbitrary states at all other sites \( x_{k}+x\B_{k} \), with the following
restriction in case of a finite space size \( N \).
Let \( \xi^{k}_{k} \) be non-vacant only for \( k\le K_{0} = K_{0}(N) \) and
\( \xi^{K_{0} }_{K_{0} }(x_{K_{0} }+z) \) non-vacant only if
 \begin{equation}\label{eq:xi^{K}.defined}
  0 \le z\B _{K_{0} } < N.
 \end{equation}
 For \( 1\le i < k\) let 
 \begin{equation}\label{eq:xi-i-k}
 \xi^{i}_{k} = \fg_{i*}(\fg_{(i+1)*}(\cdots \fg_{(k-1)*}(\xi^{k}_{k})\cdots)).
 \end{equation}
  We have
 \begin{equation}\label{eq:fitted.ext}
  \xi^{k}_{k+1}(x_{k})=\fg_{k*}(\xi^{k+1}_{k+1}(x_{k+1}))(1)
 = \xi^{k}_{k}(x_{k})
 \end{equation}
  where the first equation comes by definition, the second one by
fittedness.
  The encoding conserves this relation, so the partial configuration
 \( \xi^{i}_{k+1}(x_{k+1}+\lint{0}{\B _{k+1}}) \) is an
extension of \( \xi^{i}_{k}(x_{k}+\lint{0}{\B _{k}}) \).
  Therefore the limit \( \xi^{i}=\lim_{k} \xi^{i}_{k} \) exists for each \( i \).
The limit extends over the whole set of integer sites.
\end{proof}

  Though the code configuration
\( \xi^{1} \) above is obtained by an infinite process of encoding,
no infinite process of decoding is needed to yield a single
configuration from it: at the \( k \)-th stage of the decoding, we get a
configuration \( \xi^{k} \) with body size \( \B _{k} \).

Here is how finite and infinite sequences will be encoded.
For the following two definitions,
let \( \hierCode \) be a code complex as in~\eqref{eq:H}.

\begin{definition}[Encoding of infinite sequences]\label{def:encoding-infinite}
  Let \( \rho\in (\states_{1}.\F_{1})^{\bbZ} \) be an infinite sequence,
  we define its encoding \( \psi_{*}(\rho) \) as follows.
  Define for each \( k \) the following segment of \( \rho \)
using the notation of Definition~\ref{def:coord-hier}:
\begin{align*}
   \rho_{k}&= (\rho(x_{k}),\rho(x_{k}+1),\dots,\dots,\rho(x_{k}+\B_{k}-1)).
\end{align*}
The sequence \( s_{k}=\gamma_{k}(\rho_{k}) \), \( k=1,2,\dots \) is fitted,
allowing via Proposition~\ref{propo:fitted}
to encode \( \rho \) into an infinite code configuration
\begin{align*}
 \xi=\xi^{1}=\psi_{*}(\rho) .
\end{align*}
\end{definition}

This construction leaves plenty of freedom for the codes \( \fg_{k} \).
Indeed, the tracks of a codeword \( w=\fg_{k*}(s_{k+1}) \) 
other than \( \F_{k} \) can store all kinds of
structural and error-correcting information about \( s_{k+1} \),
and it is only its symbol \( w(1) \)
that is restricted to \( \gamma_{k}(w(1).\F_{1})=\gamma_{k}(\rho_{k})=s_{k} \).

A sequence \( \rho \) of length \( n<\infty \) will be encoded into a cell of an appropriate level.
Recall the notation \( x_{k} \) from Definition~\ref{def:coord-hier}.

\begin{definition}[Encoding of finite sequences]\label{def:encoding-finite}
  Given a code complex \( \hierCode \) as above and
 a string \( \rho\in\Sigma_{0}^{n} \) we define its encoding
 \( \xi=\psi_{*}(\rho)\in\states^{\infty} \) which has the property that
 except for an interval of size \( \le 2 n \), \( \xi(x)=q_{\states} \) with
 \( q_{\states} \) as in Definition~\ref{def:Init}. 
  Let \( k=K_{0}(n) \) as in~\eqref{eq:K0}.
  Create the string 
\begin{align*}
 \rho'=*^{-x_{k}}\rho *^{l}
\end{align*}
where \( -x _{k} \) stars come at the front and the \( l \) stars at the back are the fewest
to make \( |\rho'|=m\B_{k} \) an integer multiple of \( \B_{k} \).
Now create a string \( r_{k}\in \states_{k}^{m} \) as follows.
First let \( r'_{k}=q_{\states}^{m} \), and then obtain \( r_{k} \) by replacing the track
\( r'_{k}.\F_{k} \) in it with \( \rho' \).
Define the string \( r_{k-1}=\fg_{i*}(r_{k})\in\states_{k-1}^{m\Q_{k-1}} \).
Similarly, for each \( 1\le i<k \) let \( r_{i}=\fg_{(i+1)*}(r_{i+1}) \).
By definition of the codes \( \fg_{i*} \),
string \( r_{i} \) consists of \( m\Q_{k-1}\Q_{k-2}\cdots\Q_{i} \) symbols of \( \states_{i} \).
Then \( r_{1} \) consists of \( m\B_{k} \) symbols of \( \states=\states_{1} \).
Finally obtain \( \xi \) as follows: \( \xi(x)=r(x-x_{k}) \) for \( 0\le x-x_{k}<m\B_{k} \) and \( q_{\states} \)
otherwise.
\end{definition}
Now the symbols of \( \rho \) occupy the places \( \xi(x).\F_{1} \) for \( x\in\lint{0}{n} \),
within the string \( r_{1} \) representing a number \( m\le \Q_{k}+1 \) of \( k \)-cells.
Outside these \( k \)-cells the configuration consists of 1-cells with one and the same state \( q_{\states} \).

\subsection{Amplifiers}

An amplifier is a hierarchical code whose member codes are also simulations.
Recall general simulations in Definition~\ref{def:general-simulation}.

\begin{definition}[Amplifier]\label{def:amp}
  Suppose that a sequence \( M_{1},M_{2},\dots \) of abstract media is given
along with simulations \( \Phi_{1},\Phi_{2},\dots \) such that \( \Phi_{k} \) is a
simulation of \( M_{k+1} \) by \( M_{k} \).
Such a system will be called an \df{amplifier}\index{amplifier}.
\end{definition}

Amplifiers are somewhat analogous to renormalization 
groups\index{renormalization} in statistical physics.
So far, we have not seen in the present
work any nontrivial example of simulation other
than between deterministic cellular automata, so the idea of an amplifier
seems far-fetched at this moment.

We want the simulations in our amplifiers to have all the good properties
introduced earlier:

\begin{definition}[Trickle-down and initial stability for amplifiers]\label{def:amp-trickle}
  Suppose that an amplifier \( \Amp=(M_{k},\Phi_{k})_{k\ge 1} \) is given along
  with the sequences of parameters \( \D_{k} \),
  \( \eps''_{k} \), \( \eps_{k} \), and fields \( \F_{k} \).
  Assume that each \( \fg_{k} \) aggregates field \( \F_{k} \) into \( \F_{k+1} \),
  as in Definition~\ref{def:aggreg-field}.
 \begin{enumerate}

 \item \( \Amp \) has the 
 \df{trickle-down}\index{amplifier!trickle-down} property if for
 each \( k \), the simulation \( \Phi_{k} \) has the trickle-down
property of \( \F_{k+1} \) to \( \F_{k} \) with 
parameters  \( \D_{k} \), \( \D_{k+1} \) and \( \eps''_{k} \),
as in Definition~\ref{def:trickle-down-1}.

 \item Suppose that the sequence \( \eps_{k}>0 \), \( k=1,2,\dots \) is given 
with \( \sum_{k}\eps_{k}<1/6 \), and also a configuration \( \xi \) of \( M_{1} \)
that is a code configuration of the code sequence \( (\fg_{k}) \)
(as introduced in Definition~\ref{def:xi^k}).
We say that \( \Amp \) has the \df{initial stability property} with respect
to the initial configuration \( \xi \) and parameters \( ((\eps_{k},\T_{k}))_{k\ge 1} \)
if for each \( k \) the medium \( M_{k} \) has the initial stability
property for initial configuration \( \xi^{k} \) and parameters \( (\eps_{k},2\D_{k}) \).
 \end{enumerate}
  \end{definition}

The initial stability property for amplifiers is not as trivial as for lattice
media.
Consider the case when \( M_{1}=\CA_{\eps}(\trans,\bbZ,\B ,\T ) \).
As the trajectories of medium \( M_{2} \) are obtained
by decoding from trajectories of \( M_{1} \), their initial stability
is not apriori guaranteed.

\subsection{Information storage: proof from an amplifier assumption}
\label{sec:trickle-down}

The following lemma will be proved later in the paper.

 \begin{lemma}[Initially stable trickle-down amplifier]\label{lem:abstr-amp}
 \index{amplifier!initially stable}
 \index{lemma@Lemma!Initially Stable Amplifier}
We can construct the following objects.
 \begin{cjenum}
 \item
   Parameters \( \eps_{k}<\eps''_{k} \), \( \B _{1} \), \( \Q_{k} \), \( \D_{k} \)
   with \( \D_{k+1} \ge \D_{k} \) for \( k\ge 1 \), with \( \sum_{k}\eps''_{k}<1/6 \).
\item Media \( M_{k} \) with local state spaces \( \states_{k} \), with
  a distinguished field \( \F_{k} \).
 \begin{align}\label{eq:amp.M1}
   M_{1} = \CA_{\eps_{1}}(\trans_{1}, \B _{1}, \T_{1}, \bZ),
 \end{align}  
 and an amplifier \( \Amp=((M_{k}),(\Phi_{k})) \)
with simulations \( \Phi_{k}=(\fg_{k},\Phi_{k}^{*}) \)
where the codes \( \fg_{k} \) are block codes with blocksize \( \Q_{k} \) that
are aggregating field \( \F_{k} \) into field \( \F_{k+1} \) as per Definition~\ref{def:aggreg-field},
and where \( \F_{k} \) controls address \( 1 \) as in Definition~\ref{def:control}.

 \item For each infinite sequence \( \rho\in\Sigma_{0}^{*} \)
   an initial configuration \( \xi^{1}=\psi_{*}(\rho) \) of \( M_{1} \)
   as obtained in Proposition~\ref{propo:fitted} such that
   for each \( k \) the sequence \( \xi^{k} \) covers the space with non-vacant cells of \( M_{k} \).
   
\item \( \Amp \) is trickle-down and initially stable for \( \xi \), with the
  parameter sequences \( (\F_{k}) \), \( (\D_{k}) \) , \( (\eps''_{k}) \),
  \( (\eps_{k}) \).
 \end{cjenum}
 \end{lemma}

 Note that~\eqref{eq:amp.M1} is required only for \( M_{1} \); the media \( M_{k} \)
 for \( k>1 \) will be more complex objects, not simply cellular automata.
Let us use this lemma to prove Theorem~\ref{thm:remember-inf}.

  \begin{proof}[Proof of Theorem~\protect\ref{thm:remember-inf}]
    We will use the amplifier defined in the above lemma.
    Let \( \rho\in\Sigma_{0}^{\bbZ} \), and \( \xi \) as in Lemma~\ref{lem:abstr-amp}.
Let \( \eta^{1} \) be a trajectory of the medium \( M_{1} \) with \( \eta^{1}(\cdot,0)=\xi \),
and \( \eta^{k} \) be defined by the recursion \( \eta^{k+1}=\Phi^{*}_{k}(\eta^{k}) \).
Let \( (x_{1},t_{0}) \) be a space-time point.
There is a sequence of points \( x_{1},x_{2},\dots \) such that \( x_{k+1} \)
is a cell of \( \eta^{k+1}(\cdot,0) \) containing \( x_{k} \) in its body.
  There is a first \( n \) with \( t_{1}<\D_{n} \).
  Let \( \cE_{k} \) be the event that \( \eta^{k}(x_{k},t) \) is non-vacant and
 \begin{align*}
 \eta^{k}(x_{k},t).\F_{k}=\xi^{k}(x_{k}).\F_{k}  
\end{align*}
during \( \rint{t_{0}-\D_{k}}{t_{0}+\D_{k}} \).
The theorem follows from the bounds on \( \sum_{k}\eps_{k} \) and
\( \sum_{k}\eps''_{k} \) and from
 \begin{equation}\label{eq:cup-not-F.k}
 \Prob\bigparen{\neg(\cE_{1}\cap\cdots\cap\cE_{n})} \le \eps_{n} +
  \sum_{k=1}^{n-1}\eps''_{k}.
 \end{equation}
  To prove this inequality, use
 \[
   \neg(\cE_{1}\cap\cdots\cap\cE_{n}) = \neg\cE_{n}\cup
 \bigcup_{k=1}^{n-1}(\cE_{k+1}\setminus\cE_{k}).
 \]
 By the construction, \( \eta^{n}(x_{n},0).\F_{n}=\xi^{n}(x_{n}) \).
The initial stability property implies \( \Prob\bigparen{\neg\cE_{n}}\le\eps_{n} \).
Let us show
\( \Prob\bigparen{\cE_{k+1}\setminus\cE_{k}}\le\eps''_{k} \).
The assumption \( \cE_{k+1} \) says
\begin{align}\label{eq:E_{k+1}}
 \eta^{k+1}(x_{k+1},t).\F_{k+1}=\xi^{k+1}(x_{k+1}).\F_{k+1}.
\end{align}
Let \( a \) be the address of cell \( x_{k} \) of \( M_{k} \) in the colony
simulating cell \( x_{k+1} \) of \( M_{k+1} \).
By the definition of 
of \( \xi^{k} \) we have \( \fg_{k*}(\xi^{k+1}(x_{k+1}))(a)=\xi^{k}(x_{k}) \).
The assumption that \( \fg_{k} \) aggregates \( \F_{k} \) into \( \F_{k+1} \)
and \eqref{eq:E_{k+1}} implies, as shown in~\eqref{eq:aggreg-field-prop}, 
\( \fg_{k*}(\eta^{k+1}(x_{k+1}, t))(a).\F_{k}=\xi^{k}(x_{k}).\F_{k} \).
 The trickle-down property says that except with probability \( \eps''_{k} \),
we have \( \fg_{k*}(\eta^{k+1}(x_{k+1},t))(a).\F_{k}=\eta^{k}(x_{k},t).\F_{k} \).
during \( \rint{t_{0}-\D_{k}}{t_{0}+\D_{k}} \), which completes the proof.
  \end{proof}

  \subsection{Error-correcting codes}

  Let us define the parts of our block code that deal with error-correction.
  
\begin{sloppypar}
\begin{definition}[Error-correcting code]\label{def:err-code}
  A block code \( \fg=(\fg_{*}, \fg^{*}) \) with \( \fg^{*}: R^{\Q }\to S \) will be called
  is \(  t \)-\df{error-correcting} \index{code!error-correcting}%
  if for all \( x\in\Sigma \), \( y\in\Sigma^{\\Q } \) we
have \( \fg^{*}(y)=x \) whenever \( y \) differs from
\( \fg_{*}(x) \) in at most \( t \) symbols.
It is \( (\beta,t) \)-\df{burst-error-correcting}
if for all \( x\in S \), \( y\in R^{\\Q } \) we
have \( \fg^{*}(y)=x \) whenever \( y \) differs from
\( \fg_{*}(x) \) in at most \( t \) intervals of size \( \le\beta \).
For such a code, we will say that a word \( y\in R^{\\Q } \) is \( (\beta, r) \)-\df{compliant}
if it differs from a codeword of the code by at most \( r \) intervals of size \( \le\beta \).

Assume that strings \( R \) and \( S \) are in \( \{0,1\}^{*} \).
In the context of error-correcting codes, elements of \( R \) will be called the
\df{symbols} of the codewords, as opposed to \df{bits}, since each symbol may
consist of several bits.
\end{definition}
  \end{sloppypar}

\begin{example}\label{xmp:err-corr-code}
A popular kind of error-correcting code are codes \( \psi \) such that
\( \psi_{*} \) is a linear mapping when the binary strings in \( S \) and \( R^{\Q } \)
are considered vectors over the field \( \{0,1\} \).
These codes are called \df{linear codes}\index{code!linear}.
It is sufficient to consider linear codes \( \psi \) that are \df{systematic}:
for all \( s \), the first \( \norm{S} \) bits of the codeword \( \psi_{*}(s) \) are
identical to \( s \): they are called the
\df{information bits}\index{code!bits!information}.
(If a linear code is not such, it can always be rearranged to have this
property.)
In this case, the remaining bits of the codeword are called
\df{error-check bits}\index{code!bits!error check}, and they are linear
functions of \( s \).
For correcting a single error, in the tripling method outlined
in Section~\ref{sec:ftol-block}, the error check bits are simply the twice repeated
information bits.
 \end{example}

 \begin{example}\label{xmp:RS-code}\index{code!Reed-Solomon}
Here we define the linear code we will be using in later construction (a
generalization of the so-called Reed-Solomon code, see~\cite{Blahut83}).
In its application in our paper only its basic quantitative
properties are used, so the details can safely be skipped at first reading.

Let our codewords (information symbols and check symbols together) be
binary strings of length \( N l \)\glo{n:cap@\( N \)}\glo{l@\( l \)} for some \( l \),
\( N \).
The \df{symbols} of our codewords, 
binary strings of length \( l \), will be interpreted as elements of the
finite field \( \GF(2^{l}) \)\index{Galois field}\glo{gf@\( \GF(2^{l}) \)} and
thus, each binary string \( c \) of length \( N l \) will be treated as a vector
\( \tup{c(0),\dots,c(N - 1)} \) over \( \GF(2^{l}) \).
(Note that the word ``field'' is used in two different senses in the
present paper.)
Let us fix \( N \) distinct nonzero elements
\( \alpha_{i} \)\glo{a.greek@\( \alpha_{i} \)} of \( \GF(2^{l}) \) and let \( t < N/2 \) be a
positive integer.
The codewords are those vectors \( c \) that satisfy the equations
 \begin{equation}\label{eq:RS-code}
  \sum_{i = 0}^{N - 1}\alpha_{i}^{j}c(i).\F^{k} = 0\ (j = 1,\dots, 2t)
 \end{equation}
where the addition, multiplication and taking to power \( j \) are
performed in the field \( \GF(2^{l}) \).
These are \( 2t \) linear equations.
If we fix the first \( N - 2t \) elements of the vector in any way, (these
are the \df{information symbols}) the remaining \( 2t \) elements (the \df{error check
symbols}) can be chosen in a way to satisfy the equations.
This set of equations is always solvable, since its determinant is a
Vandermonde determinant.

Below, we will show a procedure for correcting any \( \nu \le t \) nonzero
errors.
This demonstrates that for the correction of error in any \( \le t \) symbols,
only \( 2t \) error-check symbols are needed.
(Robert Solovay noticed that in our applications,
with a little deterioration of efficiency, even the most brute-force
decoding algorithm is sufficient.
We will never have to correct more than some constant number \( c \) of errors,
so the complexity of decoding would never go much over \( n^{c} \) steps.)

If \( E =\tup{e_{0},\dots, e_{N - 1}} \) is the sequence of errors then \( C + E \)
will be the observed the word.
Only \( e_{i_{r }} \) are nonzero for \( r =1,\dots,\nu \).
Let \( Y_{r }=e_{i_{r }}, X_{r }=\alpha_{i_{r }} \).
We define the \df{syndrome} \( S_{j} \) for \( j=1,\dots,2t \) by
 \begin{equation}\label{eq:S_{j}}
	S_{j} = \sum_{i} (c_{i}+e_{i})\alpha_{i}^{j}
		= \sum_{i} e_{i}\alpha_{i}^{j} 	
	   	= \sum_{r } Y_{r } X_{r }^{j}
 \end{equation}
  which can clearly be computed from the codeword: it is the amount by
which the codeword violates the \( j \)-th error check equation.
We will show, using the last expression, that \( Y_{r } \) and \( X_{r } \) can be
determined using \( S_{j} \).
Define the auxiliary polynomial
 \[
	\Lambda(x)=\prod_{r }(1-xX_{r })=\sum_{s=0}^{\nu} \Lambda_{s} x^{s}
 \]
 whose roots are \( X_{r }^{-1} \).
Let us show how to find the coefficients \( \Lambda_{s} \) for \( s>0 \).
We have, for any \( r = 1,\dots, \nu \), and any \( j = 1, \dots, 2t - \nu \):
 \[
	0	= Y_{r }X_{r }^{j+\nu}\Lambda(X_{r }^{-1})		
		= \sum_{s} \Lambda_{s} Y_{r } X_{r }^{j+\nu-s}.     
 \]
  Hence, summing for \( r  \), 
 \begin{equation}\label{eq:key}
	0 = \sum_{s = 0}^{\nu} \Lambda_{s}(\sum_{r } Y_{r } X_{r }^{j+\nu-s})
	= \sum_{s = 0}^{\nu} \Lambda_{s} S_{j+\nu-s}\ (j = 1,\dots, 2t - \nu)
 \end{equation}
  hence using \( \Lambda_{0}=1 \),
 \( \sum_{s=1}^{\nu} \Lambda_{s} S_{j+\nu-s}=-S_{j+\nu} \).
This is a system of linear equations for \( \Lambda_{s} \) whose coefficients are
the syndroms, known to us, and whose matrix \( M_{\nu} \) is nonsingular.
Indeed, \( M_{\nu} = ABA^{T} \) where \( B \) is the diagonal matrix
\( \ang{Y_{r }X_{r }} \) and \( A \) is the Vandermonde matrix \( A_{j,r }=X^{j-1}_{r } \).

A decoding algorithm now works as follows.
For \( \nu = 1,2,\dots,t \), see if \( M_{\nu} \) is nonsingular, then compute
\( \Lambda(x) \) and find its roots by simple trial-and-error, computing
\( \Lambda(\alpha_{i}^{-1}) \) for all \( i \).
Then, find \( Y_{r } \) by solving the set of equations~\eqref{eq:S_{j}} and
see if the resulting corrected vector \( C \) satisfies~\eqref{eq:RS-code}.
If yes, stop.
(There is also a faster way for determining \( \Lambda(x) \), via the Euclidean
algorithm, see~\cite{Blahut83}).

To make the code completely constructive we must find an effective
representation of the field operations of \( \GF(2^{l}) \).
This finite field can be efficiently represented as the set of remainders
with respect to an irreducible polynomial of degree \( l \) over \( \GF(2) \), so
what is needed is a way to generate large irreducible polynomials.
  Now, it follows from the theory of fields that
 \[
	x^{2\cdot 3^{s}} + x^{3^{s}} +1 
 \]
  is irreducible over \( \GF<(2) \) for any \( s \).
So, the scheme works for all \( l \) of the form \( 2\cdot 3^{s} \).
 \end{example}

 A hierarchical code like in Proposition~\ref{propo:fitted} leaves some space for error-correction,
 but it needs to be economical, so we indeed need a code like Example~\ref{xmp:RS-code}.
 
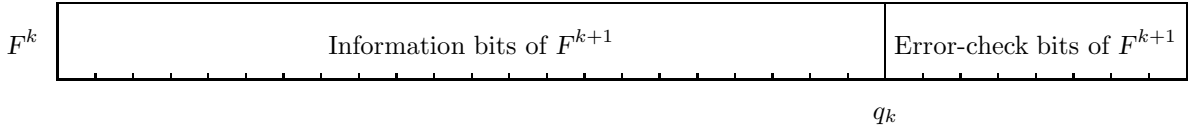
\begin{figure}[H]
 \setlength{\unitlength}{0.21mm}
 \[
   \begin{picture}(600,195)
     {\linethickness{1pt}
    \put(0,0){\framebox(600,195){}}
  }
    \put(0,85){\framebox(600,110){\( \F^{k+1} \)}}
    \put(-10,85){\makebox(0,80)[r]{\( \F^{k} \)}}

    \put(0,60){\framebox(600,25){Error-check symbols for the track \( \F^{k} \)}}
    \put(-10,60){\makebox(0,25)[r]{\( \Redun^{k} \)}}
    \put(0,25){\framebox(150,35){\( \fL^{k+1} \)}}
    \put(150,25){\framebox(100,35){\(\Redun^{k+1}\)}}
    \put(250,25){\framebox(100,35){\(\Work^{k+1}\)}}
    \put(350,25){\framebox(250,35){Error check for this track}}
    \put(-10,25){\makebox(0,35)[r]{\( \fL^{k} \)}}
    \put(0,0){\framebox(600,25){Work tracks}}
    \put(-10,0){\makebox(0,25)[r]{\( \Work^{k} \)}}
    \multiput(0,0)(20,0){30}{\line(0,1){3}}
\end{picture}
 \]
 \caption{Error-correcting code in a shared field}
 \label{fig:no-stretch}
\end{figure}

 \begin{example}\label{xmp:no-stretch}
Let us use the notation of Proposition~\ref{propo:fitted} with \( a_{k}=1 \) for all \( k \).
Our goal is to give more details of state spaces \( \states_{k} \) and
codes \( \fg_{k} \) satisfying the conditions of that Proposition.
 See Figure~\ref{fig:no-stretch}.
 The error check symbols for \( \F_{k} \) are on a track \( \Redun^{k} \).
 As the address \( a_{k}=1 \) is controlled by \( \F_{k} \),  the cell with address 1
 is not used for error checks.
 The other parts of the information needed to represent the big cell are
 on track \( \fL^{k} \).
It contains in its different segments the fields \( \fL^{k+1} \), \( \Redun^{k+1} \), \( \Work^{k+1} \)
of the big cell represented by the colony, as well as its own error check symbols.
There is also a collection of tracks shown here as a single track
\( \Work^{k} \): it contains for example \( \Addr^{k} \), \( \Mail^{k} \).
 \end{example}

In any implementation of the above example all fields other than the aggregated field \( \F^{k} \)
must be relatively quite narrow.
Recall that \( \cp_{k}=\norm{\states_{k}} \) is the number of bits in a state of \( M_{k} \).

\begin{proposition}\label{propo:small-bandwidth}
  In a hierarchical code \( (\fg_{k}) \) 
  let \( \fH_{k} \) denote any track of \( M_{k} \) that is not needed to restore information
  in a cell of \( M_{k+1} \) from the colony coding it.
  (It may contain error-checks or other administrative information like addresses.)
  Then \( \sum_{k} |\fH_{k}|/\cp_{k} < 1 \).
 \end{proposition}

 \begin{sloppypar}
   \begin{proof}
     Let \( h_{k}=|\fH_{k}|/\cp_{k} \), \( r_{k}=1-h_{k} \).
  The state of a cell of \( M_{k+1} \) must be represented by the
colony even when we exclude the track \( \fH_{k} \), so
 \begin{align*}
   \cp_{k+1}(r_{k+1}+ h_{k+1}) &\le  \Q_{k}\cp_{k} r_{k} = \cp_{k+1}r_{k},
\\    r_{k+1} + h_{k+1} &\le  r_{k},
\\    h_{k+1} &\le  r_{k} - r_{k+1},
 \end{align*}
 allowing to estimate \( \sum_{k}h_{k} \) by a telescoping sum.
\end{proof}
    \end{sloppypar}

\subsection{Major difficulties}\label{sec:trouble}

The idea of a simulation between two perturbed cellular automata is,
unfortunately, flawed in the original form: the mapping defined in the
naive way is not a simulation in the strict sense we need.
The problem is that a group of failures can destroy not only the
information but also the organization into colonies in the area where
it occurs.
This kind of event cannot therefore be mapped by the simulation into
a transient fault unless destroyed colonies ``magically recover''.
The recovery is not trivial since ``destruction'' can also mean
replacement with something else that looks locally like its
healthy neighbor colonies but is incompatible with them.
One is reminded of the biological phenomena of viruses and cancer.
Rather than give up hope let us examine the different kinds of
disruption that the faults can cause in a block simulation 
by a perturbed cellular automaton \( M_{1} \).

Let us take as our model the informally described automaton of
Section~\ref{sec:ftol-block}.
  The information in the current state of a colony can be divided into
the following parts:
 \begin{description}
  \item[Information:] an example is the content of the \( \Info \)
track.
  \item[Structure:]  the \( \Addr \) and \( \Age \) tracks.
  \item[Program:]  the \( \fld{Prog} \) track.
 \end{description}
  Less formally, ``structure'' does not represent any data for
decoding but is needed for coordinating cooperation of the colony
members.
The ``program'' determines which transition function will be
simulated.
The ``information'' determines the state of the simulated cell: it
is the ``stuff'' that the colony processes.
Disruptions are of the following kinds (or a combination of these):
 \begin{djenum}
  \item \label{i:info} Local change in the ``information''.
  \item \label{i:loc} Locally recognizable change in the ``structure''.
  \item \label{i:prog} Program change.
  \item \label{i:glob} Locally unrecognizable change in ``structure''.
 \end{djenum}

A locally recognizable structure change would be a change in the
address field.
A locally unrecognizable change would be to erase two neighbor
colonies based, say, at \( \B \Q \) and \( 2\B \Q \) and to put a new colony in the
middle of the gap of size \( 2\B \Q \) obtained this way, at position
\( 1.5\B \Q \). 
Cells within both the new colony and the remaining old colonies will
be locally consistent with their neighbors; on the boundary, the cells
have no way of deciding whether they belong to a new (and wrong)
colony or an old (and correct) one.

The only kind of disruption whose correction can be attempted along
the lines of traditional error-correcting codes and repetition is the
one of kind~\ref {i:info}: a way of its correction was indicated in
Section~\ref{sec:ftol-block}.  
The three other kinds are new; we will deal with them in different ways.

On locally recognizable changes in the structure, we will use
the method of destruction and rebuilding.
Cells that find themselves in structural conflict with their
neighbors will become vacant.
Vacant cells will eventually be restored if this can be done
fast, in a way structurally consistent with their neighbors.

To fight program changes, our solution will be that the simulation will
not use any ``program'' or, in other words, we
``hard-wire''\index{hard-wiring} the program into the transition function
of each cell.
We will not lose computational universality this way: as the hard-wired
program itself can simulate some
fixed cellular automaton on one of its tracks (which will be called the \( \Payload \) track).

To fight locally unrecognizable changes, we will ``legalize'' all
the structures brought about this way.
Consider the example where a single colony sits in a gap of size \( 2\B \Q \).
The decoding function is defined even for this configuration.
In the decoded configuration, the cell based at site 0 is followed by a
cell at site \( 1.5\B \Q \) which is followed by cells at sites \( 3\B \Q,4\B \Q \), etc.
These configurations, viewed as ``illegal'' earlier, will be ``legalized'' now;
this way they can be eliminated with their own active
participation, provided we have rules (trajectory conditions) applying to
them.
This is the main reason for the introduction of generalized cellular
automata.

The generalized cellular automaton introduced this way will be called a
\df{robust medium}.
The generalization of the notion of the medium does not weaken the
original theorem: the fault-tolerant cellular automaton that we
eventually build is a cellular automaton in the old sense.
The more general media are only needed to reason about 
structures that arise in simulations by a random process.

\section {Results for the finite space}
 \label{sec:lat}

This section overviews the main theorems concerning reliable
computation with cellular automata in discrete time.

\subsection {Relaxation time and ergodicity}
 \label{sec:toom}

Ergodicity means forgetting eventually everything about the initial
 configuration.
We quantify ``eventually'' here by the notion of a relaxation time.
 
We start  with some notation.

\begin{notation}
Let \( \sites(n) \)\glo{l.cap.s@\( \sites(n) \)}
 be an increasing sequence of finite subsets of \( \sites \) with
\( \bigcup_{n}\sites(n)=\sites \), for example
 \[
   \sites(n)=(\clint{-n}{n}\cap \bbZ)^{d}
 \] 
if \( \sites=\bbZ^{d} \).
We can view an element \( \bs \) of \( \states^{\sites(n)} \) as a vector whose
components are indexed with the elements of \( \sites(n) \).
For a measure \( \nu \) over configurations, denote%
 \glo{n:greek@\( \nu(\bs) \)}%
 \begin{equation}\label{eq:nu(s)}
  \nu(\bs)= \nu\setOf{\xi}{\xi(x) = s_{x}\text{ for all }x\in\sites(n)}.
 \end{equation}
For \( n=0 \), we have the special case
 \( \nu(s)=\nu\setOf{\xi}{\xi(0)=s} \).

If \( (\mu,\eta) \) is a random trajectory of a probabilistic cellular
automaton then let \( \mu^{t} \)\glo{m.greek@\( \mu^{t} \)} be the distribution of the
configuration \( \eta(\cdot,t) \).
  \end{notation}

 \begin{definition}[Weak convergence]\label{def:weak-conv}
A sequence \( \nu_{k} \) of measures over \( \states^{\sites} \)
 \df{weakly converges}\index{weak convergence} to measure \( \nu \) if 
\( \lim_{k}\nu_{k}(\bs)=\nu(\bs) \) for
all \( n \) and for all \( \bs\in\states^{\sites(n)} \).
 \end{definition}

In other words, the \( \nu_{k} \)-probabilities of any event determined by a finite number of
sites converge to its \( \nu \)-probability.
Let us proceed to the definition of ergodicity.
There is a linear operator \( P \)\glo{p.cap@\( P \)}
determined by the transition
function \( \bP(s,\br) \) giving \( \mu^{t+1}=P\mu^{t} \):

 \begin{definition}[Markov operator, invariance]
Assuming for simplicity, \( \sites=\bbZ \),
we define the \df{Markov operator}\index{Markov!operator}
\( P \) on measures by
 \begin{equation}\label{eq:P.def}
  (P\mu)(\bs)= \sum_{\mathbf{r}}
  \prod_{j=-n+1}^{n-1} \bP(s_{j},(r_{j-1},r_{j},r_{j+1})) \mu(\br)
 \end{equation}
 where the summation goes over all possible strings
 \( \mathbf{r}=(r_{-n},\dots,r_{n})\in\states^{2n+1} \).
Given a Markov operator, we call a measure \( \alpha \) over configurations
\df{invariant}\index{measure!invariant} if \( P\alpha=\alpha \).
  \end{definition}

It is well-known and easy to prove using standard tools (see a
reference in~\cite{Toom80}) that each continuous linear operator over
probability measures has an invariant measure.
The invariant measures describe the possible limits (in any
reasonable sense) of the distributions \( \mu^{t} \).

\begin{definition}[Ergodicity]
  A probabilistic cellular automaton is called \df{ergodic}\index{ergodic}
  if
  \begin{alphenum}
  \item\label{i:ergodic.a} It has only one invariant measure \( \nu \).
  \item\label{i:ergodic.b}
    For every possible measure \( \mu^{0} \) over configurations, \( \mu^{t}=P^{t} \mu^{0} \)
    converges weakly to \( \nu \).
  \end{alphenum}
\end{definition}

In other words, all information about the initial configuration will
be eventually lost.
A noisy cellular automaton, whenever the set of sites is finite, is a
finite Markov chain with all positive transition probabilities.
  This is ergodic by a well-known elementary theorem of probability theory.
  If the set of sites is infinite then noisiness does not imply even
ergodicity.

  \begin{remark}
  No examples are known of noisy cellular automata over an infinite space
  that satisfy~\ref{i:ergodic.a} but not~\ref{i:ergodic.b}.
  In one dimension and continuous time there is no such example, as shown in~\cite{Mountford95}.
  \end{remark}

  The first example of a non-ergodic noisy cellular automaton was given
by Toom.
(See for example \cite{Toom80}.)

\begin{definition}[Toom rule]
A deterministic cellular automaton
\( R \)\glo{rem:cap@\( R \)}\index{Toom rule} given by Toom can be defined as follows.
We start from a two-dimensional deterministic cellular automaton \( R \)
with set of states \( \{0,1\} \), by the neighborhood \( H= \)
\( \set{\pair{0}{0},\pair{0}{1},\pair{1}{0}} \).
The transition function \( \trans_R(x_{1},x_{2},x_{3}) \) is the majority of
\( x_{1},x_{2},x_{3} \).
\end{definition}

Thus, in a trajectory of the Toom rule \( R \), to obtain the next state of a cell,
take the majority among the states of its northern and eastern
neighbors and itself.
Toom showed that the rule \( R \) remembers a field (namely the whole
state of the cell) in the sense defined in Section~\ref{sec:perturb}
and is hence non-ergodic.

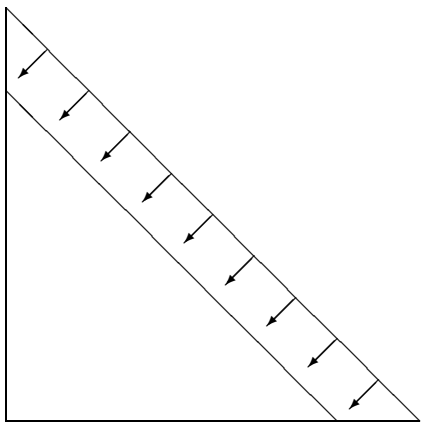
\begin{figure}
 \setlength{\unitlength}{0.55mm}
 \[
\begin{picture}(100,100)
  \put(0,0){\line(1,0){100}}
  \put(0,0){\line(0,1){100}}
  \put(100,0){\line(-1,1){100}}
  \put(80,0){\line(-1,1){80}}
  \multiput(90,10)(-10,10){9}{\vector(-1,-1){7}}
\end{picture}
 \]
 \caption{The Toom rule's effect on a large triangular island}
 \label{fig:Toom}
\end{figure}

\begin{remark}
  The notion of ergodicity was defined only for trajectories of \( R_{\eps} \)
that are also trajectories of some probabilistic cellular automaton \( R' \)
(as defined in Section~\ref{sec:PCA}) such that \( R'\subseteq R_{\eps} \),
that is that all trajectories of \( R' \) are trajectories of \( R_{\eps} \).
The difference between \( R_{\eps} \) and \( R' \) is that the local
transition probabilities of a trajectory of \( R' \) are fixed and the
same everywhere in space-time, while \( R_{\eps} \) requires them
only to be within some range.
Toom's theorem implies that no such \( R' \) is ergodic.
 \end{remark}

\paragraph{Relaxation time as a measure of information loss}
What is the relevance can results on infinite cellular automata 
for the possibilities of computation or information storage in finite systems?
We will try to answer this question here.
To stay in the context of the Toom rule, consider a (finite
or infinite) two-dimensional space \( \sites \).
We need a notion of distance for measures.

\begin{definition}[Variation distance]
Let us define the variation distance for measures \( \mu \) and \( \nu \)
are restricted to \( \states^{\sites(n)} \):%
 \glo{d@\( d^{n}(\mu, \nu) \)} 
 \[
 d^{n}(\mu,\nu) 
= \sum\evof{|\mu(\bs)-\nu(\bs)|\text{ for all } \bs\in\states^{\sites(n)}}.
 \]
\end{definition}

The maximum distance 2 means that the two measures have disjoint support.
Let us qualify distance even more as a function of the size of
the finite two-dimensional spaces that we will be considering:

\begin{definition}
For the set of sites%
 \glo{c.cap@\( \sites_{m} \)}%
 \[
 \sites_{m}=\bbZ_{m}\times\bbZ_{m}
 \]
  where \( m \) can be finite or infinite, suppose that the local transition matrix (for
simplicity, with nearest-neighbor interaction) gives rise to an
ergodic Markov process with Markov operator \( P_{m} \)\glo{p.cap@\( P_{m} \)}.
Let%
 \glo{def:cap@\( D^{n}_{m}(t) \)}
 \[
 D^{n}_{m}(t) = \bigvee_{\mu,\nu} d^{n}(P_{m}^{t}\mu, P_{m}^{t}\nu),
\]
where the dependence on the transition matrix is understood but is not shown.
\end{definition}

  It is easy to see that this is equal to 
 \[
   \bigvee_{\br,\bs} d^{n}(P_{m}^{t}\delta_{\br},P_{m}\delta_{\bs})
 \]
  for the measure \( \delta_{\bs} \) concentrated on configuration \( \bs \).
The function \( D^{n}_{m}(t) \) is monotonically decreasing in \( t \) since
for example for invariant \( \mu \),
 \[
   d^{n}_{m}(P_{m}^{t+u}\nu,\mu)= d^{n}(P_{m}^{t}(P_{m}^u\nu),\mu).
 \]
  \( P_{m} \) is ergodic if and only if for each finite configuration \( \bs \) and each measure \( \nu \) we have
\( P_{m}^{t}\nu(\bs)\to \mu(\bs) \).
By the weak compactness of the space of measures, this is equivalent
to saying that
  \[
   \lim_{t} D^{n}_{m}(t)=0
  \]
  holds for all \( n \).

  \begin{definition}[Relaxation time]\label{def:relaxation-time}
We will call the function%
 \glo{r@\( r_{m}(n,\delta) \)}%
 \[
  r_{m}(n,\delta) := \sup\setOf{t}{D^{n}_{m}(t)>\delta}<\infty
 \]
 the \df{relaxation time}\index{relaxation time}.
  \end{definition}
 Ergodicity implies \( r_{m}(n,\delta)<\infty \) for all \( \delta \).
This function is obviously increasing in \( n \) (defined only for
\( n\le (m-1)/2 \)) and decreasing in \( \delta \).
In the cases we are interested in, the order of magnitude of
\( r_{m}(n,\delta) \) as a function of \( n \) does not change fast with the
change of \( \delta \): \( r_{m}(n,0.1) \) is not much larger than \( r_{m}(n,1.9) \).
This means that once the unique invariant distribution \( \mu_{m} \) is not separated well from any
\( P_{m}^{t}\nu \) it soon becomes fairly close to all of them.
Therefore we will generally consider \( \delta \) fixed.

\paragraph{Relaxation time as a function of space size}
  If \( m<\infty \) then the medium is always ergodic.
Assume now that the medium is also ergodic for \( m=\infty \).
This has a serious implication on the relaxation times for finite \( m \):

 \begin{lemma}\label{lem:unif}
  For all \( n<(m-1)/2 \), for all \( \delta \) with \( r_{\infty}(n,\delta)<(m-1)/2-n \)
we have
  \[
   r_{m}(n,\delta)\le r_{\infty}(n,\delta).
  \]
 \end{lemma}
  This means that if the medium is ergodic for \( m=\infty \) then increasing
\( m \) in the case of \( m<\infty \) does not increase the relaxation time
significantly for any fixed \( n \): in each segment of length \( n \) of any
finite medium, information is being lost at least as fast as in the
infinite medium.
 \begin{proof}
  Take \( m,n,\delta \) satisfying the above conditions and let
\( r=r_{\infty}(n,\delta) \).
  Due to the monotonicity of \( D^{n}_{m}(t) \), it is enough to prove that 
\( D^{n}_{m}(r)\le D^{n}_{\infty}(r) \).
Take a measure \( \nu \) over configurations of \( \sites_{m} \), this will give
rise to some measure \( \nu_{\infty} \) over configurations of period \( m \) in
\( \sites_{\infty} \) in a natural way, where \( \nu_{\infty} \) is such that for all
\( n<(m-1)/2 \) and all \( \bs\in\states^{\sites(n)} \) 
we have \( \nu_{\infty}(\bs)=\nu(\bs) \).
Then, \( r<(m-1)/2-n \) implies \( 2n+1+2r<m \) and therefore via~\eqref{eq:P.def}
we have
 \[
  P_{m}^r\nu(\bs)=P_{\infty}^r\nu_{\infty}(\bs).
 \]
\end{proof}
 
We found that ergodicity of \( P_{\infty} \) implies a kind of
 \emph{uniform forgetfulness}\index{forgetful!uniformly} in the sense that
increasing the size \( m \) of the space does not increase the
relaxation time beyond \( r_{\infty}(n,\delta) \).


\paragraph{Forgetfulness: a variant of ergodicity}
Ergodicity does not express quite adequately the losing of all
information about the initial configuration in case of cellular automata,
where namely the space-time is translation-symmetric (this was noticed by
Charles Radin and Andrei Toom).

\begin{example}
Let \( \xi_{0} \) be the configuration over the one-dimensional integer
lattice that is 0 in the even places and 1 in odd places, and \( \xi_{1} \)
the one that is 1 in the even places and 0 in the odd places.
For \( b=0,1 \) let \( \mu_{b} \) be the measure concentrated on \( \xi_{b} \),
and let \( P \) be the linear operator obtained from some transition function.
Suppose that the measures \( P^{n}\mu_{0} \) converge to some measure
\( \nu_{0} \).
Then the measures \( P^{n}\mu_{1} \) converge to some measure \( \nu_{1} \).
Even if \( \nu_{1} \) is different from \( \nu_{0} \) they differ only by a
translation in space.
If the translations of \( \nu_{0} \) are the only invariant measures of
\( P \) then we would still say that in some sense, \( P \) loses all
information about the initial configuration: we might say, it loses
all ``local'' information.
Indeed, a cell has no way of knowing whether it has even or odd
coordinates.
\end{example}


\begin{definition}
  We can call the Markov operator \( P \)
 \df{strongly not forgetful}\index{forgetful!strongly not} if it has two
disjoint (weakly) closed translation-invariant sets of measures.
\end{definition}

Our non-ergodic examples will all be strongly not forgetful.

 \subsection{Information storage and computation}
 \label{sec:1dim.comp}

We recall the existence of three-dimensional reliable cellular automata.

\begin{definition}
Let us be given an arbitrary one-dimensional transition function
\( \trans \) and the integers \( N,L \).
We define the three-dimensional transition function \( \trans' \) as follows.
  The interaction neighborhood is \( H\times\{-1,0,1\} \) with the
neighborhood \( H \) defined in Section~\ref{sec:toom} above.
  The transition function \( \trans' \) says: in order to obtain your state at
time \( t+1 \), first apply the transition function \( R \) in each plane defined by
fixing the third coordinate.
  Then, apply transition function \( \trans \) on each line obtained by fixing
the first and second coordinates.
\footnote{The papers~\cite{GacsReif3dim88} and \cite{Gacs2dim89} use a
neighborhood instead of \( H \) that makes some proofs easier: the
northern, south-eastern and south-western neighbors.}

For an integer \( m \), define the space \( \sites =\bbZ_{m}^{2}\times \bbZ_{N} \).
For a trajectory \( \xi \) of \( \CA(\trans) \) on \( \bbZ_{N} \), define the
trajectory \( \zeta \) of \( \CA(\trans') \) on \( \sites \) by
 \begin{equation} \label{eq:reifCode} 
	\zeta(i,j,n,t)=\xi(n,t).
 \end{equation} 
  Thus, each cell of \( \CA(\trans) \) is repeated \( m^{2} \) times, on a whole
``plane'' (a torus for finite \( m \)) in \( \sites \).
  \end{definition}

The following has been proved in earlier work:

\begin{proposition}[G\'acs-Reif]
There are constants
\( \eps_{0},c_{1},d_{1}>0 \) such that the following holds.
For all \( N,L \), and \( m=c_{1} \log(N L) \), for any trajectory \( \xi \) of
\( \CA(\trans) \) over \( \bbZ_{N} \), if the trajectory \( \zeta \) of
\( \CA(\trans') \) is defined by~\eqref{eq:reifCode} then for any
\( \eps<\eps_{0} \), any trajectory \( (\mu,\eta) \) of \( \CA_{\eps}(\trans') \)
such that \( \eta(\cdot,0)=\zeta(\cdot,0) \) we have for all \( t \) in
\( \rint{0}{L} \) and all \( w\in \sites \)
 \[ 
  \mu\evof{\eta(w,t)\ne  \zeta(w,t)} \le d_{1}\eps .
 \]
\end{proposition}

This theorem says that in case of the medium \( \CA_{\eps}(\trans') \) and
the trajectories \( (\mu,\zeta) \), the probability of
deviation\index{deviation} can be uniformly bounded by \( d_{1}\eps \).
The trajectories \( \eta \) encode (by~\eqref{eq:reifCode}) an arbitrary
computation (for example a universal Turing machine), hence this theorem
asserts the possibility of reliable computation in three-dimensional space.
The coding is repetition \( O(\log^{2}(N L)) \) times, that is it depends on the
\df{size}\index{computation!size} \( N\cdot L \) of the computation.
The decoding is even simpler: if a plane of \( \sites \) represents a state
\( s \) of a cell of \( \CA(\trans) \) then each cell in this plane will be in
state \( s \) with large probability.
The simulation occurs in ``real time''.
The original proof of a slightly weaker version of this result used a
sparsity\index{sparsity} technique borrowed from~\cite{Gacs1dim86}.
In its current form, the theorem was proved in~\cite{BermSim88}
using an adaptation of the techniqe of~\cite{Toom80}.

Theorem~\ref{thm:1dim.nonerg} shows that one-dimensional noisy and
strongly not forgetful probability operators exist.
Section~\ref{sec:hier} gave some detail on the kind of constructs
going into the proof.
The actual proof seems to require almost the whole complexity of
the constructions of the present paper (though the continuous-time case adds
some additional nuisance to each part).
Once the basic structure (an amplifier, as asserted in
Lemma~\ref{lem:abstr-amp}) is in place, the simulations in it support
arbitrary computational actions and allow the formulation of several other
theorems.
Theorem~\ref{thm:remember-inf} asserts the possibility of storing an infinite 
string \( \rho \)---provided the space is infinite.
Now we can give more detail on the nature of the encoding \( \psi_{*} \) used
there: it will be one like in Definition~\ref{def:encoding-infinite}.
We don't give a finite version: that one can be seen as a special case of the finite version
of Theorem~\ref {thm:1dim.comp.finite} below, when the computation is a trivial
one outputting the input.
Here are some definitions needed before Theorem~\ref {thm:1dim.comp.finite}.
If \( \rho \) is a configuration in \( \Sigma^{\bbZ} \) where \( \Sigma \) is the standard
  computing alphabet, then we say that \( \rho \) is \df{supported by} interval \( I \) if
\( \rho(x)=\#\cdots\# \) or vacant for all \( x\not\in I \).
Recall \( \Init_{\psi_{*}}(\rho) \) in Definition~\ref{def:Init}.

We defined the function \( h_{1}(t,c) \) in~\eqref{eq:h_1}.

\begin{theorem}[Reliable computation in finite or infinite space]\label{thm:1dim.comp.finite}
  Let \( \trans \) be a standard computing transition function over an alphabet \( \Sigma \),
  and \( \delta>0 \) a constant.
  There is
 \begin{itemize}
  \item a transition function \( \trans' \) with a state space \( \states \) having fields \( \Input \),
\( \Output \);
\item a code \( \psi_{*}: \Sigma_{0}^{*}\to \states^{*} \) with \( |\psi_{*}(\rho)|\le 2|\rho| \),
  \item constants \( c_{1}, c_{2}>1 \)
 \end{itemize}
 such that the following holds for all sufficiently small \( \eps \).
 Let \( N \) be the size of our space, and \( \rho\in\Sigma_{0}^{*} \) with
 \( 2|\rho|<N \), \( \eps<\eps_{0} \).
 Let \( \zeta(x,t) \) be any trajectory of \( \trans \) with
 \( \zeta(x,0).\Input=\rho(x) \) for \( x\in\lint{0}{|\rho|} \) and filled with \( * \)'s
 otherwise.
 Let \( \eta \) be any trajectory over the space \( \bbZ^{N} \)
 of an \( \eps \)-perturbation of \( \trans' \), with
 \( \eta(\cdot,0)=\Init_{\psi_{*}}(\rho) \).
 For all sufficiently small \( \eps \), for all
 \( t \) for which \( \zeta(\cdot,t) \) does not use more space than \( N \),
all \( t'>t\cdot h_{1}(t,c_{1}) \), all \( x\in\bbZ \) we have
 \[
 \Pbof{\zeta(x, t).\Output \not\preceq \eta(x, t').\Output} < \delta + t'\eps^{h_{0}(N,c_{2})}.
\]
\end{theorem}

As we see the space occupied by our reliable simulation is just a constant times
more than the original space (how much larger is \( \norm{\states} \) than \( \norm{\Sigma} \)).
The time delay is more significant: we can recover the computation
result \( \zeta(x,t) \) only at times \( t'>t\cdot h_{1}(t,c_{1}) \).

The above theorem strengthens the result of~\cite{Gacs1dim86} in some ways.
\begin{itemize}
\item The need for the decoding of the computation result is eliminated: it appears
  in the same place as in \( \zeta \), due to an economical encoding.
\item The encoding of the input depends only on the input, not on the size of the computation.
  This is achieved by lifting, as indicated after Theorem~\ref{thm:1dim.nonerg}.
\end{itemize}

\section{More restrictions on media}\label{sec:med}

 Here we impose further requirements on our media
 \[
 \Med(\states, \sites, \B,\Configs, \Histories, \Trajs,  (\cA_{t})_{t\ge 0})
 \]
 from Section~\ref{sec:media}.

\subsection {Trajectories}
 \label{sec:traj}

Here, we will show in what form the set of trajectories of a medium will
be given.
The kind of conditions defining what is a trajectory in a medium will be of the form
\begin{align*}
   \g(\eta)\le b,
\end{align*}
where the function \( \g \) is an event function as in Definition~\ref{def:events}.

  \begin{example}
Let $V=\rint{a_{1}}{a_{2}}\times\rint{t}{u}$, $\B  < d = 0.2 (a_{2}-a_{1})$.
Suppose that we are given some $\br\in\states^{3}$, $s\in\states$.
Let $\g(\eta)=1$ if there is a space-time point $\pair{x}{t_{0}}$ in $V$, and an
$\eps<d$ such that $\eta(x+i\B ,t)=r_{i}$ for $i = -1, 0, 1$ during
the interval $\lint{t_{0}-\eps}{t_{0}}$ and $\eta(x,t_{0})=s$.
Thus, the event $\g(\eta)=1$ marks a certain kind of transition in the window.
  \end{example}

  The set \( \Trajs \) of a medium will be given by many such conditions.
Let us fix an arbitrary set   \glo{thm:cap.t@$\bbT$} 
\begin{align}\label{eq:cond-types}
 \bbT  
\end{align}
called the set of \df{condition types}\index{condition!type}.
  Let 
\begin{align*}
 \Rectangles  
\end{align*}
be the set of bounded space-time rectangles of the form \( \lint{a}{b}\times\rint{u}{v} \).
For a type $\alpha\in\bbT$ and a rectangle $V$, 
a \df{local condition}\index{condition!local} will be defined as a pair
\glo{a.greek@$(\alpha,V)$}%
\begin{align}\label{eq:rectangle-cond}
   (\alpha,V).   
\end{align}
Two local conditions are \df{disjoint}\index{condition!local!disjoint}
when their rectangles are.
The conditions defining a trajectory will be given by two functions:
  \begin{align}\label{eq:g-b}
   \g:\bbT\times\Rectangles\times\Histories\to \{0,1\},\quad h:\bbT\to\clint{0}{1},
  \end{align}
  where \( \g(alpha,V,\eta) \)  \glo{g@$\g(\alpha,\eta)$}%
 is for each \( \alpha, V \) an event function in the sense of Definition~\ref{def:events},
and  $b(\alpha)\in\clint{0}{1}$\glo{b@$b(\alpha)$}
is a bound: so each condition will have the form
\begin{align}\label{eq:g-b-ineq}
  \Expv \g(\alpha,V,\eta)\le b(\alpha).   
\end{align}
The  medium is then defined as
 \begin{equation}\label{eq:med}
  M=\Med(\states, \sites, \bbT,\Configs, \Histories, \{\cA_{t}\},
\g(\cdot,\cdot,\cdot), b(\cdot)),
\end{equation}
where all arguments inferrable from the context can be omitted.
An event defined by \( \g(alpha,V,\eta) \)  occurs when
\( \eta \) \emph{violates} a certain local condition.
 For technical reasons we
 we always include two special condition types: the \df{killer} type
 \( \omega, \) that does not  allow anything, and the \df{dummy} type
\( \delta \), that allows everything:
\begin{align}\label{eq:omega-delta}
 \g(\omega,V,\eta)=b(\omega)=0,\quad \g(\delta,V,\eta)=b(\delta)=1.
\end{align}
For an empty ``rectangle'' $\emptyset$ we stipulate 
 \( \g(\alpha,\emptyset,\eta) = 0  \) 
 for all $\alpha\ne \delta$, and \( \g(\delta,\emptyset,\eta)=1 \).

Generally, local conditions will only be defined for small $V$, and
in a translation-invariant\index{translation-invariance} way (in space, but
only almost so in time since the time 0 is distinguished).
Therefore one and the same pair $g,b$ can serve for all sufficiently
large finite or infinite spaces $\sites$ as well as the infinite space, just like with
deterministic cellular automata, where the same transition function
determines trajectories in finite spaces as well as the infinite space.
When we are given media $M_{1},M_{2},\dots$ then $\g_{i}(),b_{i}()$
automatically refer to $M_{i}$.

Recall the definition \( \prt \) from~\eqref{eq:proj}.

\begin{definition}[Trajectory]\label{def:traj}
  A random history $\eta$ will be called a
\df{trajectory}\index{trajectory!weak} of medium $M$ given in~\eqref{eq:med}
if for each time $u\ge 0$, for each set of disjoint local conditions
 \[
  (\alpha_{i},V_{i})_{i\le N}
 \]
  with $(\forall{i})\inf\prt V_{i}\ge u$ (where $\inf\emptyset=\infty$), 
we have
 \begin{equation}\label{eq:weak-traj-def}
   \Expv\Bigsq{\prod_{i}\g(\alpha_{i}, V_{i},\eta)\Bigm\vert \cA_{<u}}
   \le \prod_{i} b(\alpha_{i}).
 \end{equation}
\end{definition}

This says that we can multiply the probabilities of violating local
  conditions \( \Expv\g(\alpha_{i},V_{i},\eta)\le b(\alpha) \)
  in disjoint windows, in other words: these violations happen ``independently''.
  The media we are interested in have the following three types of condition:
 \begin{djenum}

  \item One type recognizes in $V$ a certain event called
``damage'', and has $b=\eps$ for some small $\eps$;

  \item One type recognizes the event that 
the transition does not occur there according to a certain transition
function.
This has $b=0$, thus prohibiting this possibility.

  \item
  One type corresponds to some ``coin tossing'' event, and has
$b=0.5+\eps'$ for a suitable $\eps'$.

 \end{djenum}

The following example will show how a discrete-time
probabilistic cellular automaton ($\PCA$) fits into this framework.

\begin{example}\label{xmp:PCA-medium}
\begin{sloppypar}
 Let us be given a probabilistic cellular automaton
$\PCA(\bP,\B ,\T )$.
For each state vector $\br\in \states^{3}$ and state $s$, let us form type
$\alpha(s,\br)$. 
The condition of this type will say that the transition rules are obeyed at times
of the form $n\T $.
Namely for cell $x$, time $t$ of the form $n\T $ with $n>0$, let 
\end{sloppypar}
 \begin{align}\label{eq:PCA-medium}
\nonumber      \g\bigl(\alpha(s,\br),\;\{x\} &\times\rint{t-\T }{t},\;\eta\bigr)
\\        &= \{ \eta(x+i\B ,t-\T /2)=r_{i}\;(i=-1,0,1),\;\eta(x,t)=s \},
\\    b(\alpha(s,r)) &= \bP(s,\br).
 \end{align}
  For all cells $x$ and rational times $u$ we define a new type
 $\beta=\beta(x,u)$ 
 by $b(\beta(x,u)) = 0$ and
 \[
  \g\bigparen{\beta(x,u),\;\{x\}\times\clint{\T \flo{u/\T }}{u},\;\eta}
   = \setof{\eta(x,u)\ne\eta(x,\T \flo{u/\T })}.
 \]
This condition says that $\eta(x,u)=\eta(x,\T \flo{u/\T })$ with probability 1.
For all other combinations $\alpha,V$ we set $\g(\alpha,V,\eta)=0$.
It is easy to see that the trajectories of $\PCA(\bP,\B ,\T )$ are
trajectories of the medium defined this way.
  \end{example}

  The new conditions introduced below are ``loosened-up'' versions of the
above ones since it is not always desirable to be as restrictive.
Instead of a bound on the probability of getting into a local state $s$, we may
just want a bound on the probability of hitting a certain set of local states
$K''$ once it is known that we were going to hit a certain bigger
set of local states $K'$.

\begin{definition}[Local state subset functions]\label{def:state-subsets}
In a probabilistic cellular automaton $\PCA(\bP,\B ,\T )$ 
for each set of states $E$ and state vector $\br=(r_{-1},r_{0},r_{1})$ let
\( \bP(E,\br)=\sum_{q\in E}\bP(q,\br) \).
Let $\bK$\glo{k.cap@$\bK$} be the set of functions $K:\states^{3}\to 2^{\states}$
such that
  \begin{equation}\label{eq:K-def}
  K(\br)\subseteq \setof{s\ne r_{0}: \bP(s,\br)>0}.
  \end{equation}
  For a space-time trajectory \( \eta \) and a \( K\in\bK \) let us write%
 \glo{k.cap@\( \ol K(x,t,\eta) \)}%
 \[
   \ol K(x,t,\eta)=K(\eta(x-\B ,t), \eta(x,t), \eta(x+\B ,t)).
 \]
  For all \( K',K''\in\bK \) with \( K'\spsq K'' \) (that is \( K'(\br)\spsq K''(\br) \)
  for all \( \br \)), let%
 \glo{c.cap@\( c(K',K'') \)}%
 \[
   c(K',K'') =\bigvee_{\br} \frac{\bP(K''(\br),\br)}{\bP(K'(\br),\br)}
 \]
  be the maximum conditional probability, over all \( \br \), of
getting into the set \( K''(\br) \) provided we get into the set \( K'(\br) \).  

In case of a continuous-time cellular automaton \( \CCA(\bR,\B ) \),
for any functions \( K,K',K''\in\bK \) with \( \bR(K,\br) = \sum_{q\in K}\bR(q,\br) \)
we define the corresponding quantity \glo{c@\( c(K',K'') \)}%
 \begin{align*}
c(K',K'')   = \bigvee_{\br}\frac{\bR(K''(\br),\br)}{\bR(K'(\br),\br)}.
 \end{align*}
\end{definition}

Let us form local conditions with the help of these concepts:

\begin{example}\label{xmp:PCA-stop}
In case of a probabilistic cellular automaton \( \PCA(\bP,\B ,\T ) \),
for each pair of such functions \( K',K''\in \bK \) with \( K'\spsq K'' \) let us form the type%
 \glo{a.greek@\( \alpha(K',K'') \)}
 \[
  \alpha =\alpha(K', K'').
 \]
  For each cell \( x \) and times \( u\le v \) of the form \( n\T  \), let us form the 
rectangle \( V=\{x\}\times\rint{u}{v} \).
  For each such pair \( \alpha,V \) let 
  \( \g(\alpha,V,\eta)=1 \) if there is a first \( t\in\rint{u}{v} \) with
 \( \eta(x,t)\in \ol K'(x,t-\T ,\eta) \) and in it, we have
 \( \eta(x,t)\in \ol K''(x,t-\T ,\eta) \).
 Also, let \( b(\alpha)=c(K',K'') \).
 The condition of type \( \alpha \) bounds the probability that as soon as
 we get into the set \( K' \) we also get into the set \( K'' \).
  \end{example}

  The proposition below is an immediate consequence of the definitions: 
  
 \begin{proposition}\label{propo:discr-stop}
  If \( M \) is the medium defined by the conditions above then all trajectories of 
 \( \PCA(\bP, \B , \T ) \) are trajectories of  \( M \).
  \end{proposition}

\begin{example}\label{xmp:CCA-stop}
In case of a continuous-time cellular automaton \( \CCA(\bR,\B ) \)
we define each \( \alpha(K',K'') \), \( V=\{x\}\times\rint{u}{v} \) and \( \g(\alpha,V,\eta) \)
as in Example~\ref{xmp:PCA-stop}, but now for all rational \( u<v \).
 \end{example}

  \begin{theorem}\label{thm:cont-stop}
  Every trajectory of \( \CCA(\bR,\B ) \) is a trajectory of the medium defined
in the above example.
  \end{theorem}
 \begin{proof}[Proof sketch]
  As in the time-discretization in
Section~\ref{sec:Markov}, for a small \( \delta>0 \), let
\( M_{\delta}=\PCA(\bP,\B ,\delta) \) with \( \bP(s,\br)=\delta\bR(s,\br) \) when
 \( s\ne r_{0} \) and \( 1-\delta\sum_{s'\ne r_{0}}\bR(s',\br) \) otherwise.
  Then a transition to the limit \( \delta\to 0 \) in
Proposition~\ref{propo:discr-stop} completes the proof. 
  The transition uses the fact that trajectories of \( M(\delta) \) converge in
distribution to trajectories of \( \CCA(\bR,\B ) \).
 \end{proof}

  Of course in both Proposition~\ref{propo:discr-stop} and in
Theorem~\ref{thm:cont-stop}, instead of the first time with \( \eta \) switching
into  \( K' \) we could have asked about the second time, or the third time, etc.

\subsection{Strong trajectories}\label{sec:strong-traj}

In one part of our construction and proof, for self-organization, in 
Section~\ref{sec:colony-birth}, we need a seemingly stronger condition on trajectories
than the one in Definition~\ref{def:traj}: namely, we must consider
sets of local conditions that are randomly chosen.
Of course, they cannot be chosen completely randomly, as then we could just choose
the places where the faults occur.
But we will allow conditions chosen in a non-anticipating way,
that is via \emph{stopping times}.
Recall the notion of a set of condition types \( \bbT \)  in~\eqref{eq:cond-types},
\eqref{eq:rectangle-cond} and~\eqref{eq:g-b}.

Here, we will specialize, for our needs,
the notion of stopping times from the theory of stochastic processes.
Let \( \cA_{t} \) be a filtration over a probability space,
as in Definition~\ref{def:filtration}.

\begin{definition}[Stopping time]\label{def:stopping-time}
A real random variable $\sigma$ with values in $\lint{0}{\infty}$ is a
\df{stopping time}\index{stopping time} if for each $t$
the event $\{\sigma\le t\}$ is in $\cA_{t}$.
Let $\cA_{\sigma}$\glo{a.cap@$\cA_{\sigma}$} be the algebra of events $A$ with
the property that $A\cap\{\sigma\le t\}\in \cA_{t}$.
\end{definition}

\begin{definition}[Constraint process]
  For some time interval $J\subset\rint{0}{\infty}$, assume that there
  is a stopping time \( \sigma\in J \), and random variables \( \alpha(\eta)\in\bbT \), 
  $V(\eta)\in\Rectangles$ \glo{v.cap@$V^{t}(\eta)$}
  measurable in $\cA_{\sigma}$, with
 \[
  \prt V^{\sigma}(\eta)\subseteq J\cap\clint{\sigma}{\infty}.
 \]
The tuple \( (J,\sigma,\alpha(\eta),V(\eta)) \) will be called a
\df{constraint process}\index{constraint process}.
  \end{definition}

 \begin{definition}[Strong trajectory]
A random history $\eta$ will be called a
\df{strong trajectory}\index{trajectory} if for each system of disjoint
constraint processes 
\begin{align*}
   (J,\sigma_{i},\alpha_{i}(\eta), V_{i}(\eta))_{i\in N}
\end{align*}
(where \( V_{i}(\eta) \) are disjoint for each \( \eta \)),
$h$ an event function over $\cA_{<\prt J}$, we have
 \begin{equation}\label{eq:traj-def}
   \Expv\bigparen{ h
     \prod_{i}\g(\alpha_{i}(\eta),V_{i}(\eta),\eta)}
   \le \Expv h\Expv\prod_{i} b(\alpha_{i}(\eta)).
 \end{equation}
 \end{definition}

 This says, just as for ordinary trajectories,
that we can multiply the probabilities of violating local
conditions in disjoint windows---but requires that this be even true if the
condition types and windows are chosen by random stopping rules.
It may seem strange that there is dependence on a random \( \alpha_{i} \) on both sides,
but note that \( \alpha_{i}(\eta) \) depends only on what happened before \( \sigma_{i} \),
while \( g(\alpha_{i}(\eta),V_{i}(\eta),\eta) \) depends on what happens in \( V_{i}(\eta) \),
whose bottom is after \( \sigma_{i} \). 
 
\begin{proposition}
  The trajectories of the cellular automaton \( \PCA(\bbP,\B, \T) \) of
  Example~\ref{xmp:PCA-stop} are strong.
\end{proposition}

\begin{proof}[Proof sketch] 
 The statement follows by a standard argument from the strong Markov
property of the Markov process \( \PCA(\bB, \B , \T ) \) (see~\cite{Neveu64}).
\end{proof}

\begin{proposition}
  Proposition~\ref{propo:discr-stop} and Theorem~\ref{thm:cont-stop}
  can be strengthened: all trajectories of 
  \( \PCA(\bB, \B , \T ) \) and \( \CCA(\bR,\B) \)  are \emph{strong} trajectories of
  the corresponding media.  
\end{proposition}

The proofs are routine.
  
\begin{figure}
 \setlength{\unitlength}{1mm}
 \[
 \begin{picture}(100, 65)(-20,0)
    \put(-10,0){\vector(0,1){60}}
    \put(-13, 55){\makebox(0,0)[r]{time}}
    \put(35, 23){\framebox(10, 10){$V_{1}^{\sigma_{1}}$}}
    \put(-11, 20){\line(1,0){2}}
    \put(-16, 20){\makebox(0,0){$\sigma_{1}$}}
    \put(55, 10){\framebox(8, 25){$V_{2}^{\sigma_{2}}$}}
    \put(-11, 8){\line(1,0){2}}
    \put(-16, 8){\makebox(0,0){$\sigma_{2}$}}
    \put(25, 5){\framebox(25, 10){$V_{3}^{\sigma_{3}}$}}
    \put(-11, 2){\line(1,0){2}}
    \put(-16, 2){\makebox(0,0){$\sigma_{3}$}}
 \end{picture}
 \]
 \caption{A system of disjoint random local conditions}
 \label{fig:loc-cond}
\end{figure}
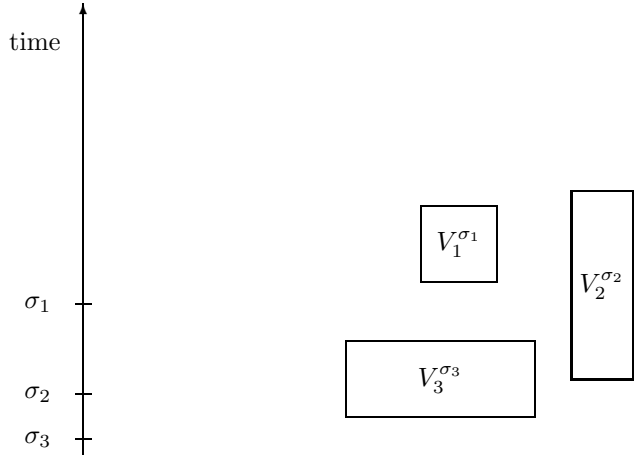

\subsection{Canonical simulations}\label{sec:canonical}

We will restrict the kind of simulations among media 
from its general form given in Definition~\ref{def:general-simulation} 
to a kind of mapping that sets up some
relations between the local conditions of the media and makes it
easier to prove it is a simulation.
\begin{definition}\label{def:canon-simul}
  Let \( M, M^{*} \) be media where the functions \( g,b \) define the constraints on \( M \),
  and \( g^{*},b^{*} \) define those for \( M^{*} \),
  and let \( \fg^{*} \) be a mapping from the
histories of \( M \) into those of \( M^{*} \): we
will write \( \eta^{*}=\fg^{*}(\eta) \).
We say that \( \fg^{*} \) is a \df{canonical simulation map}
\index{simulation!canonical}
if the following holds.
Let \( (\alpha,V) \) be any local condition for \( M^{*} \) with
\( \alpha\in\bbT^{*} \).
Then there is a system of constraint processes
\glo{baij@\( \beta_{i,j}(\alpha) \)}%
\glo{wtaijve@\( W^{t}_{i,j}(\eta) \)}%
\glo{taijve@\( \tau_{i,j} \)}%
 \begin{equation}\label{eq:weak-canonical-sys}
   \bigparen{\prt V, \tau_{i,j},\beta_{i,j}(\eta), W_{i,j}(\eta)}
 \end{equation}
 with \( i=1,\dots,n(\alpha,V) \), \( j=1,\dots,k(\alpha,V,i) \)
 with the following properties.
 For each fixed \( i \), the subsystem containing elements for that
\( i \) forms a system of disjoint constraint processes with
\( W_{i,j}(\eta)\subseteq V \), and for almost all \( \eta \) we have
 \begin{equation}\label{eq:canonical}
 \begin{aligned}
   \g^{*}(\alpha,V,\eta^{*}) &
    \le \sum_{i}\prod_{j}\g(\beta_{i,j}(\eta), W_{i,j}(\eta),\eta),
\\ b^{*}(\alpha) &\ge \sum_{i}\prod_{j} b(\beta_{i,j}(\eta)).
 \end{aligned}
\end{equation}
\end{definition}

\begin{theorem}\label{thm:canonical}
  If \( \fg^{*} \) is a canonical simulation map then for each
strong trajectory \( \eta \) of \( M \), the expression \( \fg^{*}(\eta) \) defines
a strong trajectory of \( M^{*} \). 
\end{theorem}

 \begin{Proof}
  Let \( \eta \) be a strong trajectory of \( M \) and
  \begin{align}\label{eq:the-constr-proc}
 (J, \sigma_{l}, \alpha_{l}(\eta^{*}),V_{l}(\eta^{*}))_{l\in N}  
\end{align}
  a system of disjoint random local conditions in \( M^{*} \).
Let \( h \) be an event function measurable in \( \cA_{<\inf J} \).
We must show
 \[
 \Expv \bigparen{h\prod_{l}\g^{*}_{l}(\alpha_{l}(\eta^{*}),V_{l}(\eta^{*}),\eta^{*})} 
  \le \Expv h \Expv\prod_{l} b^{*}(\alpha^{\sigma_{l}}_{l}(\eta^{*})).
 \]
  By the assumption, for each value \( (a,v) \) of the random
  \( (\alpha_{l}(\eta^{*}), V_{l}(\eta^{*})) \)
  there is a system of constraint processes
 \begin{equation}\label{eq:W-sys}
 \bigparen{\prt v,\tau_{l,i,j}(a,v),\beta_{l,i,j}(a,v,\eta), W_{l,i,j}(a,v,\eta)}
 \end{equation}
  for \( i=1,\dots,n(a,v) \), \( j=1,\dots,k(a,v,i) \), with the
  properties defined for canonical simulations, so \( \tau_{l.i,j}\in\prt v \).
  This implies
  \begin{align}\label{eq:tau-sigma}
   \tau_{l,i,j}(\alpha_{l}(\eta^{*}),V_{l}(\eta^{*}))\ge\sigma_{l}
  \end{align}
  for all \( l,i,j \).
  We will use the shorthand notation
\begin{align}\label{eq:bar-notation}
   \bar W_{l,i,j}(\eta)=W_{l,i,j}(\alpha_{l}(\eta^{*}),V_{l}(\eta^{*}),\eta),
\end{align}
and \( \bar\tau_{l,i,j} \), \( \bar \beta_{l,i,j}(\eta) \) means similarly not showing the
dependence on \( \alpha_{l},V_{l} \).
  
  \begin{claimnonumber}
    For each \( l,i,j \), the tuple
   \begin{align*}
 (J,\bar\tau_{l,i,j},\bar\beta_{l,i,j}(\eta), \bar W_{l,i,j}(\eta))  
\end{align*}
is a constraint process.
  \end{claimnonumber}
 \begin{proof}
   Let us first show that \( \bar\tau_{l,i,j} \) is a stopping time.
   Let
\begin{align*}
   E_{a,v}\eqv (\alpha_{l}(\eta^{*}),V_{l}(\eta^{*}))=(a,v).
\end{align*}
We have
   \begin{align}\label{eq:tau-union}
     \bar\tau_{l,i,j}\le t \eqv \bigvee_{a,v}\tau_{l,i,j}(a,v)\le t\land E_{a,v}.
   \end{align}
From~\eqref{eq:tau-sigma} we get that
\( \bar\tau_{l,i,j}\le t \) implies \( \sigma_{l}\le t \), therefore 
\begin{align*}
\bar\tau_{l,i,j}\le t\land E_{a,v}\eqv \bar\tau_{l,i,j}\le t\land\sigma_{l}\le t\land E_{a,v}.
\end{align*}
   By the definition of stopping times, \( \{\tau_{l,i,j}(a,v)\le t\}\in\cA_{t} \).
   Because~\eqref{eq:the-constr-proc} is a constraint process for each \( l \),
we have also 
\begin{align}\label{eq:E_{a,v}}
 \{\sigma_{l}\le t\land E_{a,v}\}\in\cA_{t},
\end{align}
hence the countable union~\eqref{eq:tau-union} is also in \( \cA_{t} \).
Similarly, for all \( (b, w)\in \bbT\times\Rectangles \), the event  
   \begin{align}\label{eq:W-event}
  \bar\tau_{l,i,j}\le t \land (\bar\beta_{l,i,j}(\eta),\bar W_{l,i,j}(\eta))=(b,w)
   \end{align}
   can be written as
\begin{align*}
  \bigvee_{a,v}\bigparen{\sigma_{l}\le t \land E_{a,v}\land \bar\tau_{l,i,j}\le t
  \land (\beta_{l,i,j}(a,v,\eta),W_{l,i,j}(a,v,\eta))=(b,w)}.
\end{align*}
We have seen~\eqref{eq:E_{a,v}}, and because~\eqref{eq:W-sys} is a system of constraint
processes, also 
\begin{align*}
 \{\bar\tau_{l,i,j}\le t\land (\beta(a,v,\eta),W_{l,i,j}(a,v,\eta))=(b,w)\}\in\cA_{t} ,
\end{align*}
therefore~\eqref{eq:W-event} is in \( \cA_{t} \).
 \end{proof}

  By~\eqref{eq:canonical} we have, with \( \bar\alpha_{l}=\alpha_{l}(\eta*) \),
 \( \bar V_{l}=V_{l}(\eta^{*}) \) and the notation in~\eqref{eq:bar-notation}
 \begin{equation}\label{eq:weak-canonical-appl}
 \g^{*}(\bar\alpha_{l},\bar V_{l},\eta^{*}) \le
   \sum_{i=1}^{n}\prod_{j=1}^{k} \g(\bar\beta_{l,i,j}(\eta),\eta),\bar W_{l,i,j}(\eta),\eta).
\end{equation}
Here we did not make \( n \) dependent on \( \bar\alpha_{l} \), \( \bar V_{l} \),
nor \( k \) dependent on these and \( i \),
setting them to the maxima of all the possible values.
Recall the special types \( \omega \), \( \delta \) defined in~\eqref{eq:omega-delta}.
Whenever \( i>n(\bar\alpha_{l},\bar V_{l}) \), we set
\( \beta_{l,i,j}=\omega \), meaning that the \( i \)th additive term is omitted.
And whenever \( j>k(\bar\alpha_{l},\bar V_{l},i) \), we set
\( \beta_{l,i,j}=\delta \), meaning that the \( j \)th multiplicative factor is omitted.
Now
 \[
  \prod_{l} b^{*}(\bar\alpha_{l})\ge
     \prod_{l}\sum_{i}\prod_{j}b(\bar\beta_{l,i,j}(\eta))
 = \sum_{i(\cdot)}\prod_{l,j} b(\bar\beta_{l,i(l),j}(\eta))
 \]
 where the last sum is over all functions \( i(\cdot) \) with \( i(l)\in[1,n] \).
  Similarly, by~\eqref{eq:weak-canonical-appl}
 \[
 \g^{*}(\bar\alpha_{l},\bar V_{l},\eta^{*}) \le 
   \sum_{i(\cdot)}\prod_{l,j}\g(\bar\beta_{l,i,j}(\eta),\bar W_{l,i,j}(\eta),\eta).
 \]
  Therefore it is sufficient to show for a fixed function \( i(\cdot) \),
 \[
  \Expv \bigparen{h\prod_{l,j} g(\bar\beta_{l,i(l),j}(\eta),\bar W_{l,i(l),j}(\eta),\eta)}
	\le \Expv h\Expv\prod_{l,j} b(\bar\beta_{l,i(l),j}(\eta)).
 \]
 For any \( l,i \), the rectangles \( \bar W_{l,i,j}(\eta) \) defined in~\eqref{eq:bar-notation}
 are all disjoint, and belong to \( V_{l}(\eta^{*}) \).
  Also the sets \( V_{l}(\eta^{*}) \) are disjoint, therefore all rectangles
 \( \bar W_{l,i(l),j}(\eta) \) are disjoint.
 The Claim has shown that the system
\[
  \bigparen{J,\bar\tau_{l,i(l),j},\bar\beta_{l,i(l),j}(\eta), \bar W_{l,i(l),j}(\eta)}_{l\in N}
\]
is a system of disjoint constraint processes for \( M \),  and hence the proof is finished by
the strong trajectory property of \( \eta \).
\end{Proof}

From for ease of reading, we omit the qualifier ``strong'', but all trajectories we will deal
with are strong trajectories.

Here is a trivial but illustrative example.

\begin{example} Let \( M_{1}\) be the \( \eps \)-perturbation \( \CA_{\eps}(\trans,\B_{1},\T_{1}) \)
  of the cellular automaton with \( \states_{1}=\{0,1\} \) and
the trivial transition function  \( \trans(x,y,z)=0 \) for all \( x,y,z \).
Then \( \eta \) is a trajectory if for \( t>0 \), \( \eta(x,t) \) are 
random variables such that for any \( u \) and any finite set \( S \) of space-time points \( (x,t) \)
with \( t\ge u \), the conditional probability that
\( \eta(x,t)=1 \) in all \( (x,t)\in S \) no matter how the \( \eta(x',t') \) are fixed for all \( t'<t \),
is bounded by \( \eps^{|S|} \).
The set of constraint types consists of a single element, \( \bbT_{1}=\{\alpha_{1} \} \)
with \( b_{1}(\alpha_{1})=\eps \).
For the event function we have \( g_{1}(\alpha_{1},V,\eta)=1 \) if and only if
for some \( (x,t) \) of the form \( (m\B_{1}, n\T_{1}) \) with integers \( m,n \) and \( n>0 \),
\( V=\lint{x}{x+\B_{1}}\times\rint{t-T_{1}}{t} \), and \( \eta(x,t)=1 \).

Let  \( M_{2}\)  be  a medium with \( \states_{2}=\{0,1\} \).
For integers \( \Q_{1},\U_{1}>0 \) let \( \B_{2}=\Q_{1}\B_{1} \), \( \T_{2}=\U_{1}\T_{1} \).
The histories will be of the form \( \eta(x,t) \) with \( x=m\B_{2} \), \( t=n\T_{2} \) for \( n\ge 0 \).
The set of constraint types consists of a single element, \( \bbT_{2}=\{\alpha_{2} \} \)
with \( b_{2}(\alpha_{2})=3\Q_{1}\U_{1}\eps^{2} \).
For the event function we have \( g_{2}(\alpha_{2},V,\eta)=1 \)
if and only if for some \( (x,t) \) of the form \( (n\B_{1}, n\T_{1}) \)
with integers \( m,n \) and \( n>0 \),
\( V=\lint{x-\B_{2}}{x+2\B_{2}}\times\rint{t-2\T_{2}}{t} \),
and \( \eta(x,t)=1 \).

Let us define the simulation \( \fg^{*} \) as follows.
If \( \eta^{*}=\fg^{*}(\eta) \) then \( \eta^{*}(x,t)=1 \) if and only if
\( (x,t) \) has the form \( (m\B_{2}, n\T_{2}) \) with integers \( m,n \) and \( n>0 \),
and in \( V=\lint{x-\B_{2}}{x+2\B_{2}}\times\rint{t-2\T_{2}}{t} \)
there are at least two elements \( (x_{1},t_{1})\ne (x_{2},t_{2}) \)  with \( \eta(x_{i},t_{i})=1 \).
To show that this is a simulation we define the system of local conditions
\[
   \bigparen{\beta_{i,j}(\alpha_{2},V), W_{i,j}(\alpha,V)}
 \]
 as follows.
 These are nonempty if and only if for some \( (x,t) \)
 of the form \( (m\B_{2}, n\T_{2}) \) with integer \( m,n \) and \( n>0 \),
 \( V=\lint{x-\B_{2}}{x+2\B_{2}}\times\rint{t-2\T_{2}}{t} \).
 In this case, let \( P=\{p_{1},\dots,p_{N}\} \) be a list of all distinct
 pairs \( \{ (x_{1},t_{1}), (x_{2},t_{2})\} \) in \( V \), where \( N=\binom{6\Q_{1}\U_{1}}{2} \).
 For \( p_{i}=\{ (x_{i,1},t_{i,1}), (x_{i,2},t_{i,2})\} \),
 let \( \beta_{i,j}(\alpha_{2},V)=\alpha_{1} \) for all \( i,j \),  and 
\begin{align*}
 W_{i,j}(\alpha_{2},V)=\lint{x_{i,j}}{x_{i,j}+\B_{1}}\times\rint{t_{i,j}-\T_{1}}{\T_{i,j}}.
\end{align*}
Then \( \g_{2}(\alpha_{2},V,\eta^{*}) \) is the event that 
\( \eta(x,t)=1 \) for at least two pairs \( (x,t)\in V \),
and this is upper-bounded by 
\begin{align*}
 \sum_{i}\prod_{j}g_{1}(\alpha_{1},W_{i,j}(\alpha_{2},V))
\end{align*}
as required.
Also 
\begin{align*}
 b_{2}(\alpha_{2})=3\Q_{1}\U_{1}\eps^{2}>N\eps^{2}
=\sum_{i=1}^{N}\prod_{j=1}^{2}g_{1}(\alpha_{1},W_{i,j}(\alpha_{2},V))
\end{align*}
as required. 
\end{example}

\paragraph{Frequency of switching times}
  Adding certain conditions does not really change some media.

  \begin{sloppypar}
  \begin{definition}[Injective canonical simulation]\label{def:inj-canon-simul}
  We will call a canonical simulation
\df{injective}\index{simulation!canonical!injective} if the mapping
\( \fg^{*} \) is the identity, and the functions
\( b_{2}(\alpha),\g_{2}(\alpha,\cdot,\cdot) \) differ from \( b_{1},\g_{1} \) only in
that they are defined for some additional types \( \alpha \), too.    
  \end{definition}    
  \end{sloppypar}

Proving a certain probability bound for the behavior of \( \eta_{1} \) within
some window does not create a canonical simulation yet: the bound must be
formulated in such a way that this type of probability bound can be
multiplied with itself and all the other bounds over a system of disjoint
windows.

We will illustrate canonical simulations on the the following example.
The simulating medium \( M_{1} \) is the \( \CCA \) of Example~\ref{xmp:CCA-stop}.
The conditions to add impose lower and upper bounds on the frequency of
switching times of \( M_{1} \).
Recall the set of functions \( \bK \) in Definition~\ref{def:state-subsets}.

\begin{definition}[Subset-function switching times]
For some function \( K\in\bK \) let%
 \glo{a@\( a_{-1}(K), a_{1}(K) \)}%
 \begin{alignat*}{3}
        &a_{-1}(K) &&= \bigwedge_{\br}\bR(K(\br),\br),
\\   &a_{1}(K)   &&= \bigvee_{\br}\bR(K(\br),\br).
 \end{alignat*}
  Let us call a time \( u \) of the history \( \eta \) a
 \df{\( K \)-switching-time}\index{switch!time} of a cell \( x \) if for all
\( t<u \) sufficiently close to \( u \) we have 
\( \eta(x,t)\notin \ol K(x,t,\eta) \) and \( \eta(x,u)\in \ol K(x,t,\eta) \), 
that is we have just jumped into the target set determined by \( K(\cdot) \).
For some rational constant \( D>0 \), integers \( k\ge 0 \) and \( j=-1,1 \) let%
 \glo{a.greek@\( \alpha(K,D,k,j) \)}%
 \begin{equation}\label{eq:switch-count}
  \g\bigparen{\alpha(K,D,k,j),\;\{x\}\times\rint{u}{u + D},\;\eta}
 \end{equation}
  be the event function of type \( \alpha(K,D,k,j) \) for \( j=-1 \) [\( j=1 \)] for the
event that site \( x \) has at most [at least] \( k \) times during \( \rint{u}{u+D} \),
that are \( K \)-switching times.  
\end{definition}

The statement below estimates the probability of having a number of switches
significantly below the lower bound on their expected number, or above the upper bound.

\begin{proposition}\label{propo:time-bounds}
Let us define, with \( a_{j}=a_{j}(K) \) and every possible \( h>0 \):
\begin{align}\label{eq:1-switch-count}
 b(\alpha(K, D, a_{j}D+j h\sqrt{D}, j)) &= 1\land \frac{a_{j}}{h^{2}},
 \end{align}
 except that in case \( j = 1 \) and \( a_{1}D+h\sqrt D= 1 \) we set it to
 \(  a_{1}D \) when the latter is smaller.
Let \(  b(\alpha(K,D,k,j))=1 \) for all other \( k \).
For each transition rate matrix \( \bR \) and function \( K\in\bK \),
the trajectories of \( \CCA(\bR,\B ) \) are also trajectories of the medium obtained
by adding all these conditions defined by \( \alpha(K,D,\dots) \), \( b(\alpha(\dots)) \).
\end{proposition}

\begin{proof}
  The number of \( K \)-switching times during \( \rint{u}{u+D} \) is dominated by a
  Poisson random variable \( X_{1} \) with parameter \( a_{1}D \).
  We use approximation as in Definition~\ref{def:discr-approx}.
  Thus, for a sufficiently small \( \delta>0 \), let \( M_{\delta}=\PCA(\bP,\B ,\delta) \) with
  \( \bP(s,\br)=\bR(s,\br) \delta \) when \( s\ne r_{0} \) and
  \( 1-\delta\sum_{s'\ne r_{0}}\bR(s',\br) \) otherwise.
  If \( A(D,\delta) \) is the random variable counting the number of switches in \( M_{\delta} \)
  during \( \rint{u}{u+D} \) then it is clearly dominated by the random variable \( B (D,\delta) \)
  that is the number of switches when a switch takes place with probability \( a_{1}\delta\) every time.
  As \( \delta\to 0 \) the distribution of this variable \( B(D,\delta) \) converges to that of \( X_{1} \), which then
  dominates the number of \( K \)-switching times.
  Similarly, the number of \( K \)-switching times dominates
  a Poisson random variable \( X_{-1} \) with parameter \( a_{-1}D \).
    The mean and variance of \( X_{1} \) is \( a_{1}D \), hence by the Chebyshev
  inequality
  \begin{align*}
  \Pbof{X_{1}\ge a_{1}D + h \sqrt D}\le \frac{a_{1}D}{h^{2}D}= \frac{a_{1}}{h^{2}}.
  \end{align*}
  Similarly, 
\( \Pbof{X_{-1}\le a_{-1}D - h\sqrt D}\le \frac{a_{1}}{h^{2}} \).
  In the special case when \( a_{1}D+h\sqrt D = 1 \) we have the simpler bound
  \( \Pbof{X_{1}\ge 1} = 1-e^{-a_{1}D}< a_{1}D \).
 \end{proof}

\subsection{Primitive variable-period media} \label{sec:asynch}

The media defined here will
serve as the bottom level of our amplifier.

\begin{definition}[Dwell period lower bound]
  We say that the number \( \Tl>0 \)\glo{thm:cap@\( \Tl,\Tu \)} is a
 \df{dwell period lower bound}\index{period!dwell!lower bound} of a
history \( \eta \) if no dwell period of \( \eta \) is shorter
than \( \Tl \).  
\end{definition}

A continuous-time cellular automaton has no such lower bound.
In the amplifier \( M_{0},M_{1},\dots \) we will construct eventually,
all abstract
media \( M_{1}, M_{2},\dots \) will have a dwell period lower bound, only
the continuous-time probabilistic cellular automaton \( M_{0} \) will not have one.
\( M_{1} \) will be a so-called
 \df{primitive variable-period medium}\index{medium!variable-period!primitive}:
these can be considered the continuous-time extension of the notion of an
\( \eps \)-perturbation of a deterministic cellular automaton with
coin-tossing.
On the other hand, \( M_{2}, M_{3},\dots \) will only fit into a more
general framework (non-adjacent cells).

\begin{definition}[Primitive variable-period medium]
The medium \glo{primvar@\( \Prim-var \)}%
 \[
   \Prim-var(\trans, \B , \Tl, \Tu, \eps).
 \]
is defined as follows.
The set \( \states \) of states is implicit in the transition function therefore
from now on, it will be omitted.
We have dwell period lower and upper bounds \( \Tl\le\Tu \) and
a failure probability bound \( \eps>0 \).
The local state, as in Example~\ref{xmp:rand-gen}, is a record with two
fields, \( \Det \) and \( \Rand \) where \( \Rand \) consists of a single bit.
To simplify the upper-bounding of dwell periods, we assume that the
transition function has the property
 \begin{equation}\label{eq:prim-time-mark}
  \trans(r_{-1},r_{0},r_{1}).\Det\ne r_{0}.\Det,
\end{equation}
that is the deterministic part of a cell will change at every transition.
In a more general setting later this will be called ``time marking'', see
Condition~\ref{cond:time-mark}.
For a history \( \eta \), site \( x \) let \( a \ge 0 \) be a rational number.
If \( a>0 \) then
let \( \sigma_{1}=\sigma_{1}(x,a,\eta) \), and \( \sigma_{2} \) be the first and second
switching times \( t>a \) of \( \eta \) in \( x \).
If \( a=0 \) then let \( \sigma_{1}=0 \) and let \( \sigma_{2} \) be the first switching time
\( >0 \).
Let us list the different types of local condition.
 \begin{cjenum}

  \item\label{i:prim-var.dwell-lb}
  Conditions of type \( \alpha(\var{dw-p-lb}) \) imposing \( \Tl \) as a
dwell period lower bound.
  We have \(  b(\alpha(\var{dw-p-lb}))=\eps \), and 
 \[
   \g\bigparen{\alpha(\var{dw-p-lb}),\; \{x\}\times\rint{a}{a+\Tu},\; \eta}
 \]
  is the event that \( \eta(x,t) \) has two switching times
closer than \( \Tl \) to each other during \( \rint{a}{a+\Tl} \).

  \item\label{i:prim-var.dwell-ub}
  Conditions of type \( \alpha(\var{dw-p-ub}) \) imposing \( \Tu \) as a
dwell period upper bound.
  We have \(  b(\alpha(\var{dw-p-ub}))=\eps \), and 
 \[
   \g\bigparen{\alpha(\var{dw-p-ub}),\;\{x\}\times\rint{a}{a+\Tu},\; \eta}
 \]
  is the event that \( \eta(x,t) \) has no switching times
in \( \rint{a}{a+\Tu} \).

  \item\label{i:prim-var.trans}
  Conditions of type \( \alpha(\var{comput}) \) saying that the new value of the
\( \Det \) field at a switching time obeys the transition function.
  We have \(  b(\alpha(\var{comput})) = \eps \), and
 \begin{multline*}
 \g\bigparen{\alpha(\var{comput}),\; \{x\}\times\rint{a}{a+2\Tu},\; \eta}
 \\ = \{\eta(x,\sigma_{2}).\Det\not\in
  \setof{\ol\trans(\eta,x,t,\B ).\Det : t\in\rint{\sigma_{1}}{\sigma_{2}-\Tl/2}}\}.
 \end{multline*}
  Here we used the notation~\eqref{eq:trans-B}.

 \item\label{i:prim-var.coin}
  Conditions of type \( \alpha(\var{coin},j) \) for \( j=0,1 \) saying that the new
value of the \( \Rand \) field at a switching time is obtained nearly by a
fresh coin-toss:
 \begin{alignat*}{3}
     & b(\alpha(\var{coin},j)) &&= 0.5 + \eps,
\\ &\g\bigparen{\alpha(\var{coin},j),\; \{x\}\times\clint{a}{a+2\Tu},\; \eta} &&=
 \{\eta(x,\sigma_{2}).\Rand = j\}.
 \end{alignat*}
 When there is no \( \sigma_{2} \) then \( \g(\alpha(\var{coin},j),\cdot,\cdot)=0 \), so
this condition has no effect.

 \end{cjenum}
\end{definition}

Condition~\ref{i:prim-var.trans} says that unless an exceptional event
occurred, whenever a state transition takes place at the end of a dwell
period \( \lint{\sigma_{1}}{\sigma_{2}} \) it occurs according to a transition made on the
basis of the states of the three neighbor cells and the random bit at some
time in the \df{observation interval} \( \rint{\sigma_{1}}{\sigma_{2}-\Tl/2} \).
  Since the observation interval does not include the times too close to
\( \sigma_{2} \), information cannot propagate arbitrarily fast.

 \begin{example}
  All trajectories of the ordinary deterministic cellular automaton
\( \CA(\trans,\B ,\T ) \) are also strong trajectories of
\( \Prim-var(\trans,\B ,\T ,\T ,0) \).
 \end{example}

\begin{theorem}[Simulation by \( \CCA \)]\index{theorem@Theorem!\( \CCA \)-simulation}
 \label{thm:sim-by-cont-pca}
  For medium
 \[
  M=\Prim-var(\trans,1,\Tl,1,\eps)
 \]
  with~\eqref{eq:prim-time-mark} and \( \Tl<1 \) there is a noisy
transition rate \( \bR(s,\br) \) over some state space \( \states' \) and a function
\( \pi:\states'\to\states \) such that \( \pi(\eta(x,t)) \) is a strong trajectory of \( M \) for each
trajectory \( \eta \) of \( \CCA(\bR,1) \).
 \end{theorem}
\begin{proof}
  The proof for trajectories is a straightforward application of
  Proposition~\ref{propo:time-bounds}, and could be modified with some
  formalism to hold also for strong trajectories.
For an appropriate \( n \) to be chosen later, let
\( \states'=\states^3\times\{0,\dots,3n-1\} \).
  If \( s'=(\bs,i) \) where \( \bs=(s_{-1},s_{0},s_{1}) \) then let
\( \pi(s')=s_{0} \).
  Let us define the transition rates \( \bR_{n}(\cdot,\cdot) \) with
\( \sum_{s'}\bR_{n}(s',\br)=3n \).
  Let \( \bK=\bK_{n} \) be the set of functions as defined in~\eqref{eq:K-def}
belonging to this set of states, assuming that all transition rates are positive.
Choose the function \( K\in\bK \) as follows.
For \( \br\in(\states')^{3} \) where \( r_{0}=(\bq,i) \), where \( \bq=(q_{-1},q_{0},q_{1}) \),
let \( K(\br) \), be the set of those \( (\bar\bq,i+1) \) for which \( \bar\bq \) is
given as follows.
  If \( i\ne n-1,3n-1 \) then \( \bar\bq= \bq \).
  If \( i=n-1 \) then \( \bar q_{j}=\pi(r_{j}) \), \( j=-1,0,1 \).
  Suppose \( i=3n-1 \): then \( \bar q_{j}= q_{j} \) for \( j=-1,1 \) and
 \( \bar q_{0}.\Det=\trans(\bq) \).
  This completes the definition of \( K \) which is a one-element set whenever
\( i\ne 3n-1 \); in the latter case it is a two-element set, with the two elements
differing by the value of \( \bar q_{0}.\Rand \).

  Now \( \bR_{n} \) is given as follows.
\( \bR_{n}(s,\br)=0 \) for all \( s\ne K(\br) \) and
\( \bR(K(\br),\br)=3n \).
  This determines \( \bR_{n}(s,\br) \) for all cases when \( K(\br) \) is a
one-element set.
  If \( K(\br) \) is a two-element set then the transition rates into each of
these elements are equal.
  This completes the definition of \( \bR_{n} \).
  Thus, in a trajectory \( \eta \) of \( \bR_{n} \), each dwell period of
  \( \pi(\eta) \) consists of \( 3n \) small dwell periods of \( \eta \).
  There is a change only in  the \( n \)th and \( 3n \)th transition.
  In the \( n \)-th transition, cell \( x \) learns the current represented
state \( \pi(\eta(x+j,t)) \) of its two neighbors.
  In the \( 3n \)-th transition, it switches its own represented state
according to the transition probabilities corresponding to the triple
of \( \states \) states it currently knows about.
  The actual transition rates \( \bR \) differ from \( \bR_{n} \) in an arbitrary
way by some small amount \( <\eps' \).
  Let us show that we can choose first \( n \) and then \( \eps' \) in such a way
that the process \( \pi(\eta) \) becomes a trajectory of \( M \).

  Let \( K_{0}(\br)=\states\setminus(K(\br)\cup\{r_{0}\}) \).
  A trajectory \( \eta \) of \( \bR_{n} \), in any switching point \( (x,t) \) makes a
switch into the set \( K(x,t-,\eta) \).
If a trajectory \( \eta \) of \( \bR \) in a switching point \( (x,t) \)
does not switch into \( K(x,t-,\eta) \) we will this a \df{faulty switch}.
Then it switches into
\( K_{0}(x,t-,\eta) \), but with a rate smaller than \( \eps'|\states| \).
  According to the part of Proposition~\ref{propo:time-bounds}
referring to bounds of the kind~\eqref{eq:1-switch-count}, by making
\( \eps' \) small enough we can always achieve that the following condition
type \( \alpha(\var{fault}) \) can be added to the others:
 \[
   b(\alpha(\var{fault})) = c\eps,
 \]
  where \( c \) is any positive constant chosen in advance and
 \[
   \g\bigparen{\alpha(\var{fault}),\;\{x\}\times\rint{a}{a+2},\; \eta}
 \]
  is the event that a faulty switch occurs in \( x \) during \( \rint{a}{a+2} \).
  Therefore in all the local conditions for the simulated \( \Prim-var() \),
when bounding the probabilities in local conditions, we can confine
ourselves to the case that no faulty switch occurs in \( x \) during 
\( \rint{a}{a+2} \).

  Consider the conditions in~\ref{i:prim-var.dwell-lb}.
  Since we can exclude faulty switches in the time period \( \rint{a}{a+1} \), a
trajectory \( \eta(x,t) \) of \( \bR \) will be such that when denoting
 \( \eta(x,t)=(\br(x,t),i(x,t)) \), then \( i(x,t) \) increases by 1 modulo \( 3n \) in
every switch.
  Therefore there are exactly \( 3n \) switching times of \( \eta \) between any
two switching times of \( \pi(\eta) \).
  Let \( d=(1-\Tl)/2 \).
  If \( \pi(\eta) \) has two switching times \( \sigma_{1}<\sigma_{2} \) during
\( \rint{a}{a+1} \) that are closer then \( \Tl \) then either they are in \( \rint{a}{a+1-d} \) or
in \( \rint{a+d}{a+1} \).
  Consider the first case, the second case is similar.
  Then \( \eta \) has at least \( 3n \) switching times during \( \rint{a}{a+1-d} \).
  By choosing \( n \) large enough and noting that the rate of \( \bR \) is
(arbitrarily close to) \( 3n \), we can use Proposition~\ref{propo:time-bounds} to
bound the probability of this event by \( \eps/3 \).
  The probability for the other case can also bounded by \( \eps/3 \) and the
probability of a faulty switch can also bounded by \( \eps/3 \).
  We thus get the bound \( \eps \) in the conditions~\ref{i:prim-var.dwell-lb}
by canonical simulation.
A similar reasoning handles condition~\ref{i:prim-var.dwell-ub}.

  Now consider condition~\ref{i:prim-var.trans}.
  Again, we can assume that no faulty switch occurs during \( \rint{a}{a+2} \).
  There are exactly \( 3n \) switching times of \( \eta \) between \( \sigma_{1} \) and
\( \sigma_{2} \).
  According to the definition of \( \bR_{n} \), the transition will occur
according to \( \trans \), with the neighbor values that were read at the \( n \)th
switching time.
  Now an argument similar to the proof of~\ref{i:prim-var.dwell-lb} shows
that this \( n \)th switching time is in \( \rint{\sigma_{1}}{\sigma_{2}-\Tl/2} \).

  Finally, the possibility of adding condition~\ref{i:prim-var.coin} on
\( \pi(\eta) \) is a straightforward consequence of Theorem~\ref{thm:cont-stop}.
  Indeed, let \( K'(\br)=\states\setminus \{r_{0}\} \) and let \( K''(\br) \) be the set of
elements \( (\bq,0)\in K'(\br) \).
  Thus, jumping into \( K'' \) means making a switch in \( \pi(\eta) \).
  At that switch, the transition rates into the possibility \( q_{0}.\Rand=1 \)
and \( q_{0}.\Rand=0 \) are within \( \eps' \) of \( 1/2 \).
  Therefore Theorem~\ref{thm:cont-stop} is applicable.
 \end{proof}

\paragraph{Generalizing the results to continuous time}

  The first publication showing the possibility of reliable computation
with a continuous-time medium (in two dimensions) is~\cite{Wang90}.
  Here, we formulate the new results for variable-period
one-dimensional information storage and computation.
  The following theorems generalize Theorems~\ref{thm:1dim.nonerg},
\ref{thm:remember-inf} and~\ref{thm:1dim.comp} allowing a
variable-period medium \( \Prim-var(\trans',1,\Tl,1,\eps) \) with \( \Tl = 0.5 \)
in place of the perturbed discrete-time automaton \( \CA_{\eps}(\trans') \).
  See Corollary~\ref{crl:ipc} below for the version for interacting particle systems.

  \begin{theorem}\label{thm:1dim.nonerg.var}
    Theorem~\ref{thm:1dim.nonerg} also holds 
    when we allow a perturbed primitive variable-period medium
    \( \Prim-var(\trans',1, 0.5, 1,\eps) \) with strong trajectories
    instead of a  perturbed discrete-time cellular automaton.
  \end{theorem}

 \begin{theorem} \label{thm:1dim.stor.var}
   Theorem~\ref{thm:remember-inf}
also holds (with the same encoding \( \psi_{*} \)) and strong trajectories
when we allow a primitive variable-period perturbation instead of a discrete-time perturbation.
  \end{theorem}

  \begin{theorem} \label{thm:1dim.comp.var}
    Theorem~\ref{thm:1dim.comp.finite} holds (with the same encoding \( \psi_{*} \))
    also with \( \Prim-var(\trans',1, 0.5, 1,\eps) \) and strong trajectories
    in place of \( \CA_{\eps}(\trans') \).
  \end{theorem}

  Theorem~\ref{thm:sim-by-cont-pca} implies the following.

\begin{corollary}\label{crl:ipc}
  In Theorems~\ref{thm:1dim.nonerg.var}, \ref{thm:1dim.stor.var},
\ref{thm:1dim.comp.var}, we can replace \( \trans' \) with a rate matrix \( \bR \),
and \( \Prim-var(\trans',1,\Tl,1,\eps) \) with a \( \CCA \) with rates coming from
an arbitrary \( \eps \)-perturbation of \( \bR \).
  This proves, in particular, Theorem~\ref{thm:1dim.nonerg.cont}.
 \end{corollary}

 In what follows when we talk about trajectories we will always mean strong trajectories.

\section{Synchronization}
 \label{sec:sync}

The random nature of the switching times of a variable-period medium is a
tame kind of nondeterminism; any deterministic cellular automaton can be
simulated by a variable-period medium.
To prove this we first introduce an auxiliary concept.

\begin{definition}[Totally asynchronous cellular automata]
We define the \df{totally asynchronous cellular
automaton}\index{cellular automaton!totally asynchronous}%
 \glo{aca@$\ACA$}%
 \[
  \ACA(\trans)=\ACA(\trans,1,1)
 \]
  associated with transition function $\trans$ as follows: $\eta$ is a
trajectory if for all $x,t$ we have either $\eta(x,t+1)=\eta(x,t)$ or the
usual
 \[
   \eta(x,t+1)=\trans(\eta(x-1,t),\eta(x,t),\eta(x+1,t)).
 \]
 \end{definition}

To analyze synchronization, some more concepts are needed.

 \begin{definition}
    A site $x$ is \df{free}\index{site!free} in a configuration $\xi$ if
$\trans(\xi(x-1),\xi(x),\xi(x+1))\ne \xi(x)$.
The set of free sites will be denoted by $L(\xi)$\glo{lem:cap@$L(\xi)$}. 
For a space configuration $\xi$ and a set $E$ of sites, let us define the
new configuration $\trans(\xi,E)$\glo{tr@$\trans(\xi,E)$} by
 \[
  \trans(\xi,E)(x) = \begin{cases}
    \trans(\xi(x-1),\xi(x),\xi(x+1))    &\text{if $x\in E$}
\\  \xi(x)                              &\text{otherwise.}
  \end{cases}
 \]
\end{definition}

  Now we can express the condition that $\eta$ is a trajectory of
$\ACA(\trans)$ by saying that for every $t$ there is a set $U$ with
 \begin{equation}\label{eq:asynch.traj}
 \eta(\cdot,t+1)=\trans(\eta(\cdot,t),U).
 \end{equation}

 \begin{definition}[Update set]
  Let the \df{update set}\index{update!set}%
 \glo{u.cap@$U(t,\eta)$}%
 \begin{equation}\label{eq:asynch-sim.ute}
  U(t,\eta)
 \end{equation}
 be the set of sites $x$ with $\eta(x,t+1)\ne \eta(x,t)$.
    \end{definition}
  
The initial configuration and the update sets $U(t,\eta)$ determine
the trajectory $\eta$.

\begin{notation}[Indicator function]
For any set $A$, let us use the indicator function%
 \glo{h.greek@$\chi(x,A)$}%
 \[
   \chi(x,A)=\begin{cases}
   1 &\text{if $x\in A$,}
\\ 0 &\text{otherwise.}
             \end{cases}
 \]
\end{notation}

\begin{definition}[Effective age]
  For given history \( \eta \), we define the function
\( \tau(x,t)=\tau(x,t,\eta) \)\glo{thm:greek@\( \tau(x,t) \)} as follows:
 \begin{equation}\label{eq:tau-x-t}
 \begin{alignedat}{3}
       &\tau(x,0)   &&= 0,
\\    &\tau(x,t+1) &&= \tau(x,t)+ \chi(x,U(t,\eta)).
 \end{alignedat}
 \end{equation}
  We can call \( \tau(x,t) \) the \df{effective age}\index{age!effective} of
\( x \) in \( \eta \) at time \( t \): this is the number of effective updatings that
\( x \) underwent until time \( t \).
  \end{definition}

  \begin{definition}[Invariant histories]\label{def:invar-hist}
Given a transition function \( \trans \) and an initial configuration \( \xi \),
we say that the function \df{has invariant histories}\index{invariant
histories} if there is a function \( \zeta(x,u)=\zeta(x,u,\xi) \) such that for all
trajectories \( \eta(x,t) \) of \( \ACA(\trans) \) with \( \eta(\cdot,0)=\xi \) we have
 \begin{equation}\label{eq:t.commut}
   \eta(x,t) = \zeta(x,\tau(x,t,\eta),\xi).
 \end{equation}
    \end{definition}

  This means that after eliminating repetitions, the sequence
\( \zeta(x,1),\zeta(x,2),\dots \) of values that a site \( x \) will go through during
some trajectory, does not depend on the update sets, only on the initial
configuration (except that the sequence may be finite if there is only a
finite number of successful updates).
The update sets influence only the delays in going through this
sequence.
The following notation will be useful:

 \begin{notation}
Denote
 \[
   \trans(\xi,E,F)=\trans(\trans(\xi,E),F).
 \]
 \end{notation}

 \begin{definition}[Commutative transition]\label{def:commutative}
  We call a transition function \( \trans \)
\df{commutative}\index{transition function!commutative} if for all
configurations \( \xi \) and all distinct pairs \( x,y\in L(\xi) \) we have
\( \trans(\xi,\{x\},\{y\}) = \trans(\xi,\{y\},\{x\}) \).
    \end{definition}

The paper~\cite{GacsCommut95} proves the theorem that if a transition
function is commutative then it has invariant histories.
In Theorem~\ref{thm:asynch-sim} below, we will give a simple example
of a universal commutative transition function.
For that example, the theorem can be proved much easier.

\begin{notation}
We will denote the smallest absolute-value remainders%
 \glo{amod@\( \amod \)}%
 \begin{equation}\label{eq:amod}
  b \amod m
 \end{equation}
 with respect to a positive integer \( m>2 \), defined by the requirement
\( -m/2< b \amod m \le m/2 \).
  \end{notation}

\begin{theorem}[Commutative Simulation]\label{thm:asynch-sim}
 \index{theorem@Theorem!Asynchronous Simulation}
 Let \( \trans_{2} \) be an arbitrary transition function with state space
\( \states_{2} \).
Then there is a commutative transition function \( \trans_{1} \) with state
space \( \states_{1}=\states_{2}\times R \) (for an appropriate finite set \( R \)) with
the following property.
Each state \( s\in \states_{1} \) can be represented as \( (s.\F,s.\G) \) where
\( s.\F\in \states_{2} \), \( s.\G\in R \).
Let \( \xi_{2} \) be an arbitrary configuration of \( \states_{2} \) and let \( \xi_{1} \)
be a configuration of \( \states_{1} \) such that for all \( x \) we have
\( \xi_{1}(x).\F=\xi_{2}(x) \), \( \xi_{1}(x).\G=0\cdots 0\in R \).
Then for the trajectory \( \eta_{1} \) of \( \CA(\trans_{1}) \), with initial
configuration \( \xi_{1} \), the function \( \eta_{1}(x,t).\F \) is a trajectory of
\( \CA(\trans_{2}) \).
Moreover, in \( \eta_{1} \), the state of each cell changes in each step.
 \end{theorem}

In other words, the function \( \trans_{1} \) behaves in its field \( \F \) just
like the arbitrary transition function \( \trans_{2} \), but it also supports
asynchronous updating.

\begin{proof}
  Let \( \U >2 \) be a positive integer and%
 \glo{cur@\( \Cur \)}%
 \glo{prev@\( \Prev \)}%
 \[
  \Cur,\Prev,\Age
 \]
  be three fields of the states of \( \states_{1} \), where \( \F=\Cur \),
\( \G=(\Prev,\Age) \).
The field \( \Age \) represents numbers mod \( \U  \).
It will be used to keep track of the time of the simulated cells mod
\( \U  \), while \( \Prev \) holds the value of \( \Cur \) for the previous value of
\( \Age \).

Define \( s' = \trans_{1}(s_{-1},s_{0},s_{1}) \).
If there is a \( j\in\{-1,1\} \) with
\( (s_{j}.\Age-s_{0}.\Age) \amod \U  < 0 \) (that is some neighbor lags behind)
then \( s'=s_{0} \) that is there is no effect.
Otherwise, let \( r_{0}=s_{0}.\Cur \), and for \( j=-1,1 \), let \( r_{j} \) be equal
to \( s_{j}.\Cur \) if \( s_{j}.\Age=s_{0}.\Age \), and \( s_{j}.\Prev \) otherwise.
 \begin{alignat*}{3}
       &s'.\Cur  &&= \trans_{2}(r_{-1},r_{0},r_{1}),
 \\  &s'.\Prev &&= s_{0}.\Cur,
 \\  &s'.\Age  &&= s_{0}.\Age + 1 \bmod \U .
 \end{alignat*}
  Thus, we use the \( \Cur \) and \( \Prev \) fields of the neighbors
according to their meaning and update the three fields according to
their meaning.
It is easy to check that this transition function simulates \( \trans_{2} \)
in the \( \Cur \) field if we start it by putting 0 into all other fields.

Let us check that \( \trans_{1} \) is commutative.
If two neighbors \( x,x+1 \) are both allowed to update then neither of them
is behind the other modulo \( \U  \), hence they both have the same \( \Age \) field.
Suppose that \( x \) updates before \( x+1 \).
In this case, \( x \) will use the \( \Cur \) field of \( x+1 \) for
updating and put its own \( \Cur \) field into \( \Prev \).
Next, since now \( x \) is ``ahead'' according to \( \Age \), cell \( x+1 \)
will use the \( \Prev \) field of \( x \) for updating: this was the \( \Cur \)
field of \( x \) before.
Therefore the effect of consecutive updating is the same as that of
simultaneous updating.%
\end{proof}

\begin{definition}[Marching soldiers]
  The commutative medium of the above proof will be called the
 \df{marching soldiers}\index{marching soldiers} scheme.
  \end{definition}

The name comes from the similarity of its handling
of the \( \Age \) field to a chain of soldiers marching ahead in
which two neighbors do not want to be separated by more than one step. 
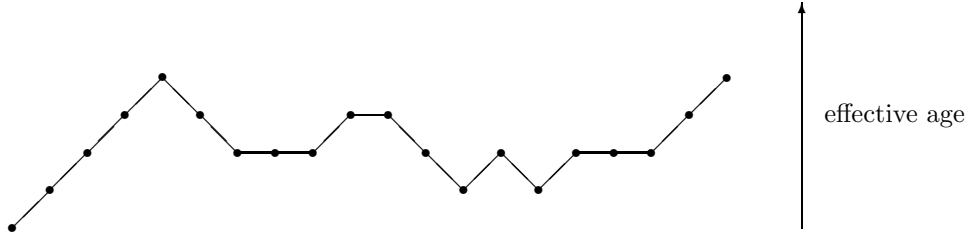
\begin{figure}
 \setlength{\unitlength}{1mm}
 \[
 \begin{picture}(100, 35)(15, 0)
    \put(0,0){\circle*{1}}
    \put(0,0){\line(1,1){5}}
    \put(5,5){\circle*{1}}
    \put(5,5){\line(1,1){5}}
    \put(10,10){\circle*{1}}
    \put(10,10){\line(1,1){5}}
    \put(15,15){\circle*{1}}
    \put(15,15){\line(1,1){5}}
    \put(20,20){\circle*{1}}
    \put(20,20){\line(1,-1){5}}
    \put(25,15){\circle*{1}}
    \put(25,15){\line(1,-1){5}}
    \put(30,10){\circle*{1}}
    \put(30,10){\line(1, 0){5}}
    \put(35,10){\circle*{1}}
    \put(35,10){\line(1, 0){5}}
    \put(40,10){\circle*{1}}
    \put(40,10){\line(1, 1){5}}
    \put(45,15){\circle*{1}}
    \put(45,15){\line(1, 0){5}}
    \put(50,15){\circle*{1}}
    \put(50,15){\line(1, -1){5}}
    \put(55,10){\circle*{1}}
    \put(55,10){\line(1, -1){5}}
    \put(60, 5){\circle*{1}}
    \put(60, 5){\line(1,  1){5}}
    \put(65,10){\circle*{1}}
    \put(65,10){\line(1, -1){5}}
    \put(70, 5){\circle*{1}}
    \put(70, 5){\line(1,  1){5}}
    \put(75, 10){\circle*{1}}
    \put(75, 10){\line(1, 0){5}}
    \put(80, 10){\circle*{1}}
    \put(80, 10){\line(1, 0){5}}
    \put(85, 10){\circle*{1}}
    \put(85, 10){\line(1, 1){5}}
    \put(90, 15){\circle*{1}}
    \put(90, 15){\line(1, 1){5}}
    \put(95, 20){\circle*{1}}
    
    \put(105, 0){\vector(0,1){30}}
    \put(108, 15){\makebox(0,0)[l]{effective age}}
 \end{picture}
 \]
 \caption[The Marching Soldiers scheme]{The Marching Soldiers scheme.
  The effective age of neighbor sites differs by at most 1.}
 \label{fig:marching}
\end{figure}

In typical cases of asynchronous computation, there are more efficient
ways to build a commutative transition function than to store the whole
previous state in the \( \Prev \) field.
Indeed, the transition function typically has a bandwidth (amount of communication
with a neighbor in one step) smaller than \( \nm{\states} \).

\begin{corollary}[Variable-period simulation]\label{crl:var-sim}
  For every deterministic transition function \( \trans_{2} \) over some
state-space \( \states_{2} \), there is a set of states \( \states_{1} \), a transition
function \( \trans_{1} \)
over \( \states_{1} \), and a code that for any values \( \Tl_{1}\le\Tu_{1} \),
is a simulation of \( \CA(\trans_{2}) \) by
\( \Prim-var(\trans_{1},1,\Tl_{1},\Tu_{1},0) \).
 \end{corollary}

\begin{proof}
  Let \( \trans_{1} \) be the commutative transition function given by
Theorem~\ref{thm:asynch-sim}, with the fields \( \F,\G \).
Let \( \xi_{2} \) be an arbitrary configuration of \( \states_{2} \) and let \( \xi_{1} \)
be a configuration of \( \states_{1} \) defined in the statement of the same
theorem. 
Let \( \eta_{1} \) be a trajectory of
\( \Prim-var(\trans_{1},1,\Tl_{1},\Tu_{1},0) \), with the starting configuration
\( \xi_{1} \).

An update set \( U(t,\eta_{1}) \) similar to~\eqref{eq:asynch-sim.ute} can be
defined now for the trajectory \( \eta_{1} \) as follows: \( x \) is in
\( U(t,\eta_{1}) \) iff \( t \) is a switching time of \( \eta_{1} \).
  Similarly, \( \tau(x,t,\eta_{1}) \) can be defined as in~\eqref{eq:tau-x-t}: 
 \begin{alignat*}{3}
        &\tau(x,0)   &&= 0,
\\    &\tau(x,t)   &&= \tau(x,t-)+ \chi(x,U(t,\eta_{1})).
 \end{alignat*}
 With these, let us define
 \begin{alignat*}{3}
      &\sigma(x,s,\xi)         &&= \bigwedge\setof{t:\tau(x,t,\eta')=s},
\\  &\eta_{2}(x,s,\xi) &&= \eta_{1}(x,\sigma(x,s)).\F
 \end{alignat*}
  By the cited theorem, \( \eta_{2} \) is a trajectory of \( \CA(\trans_{2}) \).%
\end{proof}

  The simulation in this theorem is not a local
one\index{simulation!non-local} in the sense defined in 
Section~\ref{sec:gen-sim} since it changes the time scale.
For an analysis of such simulations, see~\cite{BermSim88}.

 \section{Some simulations}\label{sec:sim}

In section, we build up the technique we will use in defining cellular automata
and simulations.

\subsection {Functions defined by programs}
 \label{sec:func-prog}

Recall the definition of a standard computing transition function
and medium as introduced in Section~\ref{sec:1dim.comp}.

\begin{definition}[Computation result]
For a standard computing medium with transition function \( \trans \), integers \( s \)
and \( t \) 
and string \( X \) consider a trajectory \( \eta \) of \( \CA(\trans) \) over the rectangle
\( \clint{0}{s}\times\clint{0}{t} \) with an initial configuration in which
\( \eta(0,0)=\eta(s,0)=\#\cdots \# \), further \( X \) is on the input track on the interval
\( \clint{1}{s-1} \) (padded with \( * \)'s to fill up \( \clint{1}{s-1} \)), \( * \)'s on the output
track and 0's on the other tracks.
This defines a trajectory \( \eta \) since the \( \# \)'s on the input field
in cells 0 and \( s \) imply that the cells outside the interval \( \clint{0}{s} \)
will have no effect.
Assume that, at time \( t \), the \( \Output \) track has no \( * \) on the interval
\( \clint{1}{s-1} \).
The string \( w \) on the \( \Output \) track on \( \clint{1}{s-1} \) will be called the
\df{result}\index{computation!result} of the computation, denoted by%
 \glo{tr@\( \trans(X; s, t) \)}%
 \[
  w= \trans(X; s, t).
 \]
 \end{definition}

Since the output is monotonic in standard computing media, the result will not
change at later times.

\begin{figure}
  \setlength{\unitlength}{1mm}
  \[
   \begin{picture}(100, 110)(0,-10)
  \put(40, 0){\framebox(60, 80){\( \nearrow\uparrow\nwarrow \trans \)}}
  \put(45,0){\line(0, 1){80}}
  \put(95,0){\line(0, 1){80}}

  \put(45,0){\makebox(50, 5){\( X \) on \( \var{Input} \)}}
  \put(45,5){\line(1, 0){50}}

  \put(45, 75){\makebox(50, 5){\( \trans(X; s, t) \) on \( \var{Output} \)}}
  \put(45, 75){\line(1, 0){50}}

  \multiput(40, 0)(0, 5){3}{\makebox(5, 5){\#}}
  \multiput(40, 75)(0, -5){3}{\makebox(5, 5){\#}}

  \multiput(95,0)(0, 5){3}{\makebox(5, 5){\#}}
  \multiput(95, 75)(0, -5){3}{\makebox(5, 5){\#}}

  \multiput(40, 40)(55, 0){2}{\makebox(5, 5){\( \vdots \)}}

  \put(42.5, -4){\makebox(0,0)[t]{0}}
  \put(97.5, -4){\makebox(0,0)[t]{\( s \)}}

  \put(36, 2.5){\makebox(0,0)[r]{0}}
  \put(36, 77.5){\makebox(0,0)[r]{\( t \)}}

  \put(0, 30){\vector(0, 1){20}}
  \put(4, 40){\makebox(0,0)[l]{time}}

   \end{picture}
 \]
 \caption{Definition of \( \text{\rmfamily\slshape\mdseries trans}(X; s, t) \)}
 \label{fig:def-of-tr-X}
\end{figure}

Universal cellular automata were introduced in Definition~\ref{def:univ}.
We will need, however, also some bounds on the time and space used
by the universal simulation.

\begin{definition}[Efficient universality]\label{def:effic-univ}
The standard computing transition function \( \trans_{0} \) is
 \df{efficiently universal\/}\index{transition function!universal!efficiently}
if for every standard computing transition function \( \trans \), there is a
string \( \rul{prog} \) and constants \( c,d \) such that for all strings \( X \) and
positive integers \( s,t \) we have%
 \glo{prog@\( \rul{prog} \)}%
 \[
     \trans(X; s, t) = \trans_{0}(\rul{prog}\cc X; c\cdot s, d\cdot t)
 \]
  whenever the left-hand side is defined.
  \end{definition}

In other words, \( \trans_{0} \) can simulate the computation of any other
standard computing transition \( \trans \) if we prepend a ``program of
\( \trans \)''\index{transition function!program of}.
The space and time needed for this simulation can only be a constant
factor greater than those of the simulated computation.

 \begin{theorem}
 There are efficiently universal standard computing transition functions.  
 \end{theorem}

\begin{proof}[Sketch of the proof]
This theorem can be proven similarly to Theorem~\ref{thm:univ}.
The main addition to the construction is that before \( \trans_{0} \) starts
a block simulation of \( \trans \) the input will be distributed bit-by-bit to
the colonies simulating the cells of \( \trans \).
At the end, the output will be collected from these colonies.%
\end{proof}

This theorem allows some standardization.

\begin{definition}
  We fix an efficiently universal standard computing transition function%
 \glo{univ@\( \univ \)}%
 \begin{equation}\label{eq:univ}
   \univ
 \end{equation}
  for the definition of various string functions.
In view of Theorem~\ref{thm:asynch-sim}, we can (and will) also assume that
\( \univ \) is commutative as in Definition~\ref{def:commutative}.
For a string function \( f(X) \) defined over some domain \( E \) we will say
that \( \rul{prog} \) is a program for \( f \) with
\df{time complexity bound} \( t \) and \df{space complexity bound \( s \) over \( E \)} if
the relation
 \[
  \univ(\rul{prog}\cc X;s,t) = f(X)
 \]
holds.
  \end{definition}

\subsection{The rule language}\label{sec:rule-lang}

It is often convenient to define a finite function (like a
transition function of some medium, or a code) by its program (for
\( \univ \)) rather than its transition table, since for many interesting
functions, the program can be given in a more concise form.
This subsection defines a language for describing a transition function.
Our purpose is twofold.
 \begin{description}

 \item[Convenience]  We want a method that is 
more convenient than giving a complete transition table.

  \item[Self-simulation]
Our self-correcting cellular automaton will simulate another self-correcting
cellular automaton, identical to itself, or just similar to it: for simplicity,
we refer to this now as ``self-simulation''.  
Self-simulation in our case has some peculiar requirements.

Let us note that self-simulation by itself is not mysterious.
A universal cellular automaton \( U \) would use the program \( p \) of any cellular
automaton with transition function \( \trans_{p} \) to simulate it.
Of course, \( U \) itself has a program, \( q \), so it could simulate itself as well,
if we give it the program \( q \).
But it seems redundant to give \( U \) its own program, and actually it is
better not to, since the written program is exposed to errors.
How to define conveniently
a machine that simulates (a modified version of) itself without needing any program?
There are various ways of doing it, but each involves a notion of a program that
is not a simple transition table, but rather has some ability to refer to itself: see
the remark~\ref{rem:self-simul} below.
 \end{description}

\begin{remark}[Self-simulation]\label{rem:self-simul}
Let us denote the output  of a universal machine (it does not matter, whether it
is a Turing machine or cellular automaton) computing the universal  partial
recursive function \( u(p,x) \), where
\( p \) is the program, and \( x \) is the input.
Let \( \tau \) be an arbitrary computable transformation of strings into strings.
Kleene's Fixpoint Theorem (Recursion Theorem) says that there is a program \( q \)
with the property that for all \( x \) we have \( u(q,x)=u(\tau(p),x) \).
What does this theorem have to do with our problem of self-simulation?
Suppose that we have some program \( p \), and a transformation \( \tau(p) \) of \( p \)
that carries out a simulation of the action of \( p \), with possibly some
modifications and extra properties (say, error-corrections, and the change of
some parameters).
The theorem says that there is a program \( q \) which already acts as \( \tau(q) \).

One  intuitive proof of the fixpoint theorem goes as follows.
Assume that we have a programming language that has at least the  following
ingredients: 1) We can define functions (procedures) in it.
2) It contains a function \( \mathtt{U(p,x)} \) that computes \( u(p,x) \), that is it interprets
the string \( p \) as a program on input \( x \).
3) It contains a function \( \mathtt{MyText()} \),\footnote{The parentheses
\( \mathtt{()} \), used for example in the programming language C, 
are a useful reminder that even when it has no explicit arguments, the function
depends on the internal state of the machine at the time of execution.} 
which returns the string that is the text of the program.
One implementation of such a function in a modern computer uses the fact that
the program is stored at some standard place in memory: it is enough to go there
and read it out.
Let \( \mathtt{T(p)} \) be the expression of the transformation \( \tau(p) \) in our
programming language: then the fixedpoint program is
 \begin{align*}
  q=\mathtt{U(T(MyText()),x)}
 \end{align*}
(the \( \mathtt{x} \) in the program is just a symbol in the program, not something
the program \( q \) depends on).
\end{remark}

A transition table can be seen as a sequence of \df{rules}: each rule could have
the following form, for all possible triples \( (a,b,c)\in\states^{3} \):

 \begin{algo}[H]
 \If{\( r_{-1}=a \) \algAnd \( r_{0}=b \) \algAnd \( r_{1}=c \)}{
       \( s \gets d \)} 
\end{algo}

where \( d=\trans(a,b,c) \).
The most important generalization of this sort of rule involves fields.
All the fields of our state will be defined in advance---as subintervals of the
interval \( \lint{0}{\nm{\states}} \).
Also, in our rule description, 
it is convenient to refer to the current cell and its neighor
in a standardized way.
Our rules, being just a description of the transition function,
will always apply to a particular cell.

\begin{notation}\label{n:neighbors}
For a field \( \F \), if \( \eta \) is the current configuration and \( x \) is a site
then we will sometimes write
 \begin{align*}
   \F(x) \text{ for } \eta.\F(x).
 \end{align*}
  We will write%
 \glo{th.greek@\( \thg_{j}(x) \)}%
 \[
 \thg_{j}(x) = x + j\B 
 \]
 for the site \( j \) steps from site \( x \).
In a semi-formal description of a rule, \( \x \)\glo{x@\( \x \)} refers to the
cell to which it is applied (the ``current'' cell).
In the condition as well as the action of the rule, a field \( \F(\x)=\eta.\F(\x) \)
will simply be written as \( \F \).
Denote%
 \glo{f.cap@\( \F^{j} \)}%
 \[
  \F^{j}=\eta.\F(\thg_{j}(\x)).
 \]
\end{notation}

The notation \( \x \) can always be replaced with references
to fields of the current cell and its neighbors.
For example, for field \( \Addr \), instead of writing
 ``\textbf{if} \( \eta.\Addr(\x)=0 \) \textbf{and} \( \eta.\Addr(\thg_{j}(\x))=1 \)'' we will write
 ``\textbf{if} \( \Addr=0 \) \textbf{and} \( \Addr^{j}=1 \)''.
Here is an example rule of this simplest form:

\begin{algo}[H]
 \If{\( \Age>0 \) \algAnd \( \Addr<100 \)}{
       \( \Mail \gets\Mail^{-1} \) }  
\end{algo}

This means that if (at the site \( \x \) under consideration)
the \( \Age \) field's value is 0 and the \( \Addr \) field's value is less than 100
then the \( \Mail \) field of the left neighbor \( \x-\B  \) 
should be copied into the \( \Mail \) field of cell \( \x \).
We allowed in~\eqref{eq:F[i:j]} the possibility of referring to a field by simply
a pair of numbers \( i,j \).
In the rule language, \( i,j \) could themselves be specified by
some (fixed) fields.
Here is a more general definition of rules:

\begin{definition}[Rules]\label{def:rule}
An \df{elementary rule} has the following form:

\begin{algo}[H]
 \If{\( C_{1} \) \algAnd \( \dotsm \) \algAnd \( C_{k} \)}{
       \( \F \gets f(\F_{1}^{j_{1}},\dots,\F_{n}^{j_{n}} \))
  }  
\end{algo}

\begin{sloppypar}
  Here \( \F,\F_{i} \) are fields (possibly specified as bit segments like in~\eqref{eq:F[i:j]}),
 \( j_{i}\in\set{-1,0,1} \).
Each of the conditions \( C_{i} \) has the form
\( a_{1}\F_{1}^{j_{1}}+\dots+a_{n}\F_{n}^{j_{n}}\le b \)
with integer constants \( a_{i}, b \).
(When the condition is not applicable, the rule does not change the state.)
The function \( f(x_{1},\dots,x_{n})\) must be computable on \( \univ \)
in linear time.
This allows for many functions: for example, if \( x ,y \) are integers in binary notation then
\( x\cdot y \), \( \flo{x/y} \), \( x\bmod y \) are all computable in linear time by known
cellular automata algorithms, see for example~\cite{Leighton92}.
\end{sloppypar}

A \df{rule} is a sequence of elementary rules.
The transition function will be defined by a set of rules.
In case some rules contradict each other, they will be ordered by explicit \df{priority}.
Rules will sometimes be associated with certain fields.
The lowest-priority rule referring to a field is called the
\df{default}\index{rule!default}.
\end{definition}

We will introduce several shorthand notations for expressing rules:
these can always be translated into a rule according to
Definition~\ref{def:rule}.

The most important shorthand is a branching conditional:

 \begin{algo}[H]
 \If{\( C_{1} \)}{
    \( A_{1} \)}
 \ElseIf{\( C_{2} \)}{
    \( A_{2} \)}
 \ElseIf{\( C_{3} \)}{
    \( A_{3} \)}
 \Else{\( A_{4} \)}
 \end{algo}

We can combine several rules.

\begin{definition}[Combination of rules]\label{def:combination}
When a sequence of rules is combined, this is a \df{parallel}
combination: all the rules apply to the same transition,
and do not refer to consecutive times.
Sometimes we may use the notation
 \glo{"+@\( \texttt{;},\prl \)}%
 \[
  \rul{R}_{1}\prl \rul{R}_{2} \prl \dots \prl \rul{R}_{n}
 \]
  for a combination of \( n \) rules.
As a shorthand, we may also combine rules in series:
the rule
 \[
  \rul{R}_{1};\rul{R}_{2}
 \]
  asks for carrying out the rules \( \rul{R}_{1} \) and \( \rul{R}_{2} \)
consecutively.
This can be understood in terms of a field \( \Age \) which
we are always going to have and can be replaced with a conditional.

In both parallel and serial combination,  we will refer to 
the rules \( \rul{R}_{i} \) as \df{sub-rules}, when they
will always be invoked (``called'') by other rules.
Still, for a better flow of discussion we will frequently refer to sub-rules
as just ``rules''.
Sub-rules may have \df{parameters} that are specified at the time of the call.
Rule \( P(i) \) with parameter \( i \) can be viewed as a different rule for
each possible value of \( i \).
  \end{definition}

For example, whenever we write

 \begin{algo}[H]
   \( \Retrieve \);\;
   \( \Eval \)
 \end{algo}

 then this can be translated, using appropriate constants \( t_{i} \), into
the following, where \( \Retrieve \) and \( \Eval \) are sub-rules:

 \begin{algo}[H]
   \If{\( t_{1}\le\Age< t_{2} \)}{
       \( \Retrieve \)}
   \ElseIf{\( t_{2}\le\Age< t_{3} \)}{
       \( \Eval \) \cmt{where \( \Eval \) was defined in (6)}
   }
 \end{algo}

This example also shows a formal way to add comments to a rule definition,
on any line after a \( \# \).
An example for a simple parameterized sub-rule would be 
\( \rul{Read-Mail}(j) \) for \( j\in\{-1,1\} \):
 \[
 \Mail\gets\Mail^{j}
 \]
 for the field \( \Mail \).  

 \begin{definition}[Conditionals and various kinds of ``for'']
Sometimes the conditions of a conditional statement
are indexed in a regular way, then we may use the shorthand 
\( \textbf{cfor} \)\glo{cfor@\( \textbf{cfor} \)} as below:

\begin{algo}
    \cFor{\( j\in\{-1,1\} \)}{
      \If{\( \Age^{j}=0 \)}{
               \( \Age\gets 0 \)
    }
  }
\end{algo}
So this rule tests the conditions in order,
carries out the command of the first satisfied condition and skips the rest.
Several conditions that are to be tested one after the other
will be combined using the construct \textbf{cond}.
This can be seen as an alternative to the
``\textbf{if}\( \dots \)\textbf{then}\( \dots \)\textbf{elsif}\( \dots \)''
construct.

\begin{algo}[H]
\Cond{
    \cFor{\( j\in\{-1,1\} \)}{
      \If{\( \Age^{j}=0 \)}{
               \( \Age\gets 0 \)
       }
    }
    \cFor{\( j\in\{-1,1\} \)}{
      \If{\( \Addr^{j}\not\equiv\Addr+j\pmod{\U } \)}{
               \( \Kind\gets\Latent \)
       }
     }
 }
\end{algo}

Both the parallel and the consecutive combination of rules can be over an
indexed set, and then we get the constructs \textbf{pfor} and \textbf{for}:

 \begin{algo}[H]
   \pFor{\( (k,d)\in\var{Mail-ind} \)}{
      \( \Mail_{k,d}\gets \var{Mail-to-receive}(k,d) \)
     }
   \end{algo}
   
Here, \glo{pfor@\( \texttt{pfor} \)} 
the rule is carried out simultaneously for all fields
indexed by the indicated values \( k \).
\glo{for@\textbf{for}}%
On the other hand in the \( \textbf{for} \) example  the application is consecutive:

  \begin{algo}[H]
   \For{\( i=a \) \KwTo \( b \)}{
      \( \ang{\text{some rule referring to } i} \)
   }
  \end{algo}

Another use of \( \textbf{for} \) can be this:

\begin{algo}[H]
   \For{\( n \) steps of \( \Age \)}{
      \( \ang{\cdots} \)
   }
 \end{algo}
 which, started at some \( \Age=t \) will mean that the commands \( \ang{\cdots} \) will
 run under the condition \( t_{1}\le\Age<t+n \).
  \glo{repeat@\textbf{repeat}}%
The \textbf{repeat} construct could be translated into a \textbf{for}:

  \begin{algo}[H]
  \pRepeat{\( k \)}{
     \( \ang{\cdots} \)
   }
 \end{algo}
 \end{definition}

We will also have \df{functions}: these are defined by the
information available in the arguments of the transition function, and can
always be considered a shorthand notation.

\begin{definition}[The ``let'' construct]
We may use temporary constants or functions in the course of defining a
rule, as in this example:

  \glo{let@\texttt{let}}%
 \begin{algo}
  \algLet \( f(i)\gets\Addr+\Age^{i} \)
\end{algo}

\end{definition}

\begin{example}
 \glo{march@\( \March_{0} \)}%
  Here is a description of the transition function given in the proof
of Theorem~\ref{thm:asynch-sim} (Asynchronous Simulation).%

 \begin{algo}\caption{rule \( \protect\March_{0} \)}\label{alg:March0}
   \pFor{\( j\in\{-1,0,1\} \)}{
      \algLet \( r(j)\gets\Cur^{j} \) if  \( \Age^{j}=\Age \), and \( \Prev^{j} \) otherwise
    }
    \If{\( \foralls{j\in\{-1,1\}}(\Age^{j}-\Age) \amod \U  \ge 0 \)}{
       \( \Prev \gets  \Cur \)\;
       \( \Cur \gets  \trans_{2}(r(-1),r(0),r(1)) \)\;
       \( \Age \gets  \Age + 1 \bmod \U  \)
       }
 \end{algo}
\end{example}

Some primitives will always be available in the rule language.
For example, for a field \( \F \), its width \( |\F| \) can be used.
Moreover, a \df{slice} of the field as a string of bits \( \F=(f_{0},\dots,f_{n-1}) \)
is available as well as \( \F[i\frTo j] \) as indicated in~\eqref{eq:F[i:j]}.
We will always have a field \( \Addr \).
As before, a \df{colony}\index{colony} is a segment of cells with addresses growing
from 0 to \( \Q -1 \) modulo some constant \( \Q  \).
Continuing to develop shorthands,
it is useful to have notation for talking about parts of a colony:

\begin{definition}[Location]
A \df{location}\index{location} is defined by a pair \( (\F,I) \) where \( \F \) is a field,
and \( I \) is an interval of addresses in the colony, or by a constant-size union of such pairs
\( (\F_{j},I_{j}) \).
  We will denote the  track \( \F \) over interval \( I \) as%
 \glo{f.cap@\( \F(I) \)}%
 \[
   \F(I).
 \]
The name of a location will generally begin with an underscore, like
 \[
   \lInfo.
 \]
The set of addresses (typically an interval) of a location \( L \) will be denoted by \( \Space(L) \).
  \end{definition}

The example in the above definition is important enough to be fixed:

  \begin{definition}[Info location]
The location on the \( \Info \) track containing the represented string is denoted by
\( \lInfo \).\glo{info.loc@\( \lInfo \)}
  \end{definition}

As an element of the language, of course, a location is simply a triple
consisting of the field, and the numbers determining the endpoints of the
interval.
A location is meant to specify a sequence of symbols on track \( \F \).

  \begin{remark}
Occasionally, we may treat the union of two or three locations as one.
It is easy to generalize the rules dealing with locations to this case.
  \end{remark}
  
Let us be given a string \( S \) consisting of a constant number of runs
of the same symbol.
For example, \( 0^{m}1^{n} \) has one run of 0's and a run of 1's.
Let us also be given a location \( \loc \).
Then \( \Write(S,\loc) \)\glo{write@\( \Write \)}, writing the string \( S \)
to the location \( \loc \), can be written as a conditional sub-rule:

 \begin{algo}\caption{sub-rule \( \protect
\Write(0^{m}1^{n},\protect\F(\lint{a}{a+m+n})) \)}\label{alg:Write}
  \If{\( a\le\Addr < \Q \wedge(a+m) \)}{
    \( \F\gets 0 \)}
  \ElseIf{\( a+m\le\Addr < \Q \wedge(a+m+n) \)}{
    \( \F\gets 1 \)
     }
\end{algo}

\begin{definition}[Name definitions]\label{def:Param}
The rule language can contain some definitions of names for constant
strings of symbols, of the form%
 \glo{param@\( \Param \)}%
 \begin{equation}\label{eq:Param}
 \Param_{1}=s_{1},\;\Param_{2}=s_{2},\dots
 \end{equation}
  where \( s_{i} \) are some strings.
The names \( \Param_{i} \) are allowed to be used in the rules.
This has the advantage that even if we use, say, \( \Param_{1} \) twice in a rule,
this does not make the rule much longer, even if the value of the 
parameter is a long string.

Parameters will have descriptive names: say, we might write 
\( \Height \) for \( \Param_{1} \),
when \( \Param_{1} \) indicates the level in a hierarchy of simulations.
Parameters will always be used via the function 
 \begin{align*}
   \WriteParam(\Param,\loc),
 \end{align*}
which writes the string value of parameter \( \Param \) to location \( \loc \).
This can be implemented as a subrule, via the primitive access function 
 \begin{align*}
   \Param_{i}(j),
 \end{align*}
referring to the \( j \)th symbol of parameter \( i \).
There will also be a special parameter: 
 \begin{align*}
\MyRules,
 \end{align*}
whose value is the string that is the sequence of all the rules.
This will implement self-reference in the style of Remark~\ref{rem:self-simul}.
\end{definition}

 \begin{theorem}[Rule Language]\label{thm:rule-lang}
 \index{theorem@Theorem!Rule Language}
  There is a string \( \Interpr \)\glo{interpr@\( \Interpr \)} and an
integer \( \cns{interpr-coe} \)\glo{interpr@\( \cns{interpr-coe} \)} such that the
following holds.
If string \( P \) is a description of a transition rule \( \trans \) over state
set \( \states \) in the above language (along with the necessary parameters), then
the machine \( \univ \) defined in~\eqref{eq:univ} computes
\( \trans(r_{-1},r_{0},r_{1}) \) (possibly padded with \( * \)'s) from
 \[
   \Interpr\cc P\cc r_{-1}\cc r_{0}\cc r_{1}
 \]
 within computation time \( \cns{interpr-coe}(|P|+1)^{2}\nm{\states} \) and space
\( \cns{interpr-coe}(|P|+\nm{\states}+1) \).
 \end{theorem}

A detailed proof would be tedious, but routine.
Essentially, each line of a rule program is interpreted in linear time on inputs
that are some fields: substrings of a state argument \( r_{i} \).
We have \( (|P|+1) \) squared in the time upper bound since we may have to
look up some parameter repeatedly.
From now on, by a \df{rule program}\index{program!rule}
\( \rul{Trans-prog} \)\glo{trans-prog@\( \rul{Trans-prog} \)} of the transition
function \( \trans \), we will understand some string to be interpreted by
\( \Interpr \).

Later, in Section~\ref{sec:comp-rules}, we will add some other, simple features to the
\( \Interpr \) string: interpreting certain special commands.

\subsection {A basic block simulation}
 \label{sec:simp-sim}

The simulation described just demonstrates the use of the
notation and introduces some elements of the later construction in a
simple setting.
Let a transition function \( \trans_{2} \) be given.
We want to define a cellular automaton \( M_{1}=\CA(\trans_{1}) \), whose
trajectories simulate the trajectories of \( M_{2}=\CA(\trans_{2},\Q ,\U ) \),
with appropriate \( \Q ,\U  \).
Of course, there is a trivial simulation, when \( M_{1}=M_{2} \), but a 
more general scheme will be set up here.
This simulation is not one of the typical simulations by a universal
medium: the cell-size of \( M_{1} \) depends on \( M_{2} \) as in 
Example~\ref{xmp:block-code}.
The construction will be summarized in Theorem~\ref{thm:simp-sim}
below.

Along the way, we introduce some more shorthand notation
in writing the rules that can also be incorporated
into the rule language---without invalidating Theorem~\ref{thm:rule-lang} (Rule
Language).

\paragraph{Overall structure}
  The transition function \( \trans_{2}:\states_{2}^3\to\states_{2} \) to be simulated is
given by a rule program \( \rul{Trans-prog}_{2} \).
To perform a simulated state transition of \( M_{2} \), a colony of \( M_{1} \)
must do the following:
 \begin{description}
  \item[Retrieve] Retrieve the states of the represented neighbor
cells from the neighbor colonies;
  \item[Evaluate] Compute the new state using \( \trans_{2} \);
  \item[Update] Replace the old represented state with the new one.
 \end{description}

The field \( \Addr \) holds a number between 0 and \( \Q -1 \), as discussed
in Section~\ref{sec:col}.
The default operation is to keep this field unchanged.
The time steps within a work period of a colony are numbered
consecutively from 0 to \( \U -1 \).
The field \( \Age \) holds a number between 0 and \( \U -1 \) intended to be
equal to this number.
The default operation is to increase this by 1 modulo \( \U  \).
These default operations will not be overridden in the simple,
fault-free simulations.
Using these fields, each cell knows its role at the current stage of
computation.

On the track \( \Info \), each colony holds a binary string of
length \( \nm{\states_{2}} \).
For a string \( S\in\states_{1}^{\Q } \), the decoded value \( \fg^{*}(S) \) is obtained by
taking this binary string.  The encoding will be defined later.
The default operation on the information field is to leave it
unchanged.
It will be overridden only in the last, updating step.
For simplicity, let us take \( |\Info|=2 \), that is the \( \Info \) track contains
only symbols from the standard alphabet, introduced in
Definition~\ref{def:std-alphabet}. 
The field \( \Cpt \)\glo{cpt@\( \Cpt \)} will be used much of the time like the
cells of a standard computing medium, so its size is the capacity
\( \nm{\states_{\univ}} \) of the fixed efficiently universal standard
computing medium.
It has subfields \( \Input \), \( \Output \).
The field \( \Cpt.\Input \) is under the control of the rule
\( \Retrieve \), \glo{retrieve@\( \Retrieve \)}%
while the rest of \( \Cpt \) is under the
control of the rule \( \Eval \)\glo{eval@\( \Eval \)}.
 \glo{update@\( \rul{Update} \)}%
  The whole program can be summarized as follows:

 \begin{algo}\caption{Basic simulation program}\label{alg:basic-sim-program}
   \( \Retrieve \);\;
   \( \Eval \);\;
   \( \rul{Update} \)
\end{algo}

\paragraph{Mail movement}
Let \( x \) be the base of the current colony under consideration.
For \( m=-1,0,1 \), we will indicate how to write a rule%
 \glo{copy@\( \Copy \)}%
 \[
 \Copy(m,\loc_{1},\loc_{2})
 \]
 that copies, from the colony with base \( x-m\Q  \), the location
\( \loc_{1} \)\glo{loc@\( \loc \)} to location \( \loc_{2} \) of the current colony.
To describe the rule \( \Copy \) we use a framework a little more general than what
would be needed here, with a variable-time version in mind.

\glo{retrieved@\( \lRetrieved_{m} \)}
 \begin{definition}[Locations for the retrieval rule]
   The rule \( \Retrieve \) uses the locations \( \lRetrieved_{m} \) to deposit the
string it retrieves from the neighbor colony number \( m\in\{-1,0,1\} \).
 \end{definition}

 With the help of the copy rule, we can formulate retrieval now as follows:

 \begin{algo}[H]\caption{rule \( \protect\Retrieve \)}\label{alg:Retrieve-1}
  \For{\( m\in\{-1,0,1\} \)}{
     \( \Copy(m, \lInfo, \lRetrieved_{m}) \)
    }
  \end{algo}

  There will be several mail tracks, indexed by the following values:
  \begin{align}\label{eq:Mail-ind}
    \var{Mail-ind} = \setOf{(k,d)}{k\in\{-1,0,1\}, d\in\{-1,1\}}.
 \end{align}
 The meaning of the two values \( -1,1 \) of \( d \) is different for \( k\in \{-1,1\} \)
 and \( k=0 \).
 If \( k\in \{-1,1\} \) then \( d=1 \) is for a track to \emph{send to}  the destination colony
 in direction \( k \), and \( d=-1 \) is for a track to \emph{receive from} that colony.
 If \( k=0 \) then the track with \( d \in\{-1,1\} \) is for sending \emph{left} or \emph{right}
 within the present colony. 

 In an ordinary cellular automaton, we could move a string on the mail
track for a certain number of steps and then copy it to another track in a single step.
But in a variable-time medium, we cannot rely on such timing.
Therefore with each piece of information, the mail track will also carry the address of the
place it is coming from, and this will allow to know when and where to deposit it
without relying on timing.
Field \( \Mail_{k,d} \) has subfields%
 \glo{toaddr@\( \Toaddr \)}%
 \[
  \Toaddr, \Info. 
 \]
  For simplicity, let \( |\Mail_{k,d}.\Info|=|\Info| \).
For adjacent cells \( x,y \) with \( j=y-x \), and their colonies \( x^{*} \) and
\( y^{*} \) (defined from their \( \Addr \) field) we define the predicate%
 \glo{edge@\( \Edge_{j}(x) \)}%
 \[
  \Edge_{j}(x)=\begin{cases}
        0 &\text{if \( x^{*}=y^{*} \),}
\\    1 &\text{if \( x \) and \( y \) are endcells of two adjacent colonies,}
\\ \infty &\text{otherwise.}
  \end{cases}
 \]
Thus, the function \( \Edge_{-1}(x) \) depends implicitly on the current values of the
address fields in the configuration \( \eta \).
We have, \( \Edge_{-1}(x)=0 \) if \( x-1 \) is in the same colony as \( x \), it is 1 if  \( x-1 \) and
\( x \) are right and left endcells of their respective colonies, and it is \( \infty \)
otherwise (abnormal case, with ``inconsistent'' address field values).
The mail track\index{track!mail} \( \Mail_{k,d} \) of cell \( x \) will read,
from the neighbor in direction \( j=j(k,d) \), from mail track 
\glo{peer@\( \peer(k,d) \)}
\begin{align}\label{eq:peer}
   \peer(k,d) = (k,d') 
\end{align}
defined as follows.
\begin{itemize}
\item If \( |k|=0 \), \( \Edge_{-d}(x)=0 \) then \( j = -d \), \( d'=d \).
\item If \( |k|=1 \), \( \Edge_{-kd}=0 \) then \( j= -k d \), \( d'=d \).
\item If \( |k|=1 \), \( \Edge_{kd}>0 \) then \( j= k d \), \( d'=-d \).
\end{itemize}
 
In all other cases, the values \( j(k,d) \), \( \peer(k,d) \) are not defined.
In words: the mail is passed along the same mail track except that when it crosses an edge then
it goes from the sending track of the sending colony to the receiving track of the receiving colony.
Formally we define
\begin{align}
  \label{eq:mail-to-receive}
   \var{Mail-to-receive}(k,d) = \Mail_{\peer(k,d)}^{j(k,d)}.  
\end{align}
as the mail to be received into \( \Mail_{k,d} \) provided \( \peer(k,d) \) is defined,
and \( \Undef \) otherwise.
The one-step sub-rule \( \MoveMail \) gets mail from the neighbor cell:%
 \glo{move-mail@\( \MoveMail \)}

 \begin{algo}[H]\caption{sub-rule \( \protect\MoveMail(k,d) \)}\label{alg:move-mail-1}
        \( \Mail_{k,d}\gets \var{Mail-to-receive}(k,d) \) if the latter is defined.
 \end{algo}

A cell will typically copy the information to be sent into
the subfield \( \Mail_{k,d}.\Info \) and at the same time, the receiver cell's address into
\( \Mail_{k,d}.\Toaddr \).
In the copy rule here, \( j \) refers to the direction of the
sending colony as seen from the receiving colony.
The next argument is the location of origin in the sending colony, the
one following is the location to receive the information in the receiving
colony.

 \begin{algo}[H]\caption{sub-rule \( \protect
 \Copy(m,\protect\F(\lint{a}{a + n}),\protect\G(\lint{b}{b + n})) \)}\label{alg:copy-1}
\algLet \( d\gets 1, d'\gets -1 \) if \( m\ne 0 \) and \( d'=d\gets\sign(b-a) \) otherwise\;
\nl\lIf{\( \Addr\in\lint{a}{a+n} \)}{\( \Mail_{m,d}.(\Toaddr,\Info) \gets  (Addr+b-a,\F) \)}
\label{alg:Toaddr}
\pRepeat{\( 2\Q  \)}{
  \( \MoveMail(m,d) \)\;
  \If{\( \Addr\in\lint{b}{b+n} \) \algAnd 
    \( \Addr = \Mail_{m,d'}.\Toaddr \) 
  }{
    \( \G\gets \Mail_{m,d'}.\Info \)
  }
}
\end{algo}

\begin{figure}
 \setlength{\unitlength}{0.2mm}
 \[
 \begin{picture}(600,150)
    \put(0, 0){\framebox(600, 120){}}

    \put(300, 0){\line(0, 1){120}}

    \put(0,20){\line(1,0){600}}

    \put(610, 0){\makebox(0,20)[l]{\( \var{Mail} \)}}
    \put(610, 40){\makebox(0,20)[l]{\( \F \)}}
    \put(610, 60){\makebox(0,20)[l]{\( \G \)}}

    \put(30, 40){\framebox(140, 20){}}
    \put(370, 60){\framebox(140, 20){}}

    \put(30, 125){\makebox(0,0)[b]{\( a \)}}    
    \put(170, 125){\makebox(0,0)[b]{\( a+n \)}} 
    \put(370, 125){\makebox(0,0)[b]{\( b \)}}   
    \put(510, 125){\makebox(0,0)[b]{\( b+n \)}}

    \multiput(20, 0)(60,0){10}{\makebox(20,20){\( \to \)}}
 \end{picture}
 \]
 \caption{ \label{fig:Copy}\( \protect\Copy(-1,\protect\F(\lint{a}{a + n}),
                   \protect\G(\lint{b}{b + n})) \)}
\end{figure}

\begin{remark}\label{rem:variable-width-copy}
Here we assumed that the source and destination locations are on tracks of the same width,
and have the same lengths.
We can write a more general version of \( \Copy \) using a loop.
For example assume that the source location, of length \( n \),
is on a track \( \F \) that is 5 times wider than its destination location
\( \G \) of length \( 5 n \).
Then we can split the field \( \F \) into 5 sub-fields \( \F_{1},\dots,\F_{5} \).
More generally, a track that has \( m \) times the width can be
treated as into \( m \) tracks of the same width.
Here \( m \) can even be a variable computed from the widths of the fields in question,
using the primitive \( \F[i\frTo j] \) as defined in~\eqref{eq:F[i:j]}.
We can do the copying in 5 iterations, copying each sub-track separately,
or (this will be our actual choice later)
we can post each cell onto the mail track in 5 steps.
The rule \( \MoveMail \) can proceed in parallel with all these iterations.

Similarly, one can write a program for copying
a longer and narrower location into a shorter and wider one.
Also, instead of giving the locations as arguments, we may just give the name of a
field where they are described.  
In applications we will assume this more general implementation.
\end{remark}

 \begin{definition}[Indirectly given locations]\label{def:indirect-copy}
  The description \( \F_{1}(\lint{a_{1}}{a_{2}}) \) of a location can fit into a
single field \( \fld{G}_{1} \) of some cell.
If for example argument \( \loc_{1} \) of the rule
\( \Copy(m,\loc_{1},\loc_{2}) \) is given by a field \( \fld{G}_{1} \) 
this way, then we write \( \Copy(m,\fld{G}_{1}*,\loc_{2}) \).
This will only happen if \( m=0 \), that is the copying proceeds within one
colony.
It will be assumed that field \( \fld{G}_{1} \) of
each cell of the colony contains the same information \( \loc_{1} \).
Therefore the rule can be written just as above, except that
some of its parameters are read now from the field \( \fld{G}_{1} \).
 \end{definition}

\paragraph {The evaluation rule}
The rule \( \Eval \) controls the track \( \Cpt \).
Before describing it we need to define some locations.

\begin{definition}[Locations for evaluation]\label{def:eval-locations}
Here are a few locations used in the rule \( \Eval \):
\begin{bullets}
 \item Track \( \Cpt \) is used for universal computation, with sub-tracks
\( \Cpt.\Input \) and \( \Cpt.\Output \).
 \item  The location for the interpreter on the \( \Cpt.\Input \) track is
denoted by \( \fld{\_Interpr} \).
  \item The location for the program to be executed by the universal medium
    \( \univ \) is denoted by \( \fld{\_Prog} \).
  \item The locations for the parameters of the program
    are denoted by \( \lParam_{i} \), \( i=1,2,\dots \), as needed,
  \item The location of the arguments of the function to be executed
    on the \( \Cpt \) track are \( \lArg_{m} \) for \( m\in\{-1,0,1\} \).
 \glo{arg@\( \lArg_{m} \)}%
  \item The location of the output of the computation on the \( \Cpt.\Output \) track
    is \( \lSimOutput \)\glo{sim-output@\( \lSimOutput \)}.
We assume it to be the same interval as that of location \( \lInfo \).
\end{bullets}  
\end{definition}

We will assume that the program to be simulated is written as a string
value of parameter \( \var{Trans-prog}_{2} \) of our program.
The first steps of \( \Eval \) get the arguments, then
write the interpreter (a constant string),
followed by writing the program \( \rul{Trans-prog}_{2} \)
of the transition function to be simulated:
 \glo{prog@\( \fld{\_Prog} \)}%
 \glo{interpr@\( \fld{\_Interpr} \)}
 \begin{algo}[H]
     \( \Write(\Interpr, \fld{\_Interpr}) \);\;
     \( \WriteParam(\var{Trans-prog}_{2}, \fld{\_Prog}) \);\;
 \end{algo}
But if \( \trans_{1}=\trans_{2} \) (self-simulation) is desired, then write
\( \MyRules \) in place of \( \var{Trans-prog}_{2} \) above.
In both cases, follow with \( \WriteParam(\Param_{i}, \lParam_{i}) \) and
copying \( \lRetrieved_{m} \) to \( \lArg_{m} \) for each \( m \).
Then the rule
\begin{align}\label{eq:Initialize}
   \Initialize
 \end{align}
writes \( *\cdots * \) to the \( \Output \) track
and the rest of the cells (including the endcells) of the \( \Cpt.\Input \)
track, and 0's to the track \( \Cpt\setminus(\Cpt.\Input\cup\Cpt.\Output) \).
Then the rule
\begin{align*}
   \rul{Interpret}
 \end{align*}
applies, for a sufficient number of steps, the transition function \( \univ \)
to the \( \Cpt \) track.
According to Theorem~\ref{thm:rule-lang} (Rule Language), the computation
finishes in
 \[
  \cns{interpr-coe}(|\rul{Trans-prog}_{2}|+1)^{2}\nm{\states_{2}}
 \]
 steps, so this number of iterations is sufficient.
 So the whole rule \( \Eval \) is (for the case of self-simulation)
 in Algorithm~\ref{alg:simp-sim.eval}, for \( N \) parameters:

 \begin{algo}\caption{sub-rule \Eval}\label{alg:simp-sim.eval}
     \( \Write(\Interpr, \fld{\_Interpr}) \);\;
     \( \WriteParam(\MyRules, \fld{\_Prog}) \);\;
     \For{\( i=1 \) \KwTo \( N \)}{
       \( \WriteParam(\Param_{i}, \lParam_{i}) \)
     }
     \lFor{\( m = -1,0,1 \)}{\( \Copy(0,\lRetrieved_{m},\lArg_{m}) \)}
     \Initialize;\;
     \Interpret
   \end{algo}
 
  The subrule 
 \[
  \rul{Update}
 \]
in Algorithm~\ref{alg:basic-sim-program} 
copies the track \( \Cpt.\Output \) into track \( \Info \) for all addresses in 
\( \Space(\lInfo) \).


\paragraph{Summary in a theorem}
Before formulating what will be accomplished by the above program,
let us define some encodings.

\begin{definition}
  The encoding \( \fg_{*} \) of a cell state \( v \) of \( M_{2} \) into a colony of
\( M_{1} \) is defined as follows.
The string \( v \) is written into \( \lInfo \).
The \( \Cpt \) track and the mail tracks are set to all 1's.
Each \( \Age \) field is set to 0.  
For all \( i \), the \( \Addr \) field of cell \( i \) of the colony is set to \( i \).
  \end{definition}

The theorem below states the existence of the above simulation.
As a condition of this theorem, the parameters
 \begin{equation}\label{eq:sim-params}
 \rul{Trans-prog}_{2}, \nm{\states_{1}}, \nm{\states_{2}}, \Q , \U 
 \end{equation}
  will be restricted as follows.

  \begin{condition}\label{cond:simple-block-sim-ineq}
The following inequalities must be obeyed.
 \begin{description}

  \item[Cell Capacity Lower Bound]\index{lower@Lower Bound!Cell Capacity}
 \[
    \nm{\states_{1}} \ge c_{1}\cei{\log \U }+ \nm{\states_{\univ}} + c_{2}
 \]
  where \( c_{1},c_{2} \) can be easily computed from the following
consideration. 
  What we really need is 
 \( \nm{\states_{1}} \ge |\Addr|+|\Age|+|\Info|+|\Mail|+|\Cpt| \) where the following
choices can be made:
 \begin{alignat*}{3}
        &|\Info|       &&= 2,
\\   &|\Mail_{i}|  &&= |\Mail_{i}.\Toaddr| 
       + |\Mail_{i}.\Info| = \cei{\log \Q } + 3,
\\   &|\Cpt|        &&= \nm{\states_{\univ}},
\\   &|\Addr|      &&= |\Age| = \cei{\log \U }.
 \end{alignat*}

 \item[Colony Size Lower Bound]\index{lower@Lower Bound!Colony Size}
  \[
   \Q  \ge \cns{interpr-coe}(\nm{\states_{2}} + |\rul{Trans-prog}_{2}|
         + \log \U  + 1).
  \]
 With the field sizes as agreed above, this provides sufficient space in
the colony for information storage and computation.

  \item[Work Period Lower Bound]\index{lower@Lower Bound!Work Period}
 \[
   \U  \ge 3\Q +\cns{interpr-coe}(|\rul{Trans-prog}_{2}| + 1)^{2}\nm{\states_{2}}.
 \]
 With the field sizes above, this allows sufficient time for the above
program to be carried out.
 \end{description}
      \end{condition}

It is not difficult to find parameters satisfying the above inequalities
since \( \log \Q \ll \Q  \).

\begin{sloppypar}
\begin{theorem}[Basic Block Simulation] \label{thm:simp-sim}
 \index{theorem@Theorem!Basic Block Simulation}
  There are strings \( \rul{Sim-prog}_{1} \)\glo{sim-prog@\( \rul{Sim-prog} \)},
\( \rul{Sim-prog}_{0} \) such that the following holds.
If \( \rul{Trans-prog}_{2} \) with parameters \( \nm{\states_{1}} \), \( \nm{\states_{2}} \), \( \Q  \),
\( \U  \) satisfies the above inequalities, and \( \rul{Trans-prog}_{2} \)
defines a transition function \( \trans_{2} \),
then \( \rul{Sim-prog}_{1} \) is a rule program with parameters
 \begin{align*}
\rul{Trans-prog}_{2}, \nm{\states_{1}}, \nm{\states_{2}}, \Q , \U ,
 \end{align*}
defining a transition function \( \trans_{1} \) such that \( \CA(\trans_{1}) \) 
has a block simulation of \( \CA(\trans_{2},\Q ,\U ) \).

 If \( \nm{\states_{1}},\nm{\states_{2}},\Q ,\U  \) satisfy the above inequalities, then
\( \rul{Sim-prog}_{0} \) is a rule program with parameters
\( \nm{\states_{1}} \), \( \Q  \), \( \U  \),
defining a transition function \( \trans_{1} \) such that \( \CA(\trans_{1}) \) 
has a block simulation of \( \CA(\trans_{1},\Q ,\U ) \).
\end{theorem}  
\end{sloppypar}

The given construction is essentially the proof.
Its complete formalization would yield \( \rul{Sim-prog}_{1} \) and
\( \rul{Sim-prog}_{0} \) explicitly.


\section {Robust media}\index{medium!robust}
 \label{sec:rob}

This section defines a special type of one-dimensional medium called
\df{robust}.
From now on, when we talk about a medium with no qualification we
will always mean a robust medium.
The definition can be compared to the
primitive variable-period media in Section~\ref{sec:asynch}.
The main distinguishing feature is the notion of damage.

\begin{remark}\label{rem:non-adjacent}
  The paper~\cite{Gacs1dim86} used communication only between adjacent neighbors.
  It used a special mechanism for preventing accidental intrusions into an intact colony:
  growth in a zigzag pattern.
  Here this is achieved by communication between non-adjacent neighbors,
  enabling a colony of ``stronger'' cells to defend itself from
  the accidental intrusion by ``weaker'' cells.
  However, to simulate such a property on higher levels, it will also be used
  to pass information other than strength to non-adjacent neighbors.
  \end{remark}

\subsection{Damage}\index{damage}
 \label{sec:rob.damage}

In a robust medium, being in the space-time set $\Damage$ defined below
excuses the site for ``not following the rules''.

\begin{definition}[Damage]\label{def:damage}
For the media to be defined, robust media,
we introduce a special set of states $\Bad\subseteq\Vacant$\glo{bad@$\Bad$}.
For a history $\eta$ we define the \df{damage set}%
 \glo{damage@$\Damage$}%
 \[
  \Damage(\eta)=\setof{ (x,t) : \eta(x,t)\in\Bad} .
 \]
  For a space configuration, the damage is defined similarly.
For a site $x$ a time interval $I$ is \df{damage-free}\index{damage-free}
if $\eta(x,\cdot)$ is not in $\Bad$ during $I$.
  \end{definition}

Damage points can be viewed as holes in the fabric of the lawful
parts of a trajectory.
When $(x,t)\in\Damage(\eta)$ then in the neighborhood of $(x,t)$, we will
not be able to make any predictions of $\eta$, that is in some sense,
$\eta$ behaves completely ``lawlessly'' there.
In all cellular media concerned with our \emph{results} we could
require $\Damage(\eta)=\eset$.
The damage concept is necessary only in a trajectory $\eta_{2}$ of a
medium $M_{2}$ obtained by simulation from a trajectory $\eta_{1}$ of some
medium $M_{1}$.
Such a simulation typically requires some structure in $\eta_{1}$.
When noise breaks down this structure the predictability of $\eta_{2}$
suffers and this will be signalled by the occurrence of damage in
$\eta_{2}$.
However, for convenience, we will define damage even in the media used on
the lowest level, as a violation of a certain transition rule.    

\begin{definition}[Damage map]\label{def:damage-map}
We will use the following constants:
 \begin{equation}\label{eq:bdist}
  \bdist = 5,\quad \bsize = 2\bdist+5.  
 \end{equation}
  Let us define the rectangle
 \begin{align}\label{eq:Vdam}
   \Vdam=\lint{-\B /2}{\B /2}\times\rint{-\Tu/2}{\Tu/2}, 
 \end{align}
 then the corresponding rectangle $\Vdam^{*}$ is defined using $\B ^{*},\Tus$.
 Two space-time points are said to be \df{too close} if they are contained
 in a single copy (space-time translation) of \( \bsize\Vdam \).
 Two sets \( A_{1},A_{2} \) are \df{too far} if no copy of \( \bdist\Vdam^{*} \) intersects both.
A subset $A$ of the damage is called an \df{island} if it is covered by a translation of
$\bsize\Vdam$ and is too far from the rest of damage.
In a simulation $\eta^{*}=\Phi^{*}(\eta)$ between two robust media,
the set $\Damage(\eta^{*})$ is obtained after we remove all islands of
$\Damage(\eta)$.  
We call this definition of $\Damage(\eta^{*})$ the \df{damage map} of simulation
$\Phi^{*}$.
\end{definition}

According to the general definition of a medium, the set of trajectories
will be defined by a pair $ b(\cdot),\g(\cdot,\cdot,\cdot)$ where $ b(\alpha)$
give the probability bound belonging to type $\alpha$ and $g(\alpha,W,\eta)$ is
the event whose probability is bounded.
We formulate these in terms of  \df{properties}.
 The \df{Computation Property}\index{axiom@Property!Computation} constrains the kinds
of events that can occur under the condition that the damage does not
intersect a certain rectangle.
The \df{Restoration Property}\index{axiom@Property!Restoration} bounds the probability of
occurrence of damage in the middle of some window.
It depends on a parameter  $\eps$, and says that the probability of the damage is small:
even if other damage occurs at the beginning of a window it will be
eliminated with high probability.  

\begin{condition}[Restoration Property] \label{cond:restor}
  Let $ b(\alpha(\var{restor}))=\eps$.
  Further, for any pair $(x,t)$ let
 \[
   \g\bigparen{\alpha(\var{restor}),\; (x,t)+(\bdist+2)\Vdam,\; \eta}
 \]
  be the event function for the event that there is damage in
$\eta((x,t)+\Vdam)$.
 \end{condition}

The Restoration Property says that damage, that is the occasional obstacle
to applying the Computation Property, 
has small conditional probability of occurrence in the middle part $(x,t)+\Vdam$ 
of any rectangle of the form $(x,t)+\bdist\Vdam$.
The property will hold automatically on the lowest level by the
property of the medium that the transition rule is only violated with small
probability.\footnote{In the model of~\cite{Gacs2dim89}, the restoration property is weaker.
There, damage does not necessarily disappear in a short time but if
it is contained in a certain kind of triangle at the beginning of the
window then, with high probability, it is contained in a smaller
triangle of the same kind later in the window.}

Damage helps us present a hierarchical construction as a sequence of
simulations.
When a large burst of faults destroys the fabric of these nested
simulations, then $\eta^{k+1}$ cannot be explained just in terms of the
$\eta^{k}$ from which it is decoded.
The damage set will cover those lower-level details of the mess that
we are not supposed to see as well as the mechanism of its removal.

Provided $\Damage(\eta^{*})$ is not nearby, $\eta^{*}(x,t)$, 
will be essentially defined by a block code $\fg$ as
 \[
  \fg^{*}(\eta(x+\clint{0}{\Q\B -1},\;t)).
 \]
  However, in the interest of stabilization there will be some look-back in time: 
the simulation will not be memoryless.

\begin{lemma}[Simulation Damage Probability Bound]\label{lem:sim-damage}
 \index{lemma@Lemma!Simulation Damage Probability Bound}
  Let $M$, $M^{*}$ be media with parameters $\eps$, $\eps^{*}$,
with $\B ^{*}\ge \B $, $\Tus\ge \Tu$, whose local
condition system includes the Restoration Property.
Let a simulation $\Phi^{*}$ be defined between them which assigns damage
in $M^{*}$ according to the damage map of Definition~\ref{def:damage-map}.
There is a constant $c_{\Dam}$ such that if
\begin{equation}\label{eq:eps*}
  \eps^{*}\ge c_{\Dam}((\B ^{*}/\B )(\Tus/\Tu)\eps)^{2}  
\end{equation}
  then $\Phi^{*}$ is a deterministic canonical simulation map for damage, as defined 
in~\ref{sec:canonical}.
 \end{lemma}

The proof observes that small bursts of damage that are not too close
to each other behave in the medium $M$ as if they were independent,
therefore the probability of the occurrence of two such bursts can be
estimated by $O(\eps^{2})$ times the number of such pairs in a rectangle
$\Vdam^{*}$.

\begin{proof}
It is enough to show an expression of the form~\eqref{eq:canonical} for
local conditions of type $\var{restor}$ in $M^{*}$.
Let $\eta$ be a trajectory of $M$.
Consider the event that 
$\Damage(\eta^{*})$ intersects the middle part $W'=\pair{x}{t}+\Vdam^{*}$
of rectangle
\begin{align*}
 W=\pair{x}{t}+(\bdist+2)\Vdam^{*},
\end{align*}
with \( \bsize,\bdist \) defined in~\eqref{eq:bdist},
then of course $\Damage(\eta)$ intersects it, too.
In what follows, by ``damage'' we mean the set $\Damage(\eta)$.
We claim that the damage in \( W \)
cannot be covered by a copy \( W'' \) of \( \bsize\Vdam \).
Suppose namely that it can: then it can also be covered by a copy of \( W'+\bsize\Vdam\).
However, the rectangle \( W'+\bsize\Vdam \) is too far from the complement of \( W \).
Then so is \( W'' \), so it is an island and as such it would have been
deleted by the damage map, contradicting the assumption that its elements
belong to \( \Damage^{*} \).

\begin{sloppypar}
We found two damage points $p_{1},p_{2}$ in $W$ that are not too close.
Let us fix some partition of the space into copies of \( \Vdam \).
Let \( \cV \) be the set of those elements of this partition intersecting \( W \).
For each point $p$ in space-time, let $U(p)$ be the $U\in\cV$ with
$p\in U$.
Let \( U'(p)=U(p)+(\bdist+1)\Vdam \).
The fact that \( p_{1},p_{2} \) are not too close implies that the two copies
\( U'(p_{1}) \) and \( U'(p_{2}) \) of \( (\bdist+2)\Vdam \) are disjoint.
Indeed, if they intersected then some copy of \( 2(\bdist+2)\Vdam \)
would contain both, but due to~\eqref{eq:bdist} then they would be too close.

Let $E$ be the set of pairs $U_{1}, U_{2}$ in $\cV$  with \( U'_{1}\cap U'_{2}=\emptyset \).
We found
 \[
 \g^{*}(\var{restor},(\bdist+2)\Vdam^{*},\eta^{*}) 
  \le \sum_{\pair{U_{1}}{U_{2}}\in E}\prod_{j=1,2} \g\bigparen{\var{restor},U'_{j},\eta},
\]
therefore \( \Expv  \g^{*}(\var{restor},\Vdam^{*},\eta^{*}) \le |E|\eps^{2} \).
Since counting bounds $|E|$ by
$c_{\Dam} ((\B ^{*}/\B ) (\Tus/\Tu))^{2}$ for
an appropriate constant $c_{\Dam}$, by the assumption~\eqref{eq:eps*}
both conditions of a weak canonical simulation in~\eqref{eq:canonical} are satisfied.
  \end{sloppypar}
\end{proof}

\subsection{Computation}\label{sec:rob.comput}

Before giving the Computation Property, let us introduce some details of the
structure of a robust medium.
Cells sometimes have to be erased in a trajectory since
cells created by the damage may not be aligned with the original ones.
At other times, the creation of a live cell at a vacant site will be
required.
Most of the trajectory requirements of a robust medium will be
expressed by a transition function $\trans$, desribed later.
Killing and creation may be indicated by some special values of this
function.

\paragraph{Neighborhood structure}
We extend the meaning of the neighbor function $\thg_{j}(x)$\glo{th.greek@$\thg_{j}(x)$}
introduced in Notation~\ref{n:neighbors}.

\begin{notation}[Non-adjacent neighbors]
Let us fix a history $\eta$.
 For a (non-vacant) cell $x$, in direction \( j=\pm 1 \)
 there is at most one cell at distance smaller than $2\B $: if it
exists and there is no damage between it and $x$ then we denote it also by \( \thg_{j}(x,t,\eta) \).
We will omit $\eta$, and sometimes even $t$, from the arguments when
it is obvious from the context.
By convention, whenever $\thg_{j}(x,t,\eta)$ is undefined then
let \( \eta(\thg_{j}(x,t,\eta),t)=\Vac \).
\end{notation}

Having an extended
notion of neighbors, we can define the notion of transitions for robust media.

\begin{definition}[Transition of a robust medium]
  The transition function in a robust medium has the form
  \begin{align}\label{eq:rob-trans}
    \trans(\bs,\ba).
  \end{align}
 Here \( \bs = (s_{-1},s_{0},s_{1}) \), \( \ba =(a_{-1},a_{1}) \)
have the following meaning. 
 \( s_{0} \) is the state of the cell whose transition will happen, and for \( j=\pm 1 \)
 \( s_{j} \) is the state of its (not necessarily adjacent) neighbors in direction \( j \).
 Further \( a_{j}=1 \) if the neighbor in direction \( j \) is adjacent and 0 otherwise.
  The transition does not depend on \( a_{j} \) if either \( s_{j} \) or \( s_{0} \) is vacant.
  Note that \( a_{j} \) is in fact \( a_{j}(x,t,\eta) \), but we may omit the arguments
that are obvious from the context.
 In analogy with~\eqref{eq:trans}, let%
 \[
   \ol\trans(\eta,x,t)=\trans((\thg_{-1}(x),x,\thg_{1}(x)),(a_{-1}(x),a_{1}(x))).
 \]
The symbol $\eta$ will be omitted from $\trans(\eta,x,t)$ when it is
obvious from the context.
\end{definition}


\paragraph{Properties of the transition function}
The state space of a robust medium will be required to have some minimal structure,
and every transition function will be required to have certain properties.
In Definition~\ref{def:damage} it has been mentioned already that the set of states has
a subset \( \Bad \).
The value of the transition function \( \trans \) will not be defined if any of its argument is \( \Bad \)
and we could even say that this defines the set \( \Bad \).

\begin{definition}\label{def:rob}
A robust medium will be denoted by%
 \glo{rob@$\Rob$}%
 \[
 \Rob(\trans, \B ,\Tl,\Tu,\eps,\eps'). 
\]
The first bit of the state of a cell will be the field $\Rand$;  we will write
 \begin{equation*}
 \Det=\All\setminus\Rand.
 \end{equation*}
\end{definition}

The transition function, a mapping from assignments to $\states$, describes
the goal for the deterministic part $s.\Det$ of the value $s$ of the state
after transition.
We impose a number of conditions on the transition functions in robust media.

\begin{condition}[Time Marking]\label{cond:time-mark}
 \index{condition@Condition!Time Marking}
  If a cell is not vacant, then the transition function always changes it.
\end{condition}

This condition is similar to~\eqref{eq:prim-time-mark},
and is helpful for nontrivial media with a strict dwell period upper bound.

\begin{condition}\label{cond:latent}
  Non-vacant cells will have a property (typically represented by a field) of being \df{latent} or not.
  Newborn cells will be latent.
  The transition function will not make a difference between a vacant and latent neighbor.
\end{condition}

This condition will help attributing all information obtained from the neighbors
to the same moment of time even when a vacant neighbor turns to a latent one.

Transitions between vacant and non-vacant state are handled in a special way.

\begin{condition}[Cling to Life]\label{cond:cling}
   \index{condition@Condition!Cling to Life}
A cell will only be erased if it may disturb a close non-aligned neighbor:
namely, suppose that \( u_{0}\ne\Vac \) and \( \trans(\bu,\ba)=\Vac \).
   Then there is a \( j=\pm 1 \) with \( u_{j}\ne\Vac \) and \( a_{j}=0 \).
 \end{condition}

\begin{definition}[Creators, emergence]\label{def:creator-emerging}
  Recall the form of the transition function in~(\ref{eq:rob-trans}).
  We will call a pair $(j,v)\in\{-1,1\}\times\states$, a
\df{potential creator}\index{creator!potential} of a non-vacant state $s$
if $s=\trans(\bu,\ba)$ where \( u_{j}=v \), \( u_{0}= u_{-j}=\Vac \).
The expectation is that an adjacent cell with state \( s \) will be created by
the cell with state \( v \) determined by the state \( u_{j} \) of the creator
\emph{independently of} the state \( u_{-j} \).

Consider a history $\eta$, with a switching time $\sigma$ of cell
$x$ when $\eta(x,t)$ turns from vacant to non-vacant.
We call $y=x+j\B$ a \df{creator}\index{creator} of $x$ for time $\sigma$ if
$(j,\eta(y,\sigma-))$ is a potential creator of $\eta(x,\sigma)$, further
$\eta(y,t)$ is non-vacant for $t\in\opint{\sigma-\Tl/2}{\sigma+\Tl}$.
(So a creator is supposed to survive the creation.)
On the other hand, if there is no
\( t \) in \( \rint{\sigma-\Tl/2}{\sigma} \) when the
interval\( \lint{x-2\B}{x+3\B} \) intersects the body of any cell,
we will say that the cell \( x \) \df{emerged} at time \( \sigma \).
  \end{definition}

\paragraph{The computation property}
The condition called the Computation Property has a form similar to the
definition of primitive variable-period media in 
Section~\ref{sec:asynch}.

\begin{definition}[Special switching times]\label{def:special-sw-times}
  For a history $\eta$, cell $x$ and number $a \ge 0$ let
 \[
  \sigma_{1},\sigma_{2},\sigma_{0}
 \]
  be defined as follows.
  The values $\sigma_{1},\sigma_{2}$ are the first two switching times of \( x \)
  in $\rint{a}{a+2\Tu}$, but \( \sigma_{2} \) is defined only if $\eta(x,\sigma_{1})$ is non-vacant.
0 is considered a switching time if $\eta(x,0)$ is not vacant.
On the other hand, $\sigma_{0}$ is the first switching time of $x$ after
$a+\Tu$ but defined only if $\eta(x,a+\Tu)$ is vacant.
  Whenever we have an event function $\g(\alpha,W,\eta)$ in whose definition
$\sigma_{2}$ occurs, this function is always understood to have value 0 if
$\sigma_{2}$ is not defined (similarly with $\sigma_{0}$).  
\end{definition}

Referring in an event function to, say, \( \sigma_{2} \) directly is a shorthand.
We will consider the times \( \sigma _{i} \) only in places where no damage occurs,
which allows us to lower-bound all dwelling times by \( \Tl \).
Then all requirements referring to \( \sigma_{i} \) can be expressed  instead in terms of some
arbitrary sequence of rational numbers \( t_{i}=t_{0}+i\Tl \).
For example the requirement \( \sigma_{2}-\sigma_{1}<\Tu \)  can be expressed by
prohibiting for each rational \( t \), and for \( n \)  with \( \Tu/n<\Tl \), that
all \( \eta(x,t+id)=\eta(x,t) \) hold for \( i=i,\dots,n \).
The requirement \( \Pbof{\eta(x,\sigma_{2}).\Rand=j}\le 1+\eps' \) can be expressed by 
\begin{align*}
 \Pbof{\eta(x,t_{1})\ne\eta(x,t_{0}) \land\eta(x,t_{1}).\Rand=j}\le 1+\eps' .
\end{align*}
for \( t_{0}<t_{1}<t_{0}+\Tl \).

\begin{figure}
 \setlength{\unitlength}{1mm}
 \[
 \begin{picture}(60,70) 
 \put(7.5, 0){\framebox(35, 60){}}

 \put(75, 20){\vector(0,1){20}}
 \put(79, 30){\makebox(0,0)[l]{time}}

 \put(8.75, 52){\makebox(8,8){$W_{1}$}}	
 \put(0, 62){\makebox(0,0)[b]{$x-1.1 \B $}}	
 \put(20, 60){\circle*{1}}	
 \put(18.75, 0){\line(0,1){60}}	
 \put(31.25, 0){\line(0,1){60}}	
 \put(20, 52){\makebox(8,8){$W_{0}$}}	
 \put(20, 62){\makebox(0,0)[b]{$x$}}	
 \put(50, 62){\makebox(0,0)[b]{$x+2.1\B $}}	

 \put(55, 60){\makebox(0,0)[l]{$a+2\Tu$}} 
 \put(55, 35){\makebox(0,0)[l]{$\sigma_{2}$}} 
 \put(55, 30){\makebox(0,0)[l]{$\sigma_{2}-\Tl/2$}} 
 \put(55, 15){\makebox(0,0)[l]{$\sigma_{1}$}} 
 \put(55, 5){\makebox(0,0)[l]{$a$}} 
 \put(55, 0){\makebox(0,0)[l]{$a-\Tu$}} 
 \put(20, 15){\framebox(10, 20){}}
 \put(20, 30){\line(1,0){10}}

 \end{picture}
 \]
 \caption[To the computation property]{To the computation property.
The large rectangle is $W_{1}(x,a)$.
The rectangle of the same height but width $1.2 \B $ around the cell body at point
$x$ is $W_{0}(x,a)$.
The small rectangle between times $\sigma_{1}$ and $\sigma_{2}$ is the cell work
period under consideration.  Its lower part is the 
observation interval\index{period!observation}.}
 \label{fig:comput-ax}
\end{figure}
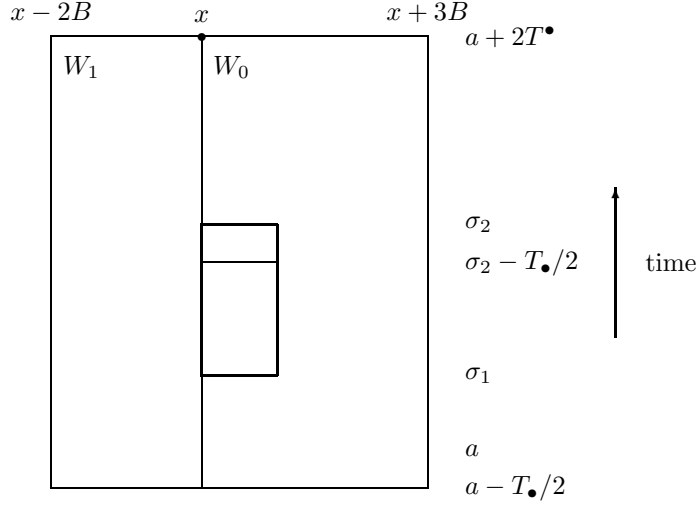

For a transition to work satisfyingly, the cell as well as its neighbors must have been
damage-free for some time before the transition.
But it is convenient to require that a transition not fail completely
even if only the cell itself is damage-free, not its neighbors.
The following space-time rectangles will play a special role.

\begin{definition}\label{def:rectangles}[Rectangles used in the computation property]
  Let%
 \glo{w.cap@\( W_{0}(x, a), W_{1}(x, a) \)}%
 \begin{align*}
       I_{0}(x)    &= \lint{x-0.1\B }{x+1.1\B }, \quad I_{1}(x) = \lint{x-2.1\B }{x+3.1\B }, 
\\   W_{j}(x,a) &= I_{j}(x)\times\rint{a-2\Tu}{a+3\Tu} \quad (j=0,1),
   \\   \f_{j}(x,a,\eta) &= \text{the event function for }
    \setof{\Damage(\eta)\cap W_{j}(x,a) = \eset} \quad (j=0,1).
 \end{align*}  
\end{definition}

Let us write
 \[
  \Damage=\Damage(\eta) ,\quad \Damage^{*}=\Damage(\eta^{*}) . 
 \]

 \begin{definition}[Affecting]\label{def:affecting}
We will say that damage \df{affects} cell \( \pair{x}{t} \) \df{directly}, if it
intersects with  \( I_{0}(x) \) at time \( t \), with \( I_{j}(x) \)
as in Definition~\ref{def:rectangles}.
It affects a cell \df{via neighbors}, if it  intersects \( I_{1}(x) \).
\end{definition}

Below we formulate Condition~\ref{cond:comput} (Computation Property),
which beside the Restoration Property is the other crucial requirement on the trajectory.
It implies that if
a cell is not affected directly during a certain time interval, then it makes a legal
transition.
It may still not make the transition required by the transition
function: for that, it must not be affected even via neighbors.
With these notions, we are ready to spell out the computation property.
Condition~\ref{cond:comput}\ref{i:comput.trans} refers to a notion called ``atomicity''
in the field of distributed computing: the fact that despite the complexity and
prolonged nature of the interaction between a colony and its neighbors, its
simulation result can be viewed as the result of an observation
of the states of the three represented big cells in a single moment.

 \begin{condition}[Computation Property]  \label{cond:comput}
  This property consists of several condition types.
For each type \( \alpha \) used here except \( \alpha(\var{rand},j) \), we 
have \(  b(\alpha)=0 \), that is the corresponding events \( \g(\alpha(\cdot),W,\eta) \) are
simply prohibited in the trajectory.
  On the other hand,
 \[
   b(\alpha(\var{rand},j))=0.5+\eps'\quad (j=0,1).
 \]
  \begin{cjenum}

 \item\label{i:comput.rand}
  Coin-tossing, for \( j=0,1 \):
 \[
   \g\bigparen{\alpha(\var{rand},j),\;W_{0}(x,a),\;\eta} = \f_{0}(x,a,\eta)
 \{\exists\sigma_{2}\land \eta(x,\sigma_{2}).\Rand=j\}.
 \]

  \item\label{i:comput.dwell-pd-bd}
  The length of dwell periods must be between \( \Tl \) and \( \Tu \):
 \[
 \g\bigparen{\var{dw-p-bd},\; W_{0}(x,a),\; \eta} = \f_{0}(x,a,\eta)
  \h(\var{dw-p-bd},x,a,\eta)
 \]
  where \( \h(\var{dw-p-bd},x,a)=1 \) if \( \eta \) has a dwell period shorter than
\( \Tl \) in \( W_{0}(x,a) \) or has a dwell period longer than \( \Tu \) there
(and, as always, 0 otherwise).

 \item\label{i:comput.legal-comp}
  Whenever \( W_{0}(x,a) \) is damage-free, the transition at
\( \sigma_{2} \) is a legal one:
 \[
   \g\bigparen{\var{legal-comp},\; W_{0}(x,a),\; \eta} 
 = \f_{0}(x, a, \eta)
    \{ \exists\sigma_{2}\land\neg\legal(\eta(x,\sigma_{2}-), \eta(x,\sigma_{2})\}.
 \]

 \item\label{i:comput.legal-birth}
   Whenever \( W_{1}(x,a) \) is damage-free and \( \sigma_{0} \) exists,
\( x \) either has a creator or is emerging (see Definition~\ref{def:creator-emerging}).
 \[
   \g\bigparen{\var{newborn},\; W_{1}(x,a),\; \eta} = \f_{1}(x, a, \eta)
 \{\exists\sigma_{0}\land\neg\h(\var{newborn}, x, a, \eta) \}
 \]
 where \( \h(\var{newborn}, x, a, \eta)=1 \) if \( x \) has a creator for time 
\( \sigma_{0} \), or \( \eta(x,\sigma_{0}) \) is emerging.
Also, the transition in \( \sigma_{0} \) is a legal one:
 \[
   \g\bigparen{\var{legal-birth},\; W_{0}(x,a),\; \eta} = \f_{0}(x, a, \eta)
    \{\exists\sigma_{0}\land\neg\legal(\eta(x,\sigma_{0}-), \eta(x,\sigma_{0}))\}.
 \]

 \item\label{i:comput.trans}
  Whenever \( W_{1}(x,a) \) is damage-free, the transition
  function applies, based on observation at some
  (unspecified) time point during the \df{observation interval} (``atomicity'')\index{atomicity}:
 \[
   \g\bigparen{\var{trans},\; W_{1}(x,a),\; \eta} = \f_{1}(x, a, \eta)
  \{\exists\sigma_{2}\land\neg\h(\var{trans}, x, a, \eta)\}
 \]
 where \( \h(\var{trans}, x, a, \eta)=1 \) if there is a
 \( t'\in \rint{\sigma_{1}}{\sigma_{2}-\Tl/4} \) with
 \[
    \eta(x,\sigma_{2}).\Det = \trans(\eta,x,t').\Det.
 \]

  \item \label{i:comput.creation}
  A cell cannot stay vacant if it has would-be creators for a
long time and has no neighbor potentially blocking the creation:
 \[
   \g\bigparen{\var{no-birth},\; W_{1}(x,a),\; \eta} = \f_{1}(x, a, \eta)
 \h(\var{no-birth}, x, a, \eta)
 \]
  where \( \h(\var{no-birth}, x, a, \eta)=1 \) if the following conditions hold:
   \begin{romanenum}
    \item \( \eta(x,t) \) is vacant for all \( t\in\clint{a}{a+3\Tu} \);
    \item for each \( t\in \clint{a}{a+2\Tu} \) there is a
\( j\in\{-1,1\} \) such that \( (j,\thg_{j}(x,t)) \) is a potential creator of some
(non-vacant) state;
\item\label{i:obstacle-c}
  there is no non-vacant \( \eta(y,t) \) with \( 0<|y-x|<\B  \), \( t\in\clint{a}{a+3\Tu} \).
   \end{romanenum}
  \end{cjenum}
\end{condition}

Our conditions do not require the state of an emerging cell to be
completely determined by the transition function (other than that it is latent,
as said in Condition~\ref{cond:latent}).
In fact as we will see, during self-organization, the state of an emerging
big cell may depend on the information in the germ cells creating the colony
representing it.

Let us take a peek ahead to self-simulation, since
at this point, condition~\ref{cond:comput}\ref{i:comput.creation} seems too
strong.
Indeed, suppose that a colony of \( \Q  \) cells of size \( \B  \) simulates a higher-order
cell of size \( \Q \B  \) in a medium \( \eta^{*} \) in such a way that these higher-order
cells must obey the conditions of a robust medium again.
Now assume again that \( \eta^{*}(x,t) \) is vacant for a long time,
and for all this time, there is a potential creator in \( x-\Q \B  \), while
there is no \( (y,t) \) with \( 0<|y-x|<\Q \B  \) with non-vacant \( \eta^{*}(y,t) \).
Can then indeed a new colony be created in \( x \)?
There could in principle be another obstacle, namely that there is a 
\( y \) with \( 0<|y-x|<\Q \B  \) such that there is a colony during all this
time starting in \( y+\Q \B  \) (that is \( \eta^{*}(y+\Q \B ,t) \) is nonvacant) and is trying to
create a colony with starting point \( y \).
The attempts of colony creation from left and right interfere.
Our solution is  simple: growth from the left will get priority over growth from
the right.
In this way, one of the conflicting attempts will succeed.

Though it is desirable to illustrate the computation property on a
number of special cases, currently, we will give only the barest minimum of these.

\begin{example}[Special cases]
To obtain a deterministic cellular automaton
as a special case of robust media,
set \( \Tl=\Tu \), \( r=1 \), \( \eps=0 \), require
that the histories have empty damage and that the
transition function not give a vacant value.

\begin{sloppypar}
The connection between primitive variable-period media and robust
media will be set up via a trivial simulation.
Let \( M_{1}=\Prim-var(\states_{1}, \sites, \trans_{1},\B ,\Tl,\Tu,\eps) \).
We make the additional assumptions that
the transition function \( \trans_{1} \) 
satisfies the time marking Condition~\ref{cond:time-mark}.
We define a simple simulation \( \Phi^{*} \) by this medium of the robust
medium
 \[
  M_{2}=\Rob(\states_{2}, \sites, \trans_{2}, \B , \Tl, \Tu, \eps, \eps, 1).
 \]
Set \( \states_{2} = \states_{1}\times \{0,1\} \), 
identify \( \states_{1} \) with \( \states_{1}\times\{0\} \), and denoting
\( \states'_{1}=\states_{1}\times\{1\} \), let
\( \Vacant_{2}=\Vacant_{1}\cup\states'_{1} \),
\( \Bad_{2}=\Bad_{1}\cup\states'_{1} \).
Thus, some new bad states are introduced by turning on a new bit which we
can call the ``bad bit''.
The transition function \( \trans_{2} \) is essentially the same as
\( \trans_{1} \): it ignores the bad bit. 
It also satisfies Condition~\ref{cond:cling} 
trivially: it never destroys or creates cells.  
\end{sloppypar}

Let \( g_{1} \) be the event function used in the definition of medium \( M_{1} \).
We give the definition of \( s=\eta^{*}(x,t)=\Phi^{*}(\eta)(x,t) \) for each
\( \eta \).
If \( x \) is not an integer multiple of \( \B  \) then \( s=\Vac \).
If \( \g_{1}\bigparen{\alpha,\; \{x\}\times\rint{t'}{t},\; \eta}=1 \) for some \( t'<t \)
and some \( \alpha\ne\var{coin} \) then \( s=(\eta(x,t),1) \) (as some condition is violated, the
bad bit is turned on), else \( s=(\eta(x,t),0) \).

It can be verified that this is indeed a simulation.
\end{example}

Peeking ahead to simulations, 
Lemma~\ref{lem:sim-damage} proves the Restoration Property
(Condition~\ref{cond:restor}) for
\( \eta^{*}=\fg^{*}(\eta) \) whenever \( \eta \) is a trajectory of \( M \).
When the computation property above is
applied to a cell \( x \) of \( \eta^{*} \), we are
given a rational number \( a \) satisfying \( \f_{1}(x,a,\eta)=1 \).
Since we have to prove this property for \( \eta^{*} \),
it will be assumed throughout the rest of the construction that
\( \Damage^{*} \) does not intersect the area we are reasoning about.
Hence all damage belongs to islands as in Definition~\ref{def:damage-map}, and for simplicity,
we will identify the damage with the islands (that is assuming it fills them out).

 \section {Amplifiers} \label{sec:amp}

In the present section, when we talk about media without qualification we
understand robust media.
Recall amplifiers from Definition~\ref{def:amp}.
We will also need a trickle-down property, as in Definition~\ref{def:amp-trickle}.
Then a lemma analogous to Lemma~\ref{lem:abstr-amp} can be formulated, and
a proof for Theorem~\ref{thm:1dim.nonerg.var} can be provided
analogously to the proof of Theorem~\ref{thm:1dim.nonerg}.

Our eventual goal is to find an amplifier $(M_{k},\Phi_{k})$,
where each medium $M_{k}$ will also simulate a
some computation with a given transition function, and damage probability bound
$\eps_{k}$, for a fast decreasing sequence $\eps_{k}$.
We introduce some parameters of the sequence
$(M_{k})$ of media that will satisfy conditions sufficient for the existence of an
amplifier.
These conditions are similar to the inequalities in the conditions
of Theorem~\ref{thm:simp-sim}.

Recall from Definition~\ref{def:rob} that
all robust media $M_{k}$ have the fields $\Det$ and
$\Rand$ whose behavior is governed by the computation and restoration
properties.

\begin{notation}\label{n:det-subfield}
Only field $\Det$ will be subdivided into subfields: therefore subfields
$\eta.\Det.\F$ will be denoted simply as $\eta.\F$, without danger of
confusion.  
\end{notation}

\begin{definition}[Amplifier complex]\label{def:amp-complex}
Consider an amplifier \( (M_{k},\Phi_{k}) \) where \( \Phi_{k}=(\fg_{k*},\Phi^{*}_{k}) \).
For a sequence of fields \( (\Payload_{k}) \) and parameters \( \D_{k},\eps''_{k}\),
our amplifiers will use a code complex
  \begin{align*}
   \hierCode = (\fg_{k},\Payload_{k}, \gamma_{k})_{k\ge 1}
  \end{align*}
as in~\eqref{eq:H},
and they will have in addition the trickle-down property as in Definition~\ref{def:amp-trickle}
for the fields \( \F_{k}=\Payload_{k} \).
We will call such an amplifier an \df{amplifier complex}.
\end{definition}

In any trajectory \( \eta \) of \( M_{1} \)  at space-time point \( (x,t) \)
the hierarchy might only be built up to a certain highest level \( k \) which may increase
as time proceeds.
The cells on this level will have the kind \( \Germ \) (on the field \( \Kind \) see later
in Section~\ref{sec:plan}), and the useful computation of $M_{k}$ will be carried out 
on their \( \Payload_{k}\) track.
This computation will be an aggregated version \( \PlTrans_{k} \)
of a certain fixed cellular automaton with transition function \( \PlTrans_{1} \).
It will have monotonic output in the sub-track \( \Output_{1} \), and correspondingly
\( \PlTrans_{k} \) will have monotonic output on \( \Output_{k} \).

\begin{condition}\label{cond:Payload-invar-hist}
We require  transition function \( \PlTrans_{1} \) to be commutative in
the sense of Definition~\ref{def:commutative}.
\end{condition}

As shown in Theorem~\ref{thm:asynch-sim}, this requirement is not strong:
the transition can still be, for example, universal; but it
makes its simulation simpler as we don't have to worry about the order of updating (for
about the effect of delays due to a temporary lack of cells at the edges).
The trickle-down property will trickle down the
latest result of \( \PlTrans_{k} \) to the lowest level.

\begin{sloppypar}
\begin{definition}[Carrying the payload]\label{def:carries}
In a robust medium $M$, recall the transition function \( t=\trans(\br,\ba) \)
from~\eqref{eq:rob-trans}. with input vector $\br=(r_{-1},r_{0},r_{1})$ and
the adjacency bits \( \ba=(a_{-1},a_{1}) \).
We will say that it 
 \df{carries}\index{transition function!carries} the payload transition
function $\PlTrans$\glo{rd-trans@$\PlTrans$} if the condition
\( t\ne\Vac \), \( r_{0}.\Kind=\Germ\),  \( \ba=(1,1) \) implies
\( t.\Payload = \PlTrans(\br.\Payload) \).
  \end{definition}
  \end{sloppypar}

Thus, the payload field is controlled by its own transition function in a germ
provided that the simulation does not command to eliminate it.
In a non-germ cell, its value will trickle down from the large cell its colony simulates.

The following concept will also be needed for the definition of amplifiers.
  Here, and always from hence, when we refer to a natural number $k$
as a binary string we always mean its standard binary representation.

\begin{definition}[Uniform program]
  The string $\var{prog}$ will be called a
 \df{uniform program (with complexity coefficient
$c$)}\index{program!uniform} for the functions $f_{k}$ if there is some
constant $c$ such that it computes on the universal computing medium
$\univ$ the value $f_{k}({\br})$ from $k,{\br}$ with space- and
time-complexities bounded by $c(\log k + |\br|)$ where $|\br|$ is the total
number of bits in the arguments $\br$.
  If a sequence $(f_{k})$ of functions has a uniform program it is called
\df{uniform}.
As a special case, a sequence of constants $c_{k}$ is uniform if it is
computable with space- and time complexities $c\log k$ from the index $k$.
  \end{definition}

Here are the parameters involved.
 The definitions
 \begin{align*}
      \states_{k} = \{0,1\}^{\cp_{k}}, \quad \B_{k}= \prod_{i<k} \Q_{i}
 \end{align*}
 are not new: \( \cp_{k} \) is the number of bits in a state in \( \states_{k} \),
 and \( \B_{k+1} \) is the size of a colony
 consisting of \( \Q_{k} \) cells of size \( \B_{k} \).
 Recall that \( \Tl_{k}<\Tu_{k} \) are the lower and upper bound on the dwell periods.
Let
 \begin{align}\label{eq:D_k}
   D_{k}=2\Tu_{k}
 \end{align}
for the parameters \( D_{k} \) in trickle-down (Definition~\ref{def:amp-trickle}). 
For simplicity, we will require
 \begin{equation}\label{eq:lg-less-2}
  \lambda:=\sup_{k}\Tu_{k}/\Tl_{k}<2.
\end{equation}
The number
\begin{align}\label{eq:U'}
 \U'_{k}
\end{align}
approximates the (variable) number of dwell periods
of \( M_{k} \) in a work period of simulation of a cell of \( M_{k+1} \),
(only approximates, as the size of dwell periods is not fixed).
For a constant 
\begin{align}\label{eq:rep}
 \rep>0  
\end{align}
to be fixed large enough later, let
\begin{align}\label{eq:Tl-Tu}
  \Tl_{k+1}/\Tl_{k}=\U'_{k}(1-\rep\Q_{k}/\U'_{k}),\quad
  \Tu_{k+1}/\Tu_{k}=\U'_{k}(1+\rep\Q_{k}/\U'_{k}).
\end{align}
We will see that \( \U'_{k}/\Q_{k} \) grows fast, and
thus both \( \Tl_{k} \) and \( \Tu_{k} \) grow almost like \( \prod_{i<k}\U'_{i} \),
only \( \Tl_{k} \) grows slightly slower, and \( \Tu_{k} \) slightly faster.
Now~\eqref{eq:lg-less-2} can be achieved
for example by setting $\U'_{k}\ge c k^{2}\Q_{k}$ for some large $c$.

Due to the aggregation property and Proposition~\ref{propo:small-bandwidth},
only a small fraction of the capacity of each cell can be used
for computation and communication during a work period.
The new parameter
\begin{align}\label{eq:bw}
   \bw_{k}
\end{align}
called \df{bandwidth rate} will approximately indicate this fraction. 
We assume that
 \begin{equation}\label{eq:Pdef}
   \nP_{k}=1/\bw_{k}
\end{equation}
and \( \bw_{k}\Q_{k-1} \) are integers, and therefore so is \( \bw_{k}\B_{k} \).
The payload transition function on level \( k \) will be an aggregated version
of \( \PlTrans_{1} \) with slowdown, as in Example~\ref{xmp:sim.aggreg}:
  \begin{align*}
   \PlTrans_{k}=\PlTrans_{1}^{\B_{k}, \bw_{k}}.
 \end{align*}
 Thus, a cell of \( M_{k} \) of body \( \B_{k} \),
 simulated by \( \B_{k} \) cells of \( M_{1} \), should contain all the
 payload of the simulating cells,
 and the payload transition function will be slowed down by a factor \( \bw_{k} \).

 The colony will be processed in chunks of \( \bw_{k}\Q_{k} \) cells,
 and we want to be able to correct a burst of constant number off cells in each of these.
 This requires that information be stored with a redundancy of at least a few cells per chunk.
 So we devote a portion
  \begin{align*}
   \red_{k}=\rep \bw_{k}/\Q_{k}
\end{align*}
of the information to redundancy.
 Let us fix some arbitrary positive initial constants:
 \begin{align*}
      \Tl_{1}\le &\Tu_{1}\le 1.5\Tl_{1},
\\   \eps'_{1} &<1/\rep^{2},
\\      \eps &=\eps_{1}.
 \end{align*}
The interpretation of the parameters  \glo{eq:greek@$\eps_{k},\eps'_{k},\eps''_{k}$}%
\( \eps_{k},\eps'_{k},\eps''_{k} \) has been seen before.
In the definition of robust media in
Section~\ref{sec:rob}, \( \eps_{k} \) bounds probability of occurrence of damage,
and \( 1/2+\eps'_{k} \) bound the probability of an outcome of the random coin.
The bound \( \eps''_{k} \) has been seen in connection with trickle-down in
Definition~\ref{def:amp-trickle}.
Here are the recursive definitions:
\begin{equation}
  \label{eq:eps-defs}
  \begin{aligned}
  \eps_{k+1}    &= c_{\Dam}(\Q_{k}\U'_{k}\eps_{k})^{2} ,   
\\   \eps''_{k}     &= 4\Q_{k}\U'_{k}\eps_{k}, 
\\   \eps'_{k}      &= \eps'_{1}+\sum_{i=1}^{k-1}\eps''_{i}. 
  \end{aligned}
\end{equation}
The formula for $\eps_{k+1}$ is natural in view of
Lemma~\ref{lem:sim-damage}, which introduces the constant $c_{\Dam}$.
The definition of $\eps'_{k+1}$ takes into account the limited
ability to simulate a coin-toss with the help of other coin-tosses.
The parameter $\eps''_{k}$ bounds
the probability that there is any $k$-level damage at all during the work
period of a colony of $k$-cells.

Let us choose
 \begin{equation}\label{eq:amp-fr}
   \begin{aligned}
     \bw_{k}   &=  1/\rep k^{2}, 
     \\    \Q_{k}     &=  \rep^{k+1},
     \\  \U'_{k}    &=  \rep\Q_{k}/\bw_{k}, 
   \\      \cp_{k} &= \B_{k}(1+\red_{k}),
   \end{aligned}
 \end{equation}
 where we omitted integer parts.
 Let us by induction, prove for sufficiently small \( \eps \):
\begin{equation}\label{eq:exp-error}
  \begin{aligned}
   \eps_{k} &\le \eps^{2^{k-2} + 2^{(k-3)/2}}.
  \end{aligned}
 \end{equation}
  For $k=1$, the statement gives $\eps\le\eps$.
  For $k>1$, using 
  the inductive assumption,
 \begin{equation*}
   \begin{aligned}
    \Q_{k}\U'_{k}&=\rep^{2k+4}k^{2}
\\   \eps_{k+1} &\le \eps^{2^{k-1} + 2^{(k-1)/2}}\cdot c_{\Dam}\rep^{4k+8}k^{4}
     \\ &= \eps^{2^{k-1} + 2^{(k-2)/2} }\cdot
          c_{\Dam} \rep^{4k+8}k^{4}\eps^{2^{(k-2)/2}(2^{1/2}-1)}.
 \end{aligned}   
 \end{equation*}
 For small enough $\eps$, the last factor is less than 1.
From here, it is easy to see that $\eps''_{k}$ also converges to 0 with similar speed.
 
The parameters introduced satisfy the requirements below, which can be compared
to the corresponding conditions for Theorem~\ref{thm:simp-sim}.

 \begin{description}

  \item[Complexity Upper Bounds]\index{upper@Upper Bound!Complexity}
All parameters in $\Frame$ are uniform
sequences with complexity coefficient $\rep$.

  \item[Capacity Lower Bound]
    Cells should be able to hold
    numbers comparable to the size of the colony and the work period within a cell,
and the colony must represent the state of the big cell with some redundancy:
  \begin{align}
\label{eq:cap-lbd-1}          \cp_{k} &\ge \rep\log \U'_{k},
\\\label{eq:cap-lbd-2} \Q_{k}\cp_{k} &\ge \cp_{k+1}(1+\red_{k+1}).
  \end{align}

\item[Bandwidth Bounds]\index{lower@Lower Bound!Bandwidth}
\begin{align}
  \label{eq:bw-lb}
  \bw_{k}\cp_{k} &\ge \rep\log \U'_{k}, 
  \\ \label{eq:bw-b2} \bw_{k}\Q_{k}\cp_{k}&> \bw_{k+1}\cp_{k+1},
  \\\label{eq:bw-b3}  \bw_{k+1}\Q_{k} &\ge \rep(1/\bw_{k}+1/\bw_{k+1}).
\end{align}
 \eqref{eq:bw-lb} says that numbers comparable to the size of the
colony and the work period should fit within the space bound $\bw_{k}\cp_{k}$.
\eqref{eq:bw-b2} says that a track of a colony of relative width \( \bw_{k} \) should be able
to hold the portion \( \bw_{k+1} \) of the state of the represented big cell, for processing.
The easily satisfiable inequality~\eqref{eq:bw-b3} (as \( \Q_{k} \) will grow much faster
than \( 1/\bw_{k} \)) will be used in Section~\ref{sec:code-decode}.

  \item[Work Period Lower Bound]\index{lower@Lower Bound!Work Period}
  There must be enough dwell periods in a work period to perform the
necessary computations of a simulation:
 \begin{equation}\label{eq:Ubd} 
   \U'_{k} \ge \rep\cdot \Q_{k}/\bw_{k}.
 \end{equation}

  \item[Error Upper Bound]\index{upper@Upper Bound!Error}
  The following bound relates colony and work-period size to
the error probability, and will help in the recursive error estimates:
 \begin{equation}\label{eq:error-ub}
  \eps_{k}^{0.2} \le 1/\rep\Q_{k}\U'_{k}.
 \end{equation}
It will be easy to satisfy in our examples, since
$\Q_{k},\U'_{k}$ will grow only exponentially and $\eps_{k}$ will decrease
super-exponentially.
 \end{description}

 \begin{lemma}\label{lem:frame-conds}
   The parameters satisfy the conditions above.
 \end{lemma}
 \begin{proof}
 Conditions~\eqref{eq:cap-lbd-1}, \eqref{eq:bw-lb}, \eqref{eq:bw-b3} and~\eqref{eq:Ubd}
 are satisfied for large 
 \( \rep \).
 For inequality~\eqref{eq:cap-lbd-2}, assuming \( \rep>2 \),  note first that
 \begin{align*}
   \frac{\red_{k+1}}{\red_{k}}=\frac{\Q_{k}\bw_{k}}{\Q_{k+1}\bw_{k+1}}
   =\frac{(k+1)^{2}}{k^{2}\rep}.
\end{align*}
  \begin{align*}
  \cp_{k+1}(1+\red_{k+1})&=\B_{k+1}(1+\red_{k+1})^{2}
    \\      &= \B_{k+1}(1+\red_{k+1}(2+\red_{k+1}))
    \\       &= \B_{k+1}(1+\red_{k}\frac{(k+1)^{2}}{k^{2}\rep}.(2+\red_{k+1}))
  \end{align*}
  The multiplier of \( \red_{k} \) on the right-hand side is smaller than 1 for large \( \rep \).
Condition~\eqref{eq:bw-b2} reduces to 
\begin{align*}
  1+\red_{k} &> (1 + \red_{k+1})k^{2}/(k+1)^{2},
\end{align*}
which is true.
    \end{proof}

 For later reference, here are explicit expressions for \( \B_{k} \) and for
 a quantity falling between $\Tl_{k}$ and $\Tu_{k}$:
 \begin{equation}\label{eq:B_k-expl}
\begin{aligned}
  \B_{k} &= \prod_{i < k} \Q_{i} = \rep^{k(k+1)/2-1},
  \\  \Tu_{k}&\eqm \prod_{i < k} \U'_{i} = \B_{k}\rep^{k-1}((k-1)!)^{2}
               = \rep^{(k+1)(k+2)/2-2}((k-1)!)^{2}.
 \end{aligned}    
 \end{equation}

\begin{lemma}[Amplifier] \label{lem:amp}\index{lemma@Lemma!Amplifier}
  For large enough $\rep$ there is a uniform amplifier complex with the
  parameters introduced above, which also obeys the following conditions:
  \begin{cjenum}

  \item  Medium 
 \begin{equation}\label{eq:M_{k}}
   M_{k}=\Rob(\trans_{k},\B_{k},\Tl_{k},\Tu_{k},\eps_{k},\eps'_{k}).
 \end{equation}
 is robust.

  \item The damage map of the simulation $\Phi_{k}$ is given as in
Definition~\ref{def:damage-map}.

\item The function $\trans_{k}$ carries the payload
  transition function $\PlTrans_{k}$ in the sense of Definition~\ref{def:carries}.

\end{cjenum}  
 \end{lemma}

 Lemma~\ref{lem:amp} implies a variable-period
generalization of Lemma~\ref{lem:abstr-amp} with $D_{k}=\Tu_{k}$,
which, via essentially the same proof as the one after
Lemma~\ref{lem:abstr-amp}, implies Theorem~\ref{thm:1dim.stor.var}
for the infinite space.
All properties but the initial stability property of an abstract
amplifier (defined in Lemma~\ref{lem:abstr-amp}) are satisfied by definition.
For the initial stability property it is sufficient to note that each
medium $M_{k}$ is a robust medium, with work period bounds $\Tl_{k},\Tu_{k}$.
Therefore, if $\eta^{k}$ is a trajectory of $M_{k}$ and $t<\Tu_{k}$ then
for each $x$ the probability that $\eta^{k}(x,t)\ne\eta^{k}(x,0)$ is less
than the probability that damage occurs in $\{x\}\times\rint{0}{\Tu_{k}}a$.
This can be bounded by the Restoration Property.

\section{Outline of the program} \label{sec:plan}

Our eventual goal is to prove Lemma~\ref{lem:amp} (Amplifier).
From now on, we fix the level \( k \) of the simulation hierarchy,
refer to medium \( M_{k} \) as \( M \) and to \( M_{k+1} \) as \( M^{*} \).
The subscript will be deleted from all parameters of \( M \), and a
superscript \( * \) will be added to all parameters of \( M^{*} \).
We will refer to cells of \( M \) as \df{small cells}, or simply \df{cells},
and to cells of \( M^{*} \) as \df{big cells}.
The program will be described in a semi-formal way; the present section
overviews it.
Later sections restrict the program more and more by giving 
some rules and conditions, and prove lemmas along the way.
The typical condition would say that certain fields can only be changed
by certain rules.
The language for describing the rules is an extension of the one given in
Section~\ref{sec:simp-sim}.
We will introduce a fair number of fields but they are all relatively
small.

Let us fix some conventions on the handling of fields and the different parts of
the transition function.
According to Section~\ref{sec:amp}, the medium \( M \) has a
special field called \( \Payload \).
Together with some other fields, it can be grouped into a field
we could call\( \Info \): in a colony, the \( \Info \) track
contains the string that encodes the state
of the represented cell via an error-correcting code.
Some other fields like \( \Addr \), will service the functioning
of the simulation program.

 \begin{definition}[Info track]\label{def:info-field}
   Field \( \Info \) consists of the subfields \( \Info.\Payload \), \( \Info.\PlRedun \),
   \( \Info.\Util \) (referred to mostly as just \( \Payload \), \( \PlRedun \), \( \Util \)).
   \glo{info.main@\( \Info.\Payload, \Info.\PlRedun, \Info.\Util \)}
The track \( \Payload \) contains the intended original information,
The track \( \PlRedun \) contains ``error check bits'' for the track \( \Payload \)
as in Example~\ref{xmp:err-corr-code},
while the \( \Util \) track contains other parts of the represented information of the
simulated colony (along with its own error checks).
 \end{definition}

In what follows we define the function \( \trans \): the
function \( \PlTrans \) is given in advance.
The field \( \Payload \) is updated via Definition~\ref{def:carries}:
thus, \( \trans \) determines the next value of \( \Payload \) only in
case of non-germ cells via trickle-down.

For the error-correcting code, the track \( \Payload\cup\PlRedun \) will be subdivided
into ``packets'' such that the parity check bits for each packet will be
computed separately.
This is done because the \( \Payload \) track will be
much wider than the work tracks, so it cannot be processed in its entirety at once on them.

Cells have an address field \( \Addr \) which determines the only \df{colony}
(\( \Q  \)-colony) to which the cell belongs.
A colony \( \cC(y) \)\glo{c.cap@\( \cC(y) \)} has base \( y \).
The \( \Age \) field of a cell, called its \df{age}, can have values in
\( \lint{0}{\U } \), where \( \U =\U'/\p_{1} \) (which we can assume to be integer without loss
of generality).
Here \( \U' \) is the parameter \( \U'_{k} \) introduced in~\eqref{eq:U'}, and \(
\p_{1} \) is a constant to be introduced in Definition~\ref{def:timing} below.

\subsection{Cell kinds}\label{sec:cell-kinds}

\begin{definition}[Adjacency]\label{def:adjacency}
Recall the transition function \( \trans(\br,\ba) \) where \( \ba=(a_{-1},a_{1}) \) says
whether the left and right neighbor is adjacent.
In our rule language we will refer to \( a_{j} \) as 
\begin{align*}
 \adjNb(j) .
\end{align*}  
\end{definition}

We will continue to understand by the word \df{cell} a
non-vacant, non-damaged site.

 \begin{definition}[Colonies]\label{def:colonies}
The values of the address field \( \Addr \) will vary in \( \lint{-\Q }{2\Q } \).
The \df{colony} of a cell \( x \) is the set \( C=\setOf{z+i\B}{i\in\lint{0}{\Q}} \)
whose starting cell is \( z = x-(\Addr(x) \bmod \Q ) \B  \).
We will also say that \( x \) \df{belongs} to colony \( C \).
The \df{originating colony} of a cell \( x \)
is the colony \( D \) whose starting cell is \( x-\Addr(x) \B  \).
We will also say that \( x \) \df{originates} at colony \( D \).
The cells whose colony is their originating colony, will be called
\df{inner cells}\index{cell!inner}, the other ones will be called
\df{outer cells}\index{cell!outer}.
Two neighbor cells will be called \df{space-consistent} if they originate at the same colony.
 \end{definition}

A certain property of cells, called their ``kind'', plays an important role in
determining their behavior.

 \begin{definition}[Cell kinds]\label{def:cell-kinds}
Cells will be of a few different \df{kinds}\index{cell!kind}, distinguished by
the field%
 \glo{kind@\( \Kind \)}%
 \[
  \Kind.
\]
The possible kinds are listed in~\eqref{eq:strength} below.
 \glo{latent@\( \Latent \)} \glo{germ@\( \Germ \)}\glo{channel@\( \Channel \)}\glo{growth@\( \Growth \)}
 \( \Member \)\glo{member@\( \Member \)} 
 \index{cell!vacant}%
 \index{cell!latent}%
 \index{cell!germ}%
 \index{cell!channel}%
 \index{cell!growth}%
 \index{cell!member}%
The kind of non-latent, non-germ, cells is determined by \( \Age \) and \( \Addr \).
Since it will be possible to determine from the age of a cell whether it
has kind \( \Channel \) or \( \Growth \),
We will sometimes use therefore just one value%
 \glo{ext@\( \Ext_{j} \)}%
 \[
  \Ext_{j}
 \]
in place of \( \Growth_{j} \) and \( \Channel_{j} \).
The kinds of cells are ordered by a property called \df{strength}%
\index{cell!strength}
as follows:
\begin{equation}\label{eq:strength}
\begin{aligned}
   \Vac &< \Latent < \Germ < \Channel_{-1}<\Channel_{1}
\\       & < \Growth_{-1} < \Growth_{1} < \Member.
\end{aligned}
\end{equation}
 If a cell needs to be created, overlapping in body a another one then if the
 new cell is stronger, the weaker one will be erased.
 This is the only way a cell can turn vacant see the later
 Condition~\ref{cond:cling} (Cling-to-Life).
Let us discuss each kind.
 \begin{itemize}
  \item
The relation \( \Kind(x,t)=\Vac \) means that there is no cell at
site \( x \) at time \( t \).
  \item
A cell will be called \df{dead}\index{cell!dead} if it is vacant or
latent, and \df{live}\index{cell!live} otherwise.
\df{Killing}\index{killing} a cell \( x \) means turning it latent. 
\df{Erasing}\index{erasing} it means turning it vacant. 

  \item The \df{member} cells are ``inner'' cells as defined above;
they are the strongest kind in order to maintain the integrity of colonies.

  \item Cells of kind \( \Ext_{j} \) are outer cells.
A \df{right outer} cell has addresses \( \ge \Q  \), and a \df{left
outer} cell has addresses \( <0 \).

  \item Left channel cells are weaker than right channel cells which are weaker than growth cells.
This way, a right-growing channel will erase the left-growing
one if it is in its way, and information need not be transmitted between two channels.
Also, left channel cells can grow only until address \( -\Q +1 \), and right ones
only until \( 2\Q -2 \): thus a channel cannot be extended to cover a full colony.

\item \df{Germ} cells exist at the highest level of the hierarchy, in the sense that
  they are not part of any colony simulating some higher-level cell.
  They also have the function to carry out a computation, the \df{payload}
  of the whole construction.
The member cells on lower levels will receive the results of this computation by a
feature of the simulation called \df{trickle-down}.
Germ cells have addresses in \( \lint{-2\Q }{3\Q } \).
 \end{itemize}
\end{definition}

In the construction leading to Theorem~\ref{thm:1dim.comp},
germ cells will play the following role.
Adjacent and ``consistent'' germ cells of the same color will 
attempt to expand and form a new colony
simulating a big germ cell: we call this ``lifting'', or ``self-organization''.
The part of the space marked with color 0 is where the actual payload computation
takes part.
The parts of space to its left and right are marked by color \( -1 \) and \( 1 \) respectively.
These pars are also populated by germ cells which also self-organize (this
needs randomization to break the translational symmetry).
The germs of the outside areas, with their appropriate colors,
will suppress possible structures that would compete with the computation area.

The transitions between different kinds are limited in advance. 

 \begin{condition}[Latent Cells]\label{cond:vac-to-lat}
 \index{condition@Condition!Latent Cells}
  A vacant cell can only turn into a latent one.
\end{condition}

The following condition helps enforce that newly created cells in the
simulation are latent.
  
 \begin{condition}\label{cond:outer-info}
Suppose that \( \fg^{*}(\eta)(x,t)\notin\Bad \), the colony with base \( x \)
at time \( t \) is full and is covered with member cells belonging to the same
work period, and is not affected by damage (as in Definition~\ref{def:affecting}). 
Then \( \fg^{*}(\eta)(x,t) \) depends only on the \( \Info \) track of this
colony via some decoding function \( \alpha^{*} \).
If the colony is covered with germ or outer cells then the state
decoded from this track is latent.
 \end{condition}

 \subsection{A colony work period}

Algorithm~\ref{alg:workperiod} describes the main stages of a colony work
period.
 \glo{extend@\( \Extend \)}%
 \glo{send@\( \Send \)}%
 \glo{compute@\( \Compute \)}%
 \glo{grow@\( \Grow \)}%
 \glo{finish@\( \Finish \)}%
 In the notation here, it sometimes does not matter much whether we call one of the rules
 here a ``sub-rule'' or not, and whether we print the semicolon,
 as the rule (like \( \Extend \)) includes the conditions under which it is applied.
 \begin{algo}[H]
 \caption{A colony work period}\label{alg:workperiod}
\PrintSemicolon
  \( \Extend\prl\Retrieve \);
\;    \( \Compute \); 
\;    \( \ProcPayload \); 
\;    \( \Grow \);
\;    Shrink;
\;   \( \Finish \)
\end{algo}
 There are some idle stages between these parts to make sure that
faults in one part have limited effect on other parts.
Rule \( \Extend \) is defined in Algorithm~\ref{alg:Extend},
 \( \Retrieve \)  in Algorithm~\ref{alg:Retrieve-2}, 
 \( \Compute \)  in Algorithm~\ref{alg:Compute}, 
 \( \ProcPayload \)  in Algorithm~\ref{alg:Proc-payload}, 
 \( \Grow \)  in Algorithm~\ref{alg:Grow},
 \( \Finish \) in Algorithm~\ref{alg:Finish}.
 Shrinking (if needed) happens via the rule \( \Decay \)
 defined in Algorithm~\ref{alg:Decay}.

 The stages are started and ended at certain specific values of the field \( \Age \),
 which runs between 0 and \( \U \), the maximum age.
 Recall for some \( k \) we are defining medium \( M_{k} \) of and amplifier.
Because of a delay parameter \( \p_{1} \) to be defined below in~\eqref{eq:p_{1}},
parameter \( \U=\U_{k} \)  is related to the estimated number \( \U'_{k} \)
of dwell periods in a work period introduced  in~\eqref{eq:U'} as follows:
 \begin{align}\label{eq:U-def}
   \U'_{k}=\p_{1}\U_{k}.
 \end{align}
Here are the ages starting or ending important stages of this program:
 \glo{compute-start@\( \computeStart \)}%
 \glo{payload-start@\( \payloadStart \)}%
 \glo{grow-start@\( \growStart \)}%
 \glo{grow-end@\( \growEnd \)}%
  \begin{equation}\label{eq:compute-start}
    \begin{aligned}
  \computeStart &= \konst\rep\nP^{*}\Q,
\\ \payloadStart &=  \computeStart+ \konst\rep\nP^{*}\Q , 
\\ \growStart &=  \payloadStart+ \rep \nP^{*}\Q, 
\\ \growEnd     &= \growStart + 6\Q \lambda, 
\\ \U              &=\growEnd+4\Q , 
    \end{aligned}
  \end{equation}
  where \( \lambda \) was defined in~\eqref{eq:lg-less-2},
\( \U \) is the number of \( \Age \) steps in the work period,
\( \konst \) is a constant to be specified later, the constant \( \rep \) was introduced in
Section~\ref{sec:amp}, and \( \nP \) was defined in~\eqref{eq:Pdef}.

\begin{sloppypar}

  \begin{definition}[Expansion and retrieval periods]\label{def:expansion-period}
The interval before age \( \computeStart \) is called the \df{retrieval period}: then
the colony retrieves information from neighbor colonies.
The constant \( \konst \) will make it substantially longer than the
transition period \( \rint{\growStart}{\U} \) at the end of the work period when the information to be
sent to neighbor colonies would change.
The age intervals \( \rint{0}{\computeStart} \) and 
\( \rint{\growStart}{\growEnd} \) for non-germ outer cells
and \( \rint{0}{\growEnd} \) for germ cells   
will be called \df{expansion periods}.
\index{period!expansion}\index{cell!expansion}
Cells in their expansion period are called \df{expansion cells}.
 \end{definition} 
\end{sloppypar}

The following short description of the major stages of a colony work
period should serve for orientation.
The part of the work period before \( \Compute \) sees the following activities
happening simultaneously: \( \Extend \), and \( \Retrieve \).

\begin{description}

 \item[Expansion]
The rule \( \Extend \) tries to extend some arms of the colony
left and right, to use in communicating with a possible non-adjacent
neighbor colony.
In direction \( j \), if it is not adjacent to another colony it will
extend an arm of cells of kind \( \Channel_{j} \).
In channel cells in the positive direction, the \( \Addr \) field
continues its values through
 \[
  \Q ,\Q +1,\dots,2\Q -1.
 \]
Similarly in channel cells in the negative direction.
Extension cells are weaker than member cells, so channel or growth cells
do not normally damage another colony.
The channels will be killed at the end of the computation.

\item[Communication]
The rule
\begin{align*}
 \Send  
\end{align*}
defined in Algorithm~\ref{alg:Send}, will be running
most of the time, sending all needed information to the neighbors.
The rule \( \Retrieve \) records the information from the neighbors.
Atomicity (Condition~\ref{cond:comput}\ref{i:comput.trans})
for the simulated medium will be guaranteed by waiting for a time when 
neither of the neighbor colonies is near the end of their work
period---when this information might change.

 \item[Computation ]
The subrule \( \Compute \) computes the output of the simulated transition
function and stores it on the track \( \Hold \).
It will doom each cell of the colony if the represented
cell is to be erased at the end of the work period.

\item[Payload processing ]
  This part carries out the part of the computation not needed for simulation but for the
  information-processing task of the cellular automaton.

 \item[Growth]
If \( \Growing_{j}=1 \) then, between values \( \growStart \) and
\( \growEnd \), the colony tries to extend an arm of length at most \( \Q  \)
in direction \( j \), making it possible to create a new colony if the encountered
area is empty.

\item[Birth]
A latent cell \( x \) turns immediately into a germ cell with address \( 0  \)
and age 0.
The \df{germ} then begins to grow, trying to fill 3 colonies until age \( \growEnd \).
A germ that fails to grow with a certain speed will be killed off.
If it succeeds then at \( \Age=\U-1 \), the germ cells turn into member cells,
thus implementing a lifting of the organization level.

 \item[Shrinking]
When they reach the end of their growth period\index{period!growth},
growth and germ cells stop producing offshoot.
In what follows, all edges whose existence is not justified by these
processes (called ``exposed'')\index{cell!exposed} will
be subject to the process \( \Decay \)\glo{decay@\( \Decay \)}.
Therefore normally, a growth either disappears before the end of the work
period or it covers a whole new colony by that time.
The rule \( \Send \) will not run during this transition time.

 \item[Finish]
Rule \( \Finish \), will called at \( \Age=\U -1 \),
reduces the addresses of growth cells (if any remain) modulo \( \Q  \).
Outer non-germ cells and inner germ cells turn into members.
Doomed cells will be killed.
Otherwise, the information from the \( \Hold \) track will be copied into the
corresponding locations on the \( \Info \) track.
\end{description}

\begin{definition}\label{def:cut}
  The application of the rule \( \Finish \) will also be called a \df{cut}.
  \end{definition}

Doomedness can be changed only in specific ways:

 \begin{condition}[Dooming]\label{lem:dooming}\index{condition@Condition!Dooming}
Only the following rules can change the field \( \Doomed \) without killing the cell:
 \begin{djenum}

  \item The sub-rules \( \Adapt \) and \( \Heal \) (defined later) can
create doomed or non-doomed cells.

  \item At \( \Age=\payloadStart \), we possibly doom each cell.

 \end{djenum}
\end{condition}

 \subsection{Timing}\label{sec:timing}

Different actions of the program will be carried out with different speeds;
this will be achieved as follows.
Some rules \( r \) have a \df{delay parameter} \( p \)
and a timing variable \( \Wait_{r}\ge 0 \), used as follows:

\begin{algo}[H]
  \lIf{\( \Wait_{r}=0 \)}{\( \Wait_{r}\gets p \)}
  \lElseIf{\( \Wait_{r}\ge p \)}{\( r \)}
  \lElse{\( \Wait_{r}\gets\Wait_{r}-1 \).
}
\end{algo}

A single assignment 
\begin{align*}
 \F\gets_{p}c  
\end{align*}
can also be seen as a rule, with its own delay parameter \( p \) and wait variable.

\begin{definition}\label{def:timing}
  Recall the definition of \( \lambda \) in~\eqref{eq:lg-less-2}.
We introduce some delay constants \( \p_{0} < \p_{1} < \p_{2} < \p_{3} \)
 \glo{p@\( \p_{i} \)}
  whose value will be defined below, further \glo{thm:greek@\( \tau_{i} \)}
 \begin{equation}\label{eq:tau_{i}}
  \tau_{i}  = (\p_{i}+1)\Tu.
 \end{equation}
 \end{definition}
Here is a summary of the roles of the delays:
 \begin{description}
 \item[\( \p_{0} \)] Default and healing;
 \item[\( \p_{1} \)] Computation;
 \item[\( \p_{2} \)] Decay;
 \item[\( \p_{3} \)] Growth.
 \end{description}
 For defining the delay constants let
\begin{align}
\label{eq:heal-span}
            \healSpan       &=2\numDirAff+1, 
\end{align}
denote the reach of the healing operation.
Let
\begin{align}
  \label{eq:p_{0}}
               \p_{0}                  &= 5\lambda,
\\ \label{eq:p_{1}}
            \p_{1}                        &=  4 \lambda \p_{0},
\\ \label{eq:p_{2}}
           \p_{2}                        &=  2 \healSpan\lambda \p_{0}, 
\\ \label{eq:p_{3}}
        \p_{3}                        &= (3\lambda+1) \p_{2}.
\end{align}

\subsection{Plan of the rest of the proof}\label{sec:plan-proof}

  In order to preserve intelligibility and modularity, rules and conditions
belonging to the program will only be introduced as they are needed to
prove some property.

Let us outline how local repairs will be made to deal with obstacles of
functioning; details follow in later sections.
The rule \( \Purge \)\glo{purge@\( \Purge \)} eliminates isolated cells, or ``exposed''
cells that are not important to heal.
The rule \( \Heal \)\glo{heal@\( \Heal \)} repairs a small hole of colony cells (typically after
\( \Purge \) eliminated ``garbage'' from it.).
Germs, and arms of communication or growth will not be healed.
An unrepaired hole will be slowly enlarged by the rule
\( \Decay \), to eventually eliminate partial colonies.

An island can destroy or alter information represented on its space
projection, therefore the information represented on the \( \Info \) track
will be a redundant, error-correcting code.
It will be decoded before computation and encoded after it.
The damage can also disrupt the computation itself, therefore the
decoding-computation-encoding sequence will be repeated several times.
The result will be temporarily stored on the track \( \Hold \) before the
final step of the work period commits to it.

The crucial Lemma~\ref{lem:attrib} (Attribution) says that soon after the
disappearence of the big damage \( \Damage^{*} \), all live non-germ cells not
immediately arising from \( \Damage \) can be attributed (via a path of
ancestors) to some nearby colonies (all disjoint from each other).
This lemma enables us to reason about the process over this area in
terms of big cells.
It helps us for example to see that a big cell can grow a
neighbor if no other cell is nearby.
In terms of colonies, this means that if no other colony is nearby, then
a colony can grow and create a neighbor colony.
This is not obvious since there could be ``debris'': earlier damage could
have left bits and pieces of larger colonies which are hard to override
locally.
The Attribution Lemma will guarantee that those bits and pieces are not
there anymore at the time when they could be an obstacle: whatever is
there is attributable to a big cell.

The Attribution Lemma also plays a role in healing.
Most local healing is performed by the rule \( \Heal \).
However, if the damage occurs at the end of some colony \( \cC \) then it
is possible in principle that foreign material introduced by damage is
connected to something large outside.
The Attribution Lemma will imply that the foreign matter is
weaker (extension or germ cells), and
can therefore be swept away by the regrowth of the member cells of \( \cC \).

The idea of the proof of the Attribution Lemma is the following.
Suppose that \( (x_{0},t_{0}) \) is a cell whose origin we want to trace.
We will be able to follow a ``steep'' path \( (x_{i},t_{i}) \) of ``ancestors''
backwards in time until time \( t_{n}=t_{0}-m\Q  \) with some large coefficient \( m \).
Lemma~\ref{lem:skirting} (Skirting)\index{lemma@Lemma!Skirting} shows that
it is possible to lead a path around an island of damage.
The attribution consists of showing that \( (x_{n},t_{n}) \) belongs to
a domain covering a whole colony.
To prove this, we will show that the decay rule, which eliminates
partial colonies, would eventually cut through the steep path unless the
latter ends in such a domain.
In actual order, the proof proceeds as follows:
 \begin{enumerate}[label = {\textendash}]

  \item
Some of the simpler killing and creating rules and conditions will
be introduced, and some lemmas will be proved that support the
reasoning about paths and domains.

  \item
We prove the Skirting Lemma.
Lemma~\ref{lem:rn-gap} (Running Gap)\index{lemma@Lemma!Running Gap} says
that if a gap is large enough then the decay process propagates it fast,
even in the presence of some damage.

  \item
Lemma~\ref{lem:bad-gap-infr} (Bad Gap Inference)
\index{lemma@Lemma!Bad Gap Inference} shows that (under certain conditions
and in the absence of damage), if there is a gap at all then it is large
enough in the above sense.
 
  \item
The above lemmas are used to prove the Attribution Lemma.

 \end{enumerate}

Here is a summary of the rest of the proof.
 \begin{itemize}

  \item
We define those computation rules not dependent on communication
with neighbor colonies.

  \item
Lemma \ref{lem:legality} (Legality)\index{lemma@Lemma!Legality} shows that
the computation terminates gracefully independently of the success of
communication.

  \item
The development of colony \( \cC \) will be followed forwards to the present
in Lemma \ref{lem:present-attrib} (Present Attribution)
 \index{lemma@Lemma!Present Attribution}.

  \item
Finally, the retrieval rules will be defined and the remaining part of
the Computation Property will be proved.
 \end{itemize}

 \section{Local consistency}\label{sec:consistency}

 Functional obstacles created by an island need to be recognized; 
 we develop the tools for this in the present section.
 
 \subsection{Local maintenance}\label{sec:loc-maint}

Information in a colony about the represented big cell will be corrected using decodings.
But on certain tracks, inconsistency be corrected almost instantaneously.

 \begin{definition}[Locally maintained fields]\label{def:loc-maint}
Certain fields, called
\df{locally maintained}\index{field!locally maintained}, will be
kept constant over the colony, for most of the work period.
Each such field \( \F \) has a \df{default value}.
The updating will of such a field always happen at a specific age
called the \df{update age}\index{update!age}.
Let
\begin{align}\label{eq:n-loc-maint}
   \nLocMaint
\end{align}
be the (constant) number of update ages for locally maintained fields.
A cell's locally maintained field \( \F \) is said to be 
\df{stable}\index{stable} if its age is at distance \( \ge \numDirAff \) from the
update age, where this constant is defined below in~\eqref{eq:split-n}.
An age \( n \) is called \df{stable} for \( \F \) if \( \F \) is stable at age \( n \).
It is called simply \df{stable} if it is stable for all locally maintained fields,
and is also at distance \( \ge \numDirAff \) from ages 0,
\( \computeStart \) and \( \growEnd \).
\end{definition}

From the definitions above it follows that within a work period, all unstable
ages are covered by \( \nLocMaint \) intervals of size \( 2\numDirAff \) each.

 \begin{example}\label{xmp:loc-maint}
Here are some examples of locally maintained fields.
 \begin{itemize}

   \item The Boolean track \( \Doomed \)\glo{doomed@\( \Doomed \)} will be set 
to 1 (true) if the represented big
cell must be removed (the site to become vacant), and 0 otherwise.
Its default value is 0 (false).
Its update age is \( \payloadStart \) (defined in~\eqref{eq:compute-start}).
Cells with \( \Doomed=1 \) are called \df{doomed}\index{cell!doomed}.

  \item The track%
 \glo{growing@\( \Growing_{j} \)}%
 \[
  \Growing_{j}\in\{0,1\},\; j\in\{-1,1\}
\]
signifies the collective decision of the colony to
grow a new neighbor colony in direction \( j \).
Its update age is \( \payloadStart \).
The field \glo{creating@\( \Creating_{j} \)}%
 \[
  \Creating_{j}\in\{0,1\},\; j\in\{-1,1\}
\]
will be used to control the creation of a new neighbor.
Its default value, in a vacant or new latent cell, is 0.
\( \Creating_{j}=1 \) makes a cell a potential creator in direction \( j \)
in the sense of Definition~\ref{def:creator-emerging}.
The value of \( \Creating_{j} \) of a big cell will be broadcast into the
locally maintained track \( \Growing_{j} \) of the cells of its representing colony.
\item The track \( \PlCommands \) will contain commands for the processing
  of payload: see Section~\ref{sec:comp}.
 \end{itemize}
 \end{example}

 \begin{condition}[Locally maintained fields]\label{cond:loc-maint}
Age 0 is an update age for all locally maintained fields other than \( \Doomed \).
If \( \Doomed=1 \) then \( \Growing_{j}=0 \) for \( j\in\{-1,1\} \).
 \end{condition}
  
 \subsection {Fitting neighbors}

 \paragraph{Color}
 We introduced the \( \Color \) field in Definition~\ref{def:std-comp}.
 It distinguishes between germ cells, but 
 can also serve as the only information conserved by the non-ergodic medium.
 Even latent cells have color; when a cell becomes latent its color does not change.

 \begin{definition}\label{def:fit}
   Two neighbor cells \( x<y \) with colors \( c,d \) will be expected to have \df{fitting} colors.
   This is expressed by the relation \( \Fit_{1}(c,d) \) which is 1 if the colors fit and 0 otherwise.
   Let
   \begin{align*}
   \Fit_{-1}(c,d)=\Fit_{1}(d,c).
\end{align*}
We will use two examples of the fitting relation:
\begin{djenum}
  \item\label{i:fit.sorg}
 A color can be an arbitrary symbol from a finite alphabet, and \( \Fit_{j}(c,d)=1 \)
 if and only if \( c=d \).
\item\label{i:fit.comp}
  A color is an element of \( \{-1,0,1\} \), and \( \Fit_{1}(c,d)=1 \) if \( c\le d\le c+1 \).
In what follows we concentrate on the second case, as the other case is treated analogously, only
simpler.
The function of color is in this case is the following: the area of space where payload computation is happening is marked by  color 0, while the space to its left and right by
colors \( -1 \) and \( 1 \) respectively.
\end{djenum}
 \end{definition}

Variant~\ref{i:fit.sorg}
will be used in the proof of Theorems~\ref{thm:1dim.nonerg}, \ref{thm:1dim.nonerg.cont}
and~\ref{thm:1dim.nonerg.var}.
Variant~\ref{i:fit.comp}
will be used in the proof of Theorems~\ref{thm:1dim.comp}, \ref{thm:1dim.comp.var}.

\paragraph{Siblings}
The basic structural soundness of a colony is expressed by some
local consistency conditions.
Space-consistency was introduced in Definition~\ref{def:colonies}.
Time consistency will also be required, but in a
continuous-time cellular automaton we cannot require all cells
to have the same age.
Requirement~\eqref{eq:same-work-period} says that the age within  an  extended colony
is highest at the position of address \( \flo{\Q /2} \) and is non-increasing (with age
difference \( \le 1 \) between neighbors) as we move away from it.

 \begin{definition}[Time consistency]\label{def:time-consistency}
Two cells \( x \) and \( y=x\pm \B  \) with \( |\Addr(x)-\flo{\Q /2}|<|\Addr(y)-\flo{\Q /2}| \)
belong to the \df{same work period}\index{period!work} if
\begin{equation}\label{eq:same-work-period}
  0 \le \Age(x)-\Age(y) \le 1.  
\end{equation}
They \df{straddle a work period boundary}
\index{period!work!boundary} if \( \Age(x)=0 \), \( \Age(y)=\U -1 \).
If one of these cases holds then we will say they are \df{time-consistent}.
\end{definition}

The most important consistency requirement, in the concept of
siblings, is a combination three kinds of consistency: of space- and time-consistency,
and agreement in stable locally maintained fields.

 \begin{definition}[Siblings]\label{def:siblings}
Two cells \( x,x+j\B  \) for \( j\in\{-1,1\} \)
are \df{siblings}\index{sibling} if
one of the following properties holds. 
 \begin{djenum}

   \item  They belong to the same work period, originate at the same colony
(see the definition of originating colony at the beginning of
Section~\ref{sec:cell-kinds}), either both are germ cells or neither of them is.
If \( x \) is a colony cell and \( x+j\B  \) is a growth cell then \( \Growing_{j}(x)=1 \) is
also required.

  \item They belong to the same colony (see the definition of the colony
of a cell at the beginning of Section~\ref{sec:cell-kinds}), 
and straddle a work period boundary.

\item They have fitting colors.

\item If they are germ cells with \( \GermGrowing(\Age) \) then
  they agree in their \( \GermSize \) and \( \Dominant \) fields
  (see Section~\ref{sec:germ-growth}).
\end{djenum}
\end{definition}

The sibling relation, unlike in biology, is not transitive.
It would be more appropriate, but more awkward,
to use the term ``half-sibling''.

\begin{definition}[Domains]\label{def:domain}
  \begin{equation}
    \label{eq:maxDepth}
 \maxDepth = \healSpan    
  \end{equation}
with value defined in~\eqref{eq:heal-span}.
An interval of cells in which the adjacent cells are siblings will be
called a \df{domain}\index{domain}.
A domain of size \( n \) will also be called a
\df{\( n \)-support}\index{support@support!\( n \)-} of its members.
For integer \( k>0 \), let us call two cells \( x,x+k\B  \)
\df{relatives}\index{relatives} if we can change the states of cells
\( x+j\B  \) for \( j=1,\dots,k-1 \) in such a way that the cells
\( x,x+\B ,\dots,x+k\B  \) become a domain.
A colony with starting cell \( x \) will be called
\df{full}\index{colony!full} if it is covered by a domain in such a way
that \( \Addr(x+i\B )\equiv i\pmod \Q  \).
\end{definition}

So if a color is an element of \( \{-1,0,1\} \), and \( \Fit_{1}(c,d)=1 \) if \( c\le d \) then
a domain is either has just one color, or starts with a \( -1 \)'s, continues with \( 0 \)'s and ends
with \( 1's \).
All of these can be empty, except that if both \( -1 \) and 1 is present
then some 0's must also be there.

 For a cell \( x \),  the distance of the boundaries of a domain containing
 \( x \) is estimated by the concept of depth.
For the purpose of depth, we only consider domains within a single colony.

 \begin{definition}[Depth]
Let us call the \df{left depth} of a cell \( x \) the size of the
largest domain of the form \( \{x-i \B , x-(i-1)\B ,\dots,x\} \)
within a single colony containing \( x \).
The right depth is defined similarly.
For example, \( \depth_{-1}(x)>1 \) \glo{depth@\( \depth \)}
just means that \( x \) has a left sibling.

For a locally maintained field \( \F \) we define the \df{left depth}
\( \depth_{j,\F}(x) \) \df{for \( \F \)} of a cell \( x \) the size of the
largest interval of the form \( \{x-i \B , x-(i-1)\B ,\dots,x\} \)
within a single colony containing \( x \) and having constant value of \( \F \).

A cell will always try to keep track of its depths, up to the maximum
possible value.
It will use the fields \glo{Depth@\( \Depth \)}
 \begin{align*}
   \Depth_{j}, \Depth_{j,\F}
   \in\{1,\dots, \maxDepth\},\  j\in\{-1,1\}.
 \end{align*}
 So if for example left depth is larger than \( \maxDepth \) then the cell may
 only see \( \Depth_{-1}=\maxDepth \).
\end{definition}

The condition \( \depth_{j}=1 \) means that in direction \( j \) either there is a
colony boundary or a domain boundary.
  On the other hand the condition \( \Depth_{j}=1 \) is a statement about the field \( \Depth_{j} \),
  maintained by observing \( \depth_{j} \) using the rule \( \rul{Watch-depth} \) 
  in Algorithm~\ref{alg:Watch-depth}.
  Let
  \begin{align}\label{eq:numDirAff}
   \numDirAff       &= 5
  \end{align}
  denote the number of cells directly affected by damage.

 \begin{algo}
 \caption{rule \( \protect\rul{Watch-depth} \)}\label{alg:Watch-depth}
 \pFor{\( j\in\{-1,1\} \)}{
   \lIf{\( \depth_{j}=1 \) \algOr \( \Addr=e_{j} \)}{\( \Depth_{j} \gets_{1}  1 \)}
   \lElse{\( \Depth_{j} \gets_{1}  \maxDepth\land(\Depth_{j}^{j}+1) \)}
   \ForAll{locally maintained fields \( \F \)}{
     \lIf{\( \depth_{j,\F}=1 \) \algOr \( \Addr=e_{j} \)}{\( \Depth_{j,\F} \gets_{1}  1 \)}
     \lElse{\( \Depth_{j.\F} \gets_{1}  \maxDepth\land(\Depth_{j,\F}^{j}+1) \)}
   }
 }
\end{algo}

Rule \( \LocMaintain \) locally maintains all locally maintained variables.

 \begin{algo}
 \caption{rule \( \protect\LocMaintain \)}\label{alg:Loc-maintain}
 \pFor{\( j\in\{-1,1\} \)}{
   \ForAll{locally maintained fields \( \F \)}{
     \If{\( \Age \) is stable for \( \F \) \algAnd \( \Depth_{j,\F}=\maxDepth \)
       \algAnd \( \Depth_{-j,\F}<\maxDepth \)}{
       \( \F\gets\F^{j} \)}     
     }
 }
\end{algo}

\paragraph{Age updating}
The rule for updating age is similar to the ``marching
soldiers''\index{marching soldiers} rule for updating \( \Age \) in 
Section~\ref{sec:sync}.
In Definition~\ref{def:time-consistency}, we imposed some extra
order on the age of all cells in the same extended colony: it must be
non-increasing as they become more distant from
the center cell of the originating colony.
While healing some inconsistency, age updating is paused with
the help of a field called \( \Frozen \).
This variable will also point to the place of  nearby inconsistency that caused the freezing,
In pointing to the position of inconsistency, playing a
role somewhat similar to the fields \( \Depth_{j} \).

 \begin{definition}\label{def:frozen}
There is 
a field %
 \glo{frozen@\( \Frozen \)}\index{cell!frozen}%
 \begin{align*}
  \Frozen\in\clint{-\numDirAff}{\numDirAff} \cap\bbZ
 \end{align*}
with the property that when \( \Frozen\ne 0 \) then \( \Age \) will not be
changed.
A cell with \( \Frozen\ne 0 \) is called \df{frozen}.
 \end{definition}

Here is the basic updating rule for age:%
 \glo{march@\( \March \)}%
 \begin{algo}\caption{sub-rule \( \protect\March \), delay \( \p_{1} \)}\label{alg:March}
  \If{\( \Frozen=0 \) \algAnd \( \Age<\U -1 \) \algAnd \;
    \nl the increase of \( \Age \) does not break the sibling relation 
    with any neighbor\label{alg:March.sibling}
  } 
  {
    {\( \Age \gets \Age + 1 \)}
   }
 \end{algo}



The main colony-organizing variables, address and age, can be changed
only by very few rules.
 
  \begin{condition}[Address and Age]\label{cond:addr-age}
 \index{condition@Condition!Address and Age}
 \begin{cjenum}

  \item
  Only \( \Finish \) can change \( \Addr \) of a live cell.

  \item
  Only \( \March \), \( \HealSync \)
and \( \Finish \) can change \( \Age \) of a live cell.
Of these, \( \HealSync \) can do it when the cell is frozen and
the others when it is not.
 \end{cjenum}
  \end{condition}

 \subsection{Edges}\label{sec:edges}

 \begin{definition}[Edges]\label{def:edges}
Suppose that cell \( x \) has no siblings in direction \( j \), that is it has
\( \depth_{j}(x)=1 \).
It will be called a \df{protected edge}\index{edge!protected} in that
direction if it is some legitimate boundary in that direction, in the sense
described below; otherwise, it will be
called an \df{exposed edge}\index{edge!exposed}.
Recall expansion periods and cells in Definition~\ref{def:expansion-period}.
Essentially, growth edges are protected during the expansion period or
if the growth reached a colony end.
Here is a list of the types of protected edges, for non-germ cells.
For germ cells, the definition will be given in Section~\ref{sec:germ-growth}.
  \begin{description}

   \item[Member] Colony end-cell towards \( j \).

   \item[Expansion cell]
     During the expansion period, while \( j \) is the direction of expansion.
     After the expansion period,
     if the cell is a growth cell and a colony end-cell towards \( j \).
     

    \end{description}
In a rule, the condition  \glo{xposed@\( \Xposed_{j} \)}%
    \[
  \Xposed_{j}.
 \]
 means that the cell \( \x \) applying the rule is an exposed edge in direction \( j \). 
\end{definition}

An exposed edge is the sign of defect, or a call to eliminate an
extension of a colony or a colony; it may be killed---fast by the purge rule,
or slow by the decay rule.

 \begin{lemma}\label{lem:expose}
If a left exposed edge dies and its right sibling was not a
colony endcell, then this neighbor becomes a left exposed edge.
The same holds if we replace left with right.
 \end{lemma}
 \begin{proof}
The one case when this is not obvious is when the exposed edge is a left outer cell
past its expansion period but its neighbor may not be.
But the address-dependent nature of expansion periods as introduced
in Definition~\ref{def:expansion-period} forces then
the right sibling also to be past its expansion period.
\end{proof}

The following notion will be needed much later:

\begin{definition}[Multi-domain]\label{def:multi-domain}
A \df{multi-domain}\index{domain!multi-} is one of the following kinds of set:
\begin{djenum}
  \item A domain.
  \item The union of some adjacent domains meeting in protected colony end-cells.
\end{djenum}
\end{definition}

From the above definitions it is clear that only a cut (as in Definition~\ref{def:cut})
can turn a domain into a multi-domain.
The following lemma is an immediate consequence of the definition of protected edges.
Its technical exceptions relate to the part of the program, detailed later, governing the growth
of germs.
A germ will start from a ``leading'' germ cell with address \( \Q/2 \),
and needs to grow to at least
a size \( 3 \) without seeing live neighbors in order to stay alive.

\begin{lemma}\label{lem:short-exposed}
  If some maximal multi-domain has size \( <\Q \), then the following holds:
  \begin{cjenum}
    \item One of its edges is exposed, with the following 
      exception: it consists of germ cells in their expansion period,
      contains a cell with address \( \Q/2 \), and the
size is either \( >3 \), or 3; in the latter case
the adjacent neighbor cells of the domain are not live.
  \item
If the size is \( \le \Q-d \) with \( d>0 \), then at least one the exposed edges
cannot become protected by the death of an outside neighbor, or 
by shrinking, or by growing by fewer than \( d \) steps.
The only exception is a germ of size 3 of the above kind: it
may become protected if a neighbor dies.
  \end{cjenum}
\end{lemma}

 \section{Killing and creation}\label{sec:kill}

 In order to eliminate local functional obstacles, information in some cells will need to be
 erased and then replaced with one in harmony with one of its neighbors.

\subsection{Killing}

As said above, a cell will be ``killed'' by making it latent.
It will advertise its death in advance:

 \begin{definition}
We will have some one-bit fields
 \glo{dying@\( \Dying \)}%
 \[
  \Dying_{j}, \ j\in\{-1,0,1\},
 \]
with default value 0.
We will say that a cell \( x \) is \df{dying} if \( \Dying_{0}(x)=1 \).
 \end{definition}

In fact, \( \Dying_{0}=1 \) announces that the cell will die soon.
The fields \( \Dying_{j} \) for \( j\in\{-1,1\} \) try to keep track of
whether a neighbor is dying, which will allow a cell to see even whether
a second neighbor is dying.
They are maintained by rule \( \rul{Watch-dying} \) in Algorithm~\ref{alg:Watch-dying} 
which acts with the minimal delay.

\begin{algo}[H]
 \caption{rule \( \protect\rul{Watch-dying} \)}\label{alg:Watch-dying}
  \pFor{\( j\in\{-1,1\} \)}{
     \lIf{\( \Kind^{j}\ne\Vac \)}{\( \Dying_{j} \gets_{1} \Dying^{j}_{0} \)}
  }
 \end{algo}

A program will kill a cell via the
rule \( \Die(\p) \) of Algorithm~\ref{alg:Die}: the argument determines the delay.
\glo{die@\( \Die \)}%
It actually takes \( 2\p \) consecutive applications of the rule \( \Die(\p) \)
to kill the cell: \( \p \)
applications to set \( \Dying_{0}\gets 1 \) and then \( \p \) more to kill it.
The rule \( \Create \) of Algorithm~\ref{alg:Create}
will disqualify a dying cell from creating a live neighbor.

 \begin{algo}
  \caption{sub-rule \( \protect\Die(\p) \)}\label{alg:Die}
  \( \Dying_{0}\gets_{\p} 1 \); \;
  \lIf{\( \Dying_{0}=1 \)}{\( \Kind\gets_{\p}\Latent \)}
\end{algo}

Algorithm~\ref{alg:Purge} (\( \Purge \)) kills exposed cells that
are not to be healed: either because they belong to small maximal 
domains thus are isolated), or expansion arms, or to a doomed area near the end of the
work period.

  \begin{algo}
  \caption{rule \( \protect\Purge \)}\label{alg:Purge}
   \cFor{\( j\in\{-1,1\} \)}{
      \If{\( \Xposed_{j} \) \algAnd 
        (\( \Depth_{-j}\le \numDirAff \) \algOr \( \Kind\ne\Member \) 
        \; \algOr (\( \Doomed \) \algAnd \( \Age>\U -2\Q \)))
      }{\( \Die(\p_{1}) \)}
   }
 \end{algo}

  \subsection{Birth, creation, adaptation}\label{sec:creation}

Rule \( \Birth \) in Algorithm~\ref{alg:Birth}
is just enforcing the condition that newborn cells are latent.
In simulation, birth will be implemented for big cells when germ cells succeed in
creating a new colony.

Rule \( \Create \) in Algorithm~\ref{alg:Create}
controls the values of the field \( \Creating_{j} \) in order
to avoid creating overlapping cells.
A new latent cell is born with \( \Creating_{j}=0 \), turning it on
with a delay \( \p_{0} \). 
Part~\eqref{alg:create.vacant} turns off \( \Creating_{j}(x) \) if 
the neighbor \( x+j\B  \) is already there.
The last part of the rule is just a constraint on the kind of cell that can be created,
namely only a latent one, thereby enforcing Condition~\ref{cond:latent}.
The cell \( \x \) is not really there to ``apply''
this part: in the simulation it is observed by a creator neighbor \( \thg_{-j}(\x) \).

We do not have to worry about other strange rules for non-cells:

\begin{condition}[Birth]\label{cond:birth}
 \( \Create \) and \( \Birth \) are the only rules applicable to a vacant cell.
 \end{condition}

The birth condition implies that in all cases different from the
ones listed in the rules \( \Create \) or \( \Birth \), the site is required by
the transition function to remain vacant.
Condition~\ref{cond:comput}\ref{i:comput.creation}
allows for the creation to be blocked by a cell whose body
intersects the cell to be created.


 \glo{birth@\( \Birth \)}%
  \begin{algo}\caption{rule \( \protect\Birth \)}\label{alg:Birth}
  \If{\( \Kind=\Vac \) \algAnd the two adjacent neighbors are vacant}{
     \( \Kind\gets \Latent \)
  }
\end{algo}

\glo{create@\( \Create \)}%
\begin{algo}\caption{rule \( \protect\Create \)}\label{alg:Create}
  \cFor{\( j\in\{-1,1\} \)}{
    \nl\lIf{ 
      \( \adjNb(j)=1 \) \label{alg:create.vacant}}{
      \( \Creating_{j} \gets_{1}  0 \)}
    \lElse{\( \Creating_{j} \gets_{\p_{0}}  1 \)}
    \lIf{\( \Kind=\Vac \) \algAnd \( \Creating_{j}^{-j}=1 \)}{
      \( \Kind\gets_{1} \Latent \)
    }
  } 
\end{algo}

Consider a cell \( x \) and its left nonadjacent neighbor
\( y \) that may want to create a cell in \( y+\B  \), overlapping the body of \( x \).
Whether \( x \) will be erased is decided not by whether \( y \) is
stronger than \( x \) but by whether the new cell \( y+\B  \) would be stronger
than \( x \).
This distinction matters when a colony attempts to create an outer cell that
would intrude into another colony.
As the created cell would be weaker
than the member cells of the other colony with whom it is competing
for space, this will not happen.
Let us introduce some auxiliary notation, useful for many rules.

 \begin{notation}
For a relation \( R \), let%
 \glo{"<@\( a\leor{R} b \)}%
 \[
  a\leor{R} b
 \]
  mean ``\( a<b \) or (\( a=b \) and \( R \) holds)''.
\end{notation}

Here is the inequality to be used when deciding whether a non-germ cell \( x \)
can be overwritten by a neighbor.
It says that the neighbor must be non-dying and backed up by a non-dying sibling, 
the kind of \( x \) must be dominated by the kind to be created, and if a neighbor on
the other side also wants to overwrite, the intended kind from that side must also be
weaker.
Ties are decided by preferring a creation towards the left.
The competing growth of germs will be regulated later in Section~\ref{sec:germ-growth}.

\begin{definition}\label{def:Kind_j}
  Let \( j\in\{-1,1\} \).
 The field \( \Kind_{j} \) shows the kind of the cell intended to be created
  in the adjacent neighbor in direction \( j \).
  Let
  \begin{equation}\label{eq:NonGermLess}
  \begin{aligned}
    \NonGermLess_{j}&\eqv Kind_{-j}^{j}\geor{j=1}\Kind
     \land \Depth_{j}^{j}>1  \land  \Dying^{j}_{j}=0
    \\ &\land (\Kind^{-j}=\Vacant\lor \Kind_{-j}^{j}\geor{j=1}\Kind_{j}^{-j}).   
\end{aligned}
  \end{equation}
  \end{definition}

For \( j\in\{-1,1\} \), rule~\ref{alg:Adapt}, \( \Adapt(j) \),
makes cell \( \x \) a sibling of its neighbor in direction \( j \).
Sub-rule~\eqref{alg:adapt.create} is not really a rule that a cell can execute, since
it is referring to the case when it is vacant.
Rather, it shows a condition under which the neighbor cell in direction \( j \) is
a potential creator, as in Definition~\ref{def:creator-emerging}.
Sub-rule~\eqref{alg:adapt.kill} kills a cell if it is in the way of the above creation,  
sub-rule~\eqref{alg:adapt.sacrif} erases a latent cell fast.

\begin{algo}\caption{sub-rule \( \protect\Adapt(j) \) }\label{alg:Adapt}
  \If{\( \Dying_{0}^{j}=0 \) \algAnd \( \Creating^{j}_{-j}=1 \)}{
    \nl \lIf{\( \Kind=\Vacant \) \algAnd \( \adjNb(j) \) \label{alg:adapt.create}}
    {\( \Kind\gets\Latent \)}
\nl    \lElseIf{\( \Kind\ne\Latent \)  \label{alg:adapt.kill} }{\( \Die(\p_{1}) \)}
\nl\lElseIf{\label{alg:adapt.sacrif}\( \neg\adjNb(j) \)
    }{\( \Kind\gets_{1} \Vac \)}
    \lElse{make \( \x \) a sibling of \( \x^{j} \), with the same color 
    }
  }
  \end{algo}

 \begin{definition}[Parent]\label{def:parent}
If a latent cell \( x \) came to life by the rule \( \Adapt(j) \)
then we will say that it has been \df{animated}.  
The neighbor in direction \( j \) will be called its \df{parent}\index{cell!parent}.
(The latent cell could have been created earlier by a neighbor, 
which we could call a ``creating parent''.)
In case the same result could also have arisen using the neighbor on the
other side then the parent cell is defined as
the one closer to the center of its colony (or to the left if the animated cell
is the center).
 \end{definition}

The following lemma is immediate from the definition of the adaptation rule.

\begin{lemma}
Suppose that
 \begin{cjenum}

  \item A cell \( x \) has just been adapted at time \( t \) to a non-germ
neighbor \( y=\thg_{j}(x) \) (this rule being a
possible explanation for its becoming live);

  \item Rectangle 
\( (x+\lint{-2.1\B }{3.1\B })\times(t+ \rint{-3\Tu}{0}) \) 
is damage-free;

  \item There is no colony-boundary between \( y \) and its sibling required
by the rule.

 \end{cjenum}
Then \( x \) and \( y \) stay siblings until after \( t \).
 \end{lemma}

  \begin{proof}
The adaptation \( y \) to be non-dying.
Due to the minimum delay \( \p_{0} \) in dying which they did not even
begin, these cells remain live till after \( t \).
Since there is no colony boundary between them, a cut (as in Definition~\ref{def:cut})
will not break the sibling relation of these cells either.
\end{proof}

The following lemma shows that the rule \( \Create \)
indeed succeeds in creating a new cell.
Here, cell \( x-\B  \) will create a cell at site \( x \).
Creation from the right is analogous.
Recall the notation \( \tau_{i} \) from Definition~\ref{def:timing}.

\begin{lemma}[Creation]\label{lem:creation}\index{lemma@Lemma!Creation}
Assume the following, with \( J=\rint{t_{0}}{t_{0} + \tau_{0} + 11\Tu} \):
 \begin{cjenum}
  \item Rectangle 
\( \lint{x-4.1\B }{x+3.1\B }\times J \) 
is damage-free;
  \item We have \( \eta(x-\B ,t).\Dying_{0}=\eta(x-\B ,t).\Dying_{-1}=0 \), and
 \[
 \eta(x-\B ,t).\Kind_{1} \ge \eta(y,t).\Kind\vee \eta(y,t).\Kind_{-1}
 \]
  for all live cells \( \pair{y}{t} \) in the rectangle \( \lint{x}{x+2\B }\times J \).
 \end{cjenum}
Then \( \eta(x,t) \) is non-vacant for some \( t \) in \( J \).
\end{lemma}
  \begin{proof}
Assume, on the contrary, that \( x \) is vacant during all of \( J \), and
we will arrive at a contradiction.
The conditions imply that noise does not affect cell \( x-\B  \), so it can obey its transition rule.
Hence part~\eqref{alg:create.vacant} of the rule \( \Create \) sets 
\( \eta(x-\B ,t).\Creating_{1}=1 \) at some \( t_{1}\le t_{0}+2\Tu \), and this
stays so while \( x \) is vacant.
Condition~\ref{cond:comput}\ref{i:comput.creation} requires \( x \) to
become a cell by time \( t_{1}+2\Tu \) unless the body of some
cell overlaps the body of \( x \).
Suppose therefore that cell \( y \) overlaps the body of \( x \) before time
\( t_{2} = t_{1}+2\Tu\le t_{0}+4\Tu \).
Then \( y>x \), since in case \( y<x \),  cell \( y \) would overlap cell
\( x-\B  \), which is assumed to be there during all this time interval.
This \( y \) disappears by time \( t_{3}=t_{2}+\tau_{0}+\Tu=t_{1}+\tau_{0}+5\Tu \).
Indeed, part~\eqref{alg:adapt.kill} of \( \Adapt \) makes \( y \) latent within
\( \tau_{0} \), then part~\eqref{alg:adapt.sacrif} erases \( y \) within \( \Tu \).

If no similar obstacle cell \( y' \) arises before time
\( t_{4}=t_{3}+2\Tu = t_{1}+\tau_{0}+7\Tu \)
then again \( x \) would be created, so assume one appears.
This can only happen if \( y'+\B  \) creates it (and hence is latent).
After it does, rule \( \Create \) of Algorithm~\ref{alg:Create}
turns off \( \Creating_{-1} \), turning which on again takes time \( \ge\p_{0}\Tl \).
Now \( \Adapt \) will erase \( y' \) by time
\( t_{5}=t_{4}+2\Tu= t_{0}+\tau_{0}+9\Tu \).
No other
cell with \( \Creating_{-1}=1 \) can arise by time \( t_{6}+2\Tu=t_{0}+\tau_{0}+11\Tu \);
partly for the same reason and partly because if \( y'+B \) is erased and a new latent cell
occupies its place then its \( \Creating_{-1} \) is turned on only by a delay of \( \p_{0}\Tl \),
as said in rule \( \Create \) of Algorithm~\ref{alg:Create}.
As a result no new obstacle \( y'' \) can appear soon after \( y' \) has been erased.

Condition~\ref{cond:comput}\ref{i:comput.legal-birth}, allows only one other way for
a cell to appear, namely ``out of nothing''; but this
requires the cell to have no neighbors with a body within distance \( 2\B \) of its body, and as
\( x-\B  \) would be such a neighbor, this cannot happen.

  \end{proof} 

\subsection{Growth}\label{sec:growth}

Growth cells as well as germ cells are born in the attempt to create a colony that encodes
a latent big cell.
Let us spell this out as a condition:

\begin{condition}\label{cond:growth-latent}
  If a cell has kind \( \Growth \) or \( \Germ \) then its \( \Info \) field's value is
  the symbol at \( \Addr\bmod \Q \) in the code of a latent big cell.
  (If, due to a fault, it does not have that value then this would be immediately corrected in the next step.)
\end{condition}

Due to this condition, we do not have to set the \( \Info \) field explicitly in the \( \Grow.\passive \)
and \( \GermGrow.\passive \) rules below.
Here are the rules of growth for non-germ cells.
For germ cells, see Section~\ref{sec:germ-growth}.
Rules~\ref{alg:Extend} (\( \Extend \)) and~\ref{alg:Grow} (\( \Grow \))
both rely on the sub-rule~\ref{alg:Grow.active},
\( \Grow.\rActive \) to show the intent to adapt a neighbor.
The adaptation will be performed by sub-rule~\ref{alg:Grow.passive}, \( \Grow.\passive \)
calling the rule \( \Adapt \).
These rules use some fields introduced in Definition~\ref{def:Kind_j},
and the following definition of end-cells:

\begin{definition}[Endcells]
Let%
 \glo{e@\( \e_{j} \)}%
 \begin{equation}\label{eq:e_j}
  \e_{-1}=0,\ \e_{1}=\Q -1.
 \end{equation}
A cell \( x \) is a \df{colony endcell in direction \( j \)}\index{colony!endcell} 
if \( \Addr(x)\equiv \e_{j} \pmod \Q  \).  
\end{definition}

Rule \( \Extend \) in Algorithm~\ref{alg:Extend} serves to extend the channel in
direction \( j \).
Rule \( \Grow \) in Algorithm~\ref{alg:Grow} depends on the fields
\( \Growing_{j} \) mentioned in Example~\ref{xmp:loc-maint}.
Rule \( \LocMaintain \) of Algorithm~\ref{alg:Loc-maintain}
keeps the locally maintained 
field \( \Growing_{j} \) constant throughout the extended colony.
We will see later that the computation rule sets \( \Growing_{j}=1 \) in
all cells of the colony iff the field \( \Creating_{j} \) in the big cell
represented by the colony has value 1.
Otherwise, it will be 0 in all cells.

 \begin{algo}
  \caption{sub-rule \( \protect\Grow.\rActive(j) \)}\label{alg:Grow.active}
  \If{\( \depth_{-j}>1 \) \algAnd ((\( \Kind=\Ext_{j} \) \algAnd \( \Addr\bmod \Q\ne e_{j} \)) \algOr \;
    (\( \Addr = \e_{j} \) \algAnd \( \Kind=\Member \)))}
     {\( \Kind_{j}\gets  \Ext_{j} \)}
     \lElseIf{\( \Kind\ne\Germ \)}{\( \Kind_{j}\gets  \Latent \)}
 \end{algo}

\begin{algo}
  \caption{sub-rule \( \protect\Grow.\passive(j) \), delay \( \p_{3} \)}\label{alg:Grow.passive}
  \If{\( \NonGermLess_{-j} \) \algAnd \( \Kind_{j}^{-j}=\Ext_{j} \) }
  {\( \Adapt(-j) \) }
 \end{algo}

 \begin{algo}
  \caption{sub-rule \( \protect\Extend \)}\label{alg:Extend}
   \pFor{\( j\in\{-1,1\} \)}{
      \If{\( \Kind=\Latent \) \algAnd \( \Age^{-j}\in\rint{0}{\computeStart} \)}{
         \( \Grow.\passive(j) \)
      }
      \lElseIf{\( \Age\in\rint{0}{\computeStart} \) \algAnd \( \Kind\ne\Germ \)}{\( \Grow.\rActive(j) \)}
    }
 \end{algo}

 \begin{algo}[H]
  \caption{sub-rule \( \protect\Grow \)}\label{alg:Grow}
   \pFor{\( j\in\{-1,1\} \)}{
     \If{\( \Age\in\rint{\growStart}{\growEnd} \)
         \algAnd \( \Growing_{j}=1 \)}{
         \( \Grow.\rActive(j) \)
      }
      \ElseIf{\( \Kind=\Latent \) \algAnd 
         \( \Age^{-j}\in\rint{\growStart}{\growEnd} \)
         \algAnd \( \Growing_{j}^{-j}=1 \)}{
         \( \Grow.\passive(j) \)
       }
     }
 \end{algo}

 \subsection{Germ growth}\label{sec:germ-growth}

A latent cell turns immediately
into a germ cell with \( \Age=0 \), \( \Addr= \Q/2 \) (which we can assume to be an integer).
A germ is a domain of germ cells; its goal is to grow into 5 colonies, and then turn
the middle one into a big cell.
It tries to grow left and right, to addresses
\begin{align}\label{eq:e'_j}
 e'_{-1}=-2\Q  \text{ and } e'_{1}=3\Q-1
\end{align}
respectively, with a notation analogous to~\eqref{eq:e_j}.
A growing germ cannot
intrude into an extended colony, as the cells of the latter are stronger,
and a growing colony can destroy germ cells in its way.
Also, the germ growth rule arbitrates between germs that are in the way of each others' growth.
See Algorithm~\ref{alg:Germ-grow.passive}.

The germ work period in Algorithm~\ref{alg:germ-work-period}
consists of a part \( \GermGrow \) when it tries to grow,
followed by a part when it will decay if it did not reach its growth goal.
So the edge of a germ cell is \df{protected} if it is pointing away from the cell
with address \( \Q/2 \), and has either age \( <\growEnd \) or reached
address \( e'_{j} \) in one of the directions \( j \).
The work period ends with the procedure \( \Lift \), to be described in Section~\ref{sec:lifting},
that lifts the computing functionality to the new level of simulation by adding a new
level of error correction to the payload.
The total number \( \U \) of age-increasing steps of a germ work period is the same
as that of a colony work period.

For convenience of notation, in what follows
properties and fields related to germs will be marked with ``G-''.

\begin{algo}\caption{A germ work period}\label{alg:germ-work-period}
  \( \GermGrow \);\;
  Shrink;\; 
  \( \Lift \)
\end{algo}

The growth of the germ happens in \df{stages}.
The field \( \StageStart \) will show the start age of the current stage,
lasting until 
 \begin{align}\label{eq:stage-end}
 \StageEnd=\StageStart+\rep\cdot\GermSize,
 \end{align}
where the field \( \GermSize \) shows the size of the germ at that time, and
the constant \( \rep \) was introduced in~\eqref{eq:rep}.
When \( \Age=\StageEnd \) then cells will set \( \StageStart\gets\StageEnd \)
and \( \StageEnd \) according to~\eqref{eq:stage-end}.
Each stage will be divided into two parts: a \df{computation} part and a \df{growth} part.
The growth part starts at age \(    \StageStart + 0.1\rep\cdot\GermSize \), and lasts until
the end of the stage.
Let
\begin{align*}
   \GermGrowing(n)\in\{0,1\}
\end{align*}
be 1 if \( n \) is in the growing part of one of the stages, and 0 otherwise.
During the computation part, the germ's cells will be informed of the addresses of the
of the germ's edges
by the rule \( \PassGermSize \) of Algorithm~\ref{alg:Pass-germ-size}.
At the end of this part, a single step
\begin{align*}
 \GermSize\gets\GermEdge_{1}-\GermEdge_{-1}  
\end{align*}
records the germ size in each cell of the germ.
We don't want a germ to grow much on one side if its growth is blocked
(for whatever reason) on the other side.
To achieve this, at the beginning of the work period, the central cell
computes 
\begin{align*}
   \GrowDir\gets\argmin\{\GermEdge_{-1},\GermEdge_{1}\},
\end{align*}
which, in case \( \GermEdge_{-1}=\GermEdge_{1} \) chooses \( 1 \).
This value will also be propagated by \( \PassGermSize \), 
and growth will only be attempted in direction \( \GrowDir \).
Also, as seen in the rule \( \GermGrow.\rActive \)
of Algorithm~\ref{alg:Germ-grow.active}, the growing side of the germ will be
allowed to grow only to a distance \( \GermSize \) from its center.

\begin{algo}\caption{rule \( \protect\PassGermSize \)}\label{alg:Pass-germ-size}
  \If{\( \Kind=\Germ \) }{
    \If{ \algNot \( \GermGrowing(\Age) \)}{   
      \cFor{\( j\in\{-1,1\} \)}{
        \lIf{ \( \Edge_{j}= 0 \) }{ \( \GermEdge_{j}\gets_{1}\GermEdge_{j}^{j} \) }
        \lElse{ \( \GermEdge_{j}\gets\Addr \) }
        \lIf{\( \sign(\Addr-\Q/2)=j \)}{\( \GrowDir\gets\GrowDir^{-j} \)}
      }
    }
  }
\end{algo}

The growth of germs over other germs has some color restrictions, which we will
express with the help of the relation \( \Fit'_{j}(c,d) \).
In case of the variant~\ref{i:fit.sorg} of Definition~\ref{def:fit},
only germs of the same color are not prohibited from overriding each other, so 
simply \( \Fit'_{j}(c,d)=\Fit_{j}(c,d) \).
In case of the variant~\ref{i:fit.comp},
also germs with color 0 are not prohibited to override fitting germs.
Formally, 
\begin{align*}
  \Fit'_{j}(c,d)=1 \eqv \Fit_{j}(c,d) \land ( c=d \lor d=0).
\end{align*}
A germ can grow (subject to the above constraints) over another, if it is significantly
larger.
By ``significantly larger'', we mean larger by a factor
\begin{align}\label{eq:germ-signif-larger}
  \grFac = 5/4.   
\end{align}
If they are close in size, that is they differ by a factor \( \le\rho \) then the random choice
described below will arbitrate between them.
However, let us see that even the smaller one, if the decision favors it, will grow to become
significantly larger than the larger one.
So suppose that germ \( G_{1} \) has size \( n \), germ \( G_{2} \)
has size \( \rho n \), and germ \( G_{1} \) is allowed to grow.
Its growth can be stopped by another germ \( G'_{1} \) growing in the opposite direction.
However, one of these will cover at least half of \( G_{2} \).
So the germ that grew will have size at least \( n(\rho/2+1)>\rho^{2}n \).

The randomized decision works as follows.
Both germs look at a field
\begin{align*}
   \Dominant\in\{0,1\},
\end{align*}
which will be set randomly by the cell with address \( \Q/2 \) at each stage start,
and then propagated by the rule \( \PassDominant \) of Algorithm~\ref{alg:Pass-dominant}.
Between germs of nearly the same size, the dominant one will prevail.
After the computation part of the stage, the values of \( \GermSize \) and \( \Dominant \)
are supposed to be constant throughout the germ: this is part of the requirement
for two neighbor cells to be siblings.
Formally, the relation \( \Fit'_{j}(c,d) \) and
\begin{align*}
\GermLess_{j}(a,b)   
\end{align*}
will be used by the rule \( \GermGrow.\passive \) of Algorithm~\ref{alg:Germ-grow.passive}.
The meaning of \( \GermLess_{j}(a,b) \) is that the germ cell with state \( b \) is
in direction \( j \) from the germ cell with state \( a \), and will be able to grow
over it or over the space that \( a \) wants to grow over.
Here are the rules to determine \( \GermLess_{j}(a,b) \).
The germ of \( b \) can definitely grow over the one of \( a \) if both
are in their growing ages, and it is it is either significantly larger or has reached the end
of its growth on the other side.
A remaining case is when both germs reached their growth end on the other side:
then the germ of \( b \) can grow over the one of \( a \)
if it is dominant while \( a \) is not.
The other remaining case is when neither germ reached its growth end on the other side,
and their sizes are also comparable: then again
the germ of \( b \)  can grow over the one of \( a \)
if it is dominant while \( a \) is not.
Formally, here is what determines whether \( \GermLess_{j}(a,b) \) holds:
\begin{sloppypar}
\begin{enumerate}
\item We always need \( \Kind(b)=\Germ \), \( \GermGrowing(b.\Age)=1 \).
\item \( a.\Kind < \Germ \) is enough,
  otherwise \( a.\Kind=\Germ \) and \( \GermGrowing(a.\Age)=1 \) is necessary.
\item If \( a.\GermEdge_{-j}\ne e'_{-j} \) then
  \( b.\GermEdge_{j}= e'_{j} \) or \( \grFac  a.\GermSize < b.\GermSize \) is enough.
\item \( a.\Dominant=0 \) and \( b.\Dominant=1 \) is enough if either
\( a.\GermEdge_{-j}= e'_{-j} \) and \( b.\GermEdge_{j}= e'_{j} \)
or \( a.\GermEdge_{-j}\ne e'_{-j} \) and \( a.\GermSize\le\grFac b.\GermSize \).
\end{enumerate}
  \end{sloppypar}
  \medskip
  
  Germ growth is governed by rule \( \GermGrow \) of Algorithm~\ref{alg:Germ-grow}.
  The ages in which the edge of a germ can expand will be restricted only to part of the
  growth period, to make sure that it can finish the growth it started:
\begin{align*}
  \GermGrowing&.\vActive(\Age,\Addr) \eqv \GermGrowing(\Age)
\\ &\land\Age<\StageEnd-0.1\rep\GermSize + 5\lambda\tau_{2}|\Addr-\Q/2|.
\end{align*}
Similarly, the ages in which the edge of a germ can be killed by the expansion of another
germ will be restricted only to part of the
growth period, to make sure that the whole germ can be killed before the end:
\begin{align*}
  \GermGrowing&.\vPassive(\Age,\Addr) \eqv \GermGrowing(\Age)
\\ &\land \Age<\StageEnd-5\lambda\tau_{2}|\Addr-\Q/2|.
\end{align*}

\begin{algo}[H]\caption{rule \( \PassDominant \)}\label{alg:Pass-dominant}
  \If{\( \Kind=\Germ \) }{
    \lIf{ \( \Addr=\Q/2 \) \algAnd \(  \Age = \StageStart \)}
    { \( \Dominant=\Rand \) }
    \lElse{ \algLet \( j=\sign(\Addr-\Q/2) \)\;
      \( \Dominant\gets_{1}\Dominant^{-j} \)
    }
  }
\end{algo}

\begin{algo}[H]\caption{rule \( \GermGrow \)}\label{alg:Germ-grow}
  \If{\( \Kind=\Germ \)}{
      \algLet \( j\gets \sign(\Addr-\Q/2 - 0.1) \)\;
      \( \GermGrow.\rActive(j) \)
    }
    \If{\( \Kind\le\Germ \) }{
      \For{\( j\in\{-1,1\} \)}{
        \lIf{ \( \Kind_{-j}^{j}=\Germ \)}
        {\( \GermGrow.\passive(j) \)}
      }
    }
  \end{algo}

  \begin{algo}[H]\caption{rule \( \GermGrow.\rActive(j) \)}\label{alg:Germ-grow.active}
    \lIf{\( \GermGrowing.\vActive(\Age,\Addr) \) \algAnd \( \Age<\growEnd \)
      \algAnd \( e'_{-1}<\Addr<e'_{1} \)
      \algAnd \( j = \GrowDir \) \algAnd \( |\Addr-\Q/2|<\GermSize \) }
    {\( \Kind_{j}\gets \Germ \) }
    \lElse{\( \Kind_{j}\gets\Latent \)}
  \end{algo}

\begin{algo}[H]\caption{rule \( \GermGrow.\passive(j) \)}\label{alg:Germ-grow.passive}
  \If{ \( \Kind=\Latent \) \algOr (\( \Kind=\Germ \)
    \algAnd \( \GermGrowing.\vPassive(\Age,\Addr) \)  \algAnd
     \( \Fit'_{j}(\Color,\Color^{j}) \)
    \algAnd \( \GermLess_{j}(\x,\x^{j}) \) \algAnd \( \GermLess_{j}(\x^{-j},\x^{j}) \) )}
  { \( \Adapt(j) \) }
  \end{algo}

\subsection{Healing rules}\label{sec:heal-rule}

Healing involves several rules.
Recall frozen cells in Definition~\ref{def:frozen}, and
that a frozen cell's age will not be advanced in the regular way
(seen in the rule~\ref{alg:March}, \( \March \)).
The default value of  \( \Frozen \) is \( 0 \); nonzero values of \( \Frozen \)
will point to some ``defects''.
This field is intended to show that in the direction indicated
by its sign, there is a time consistency
at a distance indicated by its absolute value 
\( |\Frozen|\le \numDirAff \) (the latter constant is defined in~\eqref{eq:numDirAff}).
So if \( \Frozen=1 \) then the cell has on the right
a time-inconsistent neighbor that is not a protected edge; if it is 2
then its right neighbor has this property.
The value of \( \Frozen \) will be propagated
by part~\eqref{alg:Freeze.propagate} of \( \Freeze \) of Algorithm~\ref{alg:Freeze}.
Suppose that \( -\numDirAff<\Frozen(x-\B )<0 \) 
and \( x \) is not exposed on the right:
then we will want to propagate the signal to the right
by setting \( \Frozen(x)\gets \Frozen(x-\B )-1 \).x
If at the same time \( x+\B  \) is a sibling with
\( 0<\Frozen(x+\B )<\numDirAff \), 
then \( \Frozen(x)\gets \Frozen(x+\B )+1 \) would also be suggested.
In such cases we will break the tie and propagate to the right.
The freezing process is typically started by part~\eqref{alg:Freeze.start}.
If our cell is exposed to the right then we set  \( \Frozen(x)=1 \).

The last \textbf{else} of the rule \emph{gradually} unfreezes the area after
the inconsistency was eliminated.
As expected, the field \( \Frozen \) is private to the rule \( \Freeze \):

\begin{condition}[Freeze]\label{cond:freeze}\index{condition@Condition!Freeze}
Only the rule \( \Freeze \) can change the field \( \Frozen \) (without killing
the cell).
\end{condition}

The typical result of repeated application of the \( \Freeze \) rule
is that values of the field \( \Frozen \) in consecutive cells
in the place where damage occurred will be like
 \begin{equation}\label{eq:after-freeze}
  0,\dots,0,5,4,3,2,1,*,*,-1,-2,-3,-4,-5,0,\dots,0,
 \end{equation}
(using \( \numDirAff=5 \)).
Here the \( * \)'s show latent cells.
All (weakly) exposed edges are confined between the 1 and \( -1 \) in this
picture.
One of the sides of this picture may be cut short by a protected edge.

 \glo{freeze@\( \Freeze \)}%
  \begin{algo}[H]
  \caption{rule \( \protect\Freeze \)}\label{alg:Freeze}
\Cond{
  \cFor{\( j\in\{-1,1\} \)}{  
   \nl \If{ \( \Xposed_{j} \) \algAnd the neighbor in direction \( j \) is not time-consistent \label{alg:Freeze.start}}{
      \( \Frozen \gets_{1}  j \)}
  }
   \cFor{\( j\in\{-1,1\} \)}{
     \nl\If{\( 1 \le |\Frozen^{j}| <\numDirAff \) \label{alg:Freeze.propagate}\;
       \algAnd  (\( j=-1 \) \algOr \pNot \( 1 \le |\Frozen^{-j}| < \numDirAff \))
     }{\( \Frozen \gets_{1}  \Frozen^{j}+j \)}
     }
     \lElse{\( \Frozen \gets_{p_{1}}  0 \)}
   }
 \end{algo}

Recall the notion of directly affected cells in Definition~\ref{def:affecting}.
The number of consecutive adjacent cells directly affected by the damage
rectangle is at most \( \numDirAff \), defined in~\eqref{eq:numDirAff}.
For the moment, let us call the smallest
interval encompassing these sites the \df{corruption}.
The healing rule will be able to repair a
corrupted interval of size at most \( \numDirAff \) arising in some domain.
A typical case is when the corruption occurs
inside a domain which itself is part of a colony.
In thinking of the healing rule, let us
assume that other rules (namely \( \Freeze \), \( \Purge \) and \( \Create \)) 
have achieved that the corruption is filled with 
cells aligned with its boundaries.
If our model was in discrete-time then the
correct values of \( \Addr \), \( \Age \) and the locally maintained fields could
just be inferred from either edge of the corruption and copied into its
cells.
But since the model is in continous time, the task is more complex.
First, the \( \Freeze \) rule will freeze the \( \Age \) values in
some of the cells surrounding the gap.
(Otherwise, if the \( \Age \) in one of these cells speeds
ahead, it may become impossible to bridge
the gap with continously varying \( \Age \) values.)
Then, it takes some coordination to
recreate a continous sequence of \( \Age \) values.

 \begin{remark}
It would be convenient here if our robust medium had range \( r>4 \), so cells
on both sides of the corruption can see each other.
But introducing ``action at a distance'' brings its own headaches.
 \end{remark}

The healing rule itself will carry out only a single step of the
reconstruction, but its repeated applications will close the gap.
As a further complexity, the corrupted cells can be alive, 
and therefore it might not even be clear which cell to correct.
Suppose, for example, that values of \( \Age \) in adjacent cells are
\( 0, 0, 2, 2 \).
Both \( 0,1,2,2 \) and \( 0,0,1,2 \) are possible corrections:
the healing rule will choose one.

Rule \( \Heal \) in algorithm~\ref{alg:Heal}
consists of several rules applied simultaneously.

 \begin{algo}\caption{rule \( \protect\Heal \)}\label{alg:Heal}
  \( \HealRevive \)\prl
  \( \HealSync \)
 \end{algo}

Rule~\ref{alg:Heal-revive}, \( \HealRevive \), will
revive a latent cell \( x \), making it a weak sibling to one of its neighbors.
Rule \( \HealSync \) in Algorithm~\ref{alg:Heal-sync} 
will try to adjust the \( \Age \) variables.
In rule \( \HealRevive \), assume that \( x \) is latent, and \( x-\B  \) is alive.
Then the rule may try to apply the rule \( \Adapt(-1) \)
to \( x \), adapting it to its left neighbor.
If there is a similar demand on the right-hand side, then
if one of the adaptations creates a stronger cell, that one is chosen,
otherwise the left side is preferred.
The field \( \Healing \) will count the cells revived
during healing (by numbering them), limiting this  way the size of revived
segments to \( \healSpan \).
It actually may need to revive a segment of size \( 2\numDirAff \) with (\( \numDirAff \) defined
in~\eqref{eq:numDirAff}):
if cells with \( \Addr\in\{\numDirAff,\dots, 2\numDirAff-1\} \) are changed by damage
then the purge rule may kill these cells as well as the now isolated
segment with \( \Addr\in\{0,\dots,\numDirAff-1\} \).
 
 \begin{definition}[Healing field]
We have a field called
 \begin{align*}
          \Healing &\in\{0,1,\dots,\healSpan\}.
 \end{align*}
Its default value is 0.
A cell with \( \Healing>0 \) will be called a \df{healing cell}.
The field
 \begin{align*}
   \Decaying\in\{0,1\}
 \end{align*}
will limit the healing process in time.
 \end{definition}

The rule \( \HealRevive \) is the combination of several simpler rules.
The rule \( \rul{Heal-passive} \) adapts a cell to an exposed member neighbor (making it thus 
also a member), unless the neighbor is decaying or is doomed and old.
Its part~\eqref{alg:Heal-passive.center} requires the healing to proceed
always away from a colony center (this restriction is only for proof convenience).
The end-healing will allow to kill intruder cells at the end of a colony in order to heal the
end; these intruder cells would not be necessarily isolated, so might not be eliminated by 
the \( \Purge \) rule.

\begin{condition}\label{cond:healing}
The rule \( \HealRevive \) is the only one changing the value of the field
\( \Healing \).   
\end{condition}

 \begin{algo}[H]
 \caption{sub-rule \( \protect\HealRevive \)}\label{alg:Heal-revive}
 \lIf{\( \Kind=\Latent \)}{\( \Heal.\passive \)}
 \Else{
   \( \rul{End-heal} \)\;
   \( \rul{Control-decay} \)
}
 \end{algo}

\begin{algo} \caption{sub-rule \( \protect\Heal.\passive \), delay \( \p_{1} \)}\label{alg:Heal-passive}
  \cFor{\( j\in\{-1,1\} \)}{
       \If{\( \Healing^{j}<\healSpan \) \algAnd \( \Xposed^{j}_{-j} \)
         \; \algAnd \( \Decaying^{j}=0 \) \algAnd \( \Kind^{j}=\Member \) 
         \; \algAnd\algNot(\( \Doomed^{j} \) \algAnd \( \Age^{j}>\U -2\Q \))
         \;\nl\algAnd direction \( j \) is towards the colony center of \( \x^{j} \)\label{alg:Heal-passive.center}
         \; \algAnd \( \NonGermLess_{j} \)
         }{
          \( \Adapt(j) \)\;
          \( \Healing \gets  \Healing^{j}+1 \)
       }
    }
\end{algo}

\begin{algo} \caption{sub-rule \( \protect\rul{End-heal} \)}\label{alg:End-heal}
  \cFor{\( j\in\{-1,1\} \)}{
      \If{\( \Xposed_{j} \) \algAnd \( 0<|\e_{j}-(\Addr\bmod \U ))| \le \healSpan \)
        \algAnd \( \Kind=\Member \)
         }{
            \( \Kind_{j}\gets  \Member \)
           }
   }
\end{algo}

\begin{sloppypar}
Here is the ``story'' interpreting these rules
(and helping to follow their analysis in later proofs).
Suppose healing proceeds, say, towards the left.
An exposed cell sets \( \Decaying=1 \) \emph{slowly}
in part~\eqref{alg:heal-revive-decay-1} of the rule \( \rul{Control-decay} \),
Algorithm~\ref{alg:Control-decay}, so while \( \Decaying \) is not set
it still allows the rule \( \Heal.\passive \) in Algorithm~\ref{alg:Heal-passive}
to proceed with healing.
If a healing arm created by it reaches maximum length and does not close a gap
then its exposed end eventually sets \( \Decaying=1 \), and
the rule \( \Decay \) in Algorithm~\ref{alg:Decay} will kill it.
In rule \( \rul{Control-decay} \) of Algorithm~\ref{alg:Control-decay},
as soon as a cell gets a neighbor with a larger \( \Healing \) value, it
sets \( \Decaying=1 \) \emph{fast} in part~\eqref{alg:heal-revive-decay-2}.
Thus if the neighbor with higher \( \Healing \) value
gets killed due to unsuccessful healing, it will not be revived
again by rule \( \Heal.\passive \) of Algorithm~\ref{alg:Heal-passive}.
This way eventually the whole unsuccessful healing arm will be killed.
Its first neighbor with \( \Healing=0 \) will also be killed; so
an unsuccessful healing eventually widens the gap it tried to heal.
Rule \( \rul{Control-decay} \) in Algorithm~\ref{alg:Control-decay}
also sets the fields \( \Healing \) 
and \( \Decaying \) to their default value 0 when their role has ended.
  \end{sloppypar}

\begin{algo}[H]\caption{rule \( \protect\rul{Control-decay} \)}\label{alg:Control-decay}
  \nl\lIf{\( \exist{j}\Xposed_{j} \)\label{alg:heal-revive-decay-1}}{
    \( \Decaying\gets_{\p_{2}} 1 \)}
  \nl\lElseIf{\( \exist{j}\Healing^{j}>\Healing \)\label{alg:heal-revive-decay-2}}{
    \( \Decaying\gets_{1} 1 \)}
  \Else{
    \( \Healing \gets_{1} 0 \)\;
    \( \Decaying \gets_{1} 0 \)
 }
\end{algo}

Rule \( \HealSync \) in Algorithm~\ref{alg:Heal-sync} will readjust the age of
some frozen cells, to eliminate a discontinuity.
To do this, it will try to bring \( \Age \) closer to the one
of the neighbor to which the \( \Frozen \) field points.
The conditions \( \Frozen/j>0 \) and \( \Frozen^{-j}\cdot\Frozen>0 \) say that
both \( \Frozen(x) \) and \( \Frozen(x-j\B ) \) point to the same direction \( j \). 
Condition \( \Depth_{-j}>1 \) requires 
the neighbor in direction \( (-j) \) to be a sibling.
The  condition \( \Dg(j)\cdot\Dg(-j)\ge 0 \) does not allow
an age change that would break this sibling relation.
Also, because of \( \Frozen^{-j}\cdot\Frozen>0 \), if changes happen
simultaneously in both \( x \) and \( x-j\B  \) then they happen in the same
direction, still not breaking the sibling relation.

\begin{condition}[Age decrease]\label{cond:age-decrease}
The only rule in which the \( \Age \) field may decrease is \( \HealSync \).
\end{condition}

 \begin{algo} \caption{sub-rule \( \protect\HealSync \)}\label{alg:Heal-sync}
 \lpFor{\( j\in\{-1,1\} \)}{\algLet \( \Dg(j) \gets \Age^{j}-\Age \amod \U  \)}
 \cFor{\( j\in\{-1,1\} \)}{
   \If{\( \Frozen/j>0 \) \algAnd \( \Frozen^{-j}\cdot\Frozen>0 \) \algAnd \( \Depth_{-j}>1 \)\;
     \algAnd \( \Dg(j)\ne 0 \) \algAnd \( \Dg(j)\cdot\Dg(-j)\ge 0 \)}{
       \( \Age \gets  \Age + \sign(\Dg(j)) \bmod \U  \)
   }
 }
 \end{algo}

 \begin{example}\label{xpl:heal-sync}
The typical application of \( \HealSync \) is, say, when in a domain of cells
with age values \( 1,1,1,2,3,4,5,5, \) the damage changes the values of three to 
 \begin{align*}
  1,1,1,0,7,6,5,5.
 \end{align*}
Now the cells with age values 0 and 7 become exposed towards each other, and
the \( \Frozen \) variable develops the values \( 4,3,2,1,-1,-2,-3,-4 \).
After this, \( \HealSync \) will gradually bring the ages closer to each
other.
Here is a possible synchronization history (time going downwards):
 \begin{align*}
 \begin{matrix}
      1&1&1&0&7&6&5&5
\\  1&1&0&0&6&6&5&5
\\  1&1&0&1&6&6&5&5
\\  1&1&0&1&5&6&5&5
\\  1&1&0&1&5&5&5&5
\\  1&1&0&1&4&5&5&5
\\  1&1&1&1&4&4&5&5
\\  1&1&1&2&3&4&5&5
 \end{matrix}.
 \end{align*}
As in the second line here, the synchronization
may not start in the best direction, but since the adaptation is towards
the fault in the middle, eventually the fault will disappear.

It is not hard to see that in general, if the domain to synchronize has width \( \numDirAff \) then
the synchronization will be finished within \( 2\numDirAff \) steps.
 \end{example}

Let us summarize the possible uses of animation and killing.

 \begin{condition}[Animation, killing]\label{cond:anim,kill}
 \index{condition@Condition!Animation}
 \begin{cjenum}
 \item Only the rule \( \Die(\p) \) can kill a cell.
  It is invoked always with \( \p\ge\p_{1} \).  
\item
  \begin{sloppypar}
   Only the following rules can invoke \( \Die \): \( \Adapt \), \( \Purge \), \( \Decay \). 
 \end{sloppypar}
 \item Only the rule \( \Adapt \) can erase a cell (make it vacant).
  \item\label{cond:killing.protected}
A protected edge can be killed only if a neighbor sets \( \Kind_{j}\ne\Latent \)
in its application of a growth rule or the end healing rule.

\item  Only the rule \( \Create \) can change \( \Creating_{j} \).
\item
  \begin{sloppypar}
Only the rule \( \Adapt \) can create a non-latent cell,
invoked only by \( \Heal.\passive \), \( \Grow.\passive \) and \( \GermGrow.\passive \)
(with the corresponding delay).
It is also the only rule by which a potential creator can arise in the sense of
Definition~\ref{def:creator-emerging}.
      \end{sloppypar}
  \end{cjenum}
 \end{condition}

An exposed edge is almost always a sign of past damage, the exception being when
at the end of a stage, a large group of cells needs to be killed off:

\begin{lemma}[Exposing]\label{lem:exposing}\index{lemma@Lemma!Exposing}
  A rule exposes an edge only in the following cases:
 \begin{djenum}

  \item\label{i:lem.exposing.doom}
    Member cell whose neighbor is doomed, and dies at \( \Age=0 \), that is at a cut,
    as in Definition~\ref{def:cut}.

  \item\label{i:lem.exposing.channel}
   Channel cell, \( \Age=\computeStart \).

  \item\label{i:lem.exposing.grow-end}
    Germ cell or growth, non-end-cell. \( \Age=\growEnd \), or when
    a germ cell is killed by \( \Adapt \) invoked in \( \GermGrow.\passive \).

 \end{djenum}
 \end{lemma}
 \begin{proof}
  Direct consequence of the definition of siblings and exposed edges and
Conditions~\ref{cond:addr-age} (Address and Age), \ref{cond:anim,kill}
(Killing).
\end{proof}

  \subsection{Continuity}

  Our terminology turns out to be ``incestuous'': a child cell can only be
created if it also becomes a sibling.
Recall the interval \( I_{1}(x) \) in Definition~\ref{def:rectangles}.
  
 \begin{lemma}[Parent]\label{lem:parent}\index{lemma@Lemma!Parent}
   Suppose that \( I_{1}(x) \) is damage-free during \( \rint{u-\Tu}{u} \)
   (that is \( x \) is not affected via neighbors during this time,
   as in Definition~\ref{def:affecting}),
   and the non-germ cell \( x \) becomes animated at time \( u \).
   Then there is a \( j\in\{-1,1\} \) and \( t'\in\rint{u-3\Tu}{u-\Tl/2} \) such that
the following holds:

  \begin{cjenum}

   \item\label{i:parent.gen}
\( x \) gets animated by
parent \( x+j\B  \), whose state at that time makes it a sibling of \( \pair{x}{u} \).

\item\label{i:parent.sib}
If \( I_{0}(x)\cup I_{0}(x+j\B ) \) remains damage-free during \( \rint{u}{u+\Tu} \) then
\( x \) and \( x+j\B  \) remain siblings during that time.

\end{cjenum}
\end{lemma}

  \begin{proof}
We have an application of \( \Adapt \) of Algorithm~\ref{alg:Adapt}
from some \( y=x+j\B  \).
So at the observation time \( t' \) corresponding to the switch
\( \pair{x}{u} \), we have the situation described by~\ref{i:parent.gen}.
To prove~\ref{i:parent.sib}:
according to Condition~\ref{cond:anim,kill}, the applied rule was
either \( \Heal \) or \( \rul{Grow} \), creating a sibling with \( \Age \) equal
to that of the parent.
Dying can only happen by decay or purge, 
which takes time at least \( \p_{0}\Tl\ge 2\Tu \).
The age of the child is made equal to the age of the mother.
In case of growth, the mother has time for at most change of age (with delay
\( \p_{1} \)), and this will not break the sibling relation within time \( 2\Tu \).
The healing case is trickier since the rule \( \HealSync \)
in Algorithm~\ref{alg:Heal-sync} can also decrease \( \Age \).
But given that \( x \) was just animated then if, say, it is created from the left, then
the \( \Frozen \) variables must point to the right.
Therefore the rule will not cause age change in the left neighbor as the ages are made
equal.
 \end{proof}

Normally, siblings remain siblings:

\begin{lemma}[Glue]\label{lem:glue}\index{lemma@Lemma!Glue}
Suppose that the adjacent cells \( x,x+\B  \) are siblings just before the time
\( t_{0} \), at which \( x \) breaks the sibling relation.
Suppose also that both \( x \) and \( x+\B  \) are unaffected by damage even via neighbors during
\( \rint{t_{0}-\Tu}{t_{0}} \).
Then we have the following possibilities:
 \begin{djenum}

   \item A cut as in Definition~\ref{def:cut}.

   \item Cell \( x-\B  \) was not a sibling of \( x \) at the last observation time of \( x \),
     and the switch kills \( x \) by \( \Decay \), \( \Purge \) or \( \Adapt \).

\end{djenum}
The statement also holds if we exchange left for right and \( x+\B  \) for \( x-\B  \).
 \end{lemma}

 \begin{Proof}
If a cell \( x \) breaks a sibling relation by a rule, then one of the
cases listed in the statement of the lemma holds.
This follows from the definition of siblings, 
Conditions~\ref{cond:addr-age} (Address and Age), \ref{cond:anim,kill}, 
and the healing rule.
We will show that if neither of the possibilities listed in the statement of
the lemma holds then the cells remain siblings.
  
  \begin{step+}{step:glue.new}
  Suppose that \( x \) or \( x+\B  \) were animated at some time in
  \( \rint{t_{0}-\Tu}{t_{0}} \); without loss of generality, suppose that
  the cell was \( x+\B  \).
   \end{step+}
  \begin{prooof}
  Then \( x+\B  \) will be without an age change or dying for at least \( \p_{0}\Tl \)
  time units which is longer than the whole period under consideration,
  as~\eqref{eq:p_{0}} implies \( \p_{0} \ge 5\lambda \).
If \( x \) also underwent animation during this interval then the same
is true for it, hence the two cells remain siblings.
Suppose therefore that \( x \) has been live during \( \rint{t_{0}-\Tu}{t_{0}} \).
The rule \( \Adapt \) implies that \( x+\B  \) is not a germ unless \( x \) is
a parent.
If \( x \) is a parent of \( x+\B  \), then part~\ref{i:parent.sib} of
Lemma~\ref{lem:parent} (Parent) implies
that the two cells remain siblings for at least \( \Tu \) time units; after
this, both cells have seen each other as siblings and therefore 
Condition~\ref{cond:addr-age} (Address and Age) shows that they remain siblings
until a cut or a death.

Suppose that \( x \) is not a parent of \( x+\B  \).
If it has changed its age within the last \( 4\Tu \) time units, then it
will not change the age for a long time after, and the two cells remain
siblings.
If it has not changed its age within this time then for at least
\( 2\Tu \) time units before the observation time before the animation, it was
exposed to the right (since \( x+\B  \) was latent), and therefore the rule~\ref{alg:Freeze} (\( \Freeze \))
froze it, keeping \( x \) and \( x+\B  \) siblings.
  \end{prooof} 

  \begin{step+}{step:glue.nonsib}
  Suppose now that both cells have been live during \( \rint{t_{0}-\Tu}{t_{0}} \).
  \end{step+}
 \begin{prooof}
If \( x \) changed its age within this time, then it will not change its
age soon, and therefore remains a sibling.
Suppose therefore that \( x \) does not change its age during this time.
Suppose that \( x+\B  \) became a sibling of \( x \) during this time.
Then it was not a sibling before, and the \( \Freeze \) rule must have frozen \( x \),
so \( x \) would not change its age so soon.
If \( x+\B  \) was a sibling all the time during \( \rint{t_{0}-\Tu}{t_{0}} \), 
then \( x \)
sees that \( x+\B  \) is a sibling and will not break the sibling relation.
This is guaranteed even within the rule \( \HealSync \), which may change age
in both directions.
 \end{prooof}
\end{Proof} 

\section{Gaps}\label{sec:gaps}
  
The main lemma of this section is the Running Gap Lemma, saying that if a
sufficiently large gap is found in a colony then this gap will not be
closed but will sweep through it, essentially eliminating a partial colony.
This serves as a preparation to the Attribution Lemma of the next section.

 \subsection{Paths}\label{sec:paths}
  
We will track the continuity of live areas in space-time via the notion of a path.

 \begin{definition}[Links]\label{def:links}
For times $t<u$, assume cell $x$ is live at $t$ and its body is
damage-free during $\clint{t}{u}$, and there are
no switching times in the open interval $\opint{t}{u}$.
Then we say that point $\pair{x}{t}$ is connected by a
 \df{vertical link}\index{link} to point \index{link!vertical}
$\pair{x}{u}$.
If one end of a vertical link is not a switching time, then the link
is called \df{short}. \index{link!vertical!short}
If cells $x,x+\B $ are siblings at time $t$ or $(t-)$ such that 
$\lint{x-\B }{x+2\B }$ is damage-free during $\rint{t-\Tu}{t}$,
then the points $(x,t)$, $(x+\B ,t)$ are said to be
connected by a \df{horizontal link} \index{link!horizontal}.
If point $\pair{y}{t'}$ is, according to case~\ref{i:parent.gen} of 
Lemma~\ref{lem:parent} (Parent) a parent of point
$\pair{x}{u}$, we will say that $\pair{y}{t'}$ is connected by a \df{parental}
(maternal or paternal) link to point 
$\pair{x}{u}$.
A \df{link} is a link of one of these kinds.
A link is \df{steep} or, equivalently, \df{slow} if it is a non-short
vertical link or a 
parental link.
(Since time is the second coordinate, steepness of a line is
synonymous to slowness of the  movement of a point along it.)
 \end{definition}

   By Lemma~\ref{lem:parent} (Parent),
   the parental link can be replaced by a vertical connection and a horizontal
   connection.
   The horizontal connection is to a sibling that is a successor to the parent,
   and the vertical ones continue back in time towards the parent.

\begin{definition}[Path]
A sequence $\pair{x_{0}}{t_{0}}$, $\dots$, $\pair{x_{n}}{t_{n}}$ with $t_{i}\le t_{i+1}$
such that subsequent points are connected by links, is called a
\df{path}.
A path with only steep links is \df{steep}, or \df{slow}.
A \df{backward path} is the reversed reading of a forward path,
backwards in time.
The adjective ``forward'' or ``backward'' will be omitted, when it is obvious from
the context.
For a path $P=\pair{x_{0}}{t_{0}},\dots,\pair{x_{n}}{t_{n}}$ and
$t\in\clint{t_{0}}{t_{n}}$, let%
 \glo{propo:cap@$P(t)$}%
 \[
  P(t)
 \]
 be $x_{i}$ for the smallest $i$ with $t\in\clint{t_{i}}{t_{i+1}}$.
\end{definition}

A point $\pair{x_{i}}{t_{i}}$ on a path can actually be dead,
if it has just died: indeed, it can be connected for example to a point
$\pair{x_{i+1}}{t_{i+1}}$ by a horizontal link with $t_{i}=t_{i+1}$
such that $x_{i},x_{i+1}$ are siblings at time $(t_{i}-)$.
The following lemma says that if two paths cross then they
have a common point.

 \begin{lemma}[Crossing]\label{lem:crossing}\index{lemma@Lemma!Crossing}
  Let $\pair{x_{1}}{s_{1}},\dots,\pair{x_{m}}{s_{m}}$ and
$\pair{y_{1}}{t_{1}},\dots,\pair{y_{n}}{t_{n}}$ be paths with $s_{1}=t_{1}$,
$s_{m}=t_{n}$, $x_{1}\le y_{1}$, $x_{m}\ge y_{n}$.
  Then there are $i,j$ such that $x_{i}=y_{j}$ and either
$t_{j}\in\clint{s_{i}}{s_{i+1}}$ or $s_{i}\in\clint{t_{j}}{t_{j+1}}$.
\end{lemma}

\begin{proof}
Let us call the point $\pair{x_{i}}{s_{i}}$ whose existence is asserted, the
\df{crossing point}.
Let us replace parental links with horizontal and vertical connections as described
in the remark after Definition~\ref{def:links}.
Now the two paths cannot jump over each other.
Indeed, at the time of the crossing, the cells involved have to
be at a distance \( \B  \).
At this crossing, one of the links must be horizontal;
suppose it is \( \pair{x_{i}}{x_{i+1}} \).
Then either \( \pair{x_{i+1}}{s_{i}} \) is on a vertical link between some \( t_{j} \) and \( t_{j+1} \)
or there is also some horizontal a link \( \pair{y_{j}}{y_{j+1}} = \pair{x_{i+1}}{x_{i}} \)
at time \( s_{i}=t_{j} \).
In both cases we can choose the crossing point \( x_{i+1}=y_{j} \) at time \( s_{i} \).
  \end{proof}

According to the Parent Lemma, a steep path can be continued
backwards in time until it hits some island
(as in Definition~\ref{def:damage-map}).
Moreover, occasionally we have a choice between two parents to continue to.
The lemma below says that any path started backwards not too soon
after damage (the time defined as $\purgeT$ in~\eqref{eq:purgeT})
can be diverted and continued back past any island.

\begin{definition}[Traceability]\label{def:traceable}
 A cell $\pair{x}{t}$ will be called \df{traceable} if $I_{0}(x)$ is damage-free
during $\rint{t-4\purgeT}{t}$.
\end{definition}

\begin{lemma}[Skirting]\label{lem:skirting}\index{lemma@Lemma!Skirting}
Let $\lint{a_{0}}{a_{1}}\times\rint{u_{0}}{u_{1}}$ 
denote the (only) island in the area under consideration.
Consider a traceable live point $\pair{x_{0}}{t_{0}}$.
There is a path going backwards from $\pair{x_{0}}{t_{0}}$ and ending
either in time $u_{0}$ or in a birth (in the sense of Rule \( \Birth \) of
Algorithm~\ref{alg:Birth}).
It has at most $\numDirAff+1$ non-steep links (with \( \numDirAff \) as in~\eqref{eq:numDirAff})
which, with one possible exception,
form a series of up to $\numDirAff$ consecutive horizontal links,
and can only cross a colony boundary in the inward direction, and only if
$\Age<\growEnd+\numDirAff$ on all cells involved. 

 \end{lemma}
\begin{Proof}
Let us start constructing a steep path $c_{0},\dots,c_{n}$ with
 \[
  c_{i}=\pair{x_{i}}{t_{i}}
 \] 
  backwards in time from $\pair{x_{0}}{t_{0}}$, consisting of either vertical
links, or parental links if the conditions of Lemma~\ref{lem:parent} (Parent)
are applicable.
If we get to time $u_{0}$ then we are done.
Otherwise 
let us stop just before going below $u_{1}+\Tu$, with $c_{k}$ being the last element.
Then the body of $x_{k}$ intersects $\lint{a_{0}-1.1 \B }{a_{1}+1.1 \B }$:
indeed, otherwise we could continue the path either by
a vertical or by a parental link.
The vertical link would be shorter than $\Tu$, and the parental
link would lead to a damage-free cell, so either of them would be
allowed. 

\begin{step+}{step:skirting.bridge}
There is an $i\le k$ such that
$\pair{x_{i}}{t_{i}}$ can be connected by horizontal
links that do not cross a colony boundary in the outward direction,
to a cell $y$ unaffected by the damage even via neighbors at any time (in the rectangle
under consideration).
\end{step+}
\begin{pproof}
Let us go back on the path for $k,k-1,\dots$ until the first
$i$ (counting from $k$) such that either $t_{i}>u_{1}+(\numDirAff+1)\tau_{0}$, or
$c_{i}$ is a parent of $c_{i-1}$; let us call it $i_{1}$.

Suppose first $t_{i_{1}}>u_{1}+(\numDirAff+1)\tau_{0}$.
Then there is some $i_{1}<i\le k$ such that 
$\pair{x_{i}}{t_{i}}$ can be connected by horizontal links to a directly
unaffected cell $y$. 
Indeed, if there is not, then $x_{i}=x_{k}$, $i=k,k-1,\dots,i_{1}$, since we have not
encountered a parent before.
The part of the domain containing $x_{k}$ that consists of cells directly affected by damage
has size $\le\numDirAff$.
If we cannot pass away from it via siblings within the same colony or towards the originating
colony then at least one end will be exposed and the depth will also be bounded by \( \numDirAff \).
This allows the $\Purge$ rule to gradually eliminate this part, including $x_{k}$,
leading to a contradiction.
The time this takes is at most $\numDirAff\Tu$ to propagate the depth
information, 
followed by $\numDirAff\tau_{0}$ to apply the steps of $\Purge$ until $x_{k}$ is
reached, for a total time of $\numDirAff\tau_{0}+\numDirAff\Tu$.
This shows $t_{0}\le u_{1}+(\numDirAff+1)\tau_{0}$ since \( \numDirAff<\p_{0} \).
By the definition of \( \purgeT \) this contradicts the assumption
of the traceability of  \( \pair{x_{0}}{t_{0}} \).
If the only way to get by horizontal links to a directly
unaffected cell is to cross a colony boundary in the 
inward direction with $\Age\ge\growEnd+\numDirAff$ then
there is an exposed cell in the opposite direction, at distance
$\le\numDirAff \B $, so $\Purge$ will act again in the same way.

Suppose now that $c_{i_{1}}$ is a parent of $c_{i_{1}-1}$.
Let us then go back on the path for $i=i_{1},i_{1}-1,\dots$ until the first
$i$ (counting from $i_{1}$) such that $t_{i}>t_{i_{1}-1}+(\numDirAff+1)\tau_{0}$.
There is such an $i$, and we have $x_{i}=x_{i_{1}-1}$ for all $i$, since
the newborn $x_{i_{1}-1}$ cannot be a parent sooner.
Indeed, Condition~\ref{cond:anim,kill} says that animation happens only via 
healing or growth, and even the faster healing step has a delay of \( \p_{1} \).
By the same reasoning as above, there is some $i_{1}<i\le k$ such that 
$\pair{x_{i}}{t_{i}}$ can be connected by a sequence of at most $\numDirAff$
horizontal links not crossing a colony boundary in the outward direction,
to a directly unaffected cell $y$.
\end{pproof} 

\begin{step+}{step:skirting.finish}
There is a steep path backwards in time from the point $\pair{y}{t_{i}}$ computed
above, going only either on vertical or parental links, and ending below time $u_{0}$. 
\end{step+}
\begin{pproof}
Let us construct such a path.
If it contains at most one parental link, then the body of cells it reaches
still is damage-free, hence the path can be continued.
On the other hand, it cannot contain two parental links, since
a newborn $x_{i_{1}-1}$ cannot be a parent soon.
\end{pproof} 
 \end{Proof}

\begin{definition}[Trace-back path]
A backward path is called a \df{trace-back path} if every link in it is
chosen to be a vertical or a parental link if this is possible, 
with horizontal links only as needed in Lemma~\ref{lem:skirting} (Skirting).
\end{definition}

\begin{definition}[Age progress]\label{def:age-lag}
Let $P=\tup{c_{0},\dots,c_{n}}$ be a forward path.
The \df{age progress} of  $P$ is defined as
 \begin{align*}
   \var{age-progress}(P) = \sum_{i=0}^{n-1}(\Age(c_{i+1})-\Age(c_{i}) \amod \U ).
 \end{align*}
\end{definition}

 The following lemma upper-bounds the number of colonies a trace-back path can
cross and also its age progress.

\begin{lemma}\label{lem:trace-back-bounds}
  Suppose that a trace-back path $P$ with $n$ links 
has time projection $d<\Tls/4$.
Then the following holds.
\begin{cjenum}
 \item \label{i:trace-back-bounds.colonies} $P$ passes through at most 3
   colonies, for a total space projection of $<2.1\Q \B $.
 \item\label{i:trace-back-bounds.time} $d\ge (n-\numDirAff-1)\Tl/2$. 
If the path has no parental links then the factor $1/2$ can be omitted.
 \item \label{i:trace-back-bounds.age} The age progress of $P$ is at most
   $n/\p_{1}+\numDirAff+2$. 
\end{cjenum}
\end{lemma}
\begin{proof}
Let us show that $P$ can pass through at most 3 colonies.
It moves horizontally only along parental and horizontal links.
The latter only occur in connection with the application of 
Lemma~\ref{lem:skirting} (Skirting), with at most $\numDirAff$ horizontal links per
island. 
Given the time bound on the whole path, as long as it does not cross more than
three colonies, 
at most one island will appear, hence we do not have to count with more
than one series of $\le\numDirAff$ horizontal links.

Parental links can only occur in connection with healing and growth.
Healing moves the forward path towards the edge of a colony as shown
in rule \( \Heal.\passive \) of Algorithm~\ref{alg:Heal-passive}, 
only growth can move it towards the center.
In order to move towards the center of another colony, another growth cycle must
start (except in case of germ growth, which can span three colonies).
But definition~\eqref{eq:compute-start} implies that within the time bound
$\Tls/4$ the path cannot reach from the end of one growth to the beginning
time of another.

Let us prove the time projection lower bound.
There are at most $\numDirAff$ horizontal links and 1 non-steep link, 
arising from the island.
The other at least $n-\numDirAff-1$ links are steep.
A steep link is either vertical, in which case it has size
$\ge\Tl$, or parental, in which case it has size $\ge\Tl/2$.

Let us estimate now the age progress.
By the definition of the animations in healing and growth, age does not increase
along parental links.
Along horizontal links the age can increase by 1, and it can also increase immediately
after the last horizontal link.
It can also increase by 1 at the first link (in time).
But during the remaining vertical links, age increase occurs no
sooner than every $\p_{1}$ steps, hence the total age increase is at most
$\numDirAff+2+n/\p_{1}$.
\end{proof}

\subsection{Running gaps}\label{sec:decay}

The rule $\Decay$ attempts to widen any gap that was not closed in reasonable
time. 
It erases all the cells created by $\HealRevive$, 
as well as one more cell, since they were marked
by $\Decaying$.
The delay $\p_{2}$ of decay gives a chance for the healing process to complete.
The difference between killing by $\Purge$ and killing by $\Decay$ is that
$\Purge$ kills fast, but its reach, in case of non-doomed member cells,
is only local: healing is intended after it.
On the other hand, $\Decay$ slowly eliminates the remainders of a
colony if healing fails to close a gap within a certain time limit.

   \begin{algo}
  \caption{rule $\protect\Decay$}\label{alg:Decay}
      \lIf{$\Decaying=1$ \algAnd $\exist{j\in\{-1,1\}}\Xposed_{j}$
        }{$\Die(\p_{2})$}
   \end{algo}
 
Gaps and gap paths are a tool to show that the decay rule indeed destroys
incomplete colonies.
Their definition is complicated by the consideration of germs, needed
for self-organization.
Here are some constants; 
for readability, we omit the notation $\flo{\cdot}$ for integer part.
Recall the definition $\tau_{i} = (\p_{i}+1)\Tu$ in~\eqref{eq:tau_{i}},
\( \numDirAff \) in~\eqref{eq:numDirAff} and \( \healSpan \) in~\eqref{eq:heal-span}.
 Let \glo{split-t@$\splitT$}%
\begin{align}
\nn     \islSz &= 2.5,
\\\label{eq:gap-lb}
            \gapLb &=4\healSpan+2\islSz+21,
\\\label{eq:purgeT}
            \purgeT              &= (\numDirAff+1)\tau_{0},
\\\label{eq:split-n}
          \splitN   &=(\healSpan+1)(2\gapLb+1),
\\ \label{eq:split-t}
        \splitT     &= \splitN\tau_{2}.
\end{align}                        
 
\begin{definition}[Gaps]\label{def:gap}
Consider an open interval $G=\opint{l}{r}$ 
at some time $t$, where $0 < r - l$ is divisible by $\B $.
We say that $G$ is a \df{right gap}\index{gap}
with \df{size} $|G|=r-l-\B $ if 
every traceable cell in $G$ can be space-consistent with $r$ only if it is a
germ cell. 
The \index{gap!right-age}\df{right age} of the gap $G$ is
is the upper bound of the age of these germ cells (it is 0 if the gap contains
no germ cells).
If $G$ is contained in the colony of $r$ then it is called an \df{interior gap}.
A \df{right bad gap}\index{gap!bad} is a gap of right-age
\begin{align*}
  \lambda\splitN(\healSpan/2+1) 
\end{align*}
containing an interior gap of size~$\ge\gapLb\B$. 
If the gap itself is not interior then its right edge is required to be
an exposed right outer cell.
Left gaps are defined similarly.
\end{definition}



We are interested in how a gap changes in time.

\begin{definition}[Gap path]
  Suppose that in a time interval $\clint{v_{0}}{v_{1}}$, the gap
$G(t)=\opint{l(t)}{r(t)}$ is defined for all $t$, in such a way that all cells
$r(t)$ are space-consistent with each other, and $G(t_{1})\cap
G(t_{2})\ne\emptyset$ for all $t_{1}$, $t_{2}$ with $|t_{2}-t_{1}|\le 3\Tu$.
Then $G(t)$ is called a (right) \df{gap path}, and the \df{right age} of
the gap path is the maximum of the right ages of the gaps in it.
  \end{definition}

Let us make the following simple observation.

  \begin{lemma}
If a path space-consistent with $r(t)$ has the same
time projection $\clint{v_{0}}{v_{1}}$ as the right gap path $G(t)$ and has no
germ cells younger than the right age of the gap path then it cannot cross
$G(t)$.
  \end{lemma}

The lemma below says that the decay rule causes a large enough gap
to move right rather fast.
The gap is assumed to be connected, via a series of horizontal links, to
a forward path.
This excludes irregular and unimportant cases when the gap would have to
travel through all kind of debris.

  \begin{lemma}[Running Gap]\label{lem:rn-gap}
Recall the definition of $\gapLb$ in~\eqref{eq:gap-lb}.
Let $P_{1}=\pair{x_{0}}{v_{0}},\dots,\pair{x_{n}}{v_{n}}$ be a trace-back path
(listed forwards),
let $L,k$ be positive integers with
 \begin{align*}
       L &< \Q ^{2},
\\   k &< \U  - 3 L\p_{2}\lambda.
 \end{align*}
  Assume the following:
 \begin{cjenum}

  \item $y_{0}<x_{0}$ are in the same domain at time $v_{0}$.

  \item The interval $\lint{y_{0}-2.1\B }{x_{0}+3.1\B }$ 
has been noise-free during $\rint{v_{0}-\Tu}{v_{0}}$.

  \item\label{i:rn-gap.gap-bds} $y_{0}$ is an exposed left edge at time $v_{0}$ (see 
Section~\ref{sec:edges}), and 
has been the right end of a  bad gap of right-age $k$
during $\rint{v_{0}-\Tu}{v_{0}}$.
Its age is not within $\numDirAff$ steps of any of the age values listed in
Lemma~\ref{lem:exposing} (Exposing),  with \( \numDirAff \) defined in~\eqref{eq:numDirAff}.

  \item\label{i:rn-gap.in-col}
If $y_{0}$ is a left outer cell then $P_{1}$ is in the same colony.
In this case, let $Z$ be the starting cell of the originating colony.
Otherwise, $Z=\infty$.

  \item\label{i:rn-gap.time-proj-bd}
  $v_{n}-v_{0}\le L\tau_{2}$.

  \end{cjenum}

  Then during $\rint{v_{0}}{v_{n}}$, a right gap path $G(t)=\opint{l(t)}{r(t)}$ can be
defined, with $r(v_{0})=y_{0}$,
 \begin{align}\label{eq:gap-path.fast}
     r(v_{n}) &\ge  Z\land \Paren{y_{0} + 
       \B \bigparen{\frac{|v_{n}-v_{0}-\purgeT|^{+}}{3\tau_{2}}-2\healSpan-\islSz-11}},
 \end{align}
\( \nonumber |G(t)| >0 \), with right-age $< k + 3 L\lambda\p_{2}$.
  \end{lemma}

\begin{figure}
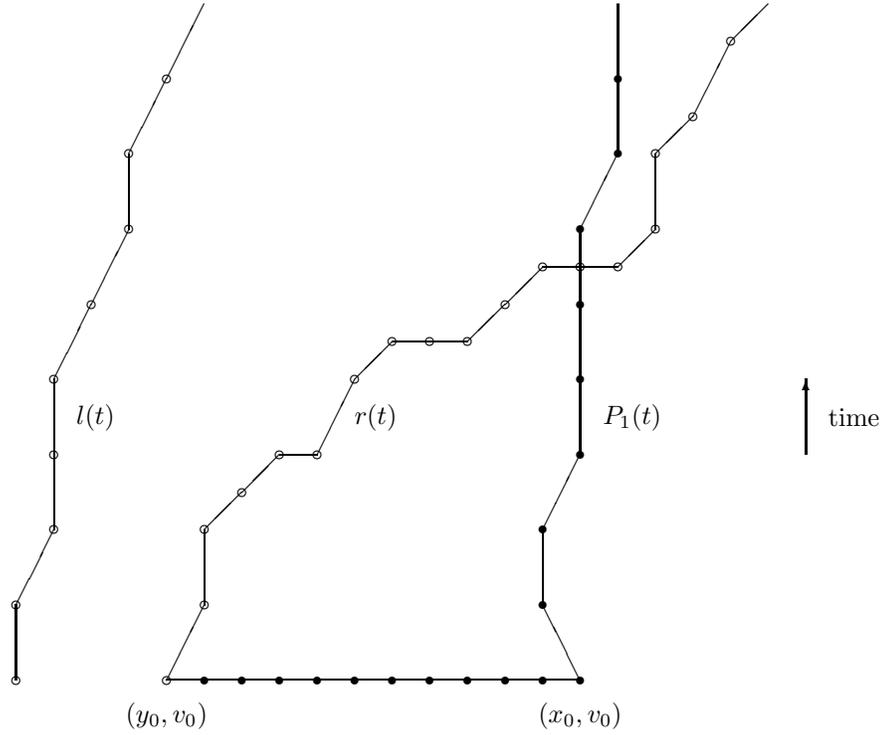

\input tex-figures/running-gap-fig.tex
 \caption{Running Gap Lemma}
 \label{fig:running}
\end{figure}

The following corollary says (in contrapositive) that if path $P_{1}$
is long, then there is no large young gap next to its beginning (say, on the
left).

\begin{align}
\label{eq:destr-time}
      \destrT  &=  10\tau_{2}, 
\\\label{eq:destr-steps}
      \destrS  &=\destrT/\p_{1}\Tl.
\end{align}
We will also use the following lower bound on \( \Q \):
\begin{align}
\label{eq:Q-req-1}
         \Q  & > 3(2\healSpan+\islSz+11)+\purgeT/\tau_{2}.
\end{align}

\begin{corollary}[Running Gap]\label{crl:rn-gap}
  Assume the conditions of the Running Gap Lemma, further that
$\pair{x_{n}}{v_{n}}$ is not a germ cell younger than $k + 6L\lambda\p_{2}$.
Then $v_{n}-v_{0} < \destrT\Q$.
  \end{corollary}
Note that the statement is meaningless unless $L>10 \Q $, since
$v_{n}-v_{0}<L\tau_{2}$ is a condition of the Running Gap lemma.
  \begin{proof}
It is easy to see that the length of the path is at most $3L\lambda\p_{2}$
(the proof is just like the proof of the lower bound on the length of path
$P_{2}$ in the proof of the Running Gap Lemma below).
Path $P_{1}$ starts to the right of the gap path $G(t)$.
Since age varies by at most 1 per link along a path 
and since $\pair{x_{n}}{v_{n}}$ is not a germ cell younger than
$k+6L\lambda\p_{2}$, no cell on $P_{1}$ is a germ cell younger than
$k+3L\lambda\p_{2}$, while according to Lemma~\ref{lem:rn-gap}, all cells in
$G(t)$ space-consistent with $r(t)$ are germ cells with such age.
Therefore $P_{1}$ never crosses the gap path.
Inequality~\eqref{eq:gap-path.fast} gives a lower bound on how
fast $r(t)$ moves right.
According to Lemma~\ref{lem:trace-back-bounds},
$P_{1}$ passes at most 3 colonies.
Edge $r(t)$ therefore can move right by at most $3 \Q $, giving
 \begin{align*}
   3 \Q  &\ge (r(v_{n})-y_{0})/\B
\ge \frac{v_{n}-v_{0}-\purgeT}{3\tau_{2}}-2\healSpan-\islSz-11,
\\ v_{n}-v_{0} &\le 9\tau_{2}\Q  + 3\tau_{2}(2\healSpan+\islSz+11)+\purgeT 
< 10\tau_{2}\Q,
 \end{align*}
 where we used~\eqref{eq:Q-req-1} and~\eqref{eq:destr-time}.
\end{proof}

 \begin{Proof}[Proof of Lemma~\ref{lem:rn-gap}]
Let $r(v_{0})=y_{0}$, and let $l(v_{0})$ be the leftmost cell such that
the gap $\opint{l(v_{0})}{r(v_{0})}$ has right-age $\le k$.
   \begin{step+}{step:rn-gap.inward}
Let $t_{1}>v_{0}$.
Assume that for all $t\le t_{1}$, a gap path $G(t)$ is defined with the
desired properties and in such a way that $\pair{r(t_{1})}{t_{1}}$ is not a germ
cell younger than $k+6L\p_{2}\lambda$, and is traceable.
Then $r(t_{1})$ is an exposed edge, space-consistent with $r(v_{0})$.
    \end{step+}
  \begin{pproof}
    \begin{step+}{step:rn-gap.inward.path}
Let us build a path.      
    \end{step+}
    \begin{prooof}
We start by building a backward trace-back path of length $m$
 \[
  P_{2}=(\pair{z_{0}}{w_{0}}=\pair{r(t_{1})}{t_{1}},\dots,\pair{z_{m}}{w_{m}})
 \]
ending at time $v_{0}$.
Lemma~\ref{lem:trace-back-bounds} and the bound~\ref{i:rn-gap.time-proj-bd}
gives
 \[
  m \le \numDirAff+1+2(v_{n}-v_{0})/\Tl \le \numDirAff+1+2L\tau_{2}/\Tl 
  <  3L\lambda\p_{2}. 
 \]
The path does not end in a birth; indeed, otherwise $\pair{r(t_{1})}{t_{1}}$ would be a germ
cell younger than $3L\p_{1}\lambda$ contrary to the assumption.
It does not enter the gap; otherwise $\pair{r(t_{1})}{t_{1}}$ would be a 
germ cell younger than $k+3L\lambda\p_{2}+3L\lambda\p_{2}$, which was excluded.
Therefore defining $x=P_{2}(v_{0})$ gives $x\ge y_{0}$.
Without loss of generality, we can suppose $x\le x_{0}$.
Otherwise, Lemma~\ref{lem:crossing} (Crossing) implies that path $P_{2}$
crosses $P_{1}$ and we can switch from $P_{2}$ to $P_{1}$ at the meeting
point.
Thus, $x$ is on the interval $\clint{y_{0}}{x_{0}}$, aligned with $x_{0}$.

Combine $P_{2}$ and the horizontal path $H$ from $\pair{x}{v_{0}}$ to
$\pair{y_{0}}{v_{0}}$ into a single chain of links.
It has at most $3\Q $ horizontal links on $H$, 
and then at most $\numDirAff+1+2(v_{n}-v_{0})/\Tl$
links on $P_{2}$, giving fewer than
 \[
   3\Q  + 2(v_{n}-v_{0})/\Tl + \numDirAff+1
 \]
 links.
 We have $r(t)\ge \B(y_{0}-2\healSpan-\islSz-11)$. 
 Indeed, the only way for \( r(t) \) to move left is via an island (can happen at most once)
 that is then connected to the right end of the gap by healing.
    \end{prooof}

    \begin{step+}{step:rn-gap.inward.exposed}
A protected left edge cannot occur on \( H+P_{2} \).      
    \end{step+}
    \begin{pproof}
      Suppose first that $\pair{y_{0}}{v_{0}}$ is a member or right outer cell,
      or a germ cell with \( \Addr>0 \).
Then the large bad gap on its left contains, by the definition of bad gaps,
a gap in the same colony or in the originating colony.
The path $P_{2}$ could cross it only by repeated parental links.
But the size of the gap is too large for using rule $\HealRevive$, since the field $\fld{Healing}$
limits the number of applications.
Instead, $P_{2}$ (when traveled forward) would first have to walk right via
parental links to
participate in the creation of a colony and then walk back via parental links through the
expansion process from that colony---for which there is not enough time due 
to~\eqref{eq:compute-start}. 

Suppose now that $\pair{y_{0}}{v_{0}}$ is a left outer cell.
By requirement~\ref{i:rn-gap.in-col}, then the path $P_{1}$ stays in the
colony of $y_{0}$, and hence so does the whole gap path.
Since the cell is exposed, its age is already past the expansion stages.
Walking right, the age of cells on the path $H$ (which remains in the same colony
as well) does not decrease
except when it crosses into another work period (this is taken care of by our Definition~\ref{def:age-lag}
of age progress that uses \( {}\amod \U \)).
Indeed, by Definition~\ref{def:time-consistency} (Time consistency), the ages of siblings
must be non-increasing as one is moving away from the center of the originating colony.
It is also non-decreasing on the trace-back part $P_{2}$, except possibly on the 
$\le\numDirAff$ horizontal links allowed by that Lemma~\ref{lem:skirting} (Skirting).
Therefore the age of $\pair{r(t_{1})}{t_{1}}$ cannot be smaller by more than
$\numDirAff$ than that of $\pair{y_{0}}{v_{0}}$.
If the age is not smaller at all then $\pair{r(t_{1})}{t_{1}}$ is also exposed, since
the number of steps on the path is not sufficient to get into an expansion stage
of the next work period.
If decreasing the age by $\numDirAff$ makes it protected, then one of
the cases in Lemma~\ref{lem:exposing} (Exposing) must have occurred within the
last $\numDirAff$ steps before the age of $\pair{y_{0}}{v_{0}}$: but this has been
excluded by the condition~\ref{i:rn-gap.gap-bds} of the lemma we are proving.

Suppose that $\pair{y_{0}}{v_{0}}$ is a germ cell with \( \Addr<0 \).
Since it is exposed to the left its age is \( \ge\growEnd \), 
and its address is \( >-\Q \).
In this case, the same argument works as for left outer cells above.
    \end{pproof}
 \end{pproof} 

A large gap cannot be closed by the healing rule, but it can be narrowed somewhat.
Suppose that for a gap $G=\opint{l}{r}$ the cells $r,r+\B ,\dots,r+k\B $ are siblings
with $\Decaying(r+i\B )=1$, $\Healing(r+i\B )=k-i$.
Then we set $r'=r+k\B$.
If these conditions are not satisfied then set $r'=r$.
We define similarly, $l'=l-k\B$ for an appropriate $k$.
The values $l',r'$ will be called the \df{adjusted edges} of the gap,
and we will write $G'=\opint{l'}{r'}$.
the \df{adjusted size} is $|G'|=r'-l'-\B $.
 
 \begin{step+}{step:rn-gap.no-dam}   
Assume that the gap $G(t)$
with the desired properties was defined up to time $f_{0}$, with adjusted size
$|G'(f_{0})|\ge 2\healSpan \B +2$ and right age $\le k$.
Assume also that in the area where we define the
path further, all cells are unaffected even via neighbors during
$\rint{f_{0}-\purgeT}{f_{1}}$. 

Then the gap path $G(t)$ can be defined further in $\rint{f_{0}}{f_{1}}$ in
such a way that as long as $r(t)$ did not reach $Z$:
 \begin{align}\label{eq:progress-right}
     (r(f_{1})-r(f_{0}))/\B             &\ge
     \frac{f_{1}-f_{0}}{3\tau_{2}}-\healSpan-1,
\\\label{eq:progress-size}
       (|G(f_{1})|-|G(f_{0})|)/\B  &\ge \frac{f_{1}-f_{0}}{3\tau_{3}} - 2\healSpan - 2.
 \end{align}
If $r(t)$ reaches $Z$ it never decreases again.
   \end{step+}
  \begin{pproof}
  Define $G(t)$ as follows.   
Suppose that it was defined up to time $t_{1}$, and let $t_{2}$ be the
next time that is a switching time of either $l(t_{1})+\B $ or $r(t_{1})-\B $ 
or $r(t)$.  
  We distinguish the following cases.
   \begin{djenum}

    \item  $t_{2}$ is a switch of $l(t_{1})+\B $.
  If the switch is an animation creating a
sibling of $l(t_{1})$, then $l(t_{2})=l(t_{1})+\B $.

    \item  $t_{2}$ is a switch of $r(t_{1})-\B $.
  If the switch is an animation creating a
sibling of $r(t_{1})$, then $r(t_{2})=r(t_{1})-\B $.

    \item  $t_{2}$ is a switch of $r(t)$.
  If $r(t_{1})$ dies, then $r(t_{2})$ is the closest cell to the right of
$r(t_{1})$ that is not a germ cell younger than $k+2(t_{2}-f_{0})/\Tl+1$.

    \item  In all other cases, we leave $G(t)$ unchanged.

    \end{djenum}

If $r(t)$ reaches $Z$, it will never decrease again.
Indeed, this happens only if left outer cells were exposed to the left.
Without loss of generality, assume that these outer cells were growth cells.
In this case, $\Age(Z)>\growEnd$.  
As $Z$ is protected to the left, the 
healing rule does not create a left neighbor of $Z$.
The only way that $r(t)=Z$ could move left is by the growth rule, but it
does not work for $\Age(Z)>\growEnd$.

Let us show
 \begin{align}\label{eq:rn-gap.no-dam.progress.spread}
       (r'(f_{1})-r'(f_{0}))/\B &\ge \Flo{\frac{f_{1}-f_{0}}{3\tau_{2}}},
 \end{align}
which will prove~\eqref{eq:progress-right}.
Here we assume that $Z$ has not been reached.
Since the conditions of part~\ref{step:rn-gap.inward} are satisfied, 
the cell $y=r(t)$ is an exposed left edge for all $t\in\rint{f_{0}}{f_{1}}$.
Hence the longest it can take to move $r'(t)$ right is to
have the rule $\HealRevive$ move $r(t)$ left by 
up to $\healSpan$ cells; then all of them die by decaying.
The $\HealRevive$ steps and the killing of the healing cells
take time at most $\tau_{0}$ each, for a total
of \( 2\healSpan\tau_{0}=\tau_{2} \),
setting the field $\Decaying$ takes time at most $\tau_{2}$, and the
final killing again time $\tau_{0}$.

Now let us show~\eqref{eq:progress-size}.
The inequality
 \begin{align}\label{eq:rn-gap.no-dam.progress.monot}
        (l'(f_{1})-l'(f_{0}))/\B  &\le \Cei{\frac{f_{1}-f_{0}}{\p_{3}\Tl}}
 \end{align}
follows from the fact that $l'(t)$ can only increase when
it is equal to $l(t)$ and increases via the the growth rule.
But the next application of $\Grow$ is away by a waiting
period\index{period!wait} of length $\ge\p_{3}\Tl$.

Let us combine~\eqref{eq:rn-gap.no-dam.progress.monot} 
and~\eqref{eq:rn-gap.no-dam.progress.spread}.
\begin{align*}
   (|G'(f_{1})|-|G'(f_{0})|)/\B  
  &\ge (f_{1}-f_{0})\Paren{\frac{1}{3\p_{2}\Tu}-\frac{1}{\p_{3}\Tl}} - 2.
    \end{align*}
Now, using the definition of \( \p_{3} \) in~\eqref{eq:p_{3}},      
    \begin{align*}
\frac{1}{3\p_{2}\Tu}-\frac{1}{\p_{3}\Tl}=\frac{\p_{3}\Tl-3 \p_{2}\Tu}{3 \p_{2} \p_{3}\Tl\Tu}
\ge \frac{\p_{3}-3 \p_{2}\lambda}{3 \p_{2} \p_{3}\Tu}
\ge \frac{1}{3\p_{3}\Tl},
 \end{align*}
hence the estimate~\eqref{eq:progress-size}.
 \end{pproof} 
 
 \begin{step+}{step:rn-gap.dam}
   Let $\lint{a_{0}}{a_{1}}\times\rint{u_{0}}{u_{1}}$ be an island with
   $u_{0}\in \rint{v_{0}}{v_{n}}$.
Then the gap path $G(t)$ with the desired properties can be defined
for $t\in \rint{u_{0}}{u_{1}+\purgeT}$.
We will have
 \begin{align*}
        r(u_{1}+\purgeT) &\ge r(u_{0})-(\islSz+9)\B ,
\\ |G(u_{1}+\purgeT)|&\ge |G(u_{0})|-(2\islSz+19)\B .
 \end{align*}
  \end{step+}
  \begin{pproof} 
Let us define $l(t)$ for $t\in\rint{u_{0}}{u_{1}+\purgeT}$.
If $a_{0}\le l(u_{0})+5.1\B $,
then we set $l(t)$ to be the first site $x$ aligned with $l(u_{0})$
to the right of $a_{1}$ and $l(u_{0})+\B $ that is unaffected even via neighbors.
This definition is justified, since growth or healing
do not have time to add any cell further to the right of $l(u_{0})+\B $.
Now $l(t)\le l(u_{0})+(5.1+\islSz+3.1)\B <l(u_{0})+\islSz+9\B $.
If $a_{0} > l(u_{0})+5.1\B $, let $l(t)=l(u_{0})+\B $ since, as noted above,
at most one cell can be added to $l(u_{0})$ by growth or healing during
$\lint{u_{0}}{u_{1}+\purgeT}$.
The definition of $r(t)$ is done symmetrically similarly (but only as a lower
bound).

Thus, the island may decrease the gap by at most
$\islSz+9\B $ on one side; with another possible decrease of size $\B $ on the
other side, this gives the size estimate of statement~\ref{step:rn-gap.dam}.
\end{pproof}

In the space-time rectangle considered, at most one island
occurs: 
indeed, \( v_{n}-v_{0}\le\Tls/2 \) follows from~\ref{i:rn-gap.time-proj-bd}
and~\eqref{eq:compute-start}.
Use the construction of part~\ref{step:rn-gap.no-dam} up to the time \( u_{0} \) of the
start of the island, then construction of part~\ref{step:rn-gap.dam}
from \( u_{0} \) to \( u_{1}+\purgeT \) (or \( v_{n} \),
whichever comes first), then possibly again the 
construction of part~\ref{step:rn-gap.no-dam} from \( u_{1}+\purgeT \) to \( v_{n} \), we
defined \( G(t) \) for all \( t\in\clint{v_{0}}{v_{n}} \).

 \begin{step+}{step:rn-gap.age} For all \( t\in\clint{v_{0}}{v_{n}} \), the gap
\( G(t) \) has right-age \( \le k+2(t-v_{0})/\Tl+1 \).
In particular, according to the assumed bound on \( v_{n}-v_{0} \),
at time \( v_{n} \), it still has age at most
 \begin{align*}
   k+2(v_{n}-v_0)/\Tl+1 \le k+2L\tau_{2}/\Tl+1 < k+3L\p_{2}\lambda,
 \end{align*}
as stated in the lemma.
 
    \end{step+}
   \begin{pproof} 
Let \( z_{1} \) be a live cell in \( G(t) \) space-consistent with \( r(t) \),
and let us build a sequence \( \pair{z_{1}}{w_{1}} \), \( \pair{z_{2}}{w_{2}} \), \( \dots \)
with \( f_{1}=w_{1}\ge w_{2}\ge\dotsm \) as follows.
Without loss of generality, assume that \( w_{1} \) is a switching time
of \( z_{1} \). 
If \( z_{1} \) is live at \( w_{1}- \) then \( z_{2}=z_{1} \) and \( w_{2} \) is the previous
switching time of \( z_{1} \).
Otherwise, either \( \pair{z_{1}}{w_{1}} \) is a newborn germ cell, in which case the
sequence ends at \( \pair{z_{1}}{w_{1}} \), or it has a parent \( \pair{x}{u} \).

If \( x\in G(u) \) then \( z_{2}=x \), \( w_{2}=u \).
Let us show that there is no other possibility.
Indeed, otherwise, \( z_{1} \) would have been 
created by animation which, by Condition~\ref{cond:anim,kill}
can only be by healing or growth.
Since \( r(t) \) is exposed, growth can only happen from the left.
We end up with one of the cases of the definition of \( G(t) \) above, and thus
with \( z_{1}\not\in G(w_{1}) \), contrary to the assumption.

By the parent construction, \( w_{2}\le w_{1}-\Tl \), and
 \( \Age(z_{1},w_{1})\le\Age(z_{2},w_{2})+1 \). 
The value \( w_{i} \) decreases at least \( \Tl/2 \)
in every step, while the age of \( \pair{z_{i}}{w_{i}} \) can decrease by at most
1.
Since it can only end in a birth or a germ cell in \( G(t) \), this
and~\ref{i:rn-gap.time-proj-bd} proves the age bound.%
\end{pproof} 

\begin{step+}{step:rn-gap.finish}
Let us lower-bound the growth of \( G(t) \) for \( t\in\clint{v_{0}}{v_{n}} \).  
\end{step+}
\begin{prooof}
Since we applied the estimate of part~\ref{step:rn-gap.no-dam} twice and the estimate
of part~\ref{step:rn-gap.dam} once, we get
 \begin{align*}
   (r(v_{n})-r(v_{0}))/\B  &\ge \frac{|v_{n}-v_{0}-\purgeT|^{+}}{3\tau_{2}}-2\healSpan-\islSz-11,
 \end{align*}
 and
\begin{align*}
    (|G(v_{n})|&-|G(v_{0})|)/\B
\\                       &\ge \frac{|v_{n}-v_{0}-\purgeT|^{+}}{3\tau_{3}}
     - 4\healSpan -2\islSz-21.
 \end{align*}
 Hence using~\eqref{eq:gap-lb}, \( |G(t)| \) never decreases below
 \( |G(v_{0})|-\gapLb\B >0 \).
  \end{prooof}
\end{Proof} 

\subsection{Non-damage gaps are large}\label{sec:gap.bad}

We now show that if healing does not succeed, then a large gap develops.
The Bad Gap Opening Lemma says that if a left exposed edge persists too
long (while possibly moving) then a right bad gap develops on its left.
This lemma will be used to prove the Bad Gap Inference Lemma,
saying that, under certain conditions, a left exposed edge
automatically has a right bad gap next to it.
Recall the definition of \( \healSpan \) in~\eqref{eq:heal-span} and \( \splitN \), \( \splitT \)
in~(\ref{eq:split-n}-\ref{eq:split-t}).

\begin{lemma}[Healing Cause]\label{lem:healing-cause}
  Let \( c_{0}=\pair{x_{0}}{t_{0}} \) be a cell with \( \Healing(c_{0})>0 \),
 \begin{alignat*}{3}
      D &= x_{0}+\lint{-(\healSpan+0.1)\B }{(\healSpan+1.1) \B },
\\   I  &= \rint{t_{0}-2\healSpan\Tu}{t_{0}}.
 \end{alignat*}
Assume that \( D\times I \) is damage-free.
Then there is a backward path with at most \( \healSpan-1 \) horizontal links,
at most 1 parental link, and at most \( 2\healSpan \) links altogether,
connecting \( c_{0} \) to an exposed cell.

A similar statement holds if \( \Frozen(c_{0})\ne 0 \), except that here the
length of the path is bounded by \( 2\numDirAff \), with \( \numDirAff \)
defined in~\eqref{eq:numDirAff}.
\end{lemma}
\begin{proof}
Let us build a path \( P=\tup{c_{0},c_{1},\dots,c_{n}} \) 
with \( c_{i}=\pair{x_{i}}{t_{i}} \) backwards with the property
that each of its cells \( c_{i} \) with \( i<n \) has \( \Healing>0 \), and \( c_{n} \) is exposed.
According to Condition~\ref{cond:healing}, \( \HealRevive \) is the only 
rule changing the value of field \( \Healing \).
It sets the value to 0 immediately unless \( x_{i} \) is exposed or has a neighbor
\( c_{i+1} \) with a higher \( \Healing \) value.
Therefore if \( \Healing(c_{i})>0 \) and it was not set to 0 in the last step,
there are the following possibilities:
\begin{bullets}
  \item There is an exposed parent.
  \item The cell itself was exposed at the observation time.
  \item There was a neighbor \( y \) with a higher \( \Healing \) value at the last observation
    time \( t'<t_{i} \) of \( x_{i} \).
\end{bullets}
In the first two cases, the path is finished with a last link.
In the last case, we add two links: a vertical link to \( c_{i+1}=\pair{x_{i}}{t'} \)
followed by a horizontal link to \( c_{i+2}=\pair{y}{t'} \).
Since this leads to an increase of \( \Healing \) value, the path can contain at
most \( \healSpan \) pairs of such links.

The statement on the field \( \Frozen \) is proved the same way, using
Condition~\ref{cond:freeze} and rule \( \Freeze \).
\end{proof}

\begin{definition}[Boundary path]\label{def:boundary-path}
Let us call the sibling relation between two cells
\df{strong}\index{sibling!strong},
if either they have been siblings for at least \( 2\Tu \) time units
or one cell is a parent of the other.

A cell is called for example a \df{weak left exposed edge}\index{edge!exposed!weak}
if it is either a left exposed edge or
would be one in case it had no left sibling, and has no strong left sibling.
A sequence \( R_{1}=\pair{y_{0}}{v_{0}},\dots, \pair{y_{m}}{v_{m}} \) of cell-time pairs
with \( v_{0}<v_{1}<\dots<v_{m} \) will be called a \df{left boundary path}\index{path!boundary}
if \( y_{i} \) is a weak left exposed edge during \( \opint{v_{i}}{v_{i+1}} \) and one of the following is true:
      \begin{djenum}
        \item Cell \( y_{i-1}=y_{i}-\B  \) dies at time \( t_{i} \) and is a strong sibling of
\( y_{i} \) at time \( t_{i}- \).
        \item Cell \( y_{i}=y_{i-1}-\B  \) is created at time \( t_{i} \) by parent \( y_{i-1} \).
      \end{djenum}
\end{definition}

 \begin{lemma}[Bad Gap Opening]\label{lem:bad-gap-opening}
 \index{lemma@Lemma!Bad Gap Opening}
 For \( m=\splitN \) let \( R_{1}= \) \( \pair{y_{0}}{v_{0}} \), \( \dots \),
 \( \pair{y_{m}}{v_{m}} \) be a left boundary
path that does not leave the colony of \( \pair{y_{m}}{v_{m}} \) on the left
(it may leave on the right).
Assume that the the segment
 \begin{align*}
 \lint{\textstyle\bigwedge_{i} y_{i} - 1.1\B }{\textstyle\bigvee_{i} y_{i}+1.1\B }
  \times\rint{v_{0}-\Tu}{v_{m}} 
 \end{align*}
is noise-free.
Then \( y_{i} \) has at most \( \p_{2}+\p_{0} \) full dwell periods during
\( \rint{v_{i}}{v_{i+1}} \), and there is a gap on the
left of \( R_{1}(v_{m}-) \) that contains a right bad gap.
The same statement holds if we interchange left and right.
  \end{lemma}
  
\begin{proof}
  Let us show that \( \pair{y_{i}}{t} \) is a left exposed edge during
\( \opint{v_{i}+2\Tu}{v_{i+1}-2\Tu} \).
Indeed, \( y_{i} \) can be a weak left exposed edge that is not a left exposed
edge only if it has a left sibling that is not strong.
Now, if \( y_{i} \) does not have a left sibling and it gets one then it
follows from Lemma~\ref{lem:glue} (Glue) (and the exclusion of cut, since the
it would not expose the edge) that the only way to lose this sibling is if
the sibling dies again, which it cannot do before making at least \( \p_{0} \)
switches, becoming a strong sibling in the meantime.
From this, it is easy to see that \( y_{i} \) can have a left sibling
only during a time interval adjacent to either \( v_{i} \) or \( v_{i+1} \).
Since the sibling stays weak these time intervals must be at most
\( 2\Tu \) long.%

Now, as soon as \( y_{i} \) is a left exposed edge, unless the healing rule creates
a left sibling, it gets \( \Decaying\gets 1 \) in \( \p_{2} \) dwell periods, and then the
decay rule kills it in further \( \p_{0} \) dwell periods.
Hence within \( \p_{2}+\p_{0} \) dwell periods, the
boundary paths moves left or right.
Now a reasoning similar to part~\eqref{eq:rn-gap.no-dam.progress.spread} of the
proof of Lemma~\ref{lem:rn-gap} (Running Gap) shows that the left movement is
associated with a healing step of rule \( \HealRevive \), and once the right
movement starts it will not stop before killing at least one cell with
\( \fld{Healing}=0 \).
Thus within each \( 2\healSpan+1 \) moves, the edge moves right.
In \( \gapLb+1 \) repetitions of this, the weak left exposed edge moves right
by this many steps.
During this time, by~\eqref{eq:p_{2}}, 
at most one growth step can occur on the left, so
a gap of width \( \gapLb\B  \) will be created which is internal by definition, since we assumed that
the boundary path does not leave the colony of \( \pair{y_{m}}{v_{m}} \) on the left.

To show that it is a bad gap, we bound the age of its germ cells.
The total number of (full or partial) dwell periods along the boundary path is
at most \( m(\p_{2}+\p_{0}+1) \).
During this time, the age of any germ cell created in place of the killed cells
can grow by at most \(  \lambda m(\p_{2}+\p_{0}+1)/\p_{1} \).
By the definition of \( \p_{i} \) in (\ref{eq:p_{0}}-\ref{eq:p_{2}}), we have
\begin{align*}
  (\p_{2}+\p_{0}+1)/\p_{1} \le\frac{2\healSpan\lambda\p_{0}+\p_{0}+1}{4\lambda\p_{0}}
  \le \healSpan/2 + 1,
\end{align*}
so the total germ age is at most \( \lambda \splitN(\healSpan/2+1) \), making the gap
by Definition~\ref{def:gap} a bad gap.
  \end{proof} 

 \begin{lemma}[Bad Gap Inference]\label{lem:bad-gap-infr}
 \index{lemma@Lemma!Bad Gap Inference}
  Let \( c_{0}=\pair{x_{0}}{t_{0}} \) be a left exposed edge,
 \begin{align*}
    D &= \lint{x_{0}-(\splitN+1.1)\B }{x_{0}+(\healSpan+1.1) \B },
\\  I &= \rint{t_{0}-\splitT}{t_{0}}.
 \end{align*}
Assume \( D\times I \) is damage-free.
  Then one of the following holds:
  \begin{djenum}

  \item\label{i:bad-gap-infr.gap}
 There is a bad gap on the left of \( c_{0} \);

\item\label{i:bad-gap-infr.exposing}
  There is a boundary path of length \( <\splitN \) leading backwards in time from \( c_{0} \)
to a cell undergoing one of the changes listed in Lemma~\ref{lem:exposing}
(Exposing).

   \end{djenum}
  The same statement holds if we replace left with right.
  \end{lemma}
  \begin{proof}
    We will distinguish a special case when \( c_{0} \) is involved in the leftward rule
    \rul{End-heal} (Algorithm~\ref{alg:End-heal})
    case of healing, calling this the (leftward) \df{end-heal} case.
    
    Let us construct a backward path \( \pair{x_{i}}{t_{i}} \), \( i=0,1,\dots,n \)
    of elements of the colony of \( c_{0} \), made up of
horizontal and vertical links such that cell \( x_{i} \) is a weak left
exposed edge during every nonempty time interval \( \opint{t_{i}}{t_{i-1}} \).
Suppose that \( \pair{x_{i}}{t_{i}} \) has already been constructed.
  \begin{romanenum}

   \item
  If \( x_{i} \) has a switching time \( t' \) immediately before \( t_{i} \) such
that \( \pair{x_{i}}{t} \) is a weak left exposed edge during \( \clint{t'}{t_{i}} \) then
let \( \pair{x_{i+1}}{t_{i+1}}=\pair{x_{i}}{t'} \).

   \item
  Otherwise, let \( t' \) be the lower bound of times \( t \) such that
\( \pair{x_{i}}{t} \) is a weak left exposed edge.
  If \( t'<t_{i} \) then let \( \pair{x_{i+1}}{t_{i+1}}=\pair{x_{i}}{t'} \).

   \item
If \( t'=t_{i} \) then at time \( t_{i} \), cell \( x_{i} \) became a weak left exposed edge without
switching, which can happen in the following ways::

\begin{varenum}{x}

\item Losing protected status via one of the cases of
  Lemma~\ref{lem:exposing} (Exposing), that is case~\ref{i:bad-gap-infr.exposing}
  of the present lemma.
  This cannot happen in the end-heal case.

\item By animation via a right parent, which by
  Condition~\ref{cond:anim,kill}, can only happen via the healing rule.
Let \( \pair{x_{i+1}}{t_{i+1}}=\pair{x_{i}+\B }{t_{i}} \).
Such an \( i \) will be called a \df{right jump} (backwards in time).

\item Losing a strong left sibling by one of the ways listed in Lemma~\ref{lem:glue} (Glue).
It is easy to check that of these, only the death of a strong left
sibling produces an exposed edge in \( x_{i} \).
Each of the rules \( \Heal \), \( \Decay \), \( \Purge \) and \( \Adapt \) that could have
killed \( x_{i}-\B  \) presupposes that this cell did not have a left sibling at
the observation time \( t'' \) of \( x_{i}-\B  \); thus, it did not have a strong left
sibling at time \( t_{i}- \).
On the other hand, as a strong left sibling of \( \pair{x_{i}}{t_{i}-} \), it has
been alive for at least \( 2\Tu \) time units.
Let \( \pair{x_{i+1}}{t_{i+1}}=\pair{x_{i}-\B }{t_{i}} \).
Such an \( i \) will be called a \df{left jump} (backwards in time).

\item\label{i:end-heal} The strong left sibling was killed by the right
  end-healing of member cells of another colony on the left (via \( \Adapt \)).
  This clearly cannot happen in the end-heal case of the present proof.
  We will see, moreover, that it does not happen at all.
  \end{varenum}
 \end{romanenum}

  The construction can stop in the following ways:
  \begin{varenum}{s}
   \item The number of jumps on the path reaches \( \splitN \).
   \item\label{i:bad-gap-infr.neighb-col}
Cell \( \pair{x_{n}}{t_{n}} \) would be a weak left exposed edge belonging to the left neighbor
colony.
This cannot happen in the end-heal case.
  \end{varenum}

By the construction, each jump is surrounded by vertical links.
For each left jump \( i \), cells \( x_{i} \) and \( x_{i+1} \) are strong siblings at
time \( t_{i}- \).
Let \( \pair{x_{i_{k}}}{t_{i_{k}}} \) for \( i_{1}<\dots<i_{m-1} \) be the points
on the backward path with the property that to
\( (i_{k},i_{k}+1) \) belongs a horizontal link.
Let us number these points forwards in time:
\( \pair{y_{0}}{v_{0}}=\pair{x_{n}}{t_{n}} \),
\( \pair{y_{1}}{v_{1}}=\pair{x_{i_{m-1}}}{t_{i_{m-1}}} \), 
\( \pair{y_{2}}{v_{2}}=\pair{x_{i_{m-2}}}{t_{i_{m-2}}} \), \( \dots \), 
\( \pair{y_{m-1}}{v_{m-1}}=\pair{x_{i_{1}}}{t_{i_{1}}} \), 
\( \pair{y_{m}}{v_{m}}=\pair{x_{0}}{t_{0}} \), 
creating a left boundary path as in Definition~\ref{def:boundary-path}.
If \( m\ge \splitN \) while we do not have case~\ref{i:bad-gap-infr.neighb-col}
then Lemma~\ref{lem:bad-gap-opening} (Bad
Gap Opening) implies a bad gap on the left of \( c_{0} \).
\begin{sloppypar}
In case~\ref{i:bad-gap-infr.neighb-col}, the death of \( \pair{x_{n}}{t_{n}} \) 
creates the weak left exposed cell \( \pair{x_{n-1}}{t_{n-1}} \), a left colony endcell.
By Definition~\ref{def:edges}, this
is possible only by the shrinking process, as described after
Algorithm~\ref{alg:workperiod}.
(It is clear from there that \( \pair{x_{n-1}}{t_{n-1}} \) is not a member cell.)
Continue the construction of the left boundary path.
Now it can be stopped in fewer than \( \splitN \) steps only by
one of the changes listed in Lemma~\ref{lem:exposing},
since \( \splitN<\Q  \); now Lemma~\ref{lem:bad-gap-opening} is applicable.
  \end{sloppypar}

Let us now show that case~\ref{i:end-heal} does not occur at all.
Assume it does and let \( y \) be the right-exposed cell trying to do the end-healing.
Then we can apply the present lemma to \( y \).
Since it is the (rightward) end-heal case, a bad gap is implied on the right of \( y \),
so \( y \) is far on the left of \( x_{i} \) and cannot possibly kill its left sibling.

 \end{proof}

\section{Attribution, progress}\label{sec:attrib}

In this section Lemma~\ref{lem:attrib} (Attribution) shows
that a non-germ cell implies a full colony nearby in the near past related
to it in one of several identifiable ways.
The main tool of the proof is a trace-back path.
The following lemma treats various possibilities for the path.
Recall the definition of $\destrT$ in~(\ref{eq:destr-time}) and
the definition \( \nLocMaint \) for the number of update ages of
locally maintained fields in~\eqref{eq:n-loc-maint} and denote
\begin{align}\label{eq:attr-time}
\attrT=\destrT + 9\tau_{1}.
\end{align}
For the lemmas below, for a time \( t_{0} \) let
\begin{equation}
  \label{eq:J-K}
\begin{aligned}
     J &= t_{0}+\rint{-\attrT\Q}{0},
  \\ K &=t_{0}+\rint{-\attrT\Q}{-\destrT\Q}.
\end{aligned}
\end{equation}

\begin{lemma}[Cover]\label{lem:cover}
Let $c_{0}=\pair{x_{0}}{t_{0}}$ be a live cell belonging to colony $\cC$ with
starting cell $z$,  
\begin{align*}
     I &=\lint{z-(\Q +1.1)\B }{z+(2\Q +1.1)\B }.
\end{align*}
For \( J \) as in~\eqref{eq:J-K} let
\begin{align*}
n_{J}&=2|J|/p_{1}\Tl+3\numDirAff 
\end{align*}
 and assume that \( \Damage^{*} \) does not intersect \( I\times J \).
Consider a trace-back path $P_{1}$ in $I\times J$ leading backwards from $c_{0}$
to a cell $c_{1}=\pair{x_{1}}{t_{1}}$.
The following holds.
\begin{cjenum}

\item\label{i:cover.stays}
  If $c_{0}$ is not an outer or germ cell younger than
$\growEnd+n_{J}$ then the path $P_{1}$ does not leave its own colony.
  
\item\label{i:cover.xposed}
Assume $t_{1}\le t_{0}-\destrT\Q$, the rectangle
 \begin{align*}
    I \times\rint{t_{1}-\splitT}{t_{1}+\splitT}
 \end{align*}
 is damage-free, \( c_{1} \) is not a doomed member cell with  \( |\Age(c_{1}) > \U-2\Q-\splitN \).
Then the maximal domain containing the cell $c_{1}$  contains the colony of \( c_{0} \) as
 well as its originating colony.
\end{cjenum}
\end{lemma}
We call this the Cover Lemma since part~\ref{i:cover.xposed} implies that the originating colony is
covered by a domain.
\begin{proof}
  \begin{enumerate}
  \item\label{i:cover.age-incr}
    The age increase along the path is at most \( n_{J}\).
Indeed, for the path length \( n \) we have
$n\le 2|J|/\Tl+2\numDirAff+2$, from part \ref{i:trace-back-bounds.time}
of Lemma~\ref{lem:trace-back-bounds} (where \( \numDirAff \) was defined in~\eqref{eq:numDirAff}).
The age increase along the path can be upper-bounded 
using Lemma~\ref{lem:trace-back-bounds} as
 \begin{align*}
     n/\p_{1} +2\numDirAff+3 &\le 2|J|/\p_{1}\Tl +2\numDirAff(1+1/\p_{1})+3+2/p_{1}.
 \end{align*}
The sum all terms but the first can be bounded by \( 3\numDirAff \).

\item Let us prove~\ref{i:cover.stays}.
According to Lemma~\ref{lem:skirting} (Skirting), the horizontal links on the path
can only cross a colony boundary in the inward direction, and only if
$\Age<\growEnd+\numDirAff$ on all cells involved. 
A parental link could be created only by the healing rule or one of the
expansion rules.
The healing rule does not act across colony boundaries, so
the later end of the
link along which the path is leaving would need to be an expansion or germ
cell. 
The expansion rules do not work for cells past age $\growEnd$,
so in each case, 
the cell at the crossing would need to have an age $\le\growEnd$.
By part~\ref{i:cover.age-incr} above, the age could not grow beyond
\( \growEnd+2|J|/\p_{1}\Tl+3\numDirAff \), so it does not reach the lower
bound required by~\ref{i:cover.stays}.

\item Let us prove~\ref{i:cover.xposed}.
  Assume that the maximal domain contains an exposed edge \( e \).
Lemma~\ref{lem:bad-gap-infr} (Bad Gap Inference)
implies that either there is a bad gap next to $e$ or
we have case~\ref{i:bad-gap-infr.exposing} of that lemma. 

In case of a bad gap,
we can apply Corollary~\ref{crl:rn-gap} (Running Gap) to the assumption 
$t_{1}\le t_{0}-\destrT\Q$, where the role of the
path \( P_{1} \) there is played by path \( P_{1} \) here.
It implies that the bad gap cannot be inside the originating colony of \( c_{0} \),
nor can it be in the same colony (as the gap path cannot cross \( P_{1} \)).

In case~\ref{i:bad-gap-infr.exposing} of Lemma~\ref{lem:bad-gap-infr}
there is a path of length $<\splitN$ leading back to an edge that becomes exposed
by one of the cases listed in Lemma~\ref{lem:exposing} (Exposing).
If it is one the cases different from~\ref{i:lem.exposing.doom} then
it is outside the originating colony and pointing away from it.
Case~\ref{i:lem.exposing.doom} does not occur; indeed, it could only occur to member
cells \( c_{1} \), but the condition
\( |\Age(c_{1}) \amod \U|>2\Q+\splitN \) implies \( |\Age(e) \amod\U|>\splitN \).
  \end{enumerate}
\end{proof}

Here is the sense in which we will be attributing a cell \( c_{0}=(x_{0},t_{0}) \)
to some colony.
Recall the definition of \( J,K \) in~(\ref{eq:J-K}).

\begin{definition}[Attribution]\label{def:attrib}
We  say that \( c_{0} \) is \df{attributed}\index{attribution} to colony $\cC$ if there is a path
$P_{1}$ going back to time $t_{0}-\attrT\Q$, and a union $E_{0}$ of 
intervals of total length $< |K|/2$ such that $P_{1}(t)\in \cC$ for $t$ in \( K \),
 and $\cC$ is covered by a domain for all times in $K\setminus E_{0}$.
\end{definition}

Suppose that another cell $\pair{x_{1}}{t_{0}}$ is attributed to some colony $\cC_{1}$;
Then even the union of their exception sets \( E_{0}\cup E_{1} \) has a total length \( <|K| \),
and therefore there is some time in \( K \) where both colonies \( C_{0} \) and \( C_{1} \)
are covered by a domain.
It follows that they are either disjoint or identical.
The following lemma again uses the intervals \( J,K \) in~\eqref{eq:J-K}.

\begin{lemma}[Attribution]\label{lem:attrib}
Assume that the live cell $c_{0}=\pair{x_{0}}{t_{0}}$ is traceable, and is not a
germ cell younger than $\growEnd+2|J|/p_{1}\Tl$.  
Then we have:
  \begin{cjenum}

  \item\label{lem:attrib.expan}
  $c_{0}$ can be attributed to its originating colony.

\item\label{lem:attrib.no-expan}
  If it is a member cell, or an outer or germ cell older than $\growEnd + |J|/p_{1}\Tl$
  then it can also be attributed to its own colony.
   \end{cjenum}
  \end{lemma}
  \begin{proof}
Build a trace-back path $P_{1}$ from \( c_{0} \) backwards in time.
Due to the absence of  $\Damage^{*}$, at most one island occurs in the area of \( P_{1} \).
Lemma~\ref{lem:cover} (Cover) implies that unless
$c_{0}$ is an outer or germ cell younger than $\growEnd+2|J|/p_{1}\Tl$, we can continue
the path $P_{1}$ all the way to time $t_{0}-\attrT\Q$ and stay in the colony
of $c_{0}$.
Let \( E \) be the set of times \( t \) of \( K \) that are either within distance \( \splitT \)
of the damage or such that \( |\Age(P_{1}(t)) \amod \U|\le 2\Q+\splitN \).
Then \( |E|\le 2\splitT + (4\Q+\splitN)\tau_{1}<|K|/2 \).
Lemma~\ref{lem:cover} implies that for all times \( t \) in \( K\setminus E \)
the cell \( P_{1}(t) \) is contained by a domain covering its originating colony.

In case~\ref{lem:attrib.no-expan}, the domain must contain the colony of \( c_{0} \),
because being past \( \growEnd \) otherwise it would have an exposed edge.
So \( c_{0} \) can be attributed to its own colony.
 \end{proof}

We will estimate the amount of progress in a given time
made by a colony in the absence of higher-order damage.

\begin{lemma}[Colony Trace-back]\label{lem:colony-trace-back}
Consider a time interval $\rint{t_{0}}{t_{1}}$.
Suppose the following, with \( I=z_{0}+\lint{-1.1\Q\B }{2.1\Q\B } \):
\begin{enumerate}
\item The rectangle \( I\times \rint{t_{0}-\splitT}{t_{0}} \) is damage-free,
  and \( \Damage^{*} \) does not intersect \( I\times\rint{t_{0}-\splitT}{t_{1}} \).
\item At time $t_{1}$, colony $\cC$, with starting cell $z_{0}$,
is covered by a maximal domain with member cells in $\cC$, and
no exposed edges.

\end{enumerate}
Assume that some trace-back path \( P \)
starts within $0.2\Q \B $ of the center at time $t_{1}$,
goes back to time $t_{0}$, and does not
come within \( 2 \Q  \) steps of a work period boundary time.
Then:
\begin{cjenum}

  \item\label{i:lem.colony-trace-back.path}
Path \( P \) moves at most to a distance $2\lambda\healSpan\splitN\B$ from where it started.
If there is no damage in \( I\times\rint{t_{0}}{t_{1}} \) then \( P \) does not move
horizontally at all.

  \item\label{i:lem.colony-trace-back.cover}
The colony $\cC$ is covered by a maximal domain without exposed edges
also at time $t_{0}$.

  \item\label{i:lem.colony-trace-back.progress}
    The age progress from $t_{0}$ to $t_{1}$ along path $P$ is at
most $\frac{t_{1}-t_{0}}{\p_{1}\Tl}+\numDirAff + 3$,
and at  least $\frac{t_{1}-t_{0}}{4\p_{1}\Tu}-6\Q$.
\end{cjenum}

\end{lemma}

The following corollary is just the inversion of part~\ref{i:lem.colony-trace-back.progress}
of the above lemma.

\begin{corollary}[Progress]\label{crl:progress}
  Under the assumptions of Lemma~\ref{lem:colony-trace-back}, suppose that
  the age progress along the path is \( n \).
  Then
  \begin{align*}
   (n-\numDirAff -3)\p_{1}\Tl \le t_{1}-t_{0} \le 4(n+6\Q)\p_{1}\Tu.
\end{align*}
\end{corollary}

\begin{Proof}[Proof of Lemma~\protect\ref{lem:colony-trace-back}]
  Let $P=(c_{1},c_{2},\dots)$ be the trace-back path mentioned in the lemma,
  listed backwards  in time, with $c_{i}=\pair{x_{i}}{u_{i}}$.
In the given time interval, at most one island intersects colony $\cC$.
Let \( K \) be the interval of times within \( \splitT \) of 
this island, then \( |K|\le 2\splitT+\islSz\Tu \).

\begin{step+}{step:colony-trace-back.progress.protected}
For times outside \( K \), 
the maximal domain containing the cell \( P(t) \) has only protected edges.
\end{step+}
\begin{pproof}
If there was an exposed edge, then Lemma~\ref{lem:bad-gap-infr} 
(Bad Gap Inference), combined with Lemma~\ref{lem:rn-gap} (Running Gap) would
imply that the colony \( \cC \) would not be covered by a domain at time \( t_{1} \).
\end{pproof} 

\begin{step+}{step:cover-trace-back.path}
Claims~\ref{i:lem.colony-trace-back.path} 
and~\ref{i:lem.colony-trace-back.cover} are true.
\end{step+}
\begin{pproof}
The island may cause at most $\numDirAff$ horizontal links
(with \( \numDirAff \) defined in~\eqref{eq:numDirAff}), all other links are
vertical or parental.
As shown in part~\ref{step:colony-trace-back.progress.protected}, the
path can have parental links only within distance $\splitT$ in time of the island.
Since they use $\Adapt(\cdot)$ with a delay of at least $\p_{1}$,
parental links are separated from each other by at least $\p_{1}$ vertical
links.
If follows then, using the definition of \( \splitT \) and \( \splitN \)
in~(\ref{eq:split-n}-\ref{eq:split-t}), that the path can move at most
\(    \numDirAff+|K|/\p_{1}\Tl \) steps away, where
 \begin{align*}
   |K|/\p_{1}\Tl=(2\splitT+\islSz)\tau_{2}\splitN/\p_{1}\Tl \le 3\lambda\healSpan\splitN/2
\end{align*}
(and will therefore stay in the colony).
Assume now that, contrary to claim~\ref{i:lem.colony-trace-back.path},
there is a first cell $c_{i}$ with a parental
link to $c_{i-1}$, farther than $\splitT$ in time from the island.
Then $x_{i}$ is at distance at most $\lambda\healSpan\splitN\B$ from $x_{1}$, and
as by assumption does not come within \( 2 \Q  \) steps of a work period boundary time,
is an exposed member cell, or growth cell with age $\ge\growEnd$, inside the the
colony. 
Lemma~\ref{lem:bad-gap-infr} (Bad Gap Inference)
implies then a bad gap which would not close until
time $t_{1}$, contrary to the assumptions.
Reference to the same lemma proves also the
statement~\ref{i:lem.colony-trace-back.cover}. 
\end{pproof} 

\begin{step+}{step:colony-trace-back.progress}
Claim~\ref{i:lem.colony-trace-back.progress} is true.
\end{step+}
\begin{pproof}
  The upper bound follows because according to
  Lemma~\ref{lem:skirting} (Skirting)
  there are at most \( \numDirAff+1 \) horizontal links in the path, while
  along the vertical part of the path, age increases at most by 1 every \( \p_{1} \) steps.
 We proceed to the lower bound.

\begin{step+}{step:colony-trace-back.progress.path}
We will build a new path, \( P'=\tup{c'_{0},\dots,c'_{m}} \) 
with \( c'_{i}=\pair{x'_{i}}{t'_{i}} \).
\end{step+}
\begin{prooofi}
Start with \( c'_{0}=c_{0} \).
Within the set \( K \) we continue  \( P' \) as a trace-back path;
consider the cases \( t'_{i}\not\in K \).
We have \( \Healing(x'_{i})=0 \) and \( \Frozen(x'_{i})=0 \)
for times near \( t'_{i} \), since otherwise
Lemma~\ref{lem:healing-cause} (Healing Cause)
would imply the nearby presence of exposed edges.  
In the absence of freezing, the age cannot decrease in time.
If \( \Age(x'_{i}) \) did not increase in time \( t'_{i} \),
then this can only have the following reasons:
\begin{djenum}
  \item\label{i:delay.wait}
 Cell \( x'_{i} \) is in the process of decreasing its \( \Wait.\March \) variable
    from \( \p_{1} \).
  \item\label{i:delay.parent} Cell \( c'_{i} \) has just been born.
  \item\label{i:delay.neighbor}
 There is a neighbor \( y \) at the last observation time \( t' \) of \( x'_{i} \) that does
 not let its age advance, by condition~\eqref{alg:March.sibling}
 of the \( \March \) rule~\ref{alg:March}.
 Thus it is either lower in age (call this a \df{delay link of the first kind})
 or is equal in age but then it is closer to the center of the originating colony (call it a
 \df{delay link of the second kind}).
\end{djenum}
In cases~\ref{i:delay.wait} and \ref{i:delay.parent}, let \( c'_{i+1} \) be
 found by regular trace-back: thus, the previous switching time of \( x'_{i} \) if
there is one, and the parent, if \( x'_{i} \) has just been born.
In case~\ref{i:delay.neighbor}, let \( c'_{i+1}=\pair{x'_{i}}{t'} \),
\( c_{i+2}=\pair{y}{t'} \).
\end{prooofi} 

We will now lower-bound age decrease and upper-bound age decrease
backwards on path \( P' \).

\begin{step+}{step:colony-trace-back.progress.inK}
Let us upper-bound age increase inside the set \( K \).
\end{step+}
\begin{prooofi}
Inside the time set \( K \), the only horizontal links
on the path \( P' \) come from Lemma~\ref{lem:skirting} (Skirting),
adding at most \( \numDirAff \) to the total age increase backward.
There is no increase on parental links, and at most one backward increase on each
vertical link (as \( \HealSync \) may decrease age).
This is a total increase of at most \( |K|/\Tl+\numDirAff \).
\end{prooofi} 

\begin{step+}{step:colony-trace-back.progress.nonK}
Let us lower-bound age decrease outside the set \( K \).
\end{step+}
\begin{prooofi}
Going to the parent via
alternative~\ref{i:delay.parent} above is always going towards the
originating colony:
since we are in a maximal domain without exposed edges, the only possible kind
of cell creation is growth.
The cells created by growth have age equal to the creator cell.
During a colony work period, there are at most \( 4\Q  \) growth steps (growing a
channel and then a growth arm to the left, and also to the right), this bounds the
number of parental links.
Let \( n_{1} \) be the number of delay links of the first kind and \( n_{2} \) that of
delay links of the second kind, \( n=n_{1}+n_{2} \).
We have
\(  n_{2}\le 2 \Q  + n_{1} \), as one can move towards the center
of the originating colony only after being moved away from it.
The total age decrease on these non-vertical links is hence at least
\(  n_{1}\ge n/2 -  \Q  \).
The total number of non-vertical links is \( m \le 4 \Q  + n \), so
\begin{align}\label{eq:n1-lb}
  n_{1}\ge (m-4 \Q )/2 - \Q  = m/2 - 3 \Q .
\end{align}
In a vertical segment of the path of time projection \( l \),
age decreases at least by 1 in every \( \p_{1} \) vertical steps.
So it decreases at least \( 0\lor(\frac{l}{\p_{1}\Tu}-1) \).
Let us call a vertical segment \df{short} if this number is \( <1 \),
that is whose length is \( \le 2\p_{1}\Tu \), else \df{long}.
On a long segment of length \( l \), the decrease is at least
\( \frac{l}{\p_{1}\Tu}-1 \ge \frac{l}{2\p_{1}\Tu} \).
Let \( L_{1} \) be the union of long vertical segments, and \( L_{2} \) that
of the short ones.
The time decrease on \( L_{1} \) is at least \( \ge |L_{1}|/2\p_{1}\Tu \).
The number of non-vertical links is at least \( m\ge |L_{2}|/2\p_{1}\Tu -1\),
and by~\eqref{eq:n1-lb} the age decrease on these links is at least
\begin{align*}
   m/2 - 3 \Q  \ge |L_{2}|/4\p_{1}\Tu -3 \Q  - 1.
 \end{align*}
So the total decrease of at least
\begin{align*}
  |L_{1}|/2\p_{1}\Tu + |L_{2}|/4\p_{1}\Tu -3 \Q  - 1
  \ge \frac{t_{1}-t_{0}-|K|}{4\p_{1}\Tu}- 3\Q  - 1.
\end{align*}
\end{prooofi} 

\begin{step+}{step:colony-trace-back.progress.sum}
Let us add up all the decreases inside and outside of \( K \).
\end{step+}
\begin{prooofi}
We get
 \begin{align*}
     &\frac{t_{1}-t_{0}-|K|}{4\p_{1}\Tu}-3\Q -1 - |K|/\Tl-\numDirAff =  \frac{t_{1}-t_{0}-|K|}{4\p_{1}\Tu}-V,
\\ V &\le 3\Q  +|K|/\Tl+\numDirAff+1 \le 4 \Q ,
 \end{align*}
where the last step follows from our assumptions about \( \Q  \).  
At time \( t_{0} \), there is no exposed edge between 
\( P(t_{0}) \) and \( P'(t_{0}) \), since this would imply a bad gap between these cells.
This could not close, and at least one of the paths would need to cross
it in order to meet, which is not possible.  
Therefore these two nodes are connected by a chain of siblings and so
their age differs by at most \( 2\Q  \).
Since \( P(t_{1})=P'(t_{1}) \), this shows that the age progress along \( P \) differs from that
of \( P' \) by at most \( 2\Q  \).
\end{prooofi} 
\end{pproof} 
\end{Proof}

\section{Healing}\label{sec:heal}

This section shows how the effect of an island will be healed.
By its definition, the healing process only tries to bridge gaps inside a colony
of member cells.
If a connecting arm gets damaged, it can just regrow from the healthy part.
If growth is thwarted by damage, and is still possible in the next work period,
then it will happen then.
Similarly, a germ killed by damage can regrow later.

The effect of an island can be delayed, in one case.  

  \begin{definition}[Successor]\label{def:successor}
Given a set \( E \) of cells at time \( t \), let us call a cell \( x' \) a
\df{successor}\index{successor} of \( E \) at time \( t'>t \) if \( \pair{x'}{t'} \) is
reachable by a forward trace-back path of cells 
from some \( \pair{x}{t} \) with \( x \) in \( E \).    
  \end{definition}

Recall the definition of 
siblings in Definition~\ref{def:siblings}, 
and the definition of 
(multi-) domains in Definition~\ref{def:multi-domain}.
Let
\begin{align}\label{eq:heal-t} 
            \healT        &= \purgeT+4\numDirAff \tau_{1},
\end{align}
 where \( \purgeT \) was defined in~\eqref{eq:purgeT}.
 \glo{heal-t@$\healT$}

  \begin{lemma}[Healing]\label{lem:heal}\index{lemma@Lemma!Healing}
Let \( \lint{a_{0}}{a_{1}}\times\rint{u_{0}}{u_{1}} \) be an island. 
Let \( \cC \) be a colony with starting cell \( z_{0} \), 
covered at time \( u_{0} \) by a 
domain of member cells, with no doomed cells of age \( >\U -3\Q  \).
Then a domain covers \( \cC \) at time \( w_{0}=u_{1}+\healT \).
  \end{lemma}

 \begin{Proof}
   The damage affects directly at most \( \numDirAff \) neighbor cells of \( \cC \)
   (with \( \numDirAff \) defined in~\eqref{eq:numDirAff}).
At time \( t \), let \( L(t) \) be the set of successors in \( \cC \) of 
cells of \( \cC \) on the left of the affected part, and \( R(t) \) the successors from
cells to the right.
(These sets are not necessarily disjoint.)
Both \( L(t) \) and \( R(t) \) are (possibly empty) domains; indeed,
the only rule that would change the domain property is the killing
of a doomed cell at age 0, but there is not enough time for this. 


\begin{step+}{step:heal.purge}
Let \( u_{2}=u_{1}+\purgeT \).
If the sets \( L(u_{2}) \) and \( R(u_{2}) \) are both nonempty
then there is no live cell between them at time \( u_{2} \)
other than possibly germ cells with age \( \le 1 \).
\end{step+}
\begin{pproof}
The only cells in question are the up to \( \numDirAff \) cells affected directly by the
damage, others could not have reborn yet after dying.
Consider the up to \( \numDirAff \) domains between \( L(u_{1}) \) and \( R(u_{1}) \);
they can shrink or merge during the time interval \( \rint{u_{1}}{u_{2}} \).
According to Lemma~\ref{lem:short-exposed}, each will still have an exposed side,
hence each will be killed by the time \( u_{2} \) by
the purge rule (while having no time yet to grow due to the delay in the
rule \( \Adapt \))---unless they are germs so new
that the purging did not have time to finish yet: in this case they have age \( \le 1 \).  
\end{pproof} 

\begin{step+}{step:heal.purge-end}
Any live cell at time \( u_{2} \) whose body intersects \( \cC \) and that
does not belong to \( L(u_{2})\cup R(u_{2}) \) cannot be a member cell.
\end{step+}
\begin{pproof}
  Assume by way of contradiction, that there is such a member cell \( x \).
  In this case clearly either \( L(u_{2}) \) or \( R(u_{2}) \) is empty; without
  loss of generality assume that \( L(u_{2}) \) is empty.
Consider a trace-back path from \( x \).
It cannot lead to \( \cC \) at time \( u_{0} \), since then according  to  part~\ref{step:heal.purge}
above, \( x \) would belong to either \( L(u_{2}) \) or \( R(u_{2}) \).
Therefore it must lead to a domain outside; given that \( R(u_{2}) \) is nonempty, this must be on the left.
Let \( y \) be the right end of this domain at time \( u_{0} \).
It must be outside \( \cC \), so it is to the left of \( x \) whose body, by assumption, intersects \( \cC \).
By Lemma~\ref{lem:skirting} (Skirting), the path has at most \( \numDirAff \) horizontal links,
and one parental link (there is no time for a complete animation delay),
so \( x \) and \( y \) are within distance \( (\numDirAff+1)\B  \) from each other.
All cells of the short chain from \( \pair{x}{u_{2}} \) to \( \pair{y}{u_{0}} \) 
have ages within a few steps of each other.
If \( \pair{y}{u_{0}} \) is a right colony end then \( x \)
(assumed to be a member cell to the left of \( R(u_{2}) \))
would have been purged away due to its small depth.
If it is not then the right edge of \( \pair{y}{u_{0}} \)  is exposed,
since it is either a member cell not at a colony end or a growth cell near the end
of the work period, and so during its shrinking phase.
But then Lemma~\ref{lem:bad-gap-infr} (Bad Gap Inference) 
would imply large gap on the right of \( \pair{y}{u_{0}} \), making the path from it
to \( \pair{x}{u_{2}} \) impossible.
(Indeed, due to the age of \( \pair{y}{u_{0}} \) the
conclusion~\ref{i:bad-gap-infr.exposing} of that lemma is not applicable.)
\end{pproof} 


\begin{step+}{step:heal.size}
If both \( L(t) \) and \( R(t) \) survive until time \( u_{2} \), then
the gap between them will have size \( \le \numDirAff \) during the whole time.
If \( L(t) \) disappears then the 
the gap from \( z_{0} \) to \( R(t) \) has size \( \le 2\numDirAff \) during the whole time.
\end{step+}
\begin{pproof}
If both \( L(t) \) and \( R(t) \) survive, then there were up to \( \numDirAff \) cells
affected by damage between them.
Only these cells can disappear.
If \( L(t) \) disappears, this is only since its size was \( \le\numDirAff \) to begin with.
\end{pproof} 

 \step{step:heal.no-gap}{If both \( L(t) \) and \( R(t) \) remain nonempty
then they will become part of a domain.}
 \begin{pproof}
As shown above, between times \( u_{2} \) and \( w_{0}  \), the gap between
\( L(t) \) and \( R(t) \) does not become larger than \( \numDirAff \).
Rule \( \HealRevive \) brings together \( L(t) \) and \( R(t) \) 
at some time \( t\le \numDirAff\tau_{1} \), at  some point \( x\in L(t) \), \( x+\B \in R(t) \).
The possible germ cells of age \( <2 \) in the gap are not an obstacle: they
are weaker than the cells that are being created, so they will be erased.
It can only happen that \( L(t) \) and \( R(t) \) are not joined if
\( |\Age(x)-\Age(x+\B )\amod \U |>1 \).
But as said in the remark at the end of Example~\ref{xpl:heal-sync},
then rule \( \HealSync \) will bring the two ages together
in time at most \( 2\numDirAff\tau_{0} \).  
This creates a single domain by time 
\( u_{2}+\numDirAff(\tau_{1}+2\tau_{0}) \). 
 \end{pproof} 

 \step{step:heal.fill}{If \( L(t) \) does not survive then \( z_{0}\in R(w_{0}) \).}
  \begin{pproof}
The left endcell of \( R(t) \) will be exposed to the left and therefore the
rule~\ref{alg:End-heal} (\( \rul{End-heal} \)) will apply to it.
Part~\ref{step:heal.purge-end} above shows that no member cell's body intersects the
colony \( \cC \).
Other possible cells that might intersect are weaker than the member cells that
end-healing creates,
and will therefore be killed off by rule \( \Adapt \) of Algorithm~\ref{alg:Adapt}.
The end healing thus proceeds in time \( 2\numDirAff\tau_{1} \).
 \end{pproof}
 \end{Proof} 

\section{Communication}\label{sec:commun}

Here, we will give the rules for information retrieval,
which eventually will be used in rules \( \Send \) and \( \Retrieve \) at the beginning of the
work period as in Algorithm~\ref{alg:workperiod}.
Compared to Section~\ref{sec:simp-sim}, there are two main differences:
\begin{Alphenum}
\item\label{i:not-adj}
  the neighbor colonies may not be adjacent;
\item\label{i:asynch} they may be in different stages of their work period.
\end{Alphenum}
To deal with issue~\ref{i:not-adj}, during the time of observation, a colony whose neighbor is
not adjacent extends an arm of channel cells.
By Definition~\ref{def:cell-kinds}, rightward channels would override leftward channels.
This way at the time of communication there will be only one channel arm between the
neighbor colonies, assuring that their distance is less than \( \Q\B \).

To deal with issue~\ref{i:asynch}, within the colony work period,
note that the \df{unsafe} time interval of the colony work period during which
information to be sent to neighbor colonies can change 
starts with the time when one of the cells is within \( 2 \Q  \)
of the beginning of the work period or the end 
of the \( \Compute \) rule.
The retrieved information from the neighbor colonies will carry the age of sending.
The information will be used only when it comes from both of these
colonies at a safe ages.
The safe interval is designed to be much longer than the unsafe one,
so a safe time will be found.
We define $\var{Mail-ind}$ as in Section~\ref{sec:simp-sim}, with some additions:
 \glo{mail-ind@$\var{Mail-ind}$}%
 \glo{k@$\ol k$}%
 \begin{align*}
\var{Mail-ind} =\setOf{(k,d)}{k\in\{-1,0,1\}, d\in\{-1,1\}}.
 \end{align*}
 Recall that we defined \( e_{-1}=0 \), \( e_{1}=\Q -1 \) in~(\ref{eq:e_j}).
The relation of cells and their originating colonies will define the
predicate \( \Edge_{j}(x) \).
 In words: \( \Edge_{j}(x)=0 \) if \( y \) is a neighbor in the same extended colony,
 1 if it is an adjacent neighbor in an adjacent colony, 1.5 if it is a neighbor in an (extended)
 neighbor colony, but only one of the two colonies at the junction can be an extended one.
Formally:
\[
  \Edge_{j}(x) = \begin{cases}
     0 &\text{if $\Addr(y)=\Addr(x)+j$,}
                   \\   1 &\text{if }y=\thg_{j}(x),\; \Addr(x)=e_{j},\; \Addr(y)=e_{-j},
                   \\ 1.5 &\text{if }\Edge_{j}(x)\not\in\{0,1\},\; y=\thg_{j}(x),
                   \\     & (\Addr(x)-e_{j})/j\ge 0,\; (\Addr(y)-e_{-j})/j\ge 0,
                    \\    & \text{but at most one of these is \( > 0 \)}.
\\ \infty &\text{otherwise.}
     \end{cases}
\]
The additional issues that copying has to deal with are the unsafe times, as mentioned above,
non-adjacent colonies,
and that the mail track may be narrower than the information in any one cell it needs to carry.
In particular, the \( \Payload \) field will consist of \df{slices} \( \Payload_{i} \) only one of which
fits into the mail field.
It will be convenient to have a separate mail field for every kind of copying from
some track \( \F \) to some track \( \G \), that is
\begin{align*}
   \Mail_{k,d,\F,\G}.
\end{align*}
However, while discussing the general rules covering all kind of mail,
we will suppress the indices \( \F,\G \).
We have three widths to consider: that of the sender, of the receiver, and of \( \Mail_{k,d}.\Info \).
The latter will always have the same width as the smaller of the other two,
and the bigger width will always be an integer multiple of the smaller one.
To deal with these issues, the field $\Mail_{k,d}$ for \( (k,d)\in\var{Mail-ind} \)
will have some additional sub-fields besides the sub-fields defined in 
Section~\ref{sec:simp-sim}.
The sub-field \( \Toaddr \) has already been mentioned,
but there will also be sub-field \( \Fromaddr \).
The sub-field 
\begin{align*}
 \Adj
\end{align*}
 will show whether the
neighbor the direction from which the mail is retrieved is adjacent.
Recall the definition of
\begin{align*}
   j(k,d),\peer(k,d)
\end{align*}
in Section~\ref{sec:simp-sim}:
$\Mail_{k,d}$ of a cell will read, from direction $j(k,d)$, the mail field with index \( \peer(k,d) \).
Recall the definition of \( \var{Mail-to-receive}(k,d) \) in~\eqref{eq:mail-to-receive}:
 \begin{align*}
   \var{Mail-to-receive}(k,d) = \Mail_{\peer(k,d)}^{j(k,d)}.   
 \end{align*}
   Now we modify this to account for missing neighbor colonies and to pass
the information on whether the mail is coming from an adjacent colony 
to the bit track \( \Mail_{k,d}.\Adj \) of the receiving colony.
\begin{itemize}
\item If  \( \Edge_{j(k,d)} = \infty \) then \( \var{Mail-to-receive}(k,d) = \# \)
(no non-adjacent neighbor colony, neighbor big cell will be considered vacant).
\item If  \( \Edge_{j(k,d)} = 1 \) then \( \var{Mail-to-receive}(k,d).\Adj = 1 \)
(adjacent neighbor colony).
\item If  \( \Edge_{j(k,d)} = 1.5 \) then \( \var{Mail-to-receive}(k,d).\Adj = 0 \)
(non-adjacent neighbor colony).
\end{itemize}

\begin{sloppypar}
To deal with asynchrony, we add a 
``hand-shaking'' condition  \( \MailFree_{k,d}(x) \):
 \glo{mail-used@$\MailFree$}%
it says that the cell \( x \) is free to rewrite $\Mail_{k,d}$ with a defined value
 if the new value is not destined for an earlier address, and either its current value is undefined or
 if it is for a later or equal address and has already been passed on.
 (This is sufficient: every piece of mail must be distinct due to (\( \Toaddr \), \( \Fromaddr \)).)
 Let \( j=j(k,d) \), where \( (k,d)=\peer(k,d') \) and note that necessarily \( j(k,d)=j(k,d') \).
 Then
 \begin{equation}\label{eq:mail-free}
   \begin{aligned}
   \MailFree_{k,d}\eqv& \sign(\Addr-\var{Mail-to-receive}(k,d).\Toaddr)\ne -j(k,d)
\\         &\land (\Mail_{k,d}\in\{\#,\Mail_{k,d'}^{-j}\} \lor \Mail_{k,d}.\Toaddr=\Addr).
   \end{aligned}
 \end{equation}
 Some situations in which mail is used will also need some error-checking, as
 its output cannot be stored in extra copies for voting.
 For this, we will check that mail has not gone
 beyond the cell with \( \Addr=\Toaddr \) (in the direction \( -j \) of mail movement),
 and that both \( \Toaddr \) and \( \Fromaddr \) are increasing by 0 or 1 from left to right,
 between a cell and the one to which mail would be forwarded. 
 Here we assume that if \( \Mail_{k,d}=\# \) then \( \Mail_{k,d}.\Toaddr \), \( \Mail_{k,d}.\Fromaddr \)
 and any expressions they enter into have the value \( \Undef \).
For \( r\in\{-1,1\} \), \( \MailLegal_{k,d}(r) \) checks whether
mail is legally formed between the cell and its neighbor in direction \( r \).
\begin{align*}
  \MailLegal_{k,d}(r)\eqv& \sign(\Addr-\Mail_{k,d}.\Toaddr)\ne r
  \\                         &\land \Mail^{r}_{k,d}.\Toaddr-\Mail_{k,d}.\Toaddr\in\{0,r,\Undef\}
\\                            &  \land \Mail^{r}_{k,d}.\Fromaddr-\Mail_{k,d}.\Fromaddr\in\{0,r,\Undef\}.
\end{align*}
Rule \( \MoveMail \) in Algorithm~\ref{alg:move-mail-2} will be more elaborate than its
example version in Algorithm~\ref{alg:move-mail-1} by checking these two conditions.
\end{sloppypar}
 
 \begin{algo}[H]\caption{sub-rule $\protect\MoveMail(k,d)$}\label{alg:move-mail-2}
     \algLet \( j\gets j(k,d) \)\;
     \lIf{\algNot \( \MailLegal_{k,d}(-j) \)}{\( \Mail_{k,d}\gets\# \)}
     \ElseIf{$\MailFree_{k,d}$ \algAnd \( \MailLegal_{k,d}(j) \)} 
     {$\Mail_{k,d}\gets \var{Mail-to-receive}(k,d)$}\;
   \end{algo}
   As \( \Mail_{k,d} \)
   has the additional indices \( \F,\G \), eventually we would write \( \MoveMail(k,d,\F,\G) \).
   The condition \( \MailLegal \) will not let pass any piece of mail whose \( \Toaddr \), \( \Fromaddr \)
is not gradually non-decreasing from left to right.
The replacement with \( \# \) 
eats up the illegal pieces and allows the rest behind it to move.
This way the damage can have two effects: one in the place where it occurred,
and another, in (or near) the \( \Toaddr \) where the mail is headed.
But it cannot affect several other destinations by substantially changing \( \Toaddr \).

\begin{sloppypar}
Mail will be posted by the rule \( \PostMail \)  and received by the rule \( \ReceiveMail \).
Both of these will be running for certain intervals of \( \Age \).
We start with defining \( \ReceiveMail(k,\F, \nS, a, b, n) \) in Algorithm~\ref{alg:receive-mail},
as it is somewhat simpler.   
Parameter \( k \) has the same role as \( k \) in \( \Mail_{k,d} \).
Parameter \( \F \) is the track \emph{into} which the mail is to be received.
Parameter \( \nS \) shows how many times the width of this field is a multiple of \( \Mail_{k,d}.\Info \);
if \( \nS>1 \) then \( \F_{i} \) are the slices with the same width as \( \Mail_{k,d}.\Info \).
If \( \nS=1 \) then let \( \F_{0}=\F \).
Parameters \( a, b \) are the starting addresses of the sender and receiver location, and \( n \) is the length
of the \emph{receiver} location.
When the width of the two locations is different then if their intervals could intersect then
this would complicate the posting and receiving rules. 
Therefore we will always arrange that \emph{when the width of the two locations is different
then their intervals are disjoint}.
\end{sloppypar}

 \begin{algo}[H] \caption{sub-rule $\protect\ReceiveMail(k,\F,\G, \nS, a, b, n)$}\label{alg:receive-mail}
   \leIf{\( k\ne 0 \)}{\algLet \( d\gets -1 \)}{\algLet \( d\gets \sign(b-a) \)}
   \algLet \( j\gets j(k,d) \), \( A \gets (\Mail_{k,d,\F,\G}^{j}.\Fromaddr-a) \)\;
   \If{\( \Addr\in \lint{b}{ b+n} \) \algAnd \( \Addr-b = \Toaddr  \)}{
   \( \G_{A\bmod \nS}\gets\Mail^{j}_{k,d,\F,\G}.\Info \)
 }
  \end{algo}

In \( \PostMail(k,\F,\G,\nS, a, b, n) \), Algorithm~\ref{alg:post-mail},
parameters \( k,a,b \) have the same role as in \( \ReceiveMail \).
Parameter \( \F \) is the track \emph{from} which the mail is to be sent: its width is
\( \nS \) times that of \( \Mail_{k,d}.\Info \), with slices \( \F_{i} \) as in \( \ReceiveMail \).
Parameter \( n \) is the length of the \emph{sender} location: the receiver location has length \( \nS n \).
It assumes that at its start the \( \Mail_{k,d} \) track on interval \( \lint{a}{a+n} \)
is covered by \( \# \).
\( \PostMail \) is somewhat more complex than \( \ReceiveMail \), as the order in which the slices
\( \F_{i} \) are posted depends on the direction in which mail will be sent.
For simplicity we first look at $\protect\PostMail(1,\F,\G,\nS, a, b. n)$, that is
the case for \( \Mail_{1,1} \): the slices of location \( \F(\lint{a}{a+n}) \) will be posted
starting with address \( a+n-1 \), going backwards towards address \( a \),
also the slices are posted one-by-one backwards starting from \( \F_{\nS-1} \).
The information will just get posted onto the \( \Mail_{1,1} \) track, relying
on \( \MoveMail \) to forward it to the right.
The posting in each cell begins only when \( \Mail_{1,1}=\# \), and except for the right
end-cell, the right neighbor must have already posted all its slices.
If a slice \( \F_{i} \) with \( i>0 \) is currently in \( \Mail_{1,1} \)
then slice \( \F_{i-1} \) is posted only if the current one 
has already been forwarded to the right neighbor.
When everything has been posted and forwarded from the leftmost cell then it posts \( \# \).
This will be forwarded by \( \MoveMail \), and eventually the whole interval \( \lint{a}{a+n} \)
will be covered by \( \# \); at this point the posting will start again from the right end.
A fault will not make the posting process stuck forever, as \( \MoveMail \) will keep moving whatever
it finds on the mail track.
For readability here as well in Algorithm~\ref{alg:post-mail} below we suppress the
indices \( \F,\G \) from \( \Mail_{k,d,\F,\G} \).

\begin{algo}
\caption{sub-rule $\protect\PostMail(1,\F, \G,\nS, a, b, n)$}\label{alg:post-mail-1}  
   \algLet \( \Addr'\gets\Addr-a \)\;
   \If{\( \Addr'\in \lint{0}{n} \)}{
      \If{\( \Mail_{1,1}=\# \)}{
        \nl\lIf{\( \Addr'=n-1 \) \algOr \( \Mail_{1,1}^{1}.\Toaddr- b = \nS(\Addr'+1) \)
          \label{alg:post-mail.first.x}
        }{\( \Mail_{1,1}.(\Info,\Fromaddr,\Toaddr)\gets(\F_{\nS-1},\Addr, b + \nS\cdot\Addr'+(\nS-1)) \)
        }
      }
      \nl\ElseIf{\( \Mail_{1,1}.\Fromaddr=\Addr \) \algAnd
        \( \Mail_{1,1}^{1} = \Mail_{1,1} \)\label{alg:post-mail.sent.x}}{
        \algLet \( i \gets (\Mail_{1,1}.\Toaddr - b)\bmod\nS \)\;          
         \lIf{\( i\ne 0\)}{
            \( \Mail_{1.1}.(\Info,\Fromaddr,\Toaddr)\gets(\F_{i-1},\Addr,\Mail_{1,1}.\Toaddr -1 ) \)
          }
          \lElseIf{ \( \Addr'= 0 \)}{\( \Mail_{1,1}\gets\# \)}
        }
      }
\end{algo}

Case~\eqref{alg:post-mail.first.x} is when the slice of the field \( \F \) to be posted first
can actually be posted. 
Condition~\eqref{alg:post-mail.sent.x} is testing whether
the posted value has been passed on (by \MoveMail): this is similar to part of the condition
\( \MailFree_{k,d} \) in~\eqref{eq:mail-free}.
For the general \( k,d \), recall the direction \( j(k,d) \) from which mail is received,
defined after~\eqref{eq:peer}.
 \glo{post-mail@$\PostMail$}%

\begin{algo} \caption{sub-rule $\protect\PostMail(k,\F, \G,\nS, a, b, n)$}\label{alg:post-mail}
   \leIf{\( k\ne 0 \)}{\algLet \( d\gets 1 \)}{\algLet \( d\gets \sign(b-a) \)}
   \algLet \( \Addr'\gets\Addr-a \), \( \delta_{j} \gets (j+1)/2 \)\;
    \If{\( \Addr'\in \lint{0}{n} \)}{
      \If{\( \Mail_{k,d}=\# \)}{
        \lIf{\( \Addr'=\delta_{-j}(n-1) \)
          \algOr \( \Mail_{k,d}^{-j}.\Toaddr- b = \nS(\Addr'-j)-\delta_{j}(\nS-1) \)
        }{\( \Mail_{k,d}.(\Info,\Fromaddr,\Toaddr)\gets
          (\F_{\delta_{-j}(\nS-1)},\Addr, b + \nS\cdot\Addr'+\delta_{-j}(\nS-1)) \)
        }
      }
      \ElseIf{\( \Mail_{k,d}.\Fromaddr=\Addr \) \algAnd \( \Mail_{k,d}^{-j} = \Mail_{k,d} \)}{
        \algLet \( i \gets (\Mail_{k,d}.\Toaddr - b)\bmod\nS \)\;
        \lIf{\( i\ne \delta_{j}(\nS-1)\)}{
          \( \Mail_{k,d}.(\Info,\Fromaddr,\Toaddr)\gets(\F_{i+j},\Addr, \Mail_{k,d}.\Toaddr +j ) \)
        }
        \lElseIf{\( \Addr' = \delta_{j}(n-1) \)}{\( \Mail_{k,d}\gets\# \)}
      }
    }
  \end{algo}

Let us estimate the time it takes the posted information to arrive, in a typical special case.

\begin{lemma}\label{lem:mail-speed}
  Consider a case of \( \PostMail(1,\F,\G,1,0,0,\Q) \) when a track \( \F \) is sent to a track \( \G \)
  of a right neighbor colony (not necessarily adjacent).
  Assume that between the sending time \( t_{0} \) and the receiving time,
  the space-time area of concern is damage-free.
  For a cell \( x \) let \( d(x) \) be the distance of \( x \) from the right end of the
  receiving range.
  Let \( x \) be a cell and \( y \) be a cell that is \( m \) steps
  (over neighbors) to the right of \( x \).
  Then at time \( t_{0}+2\Tu(m+2 d (x)) \) the information from \( x \) has reached \( y \).
\end{lemma}
\begin{proof}
  Induction on \( n= m+2 d(x) \).
  For \( i=0 \) the information did not have to travel any steps.
  So the statement holds for \( i=0 \) and hence also for \( n=0 \).
  Suppose that the statement holds for \( n-1 \), we will prove it for \( n \).

  Let \( y' \) be the left neighbor of \( y \), and let \( t'_{1} \) be its first switching time
  after \( t_{1}-2\Tu \).
  For \( y' \), distance \( m'=m-1 \), we have \( n'=m'+2 d(x) =n-1 \).
  By the inductive assumption on \( n \), information from \( x \)
  must have reached \( y' \) by the time \( t' _{1}\).
  It follows that at the decision time \( t_{1} \) of \( y \), it must have seen this information.

  If \( d(y)=0 \) then it has no right neighbor to consider and
  will receive the information of \( y' \) at time \( t_{1} \).
  Suppose that it has a right neighbor \( y'' \) to consider, and let \( t''_{1} \) be the
  first switching time of \( y'' \) after \( t_{1}-2\Tu \).
  Then \( y'' \) is at distance \( m''=m+1 \) to the right of \( x \).
  Let \( x'' \) be the right neighbor of \( x \), then \( d(x'') = d(x)-1 \),
  and for the pair \( x'',y'' \) we have \( n''=m''+2 d(x'') = n-1 \), therefore
  by time \( 2\Tu(n-1) = 2\Tu n -2\Tu \) the cell \( y'' \) has the information from \( x'' \).
  Cell \( y \) will see this and hence will be free to
  receive the information from \( y' \) at time \( t_{1} \) (if it still does not have it).
  \end{proof}

The copy command will now be defined with more care than
in the example simulation in Algorithm~\ref{alg:copy-1}.
Actually, we will only define it for copying within the colony, so from now on we will just write 
\begin{align*}
 \Copy(\loc_{1},\loc_{2}) .
\end{align*}
For retrieving information from the neighbor colonies, another mechanism will be used.
For the command
\begin{align*}
 \Copy(\F(\lint{a}{a+n_{1}}),\G(\lint{b,b+n_{2}}),
\end{align*}
we will need only the case when for some \( \bw=|\Mail_{k,d}.\Info| \)
and \( \nS\ge 1 \) we have either \( |\F|=\nS \bw \), \( |\G|=\bw \), \( n_{2}=\nS n_{1} \)
or \( |\G|=\nS \bw \), \( |\F|=\bw \), \( n_{1}=\nS n_{2} \).
In the implementation here,
we simply run the rules \( \PostMail \) and \( \ReceiveMail \) for a certain length of time.
If \( \nS\ne 1 \) then we require here the
two intervals \( \lint{a}{a+n_{1}} \), \( \lint{b}{b+n_2} \) to be disjoint.
In applications we will not insist on this condition because if
the intervals are not disjoint we could replace one copy operation with two
that satisfy the condition.

\begin{algo}\caption{sub-rule \(\Copy(\F(\lint{a}{a+n_{1}}),\G(\lint{b}{b+n_{2}}))\)}
  \label{alg:copy-2}
  \algLet \( d\gets \sign(b-a) \)\;
  \lIf{\( n_{1}>n_{2} \)}{\algLet \( \nS_{1}\gets n_{1}/n_{2},\; \nS_{2}\gets 1 \)}
  \lElse{\algLet \( \nS_{1}\gets 1,\;\nS_{2}\gets n_{2}/n_{1} \)}
  \( \Mail_{0,d,\F,\G}\gets\# \)\;
  \For{ \( \lambda\nS\Q \) steps of \( \Age \)}{
    \( \PostMail(0,\F,\G,\nS_{1},a,b,n_{1}) \)\; 
    \( \MoveMail(0,d,\F,\G) \)\; 
    \( \ReceiveMail(0,\F,\G,\nS_{2},a,b,n_{2})  \)
    }
\end{algo}
    
Two kinds of information will be retrieved before the computation.
The information that is found within the colony itself will be retrieved using the rule
\( \Copy \) of Algorithm~\ref{alg:copy-2}.
Retrieving from the neighbor colonies is more complex, as these
colonies may not be there or may be in different stages of their work period.
It will be handled in Section~\ref{sec:send-retrieve}.

 \section{Computation}\label{sec:comp}

Before defining the rules that carry out the computation of a colony, we
define some of the simple data structures that support it.

\subsection{Coding and decoding}\label{sec:code-decode}

Let us say more about the error-correcting code used.

\begin{definition}[Errors]\label{def:the-err-corr-code}
Our error-correcting code will distinguish, similarly to  
Examples~\ref{xmp:err-corr-code} and \ref{xmp:no-stretch}, information symbols
and error-check symbols.
A string $s$ of symbols that can be the argument of a
certain error-correcting decoding \df{has $d$ errors}\index{error} if,
when applying the decoding and encoding to it, the result differs from $s$
in $d$ symbols.
We will work with a code%
\glo{a.greek@$\alpha$}%
 \begin{equation}\label{eq:err-corr-code-alpha}
  \alpha
 \end{equation}
of  Example~\ref{xmp:RS-code} that can
correct errors in $t$ symbols with just $2t$  error-check symbols.
When $\loc$ is a location representing a string $s$ then we will use the
notation
 \[
   \alpha^{*}(\loc)=\alpha^{*}(s)
 \]
 for the value decoded from it.
 We expect having to deal with at most one  damage rectangle.
 It can corrupt a group of \( \numDirAff \)  
 consecutive cells (where \( \numDirAff \) was defined in~\eqref{eq:numDirAff}).
The symbols of the code will be chosen in such a way that a code symbol 
stretches over $\numDirAff$ cells.
A damage rectangle can damage therefore up to 2 symbols, for which 4 error check
symbols suffice.
  \end{definition}

In the subdivision of the state of our simulation cells into fields, we
develop here a more elaborate version of Example~\ref{xmp:no-stretch}.
Like there, we have a field \( \Work \), to be subdivided into a constant number of
fields like \( \Addr \), \( \Mail \), and so on.
In Definition~\ref{def:info-field}, we introduced a field called $\Info$
partitioned into sub-fields $\Payload$, $\PlRedun$ and \( \Util \).
We will build a uniform amplifier as in Lemma~\ref{lem:amp},
with an aggregated field \( \Payload \).
Different parts of the \( \Util \) track in the colony will
represent the fields \( \Work^{*} \), \( \PlRedun^{*} \), \( \Color^{*} \) and \( \Util^{*} \)
of the simulated big cell, and contain also the error check symbols for the \( \Util \)
track itself.
Let
\begin{align*}
   \Payload'=\Payload\cup\PlRedun
\end{align*}
denote the payload along with its error-check symbols.
Just as in Example~\ref{xmp:no-stretch}, all fields other than \( \Payload \) have relative 
width \( O(\bw) \); in particular, \( \PlRedun \) has relative width \( \bw \).
With \( \bw_{k} \) chosen as in~\eqref{eq:amp-fr}, this allows the
capacity \( \cp_{1} \)  of the base
medium \( M_{1} \) to still be constant, as the sum over all levels \( k \)
of the relative widths of all non-aggregated fields is bounded.

Our transition function \( \trans \) will have bandwidth rate \( \bw=\bw_{k} \), 
as in~\eqref{eq:bw}, and so the
simulated transition function $\trans^{*}$ has bandwidth rate $\bw^{*}=\bw_{k+1}$.
This way the amount of information needed from the neighbor colonies will be small.
Small bandwidth requires, however, another complication.
All information processing will happen on the \( \Work \) track, which has relative width \( O(\bw) \).
As at any one time, this track can only handle a small portion of the
the \( \Payload \) track of the colony, we introduce the notion of \df{packets}.

\begin{definition}[Packets]\label{def:packets}
  As the payload of a simulated cell is the aggregated payload of the cells
  of its simulated colony, we define the \df{payload capacity} as
  \begin{align*}
   \cp'_{k}=|\Payload_{k}|=\B_{k}|\Payload_{1}|.
  \end{align*}
We partition the \( \Payload \) field
into \( \nP \) sub-fields of relative width
\( \bw \) called \df{field packets} numbered \( 0,1,\dots,\nP-1 \).
The \( i \)th field packet of the \( \Payload \) field will be denoted \( \Payload_{i} \).
We also partition each track \( \F \) of whole colony into \( \nP^{*} \) \df{track packets}
of \( \bw^{*}\Q  \) cells each.
The \( j \)th track packet, with addresses in \( \lint{j\bw^{*}\Q}{(j+1)\bw^{*}\Q} \),
will be denoted as the \df{location} \( \_\F_{j} \),
In particular the \( j \)th track packet of the \( \Payload \) track is location \( \lPayload_{j} \)
\end{definition}

\begin{definition}[Payload error checks]  \label{def:PlRedun}
The error checks of track packet \( \lPayload_{j} \), will be in location
  \( \lPlRedun_{j} \) on the \( \PlRedun \) track.
Location 
\begin{align}\label{eq:lPayload'}
   \lPayload'_{j}=\lPayload_{j}\cup\lPlRedun_{j}
\end{align}
combines location \( \lPayload_{j} \) with its error check symbols for a complete code.

We consider the whole \( \Payload \) track
of each consecutive \( \numDirAff \) cells in the packet \( \lPayload_{j} \) as an
information symbol; the error check symbols of \( \lPayload_{j} \) go to
location \( \lPlRedun_{j} \).
To avoid that the damage rectangle affect both an information symbol and an error-check symbol
of the same code, we divide the track packet \( \lPayload'_{j} \) into left and right halves.
The error checks for the left half of \( \lPayload_{j} \) go into the right half of \( \lPlRedun_{j} \),
and those for the right half go into the left half of the same.
\end{definition}

The error check symbols in each half of the track packet \( \lPlRedun_{j} \) have
\( 4\numDirAff \) error check symbols of size \( <\cp \).
The number of bits in \( \lPlRedun_{j'} \)
is \( \Q\bw^{*}\bw\cp \) so these error checks will
fit in, according to condition~\eqref{eq:bw-b3} on amplifier our parameters.

\begin{figure}
 \setlength{\unitlength}{0.23mm}
 \[
   \begin{picture}(600,260)
     {
       \put(0,0){\framebox(600,240){}}
     }
    \put(0,245){\makebox(120,20){\( 0 \)}}
    \put(120,245){\makebox(120,20){\( 1 \)}}
    \put(240,245){\makebox(240,20){\( \dots \)}}
    \put(480,245){\makebox(120,20){\( \nP^{*}-1 \)}}
    \put(610,215){\makebox(0,25)[l]{0}}    
    \put(610,190){\makebox(0,25)[l]{1}}    
    \put(610,145){\makebox(0,25)[l]{\vdots}}    
    \put(610,90){\makebox(0,25)[l]{\(\nP-1\)}}    
    \put(160,215){\line(1,0){20}}    
    \put(160,190){\line(1,0){20}}    
    \put(-10,75){\makebox(0,150)[r]{\( \Payload \)}}
    \put(200,115){\makebox(300,75){\(\ddots\)}}
    \put(480,65){\line(0,1){175}}
    \put(240,65){\line(0,1){175}}
    \put(120,65){\line(0,1){175}}
    \put(120,65){\line(0,1){175}}
    \put(160,90){\line(0,1){150}}
    \put(180,90){\line(0,1){150}}
    \put(40,65){\line(0,1){25}}
    \put(160,115){\makebox(20,80){\( \vdots \)}}    
    \put(160,115){\line(1,0){20}}    
    \put(0,90){\line(1,0){600}}
    \put(0,65){\line(1,0){600}}
    \put(0,40){\line(1,0){600}}
    \put(-10,105){\makebox(0,25)[r]{\( \lPlRedun_{1} \)}}
    \put(-10,105){\vector(1,-1){30}}

    \put(-10,65){\makebox(0,25)[r]{\( \PlRedun \)}}
    \put(-10,40){\makebox(0,25)[r]{\( \Util \)}}
    \put(-10,0){\makebox(0,40)[r]{\( \Work \)}}
    \multiput(0,0)(20,0){30}{\line(0,1){3}}
\end{picture}
 \]
 \caption{Subdivision of the payload into field packets (horizontal lines) and track packets (vertical lines).}
 \label{fig:packets}
\end{figure}


The leftmost and rightmost packet of the (encoded) payload will be displayed to neighbor colonies as
tracks 
\begin{align}\label{eq:PayloadToNb}
  \PayloadToNb_{j},\; j=\pm 1
\end{align}
of width \( \bw \);
these displays will only be updated in the last step of the work period.

\paragraph{The code of the amplifier}
Let us define the code and decoding $\fg_{*},\Phi^{*}$ 
used in the simulation $\Phi$
taking medium $M=M_{k}$ into medium $M^{*}=M_{k+1}$.

\begin{definition}[Encoding]\label{def:encoding}
The encoding 
 \begin{align*}
  \fg_{*}:\states_{k+1}\to\states_{k}^{\Q _{k}}
 \end{align*}
takes the state of a cell, applies to it the 
error-correcting code $\alpha_{*}$ introduced above in
Definition~\ref{def:the-err-corr-code} and
copies the result to the $\Info$ track of a colony.
The other parts of the colony are set in such a way as they are to be
found at the beginning of a work period.
\end{definition}

There belongs a decoding function $\fg^{*}:\states_{k}^{\Q _{k}}\to\states_{k+1}$ 
to this encoding, but the actual decoding $\Phi^{*}:\eta\mapsto\eta^{*}$ 
depends on the history, not only on the current configuration \( \eta(\cdot,t) \).
One reason is that the transition of a colony from one work period to the next
one takes some time, due to asynchrony.
But we can also use the looking-back to avoid dependency on small noise in the near past.

\begin{definition}[Decoding]\label{def:decoding}
If the \( \pair{x}{t} \) is not in \( \Damage^{*} \), then define
 \[
 \eta^{*}(x,t)=\fg^{*}(x,t)=\Vac,
 \]
 unless there is a latest time \( t'\in \rint{t-2\Q \Tu}{t} \) that is a switching time
of a cell of \( \cC(x) \) such that the rectangle
\begin{align*}
 \lint{x-(\Q +1.1)\B }{x+(2\Q +1.1\B )}\times\rint{t'-\splitT}{t'}  
\end{align*}
is damage-free, and the colony is covered by a domain of 
cells belonging to the same work period
(the condition on switching just makes sure that a ``latest'' value exists).
In the latter case, let us define the 
strings \( s=\eta(\cC(x),t').\Util \) and \( p_{i}=\eta(\cC(x),t').\lPayload'_{i} \).
For the decoded value to be nonvacant, we also require that each of \( s,p_{i} \)
contains at most 4 errors in the sense that \( \alpha_{*}(\alpha^{*}(s)) \) differs
from \( s \) in at most 4 places of up to 2 consecutive symbols) and the same for each \( p_{i} \)
for \( i=0,\dots,\nP^{*}-1 \).
If this is satisfied, then \( \eta^{*}(x,t).\Payload^{c}=\alpha^{*}(s) \) while
\( \eta^{*}(x,t).\Payload \) is the concatenation of all \( \alpha^{*}(p_{i}) \).
We will say that \( \eta(\cC(x),t) \) is \df{error-free} if each of \( s,p_{i} \) are error-free.
\end{definition}

Note that if this definition gave a nonvacant, non-bad value to
\( \eta^{*}(x,t) \), then it gave the same value to all \( \eta^{*}(x,u) \) for
\( u\in\clint{t'}{t} \) where the time \( t' \) was defined above.

\paragraph{Rules for coding and decoding}
Before introducing the coding and decoding rules, a few more notions are needed.
In each rule of the present section, the operation%
 \glo{maj@\( \Maj \)}%
 \[
   \Maj_{i = 0}^{2}s_{i}
 \]
  when applied to strings \( s_{0},s_{1},s_{2} \) is the result of taking the
bit-wise majority of these three strings.
The rule
  \begin{equation}
    \label{eq:Broadcast}
  \Broadcast(\loc, \F(I))    
  \end{equation}
 takes the information found in location \( \loc \) and writes it into the
 \( \F \) field of every cell in interval \( I \).
 The default for \( I \) is the whole colony \( \lint{0}{Q} \).

We will use a general type of rule%
 \glo{check0@\( \rCheck_{0}\)}%
   \begin{equation}
     \label{eq:Check0}
  \rCheck_{0}(\var{prop},\F(I),X_{1},X_{2},\dots)     
   \end{equation}
  checking some global property of some parameters \( X_{1},X_{2},\dots \)
each of which is either an explicit string or a location.
Here it is assumed that
\( \var{prop} \) is some program for the universal computing medium
such that after we run the rule of that medium for \( \Q  \) steps the bit \( b \)
representing the outcome of the check will be broadcast into location \( \F(I) \)
(with the help of the rule \( \Broadcast \) given above).

The coding and decoding functions can be computed on the
universal medium simulated on one of the tracks.
Recall rules~\ref{alg:copy-1} (\( \Copy \)) and~\ref{alg:Write} (\( \Write \)) defined in
Section~\ref{sec:simp-sim} and again in Section~\ref{sec:commun}.

\begin{definition}
  Let%
 \glo{vacant-str@$\VacantString $}%
 \[
  \VacantString 
 \]
 be a string representing the state $\Vac$ in our code $\alpha$ as in
 Definition~\ref{def:the-err-corr-code}.  
  Let \( \var{Decode-prog} \) \glo{decode-prog@\( \rul{Decode-prog} \)}%
  be a program for the medium \( \univ \) such that after applying
 \[
  \univ(\var{Decode-prog}\cc S;\Q ,2\Q ),
 \]
  location
 \( \fld{\_Decode-output} \) \glo{decode-output@\( \fld{\_Decode-output} \)}%
on track \( \Cpt.\Output \) contains \( \alpha^{*}(S) \) (see~\eqref{eq:err-corr-code-alpha}) if
\( \alpha_{*}(\alpha^{*}(S)) \) differs from \( S \) in at most 4 groups of 2
consecutive symbols each; otherwise, it contains \( \VacantString  \).
  \end{definition}

Rule~\ref{alg:Decode} (\( \Decode \))  \glo{decode@\( \Decode \)}%
takes a string from location \( \loc_{1} \), decodes it and
copies the result to location \( \loc_{2} \).
It will always give a result, even if its input in \( \loc_{1} \) is not a code word.
(By the assumption made in Definition~\ref{def:effic-univ}, \( \univ \) is
commutative, so we do not have to worry about the order of execution in
different cells.)
\begin{algo}\caption{sub-rule \( \protect\Decode(\protect\loc_{1},
\protect\loc_{2}) \)}\label{alg:Decode}
  \( \Cpt.\Input \gets  * \);\;
  \( \Write(\var{Decode-prog}, \fld{\_Prog}) \);\;
   \( \Copy( \loc_{1}, \fld{\_Decode-arg}) \);\;
   apply \( 2\Q  \) times the rule \( \univ \)  to track \( \Cpt \);\;
   \( \Copy(\fld{\_Decode-output}, \loc_{2}) \)
 \end{algo}

The rule  \glo{encode@\( \Encode \)}%
\( \Encode(\loc_{1},\loc_{2}) \)  performs encoding in a similar way.
It will always output a code word; so applying first \( \Decode \) and then \( \Encode \)
turns every word into a code word. 
The locations \( \loc_{i} \)\glo{loc@\( \loc_{i} \)} in the encoding and decoding
rules will be allowed to be given indirectly, as explained in 
Definition~\ref{def:indirect-copy}.

In view of Condition~\ref{cond:comput}\ref{i:comput.legal-comp},
for certain  information \( X \)
retrieved from the neighbors we want to make sure that it
has the proper form of a code-word.
Decoding and then encoding it would do this in the absence of the damage rectangle.
The rule \( \Legalize \) in Algorithm~\ref{alg:legalize}
will make sure that this happens even in the presence of a damage
rectangle.
It repeats the decoding-encoding three times (using a temporary location \( \lDecoded \))
and then takes the majority of the results.
Moreover, then repeats the whole procedure again.

\begin{algo}[H]\caption{sub-rule \( \Legalize(\F(I)) \)}\label{alg:legalize}
  \pRepeat{\( 2 \)}{
    \For{\( j=0 \) \KwTo \( 2 \)}{
      \( \Decode(\F(I),\lDecoded) \);\;
      \( \Encode(\lDecoded,\Vote_{j}(I)) \)
    }
    \lIf{\( Addr\in I \)}{\( \F\gets \Maj^{2}_{j=0}\Vote_{j} \)}
  }
\end{algo}

\begin{lemma}\label{lem:legalize}
  After an application of the rule \( \Legalize(\F(I)) \), location
  \( \F(I) \) contains a word with at most 2 errors.
  Moreover, if \( \F(I) \) has at most 2 errors and decodes to some value \( X \)
  then result also decodes to \( X \).
\end{lemma}
\begin{proof}
Suppose first that \( \F(I) \) has at most \( 2 \) errors and decodes to \( X \).
Consider any one of the two repetitions;
the decoding-encoding writes \( \Vote_{j}(I) \).
If damage happens during it then the result is worthless.
But otherwise, earlier damage could add at most 2 more errors,
and since \( \Decode \) corrects 4 errors, the result in \( \Vote_{j} \)
is the code of \( X \).
Therefore the result in the first repetition differs from the code of \( X \) in
the at most 2 errors added possibly in the last majority step.
The same reasoning applies to the second repetition.

Suppose now that \( \F(I) \) is an arbitrary word.
Suppose that damage does not occur in the first repetition.
Then this part turns \( \F(I) \) into a code word (without errors),
and the previous analysis applies to the second repetition.
If damage occurs in the first repetition then
\( \F(I) \) may be again arbitrary as the result of this part.
But then the second repetition is damage-free and writes a code word
into \( \F(I) \).  
\end{proof}

\subsection{Sending and retrieval}\label{sec:send-retrieve}

In order to satisfy Condition~\ref{cond:comput}\ref{i:comput.trans},
it is necessary to find a time when the information received from both neighbor
colonies can be attributed to the same moment.
For this, the communication part of the work period will be substantially longer than
the rest.
During this time, the colony will continually post the information destined for the neighbor
colonies.
One of these pieces of information will be the age (the age of one cell would be sufficient, but
for simplicity and safety, all cells will send their age).
On the other hand, the receiving colony keeps checking (using the sender's age)
whether the sender's age of \emph{both} neighbors is from a safe part of their
work period: if yes then it stops receiving (by setting \( \Receiving\gets 0 \)), and keeps the
information already has to be used later by \( \Compute \).

Recall the colony work period in Algorithm~\ref{alg:workperiod}.
The rule \( \Send \), running almost all the time,
will try to mail information to the neighbor colonies,
running the rules \( \PostMail \) listed in Algorithm~\ref{alg:Send} simultaneously.
The retrieval part is played by two more rules running simultaneously.
The rule \( \Extend \) of Algorithm~\ref{alg:Extend}
tries to extend some arms of the colony
(consisting of cells of kind \( \Channel_{j} \))
left and right, for use in communicating with a possible non-adjacent
neighbor colony.
The rule \( \Retrieve \) in Algorithm~\ref{alg:Retrieve-2} (which is very different from the
example in~\ref{alg:Retrieve-1}) 
will try to retrieve information sent by the neighbors this way.
It records in track \( \Age_{k} \), etc.~the \( \Age \), \( \Doomed \) and \( \New \) fields sent
from the neighbor colonies in direction \( k \) (see below the explanation for \( \New \)).
For example \( \Age_{-1} \) at address \( a \) is
expecting the \( \Age \) of cell at address \( a \) of the left neighbor colony.
As the neighbor keeps sending this information, the receiving cell keeps updating it.

The main fields of the represented big cell of each colony are encoded onto its \( \Util \) track.
This entire track of the left neighbor will be sent to the \( \Retrieved_{-1} \) track of the receiving
colony.
The slices of \( \Payload_{i} \) of the represented cell \( (i=0,\dots,\nP^{*}-1) \)
are encoded into the locations \( \lPayload'_{i} \) of the colony.

\begin{algo}\caption{rule \( \Send \)}\label{alg:Send}
  \If{\( \Age<\growStart \)}{
    \For{\( k\in\{-1,1\} \)}{
      \( \PostMail(k,\Age,\Age_{-k},1,0,0,\Q) \)\; 
      \( \PostMail(k,\Doomed,\Doomed_{-k},1,0,0,\Q) \)\; 
      \( \PostMail(k,\New,\New_{-k},1,0,0,\Q) \)\; 
      \( \PostMail(k,\Util,\Retrieved_{-k},1,0,0,\Q) \)\;
    }
  }
\end{algo}

The \( \Send \) rule of Algorithm~\ref{alg:Send}
will run most of the time, and each \( \PostMail \) in it will
keep sending its information over and over.
At the beginning of the work period the field 
\begin{align*}
 \Receiving  
\end{align*}
gets its default value 1, and the 
\begin{align*}
 \AdjCol_{k}, k=\pm 1
\end{align*}
field, indicating whether there is an adjacent neighbor colony in direction \( k \), gets the default value 0.
The \( \Retrieve \) rule in Algorithm~\ref{alg:Retrieve-2} will set \( \Receiving\gets 0 \)
as soon as it finds that the information from both directions is \df{safe}:
it has not changed ``recently'' and will not change ``soon'' (precise definition is in~\eqref{eq:Safe}).
When \( \Age \) is safe any information sent from
any other cell of the colony must have already been from the same work period,
and the updating of payload has not yet been
started---so all payload information sent from this colony is consistent.
Here we use the fact from Lemma~\ref{lem:mail-speed} that the sending process itself
is reasonably fast, happening in \( O(\Q) \) time.
Let
\begin{align}\label{eq:Safe}
  \Safe= \{\Undef\}\cup\rint{\konst\cdot\Q}{\payloadStart-\konst\cdot\Q},
\end{align}
where \( \konst \) is the constant started to be used in~\eqref{eq:compute-start}.
Information from direction \( k \) has not changed lately and will not change soon
if \( \Age_{k}\in\Safe \).
Indeed, if \( \Age_{k}=\Undef \) then for a while before, the neighbor colony in direction \( k \)
was either vacant or consisted of doomed cells, and for a while after,
if it is not vacant or doomed then its \( \New \) track is filled with 1's---in all
these cases the neighbor big cell in direction \( k \) will be treated as vacant.

The definitions in~\eqref{eq:compute-start} show that, with \( \rep \) and \( \konst \)
large enough the overwhelming majority of the work period is safe.
Indeed, the safe part has age length \(   \konst\Q(2\rep\nP^{*}-1) \),
hence has a time interval of length at least
\begin{align}\label{eq:safe-part}
  \konst\Q(2\rep\nP^{*}-3)\Tl\ge   1.5\konst\rep\nP^{*}\Q\Tl
\end{align}
in which all age is safe.
Unsafe are only parts outside it, in up to two time intervals of total age length
\( \Q(\rep\nP^{*}+6\lambda+4 + \konst) \),
and hence total time length \(  \le 1.5\rep\nP^{*}\Q\Tu \) if \( \rep \) is large.

The \( \Retrieve \) rule uses the \( \Adj \) field of the one of the mails (say the
one sending the age information) to set \( \AdjCol_{k} \).
The bits \( \Doomed \), \( \New \), \( \AdjCol_{k} \) would normally be the same over the whole colony,
and would differ from the common value only at the places directly affected by
damage (the damage may have happened elsewhere, with the mail carrying the change
to its destination).

\begin{algo}\caption{rule \( \Retrieve \)}\label{alg:Retrieve-2}
  \If{\( \Receiving=1 \)}{
    \For{\( k\in\{-1,1\} \)}{
      \( \ReceiveMail(k,\Age,\Age_{-k},1,0,0,\Q) \)\; 
      \( \ReceiveMail(k,\Doomed,\Doomed_{-k},1,0,0,\Q) \)\; 
      \( \ReceiveMail(k,\Util,\Retrieved_{-k},1,0,0,\Q) \)\;
      \( \ReceiveMail(k,\PayloadToNb_{k},\PayloadFromNb_{-k},1,0,0,\Q) \)\;
      \lIf{\( \Doomed_{k}=1 \) \algOr (\( \New_{k}=1 \) \algAnd \( \Age_{k}<\konst\cdot\Q \))}
      {\( \Retrieved_{-k}\gets\VacantString(\Addr) \)}
      \( \AdjCol_{-k}\gets\Mail_{k,-1,\Age,\Age_{-k}}.\Adj \)
    }
    \lIf{\( \{\Age_{-1},\Age_{1}\}\subset\Safe \)}{\( \Receiving\gets 0 \)}
  }
\end{algo}

\subsection{Computation rules}\label{sec:comp-rules}

Recall that we are proving Lemma~\ref{lem:amp} (Amplifier) by 
building the transition function of the medium \( M_{k} \).
Given the amplifier parameters and the number \( k \) denoting the level,
let us fix where the frame information is to be found, using the scheme
with parameteres introduced in Definition~\ref{def:Param}.
(The rule evaluating the transition function assumes that the collection of
amplifier parameters are described in a parameter called \( \Frame \), 
the parameter \( k \) is described in parameter called \( \Height \).

Recall also the structure of the colony work period in algorithm~\ref{alg:workperiod}.
Here we will describe the \( \Compute \), \( \ProcPayload \) and \( \Finish \) rules.
The main simulation will interpret the program that needs to be applied
to the represented states of the big cell and its neighbors.
Given that all computation happens on the narrow \( \Work \) track, 
the rules processing the payload rely extensively on the copying rule of
Algorithm~\ref{alg:copy-2}.
The copy process, by its simplicity, is fault-tolerant: the
damage rectangle can only locally corrupt the information processed by it.

By Theorem~\ref{thm:rule-lang}, our rule language can be interpreted by a cellular automaton
on some work track.
In our example simulation in Section~\ref{sec:simp-sim} the evaluation
rule \( \Eval \) worked on a track was named \( \Cpt \), with
its result on the track \( \Cpt.\Output \).
Now some post-processing will also follow, but let us discuss first
the input information to the rule, some of which 
will be retrieved from neighbor colonies.
That retrieval is defined in Section~\ref{sec:send-retrieve}.
According to the requirement on amplifiers,
the transition function \( \trans \) carries out the computational part \( \PlTrans \)
via Definition~\ref{def:carries}.
Denote by
\begin{align*}
   \Payload^{c}=\All\setminus\Payload
\end{align*}
the complement of the payload field.
As arguments to the transition function, let \( \br=(r_{-1},r_{0},r_{1}) \), and
\( \ba=(a_{-1},a_{1}) \): the latter are
the bits showing whether the left and right neighbors are adjacent.
Recall the notations in Example~\ref{xmp:sim.aggreg}.
From \( \br \), the transition function \( \trans(\br,\ba) \) will depend
only on
\begin{align*}
  r_{0},\; r_{j}.\Payload^{c},\;  r_{j}.\Payload.\slice^{\bw,-j},\; j=-1,1.
\end{align*}
Thus, from the \( \Payload \) of the neighbors only the last segment
from the left neighbor and the first one from the right neighbor is needed.
This way, when simulating \( \trans^{*} \) by a colony then
due to~\eqref{eq:bw-b2}, all needed information from the
neighbor colonies for computing \( \trans^{*}(\br,\ba) \) will fit onto
the \( \Work \) track of the computing colony.
Indeed, for \( j=-1,1 \), for \( r_{j}.\Payload^{c} \) one needs
only the \( \Util \) track of the neighbor colony in direction \( j \).
For \( r_{-1}.\Payload.\slice^{\bw^{*},1} \) one needs only
the information in location \( \lPayload'_{\nP^{*}-1} \)
of the left neighbor colony as defined in~\eqref{eq:lPayload'}, 
and for \( r_{1}.\Payload.\slice^{\bw^{*},-1} \) one needs \( \lPayload'_{0} \)
from the right neighbor colony.
(The part of the simulating program that needs the latter information
is the rule \( \UpdatePayload \) as defined in Algorithm~\ref{alg:Update-payload}.)

\begin{sloppypar}
After retrieval, as described in Section~\ref{sec:send-retrieve}, 
the encoded argument \( r_{j}.\Payload^{c} \), \( j=\pm 1 \)
for the simulation will be on track \( \Retrieved_{j} \), while for \( r_{0}.\Payload^{c} \)
it is on track \( \Util \).
After decoding, the results will appear in some locations
\( \lArg_{m} \) for \( m\in\{-1,0,1\} \).
The bits \( a_{-1}\), \( a_{1} \) can be determined directly from tracks \( \AdjCol_{j} \),
as the \( \Retrieve \) rule of Algorithm~\ref{alg:Retrieve-2} has set them.
(Each of these two tracks contains the repetition of the same single bit, except for the places
affected by damage, so decoding can just take the majority.)
 \end{sloppypar}

 \begin{sloppypar}
The rule \( \Eval \) of Section~\ref{sec:simp-sim}
computes the transition function on the \( \Cpt \) track,
using the program \( \Interpr \) similarly to Algorithm~\ref{alg:simp-sim.eval}
in Section~\ref{sec:simp-sim}.
It uses the locations \( \lArg_{m} \) and sub-rule \( \Initialize \) mentioned in~\eqref{eq:Initialize}.
But in the current version, in Algorithm~\ref{alg:eval},
there are two more arguments, in locations \( \lAdjArg_{j} \) for \( j=\pm 1 \),
for the bits \( a_{-1},a_{1} \) in \( \trans^{*}(\br,\ba) \).
\end{sloppypar}

\begin{sloppypar}
The modification also accounts for payload processing.
The program \( \Interpr \) for the universal cellular automaton,
used in the rule \( \Interpret \),
will now be enhanced as follows, using some strings, small in size,
called \df{symbolic commands}, from the following list.
\begin{equation}\label{eq:symb-outputs}
  \begin{aligned}
 &  \prg{Update-payload},\;  \prg{Refresh-payload},
  \\ & \prg{Mail}_{k,d}.\prg{Info}\mathtt{\gets} \prg{Payload}_{i},\;
       \prg{Payload}_{i}\mathtt{\gets}\prg{Mail}_{k,d}.\prg{Info},
  \end{aligned}
\end{equation}
where \( i\in\{0,\dots,\nP^{*}-1\} \) and \( k,d \) are numbers, and the role of the
field \( \Mail_{k,d} \) is explained after~\eqref{eq:Mail-ind}.
When a condition requires it will \emph{broadcast} them (using the rule \( \Broadcast(x,\F) \)
of~\eqref{eq:Broadcast}):
write them into the field \( \Cpt.\Output.\PlCommands\) of \emph{each cell}.
In particular:
\begin{itemize}
\item
  \begin{sloppypar}
    Encountering the command \( \Payload_{i}\gets\Mail_{k,d}.\Info \) \\
    in the program
for the big cell, broadcasts \( \prg{Payload}_{i}\mathtt{\gets}\prg{Mail}_{k,d}.\prg{Info} \)
onto track \( \Cpt.\Output.\PlCommands\)  while on command \( \Mail_{k,d}.\Info\gets\Payload_{i} \)
it broadcasts \( \prg{Mail}_{k,d}.\prg{Info}\mathtt{\gets}\prg{Payload}_{i} \).    
  \end{sloppypar}
\end{itemize}
The broadcast symbolic commands will be executed by the rule \( \ProcPayload \)---eventually,
after they have been written onto a locally maintained track \( \PlCommands \).
  \end{sloppypar}

 \begin{algo}[H]\caption{sub-rule \Eval}\label{alg:eval}
     \( \Write(\Interpr, \fld{\_Interpr}) \);\;
     \( \WriteParam(\MyRules, \fld{\_Prog}) \);\;
     \lFor{\( i=1 \) \KwTo \( N \)}{
       \( \WriteParam(\Param_{i}, \lParam_{i}) \)
     }
     \Initialize;\;
     \Interpret
   \end{algo}

\begin{sloppypar}
To deal with the possible damage rectangle, 
rule \( \Compute \) of Algorithm~\ref{alg:Compute}
will call \( \Eval \) three times, storing its output
on the track \( \Vote_{r} \) for \( r=0,1,2 \).
The track \( \Hold \) will receive then the majority of these results.
A post-processing step \( \UpdateLocMaint \) \glo{update-loc-maint@\( \UpdateLocMaint \)}
will update the locally-maintained fields
 \[
 \Doomed, \Growing_{j}\ (j\in\{-1,1\}),\PlCommands.
\]
It checks via rule \glo{check-vacant@\( \rul{Check-vacant} \)}
 \( \rCheck_{0}(\var{Check-vacant}, \F(I), \loc) \)
whether the string represented in location \( \loc \) is \( \Vac^{*} \) (and broadcasts the
result into track \( \F \) on interval \( I \)).
Let \glo{creating@\( \fld{\_Creating}_{j} \)}
\begin{align*}
 \fld{\_Creating}_{j}  
\end{align*}
be the location of the field \( \Creating_{j}^{*} \) of the represented cell
on the \( \Cpt.\Output.\Util \) track.
If the represented cell becomes vacant then, as expected, \( \Creating_{j}^{*}=0 \).
(This is how we will satisfy Condition~\ref{i:comput.legal-birth} in simulation, as
Definition~\ref{def:creator-emerging} requires creator to survive until after the creation.)
The value found here will be broadcast to the locally maintained track
\( \Growing_{j} \), since when the represented cell
is allowed to be creating, say, to the right, then the representing colony should be
allowed to grow to the right in an attempt to create an adjacent new colony to the right.
As usual, the argument \( \F \) means \( \F(\lint{0}{\Q }) \).  
\end{sloppypar}

 \begin{algo}[H]\caption{sub-rule \( \protect\UpdateLocMaint(r),  r=0,1,2 \)}
\label{alg:Update-loc-maint}
\( \rCheck_{0}(\var{Check-vacant}, \Vote_{r}.\Util, \Vote_{r}.\Doomed) \);\;
 \lpFor{\( j=-1,1 \)}{
   \( \Broadcast(\fld{\_\Creating}_{j}, \Vote_{r}.\Growing_{j}) \)
 }
 \end{algo}

 The rule \( \Compute \) repeats the evaluation process 3 times:
after repetition \( r \), we encode the result into
\( \Cpt.\Output \) and copy it onto track \( \Vote_{r} \).
Finally \( \Hold \) will be obtained by majority
vote from \( \Vote_{r} \), \( r=0,1,2 \).
Recall that \( \Cpt.\Output \) has sub-tracks \( \Util \) and \( \PlCommands \).
Only the sub-track \( \Util \) needs encoding, the track \( \PlCommands \)
contains the same value repeated in each cell.
Rule \( \Compute \) calls some more rules.
The sub-rule \( \rul{Randomize} \) takes the \( \Rand \) bit of the first cell of the colony,
and places its \( \alpha_{*} \)-code into the appropriate location on the \( \Hold \) track.
This is done only once, without any attempt of error-correction.
The line~\eqref{alg:Kind*} tests whether the colony encodes a 
germ big cell that expects to become part of a (big) colony.
It is interpreted easily, when we stipulate that \( \Kind^{*} \) and \( \Addr^{*} \)
occupy a special location on the \( \Util \) track.
The last lines of the code of \( \Compute \) update the locally maintained variables \( \Doomed \),
\( \Growing_{j} \) and \( \PlCommands \) in a single step, so this step is their update time.

 \begin{algo}\caption{sub-rule \( \protect\Compute \)}\label{alg:Compute}
   \lFor{\( m\in\{-1,1\} \)}{\( \Legalize(\lRetrieved_{m}) \)}
   \( \Legalize(\Util) \);\;
   \For{\( r=0 \) to \( 2 \)}{
     \( \Decode(\Util, \lArg_{0}) \) ;\;
     \For{\( m\in\{-1,1\} \)}{
       \( \Decode(\lRetrieved_{m}, \lArg_{m}) \);\; 
       \( \Decode(\AdjCol_{m}, \lAdjArg_{m}) \)
     }
     \( \Eval \);\;
     \( \Encode(\Cpt.\Output, \Vote_{r}.\Util) \);\;
     \( \UpdateLocMaint(r) \)\;
     \nl\lIf{\( \Kind^{*}=\Germ \)  \label{alg:Kind*} \algAnd \( 0\le \Addr^{*}<\Q^{*} \)
     }
     {\( \Broadcast(\prg{Update-payload},\Cpt.\Vote_{r}.\PlCommands \)) }
     \lElse{\( \Broadcast(\prg{Refresh-payload},\Vote_{r}.\PlCommands) \) }     
   }
   \( \rul{Randomize} \);\;
   \( \Hold \gets  \Maj_{r=0}^{2}\Vote_{r} \);\;
   \( (\Doomed, \PlCommands)\gets (\Hold.\Doomed, \Hold.\PlCommands) \)\;
   \lpFor{\( j=-1,1 \)}{ \( \Growing_{j}\gets\Hold.\Growing_{j} \) }
\end{algo}

  \begin{sloppypar} 
If we had to do with a single cell and not a represented one then no payload refreshing
would be needed, but the rule \( \RefreshPayload \) (interpreting the corresponding
symbolic command), as  it applies to a represented cell,
will perform error-correction.
It may seem at first sight that some of these commands might contradict each other.
However, when the symbolic command
\( \prg{Update-payload} \) is outputted then it stands alone; indeed,
this is done when a germ cell is represented,
and germ cells are not part of any colony performing simulation.
Similarly, \( \RefreshPayload \) may seem to possibly contradict the interpretation of the
symbolic command
\( \prg{Payload}_{i}\mathtt{\gets}\prg{Mail}_{k,d}.\prg{Info} \); however, as they are both
are trying to keep (or transport) the same string, only some new fault can cause
a (local) contradiction.
  \end{sloppypar}

 The rule \( \RefreshPayload \) described in Algorithm~\ref{alg:Refresh-payload} 
uses the location \( \lToRefresh \) on the \( \ToRefresh \) track of the
work tape to which the packets of \( \Payload \)
will be spread out for processing.
To control for the occurrence of possible damage during processing,
each packet is processed three times for \( r=0,1,2 \):
after decoding and encoding, the result is stored in a location \( \lVote_{j} \)
on a track \( \Vote_{j} \).
Then the majority is copied back to the packet.

\begin{algo}\caption{sub-rule \RefreshPayload}\label{alg:Refresh-payload}
  \For{\( i=0 \) \KwTo \( \nP^{*}-1 \)}{
    \( \Copy(\lPayload'_{i},\lToRefresh) \);\;
    \For{\( r=0 \) \KwTo \( 2\)}{
      \( \Decode(\lToRefresh,\lDecoded) \);\; 
    \( \Encode(\lDecoded,\lVote_{r}) \)
  }
  \( \ToRefresh\gets\Maj_{r=0}^{2}\Vote_{r} \);\;
  \( \Copy(\lToRefresh,\lPayload'_{i}) \)
}
\end{algo}

\begin{sloppypar}
The computation of \( \PlTrans^{*}=\PlTrans_{k+1} \)
in rule \( \UpdatePayload \) as shown in Algorithm~\ref{alg:Update-payload},
takes place when \( \Kind^{*}=\Germ \) and the represented germ cell can expect to become
part of a colony (so \( 0\le\Addr^{*}<\Q^{*} \)).
This level is the highest, and as will be seen the possibility of damage on this
level can be ignored, so the rule \( \UpdatePacket \) in
Algorithm~\ref{alg:Update-packet} does not need the vote over
three repetitions used elsewhere.

Recall that \( \PlTrans_{k+1}=\PlTrans_{1}^{\B_{k+1},\bw_{k+1}} \) is an aggregated
and slowed version of \( \PlTrans_{1} \).
We compute it part-by-part, as we can only use
the narrow track \( \Work \).
So we will copy, starting from the left, new and new packets to the \( \Work \) track,
perform the updating using three neighboring ones, and copy back the result.
The actual simulation of \( \PlTrans_{1} \) will take place on a sub-track 
called \( \PlUpd \), of width \( \bw\cdot\cp \).
Recall from Definition~\ref{def:packets}
that the capacity of the medium \( \PlTrans_{1} \) is \( \cp'_{1} \).
After copying and decoding, we will have on the \( \PlUpd\) track
three consecutive locations of the same length:
\( \lPlUpd_{p} \) for \( p=-1,0,1 \),
A packet \( \lPayload_{i} \) will be distributed in this location.
The amount of information in \( \lPayload_{i} \) is \( \bw^{*}\Q\cdot\cp \),
so we need 
\begin{align}\label{eq:lPlUpd}
 \bw^{*}\Q\cdot\cp/\bw\cdot\cp=\bw^{*}\Q\nP  
\end{align}
cells as the width of the track is \( \bw\cdot\cp \)..
Each cell in this location will contain 
\begin{align*}
 S=\bw\cdot\cp/\cp'_{1}  
\end{align*}
symbols of \(  \PlTrans_{1} \),
so this is an aggregation code, and each cell can compute on this track the
aggregated transition function \( \PlTrans_{1}^{S} \) in one step.
The number of steps is equal to the number of cells in \( \lPlUpd_{0} \), as in~\eqref{eq:lPlUpd}.
At the end, only the result in \( \lPlUpd_{0} \) will be used.
Decoding, updating and encoding will be repeated three times just as in the refreshing case.
The rule \( \UpdatePacket \) assumes that the needed three
packets have already been copied to the consecutive
locations \( \lToUpdate_{p}\) for \( p=-1,0,1 \), and writes its result into location \( \lUpdated \).
In rule \( \rul{Update-payload} \),
in order to update the first and last package, the last package 
of the left and the first package of the right neighbor colonies are needed;
but they have already been retrieved by \( \Retrieve \) into \( \PayloadFromNb_{j} \),
\( j=\pm 1 \).
  \end{sloppypar}

  \begin{algo}\caption{sub-rule \(\rul{Update-packet}\)}\label{alg:Update-packet}
      \For{\( p=-1,0,1 \)}{
        \( \Decode(\lToUpdate_{p},\lPlUpd_{p}) \)
        }
        \pRepeat{\( \bw^{*}\Q\nP \)}{
          \( \PlUpd\gets \PlTrans_{1}^{S}(\PlUpd^{-1},\PlUpd,\PlUpd^{1}) \)
        }
        \( \Encode(\lPlUpd_{0},\lUpdated) \)
  \end{algo}

  \begin{algo}\caption{sub-rule \UpdatePayload}\label{alg:Update-payload}
    \( \Copy(\PayloadFromNb_{-1},\lToUpdate_{-1}) \);\; 
    \( \Copy(\lPayload'_{0},\lToUpdate_{0}) \);\; 
    \( \Copy(\lPayload'_{1},\lToUpdate_{1}) \);\; 
    \( \UpdatePacket \);\;
    \( \Copy(\lUpdated,\lPayload'_{0}) \);\; 
    \( \Copy(\lToUpdate_{0}, \lToUpdate_{-1}) \);\; 
    \( \Copy(\lToUpdate_{1}, \lToUpdate_{0}) \);\;
  \For{\( i=1 \) \KwTo \( \nP^{*}-2 \)}{
    \( \Copy(\lPayload'_{i+1}, \lToUpdate_{1}) \);\; 
    \( \UpdatePacket \);\;
    \( \Copy(\lUpdated,\lPayload'_{i}) \);\; 
    \( \Copy(\lToUpdate_{0}, \lToUpdate_{-1}) \);\; 
    \( \Copy(\lToUpdate_{1}, \lToUpdate_{0}) \);\;
  }    
   \( \Copy(\PayloadFromNb_{1}, \lToUpdate_{1}) \);\;
   \( \UpdatePacket \);\;
    \( \Copy(\lUpdated,\lPayload'_{\nP^{*}-1}) \)
  \end{algo}

\begin{algo}\caption{sub-rule \ProcPayload}\label{alg:Proc-payload}
  \lIf{\( \PlCommands \) contains \prg{Update-payload}}{\UpdatePayload}
  \lIf{\( \PlCommands \) contains \prg{Refresh-payload}}{\RefreshPayload}
  \lIf{\( \PlCommands \) contains \( \prg{Mail}_{k,d}.\prg{Info}\mathtt{\gets}\prg{Payload}_{i} \) for some \( k,d,i \)}{
    \( \Copy(\lPayload'_{i},\lMail^{*}_{k,d}.\Info) \)
  }
  \lIf{\( \PlCommands  \) contains \( \prg{Payload}\mathtt{\gets}\prg{Mail}_{k,d}.\prg{Info} \)
    for some \( k,d,i \)}{
    \( \Copy(\lMail^{*}_{k,d}.\Info,\lPayload'_{i}) \)
    }
   \For{\( k\in\{-1,1\} \)}{
     \nl\algLet \( i= (\nP-1)(k+1)/2 \)\label{alg:endslice}\;
     \( \Copy(\lPayload'_{i},\TempPayloadToNb_{k}) \)
     }
  \end{algo}

  \begin{sloppypar}
Algorithm~\ref{alg:Proc-payload} (\( \ProcPayload \)) depends on the locally maintained
track \( \PlCommands\).
We will analyze later what happens when it is damaged.
(Only the command \( \UpdatePayload \) is not controlled for damage; however, as it is
called only on the highest level, we will be able to assume no damage on that level.)
Location \( \lMail^{*}_{k,d}.\Info \) is, just as all fields of the represented cell,
on the \( \Util \) track.
In line~\eqref{alg:endslice}, \( i=0 \) if \( k=-1 \) and \( \nP-1 \)  if \( k=1 \).
As said above, now this program outputs not just the candidate
values of the fields in \( \All^{*}\setminus\Payload^{*} \) onto \( \lCpt.\Output.\Util \)
but also a couple of symbolic commands from the list~\eqref{eq:symb-outputs}
onto \( \lCpt.\Output.\PlCommands \).
In particular, every time when an instruction \( \Mail_{k,d}\gets\Payload_{i} \)
is encountered it is not interpreted by an action, only by outputting the 
corresponding symbolic command.
   \end{sloppypar}

At \( \Age=\U -1 \) 
  a final computation step takes place in the rule \( \Finish \) in Algorithm~\ref{alg:Finish}, which indeed
  is only a single step.
  The part applying to germ cells will be explained later.

  \begin{algo}[H]\caption{rule \( \Finish \)}\label{alg:Finish}
  \eIf{\( \Doomed=1 \) \algOr (\( \Kind=\Germ \) \algAnd \( \Addr\not\in\lint{0}{\Q} \)) }
  {\( \Kind\gets\Latent \)} {
    \( \Addr\gets\Addr\bmod\Q \)\;
    \( \Age\gets 0 \)\;
    \( \Receiving\gets 1 \)\;
    \eIf{\( \Kind\in\{\Growth, \Germ \} \) }{
      \( \Kind\gets\Member \)\;
      \( \New\gets 1\) }{
      \( \New\gets 0 \)\;
      \( \Util\gets\Hold \)\;
      \lFor{\( k\in\{-1,1\} \)}{\( \PayloadToNb_{k}\gets \TempPayloadToNb_{k} \)}
    }
  }
\end{algo}

\subsection{Lifting}\label{sec:lifting}

Recall that the payload computation was carried out in germ cells.
As seen in Section~\ref{sec:germ-growth}, a germ is trying to grow into an area
occupied by five colonies.
Only the middle part will become a new colony, intended to simulate a new big cell:
a latent cell on a higher level.
Before the germ work period ends, however, the new colony will carry out
a much simplified version of Algorithm~\ref{alg:Refresh-payload}
(\( \RefreshPayload \)) called \( \Lift \), creating the error checks of the new level for the payload.
Also, all cells of the new colony are supposed to be of the same color, but this color must be
recorded also in the new, represented big cell:
assume this information is (along with its error checks) in place \( \lColor^{*} \) on the track \( \Util\).

As with all germ actions, no damage rectangle is expected, so there is no repetition-and-vote.

\begin{algo}\caption{sub-rule \Lift}\label{alg:Lift}
  \For{\( i=0 \) \KwTo \( \nP^{*}-1 \)}{
    \( \Copy(\lPayload_{i},\lToEncode) \);\;
    \( \Encode(\lToEncode,\lEncoded) \);\;
    \( \Copy(\lEncoded,\lPayload'_{i}) \)
  }
  Encode into \( \lColor^{*} \) the \( \Color \) of cell with address 0.
\end{algo}

\section{The simulated medium is robust}\label{sec:sim-robust}

\subsection{Legality}

Here we will prove parts~\ref{i:comput.rand}-\ref{i:comput.legal-birth}
of Condition~\ref{cond:comput} (Computation Property) for big cells.
In the lemmas that follow, we fix a cell \( x \) and denote
\begin{align*}
 I=\lint{x-(\Q +1.1)\B }{x+(2\Q +1.1\B )}.
 \end{align*}

\begin{lemma}\label{lem:free-past}
For some time \( t \), let \( J=\rint{t-2\Tus-\splitT}{t} \).
Assume that the rectangle \( I\times J \) is damage-free, and at time \( t \)
the colony \( \cC=\cC(x) \) is covered by a domain of cells belonging to the same work period.
Then the string \( \eta(\cC,t) \) is error-free in the sense of Definition~\ref{def:decoding}.
\end{lemma}
\begin{proof}
Suppose that during the whole interval \( \rint{t-2\Tus}{t} \), the colony \( \cC \)
is covered by a domain.
In the absence of damage, this colony performs at least one complete work period
of computation without any interference from errors inside \( \cC \).
It follows then from the definition of the program (summarized in Algorithm~\ref{alg:workperiod})
that the string \( s \), as a result of such a computation, is error-free.

Suppose now that there is a first time \( t'>t-2\Tus \) such that the colony
is covered by a domain after \( t' \), but not before.
Leading back a trace-back path from any cell of \( \cC \), starting at time \( t' \)
the path is at any time in a domain covering its originating colony.
Indeed, otherwise Lemma~\ref{lem:bad-gap-infr} (Bad Gap Inference) implies
a bad gap which would not be closed, contradicting the assumption that \( \cC \)
is covered by a domain at time \( t \).
The only possibility is thus that of a growing arm that finishes covering \( \cC \) at
time \( t' \).

Now if \( t'\le t-1.5\Tus \) Lemma~\ref{lem:colony-trace-back} (Colony Trace-back)
together with the number of steps allowed between the end of growth and the
work period end
implies that a new work period starts before time \( t-\Tus \),
leaving a whole work period until time, guaranteeing an error-free colony.
If \( t'> t-1.5\Tus \) then Lemma~\ref{lem:colony-trace-back}
together with the number of steps allowed between the start and end of growth
implies that the growth started after time \( t-2\Tus \)
using the growth rule, so by the time \( t' \) the whole colony 
encodes a latent big cell in an error-free way.
No new errors arise later.
\end{proof}

 \begin{lemma}[Legality]\label{lem:legality}
  Assume that, for some number \( a \), the set \( \Damage^{*} \) does not
intersect \( \lint{x}{x+\Q \B }\times\rint{a-2\Tus}{a+2\Tus} \), 
further that \( \sigma_{1},\sigma_{2} \),
as in Definition~\ref{def:special-sw-times} (Special switching times)
are defined for cell \( x \) of medium \( M^{*} \) in the interval \( \rint{a-2\Tus}{a+2\Tus} \),
using Definition~\ref{def:decoding} (Decoding).
Then the following holds, using Definition~\ref{def:legal} (Legality):
\begin{cjenum}
  \item\label{i:legality.legal}
    \( \legal_{k+1}(\eta^{*}(x, \sigma_{2}-), \eta^{*}(x, \sigma_{2})) = 1 \), thus
    proving Condition~\ref{cond:comput}\ref{i:comput.legal-comp}.
  \item\label{i:legality.time-bounds}
    \( \Tls \le \sigma_{2} - \sigma_{1} \le \Tus \), thus proving
    Condition~\ref{cond:comput}\ref{i:comput.dwell-pd-bd}.

    \begin{sloppypar}
\item\label{i:legality.rand} The randomization property
  \( g(\alpha(\var{rand}),j, W_{0}(x,a),\eta^{*}) \), as defined
  in Condition~\ref{cond:comput}\ref{i:comput.rand}
  will hold: it can be expressed in the needed canonical simulation of
  Definition~\ref{def:canon-simul}.      
    \end{sloppypar}
\end{cjenum}
  \end{lemma}

\begin{Proof}
Let \( I=\lint{x}{x+\Q  \B } \).
Our assumption implies that damage is covered by at most one island in 
\( I\times\rint{a-2\Tus}{a+2\Tus} \).
According to Definition~\ref{def:decoding} (Decoding), there is a latest time \( v_{2} \) 
in \( \rint{\sigma_{2}-2\Q\Tu}{\sigma_{2}} \) that is a switching time
of a cell of colony \( \cC(x) \) such that \( I \) is damage-free during
\( \rint{v_{2}-\splitT}{v_{2}} \), and at time \( v_{2} \) the colony
is covered by a domain of cells belonging to the same work period.
The value \( \eta^{*}(x,\sigma_{2}) \) was decoded from the state of the colony at
time \( v_{2} \).
Let \( y_{1} \) be the cell with address \( \flo{\Q /2} \) in colony
\( \cC(x) \), then \( \pair{y_{1}}{v_{2}} \) is traceable, as in
Definition~\ref{def:traceable}.
Build a trace-back path from \( \pair{y_{1}}{v_{2}} \) 
to the latest time \( v_{1} \) in
\( \rint{\sigma_{1}-2\Q\Tu}{\sigma_{1}} \) at which the colony is again covered with
member cells belonging to the same work period,
moreover \( I \) is damage-free during \( \rint{v_{1}-\splitT}{v_{1}} \). 
There is such a time, due to the definition of \( \sigma_{1} \).
Since \( \sigma_{1},\sigma_{2} \) are defined, 
the colony \( \cC(x) \) encodes at time \( v_{1} \) with at
most 4 errors a big cell with value different from the one
encoded at time \( v_{2} \).

\begin{step+}{step:legality.free}
If \( I \) is damage-free during \( \rint{v_{1}}{v_{2}} \) then
claim~\ref{i:legality.legal} holds.
\end{step+}
\begin{pproof}
The computation runs now without faults, and as planned, computes an error-free
value.
Rule \( \Compute \) in Algorithm~\ref{alg:Compute} applies \( \rul{Legalize} \) to
the information retrieved from the neighbor colonies, so the result
computed from this information will be a legal one.
The switching time \( \sigma_{2} \) occurs just when the last cell of the colony
passes to a new work period. 
The process makes sure that either all cells of the colony are doomed (when the
computed value is vacant), or none them.
If they are doomed then the whole colony dies.
Otherwise, the rule \( \Finish \) installs the remaining changes in one step.
In both cases, the state at time \( \sigma_{2} \) will be a legal consequence of
the state at time \( \sigma_{2}- \).
\end{pproof} 

\begin{step+}{step:legality.non-free}
If \( I \) is not damage-free during \( \rint{v_{1}}{v_{2}} \) then 
claim~\ref{i:legality.legal} still holds.
\end{step+}
\begin{pproof}
Let \( t_{1}\in\rint{v_{1}}{v_{2}} \) be some time at which damage intersects \( I \),
and let \( t_{0} \) be the latest time before \( \min(t_{1}-1.1\destrT\Q, v_{1}) \) which is the
switching time of some cell, with the
property that all cells of the colony belong to the same work period.
By Lemma~\ref{lem:cover} (Cover), at time \( t_{0} \), the colony is covered by a
domain of cells.
Lemma~\ref{lem:free-past} implies that the colony encodes an error-free string.

Let us follow the development of the colony forward in time, 
starting from \( t_{0} \).
Healing succeeds, due to the Cover Lemma.
Locally maintained fields will be restored in a short time following the
occurrence of damage, via the rule \( \LocMaintain \) while other fields
changed in at most an interval of \( \numDirAff \) cells, with \( \numDirAff \)
defined in~\eqref{eq:numDirAff}.
After this, the work of the colony can be followed as before.

The time of the damage may fall into at most one of the three iterations of
any of the rules with three iterations (\( \RefreshPayload \), \( \UpdatePayload \),
\( \Legalize \),  \( \Compute \)).
The output of this one iteration in \( \Vote_{r} \)
will be outvoted by the other two values in the final majority vote.
So the only errors in \( \Hold \) can be the ones due
to their short-term damage in the final majority vote; 
these change at most \( \numDirAff \) cells.
In summary, we can make the same conclusion as
part~\ref{step:legality.free}, except that the encoded string may contain up to
2 errors (in case the island occurs in the last part of computation or later).

There is an exception to this reasoning: in the rules \( \RefreshPayload \) and
\( \UpdatePayload \), the copy command is not repeated three times, so
it is protected from the damage only in a limited
way in the rule \( \MoveMail \) of Algorithm~\ref{alg:move-mail-2}.
This protection makes sure that wrong information planted by damage 
is copied to at most one segment of size \( \numDirAff \); however, that one segment
will still be affected, so the damage rectangle can cause 2 errors in the place it occurs and 2
more in the place to which the information is carried.
But the error-correcting code, being able to correct 4 errors, will still deal with this.
\end{pproof} 

\begin{step+}{step:legality.time-lb}
We have \( \Tls \le \sigma_{2} - \sigma_{1} \).
 \end{step+}
 \begin{pproof}
Consider a cell \( y \) of the colony,
for example the one with address \( \flo{\Q /2} \), and lead a trace-back path from
time \( \sigma_{2} \) to time \( \sigma_{1} \).
Suppose it has \( n \) links: none of these are parental.
At time \( \sigma_{j}- \) all cells have ages \( \ge \U -2\Q  \), 
hence the age progress along the path is at least \( \U -2\Q  \).

By Lemma~\ref{lem:trace-back-bounds}, the age progress is 
\( \U -2\Q \le n/\p_{1}+\numDirAff+2 \), hence 
 \begin{align*}
   n\ge \p_{1}\U -\p_{1}(2\Q +\numDirAff+2).
 \end{align*}
By the same lemma, 
 \begin{align*}
\sigma_{2}-\sigma_{1} &\ge n\Tl-(\numDirAff+1)\Tl 
\\ &\ge \Tl\Paren{\p_{1}\U  -\Paren{2\p_{1}\Q +(\p_{1}+1)\numDirAff+2\p_{1}+1}}
\\ &=  \Tl\U' -\Tl\Paren{2\p_{1}\Q +(\p_{1}+1)\numDirAff+2\p_{1}+1}
   \\ &=  \Tl\U'\Paren{1-\Paren{2\p_{1}\Q +(\p_{1}+1)\numDirAff+2\p_{1}+1}/\U'}
\\ &\ge\Tl\U' (1-\rep\Q/\U') = \Tl^{*},
 \end{align*}
 where we used~\eqref{eq:U-def}, further
 \( \U'=\U'_{k} \) was defined in Section~\ref{sec:amp}.
Indeed, sufficiently large \( \rep \) will satisfy the last inequality.
\end{pproof} 

\begin{step+}{step:legality.time-ub}
We have \( \sigma_{2} - \sigma_{1}\le\Tus \).
\end{step+}
\begin{pproof}
  Lemma~\ref{lem:colony-trace-back} (Colony Trace-back) lower-bounds
the progress along a trace-back path;
then the proof is finished similarly to part~\ref{step:legality.time-lb} above.
\end{pproof} 
\begin{step+}{step:legality.rand}
 Let us prove~\ref{i:legality.rand}. 
\end{step+}
\begin{prooof}
  We need to prove essentially \( \Pbof{\eta^{*}(x,\sigma_{2}).\Rand=j}\le 1/2+eps'+\eps'' \),
  given how \( \eps'_{k} \) is defined in~\eqref{eq:eps-defs}.
  Here \( \eps'' \) bounds the probability that damage occurs in the window \( W_{0}(x,a) \) at all.
  On the other hand, when it does not occur then each cell at each time within the work period
  performs its needed action: in particular, it carries out the rule \( \Compute \) of
  Algorithm~\ref{alg:Compute}, and within it, the rule \( \rul{Randomize} \).
  This rule relies on the randomization action of a well-defined cell at a well-defined age
  of the work period, with the result being \( j \) with probability \( \le 1/2+\eps' \).
  The probability can only be increased by the probability bound that damage does occur in \( W_{0} \).
  It would be now just a tedious exercise to actually express \( g(\alpha(\var{rand}),j, W_{0}(x,a),\eta^{*}) \)
  formally as needed by a canonical simulation.
\end{prooof} 

\end{Proof}

\begin{lemma}
  Assuming \( \eta \) is a trajectory, \( \eta^{*} \) satisfies
  Condition~\ref{cond:comput}\ref{i:comput.legal-birth}.
\end{lemma}
\begin{proof}
  The proof is similar to that of the above lemma, only simpler.
There is a latest time \( v_{2} \) 
in \( \rint{\sigma_{2}-2\Q\Tu}{\sigma_{2}} \) that is a switching time
of a cell of colony \( \cC=\cC(x) \) such that \( I \) is damage-free during
\( \rint{v_{2}-\splitT}{v_{2}} \), and at time \( v_{2} \) the colony
is covered by a domain of cells belonging to the same work period.
The value \( \eta^{*}(x,\sigma_{2}) \) was decoded from the state of the colony at
time \( v_{2} \).
Let \( y_{1} \) be the cell with address \( \flo{\Q /2} \) in colony
\( \cC(x) \), then \( \pair{y_{1}}{v_{2}} \) is traceable, as in
Definition~\ref{def:traceable}.
Now the a trace-back path from \( \pair{y_{1}}{v_{2}} \) will lead to a
cell still near the middle of colony \( \cC \), where \( \cC \) is now not covered
by member cells.
According to the Cover Lemma it the cell is still in a domain without exposed edges,
so this domain will cover \( \cC \) as well as an originating colony.
Tracing back further one can see that the cells covering \( \cC \) (whether outer or germ)
were all formed by the rule \( \Grow.\passive \) or \( \GermGrow.\passive \)
of Algorithms~\ref{alg:Grow.passive}, \ref{alg:Germ-grow.passive}
so they encode a latent big cell.
This is in accordance with Condition~\ref{cond:comput}\ref{i:comput.legal-birth}, as the
rule \( \Birth \) of Algorithm~\ref{alg:Birth} requires that a cell born from a vacant one is latent.
The trace-back also shows that the germ growth succeeded only if the created colony
is at a distance \( \ge 2\Q\B \) from existing big cells, as requires the definition of emergence.
\end{proof}

\begin{lemma}\label{lem:update-loc-maint}
Under the same condition on \( \Damage^{*} \) as in Lemma~ \ref{lem:legality},
  assume that we have a full colony in which all cells have ages before
  the start of rule \( \UpdateLocMaint \) (Algorithm~\ref{alg:Update-loc-maint}).

 \begin{cjenum}

 \item\label{i:update-loc-maint.damage}
The locally maintained fields \( \Doomed \) and \( \Growing_{j} \)  will be
homogeneous with the exception of an interval of 
at most \( \numDirAff \) cells.

  \item The value of \( \Doomed \) is true (almost) everywhere 
 if and only if \( \Hold \) represents the vacant state \( \Vac^{*} \).

  \item For \( j\in\{-1,1\} \), the value of \( \Growing_{j} \) is equal (almost)
everywhere to \( \Creating_{j}^{*} \) of the cell state represented by the colony.

 \end{cjenum}
  \end{lemma}
 
    \begin{sloppypar}
  \begin{proof}
The statement~\ref{i:update-loc-maint.damage}
holds by the same reasoning as that of
part~\ref{step:legality.non-free} of the proof of Lemma~\ref{lem:legality},
because being part of the rule \( \Compute \), rule \( \UpdateLocMaint \) 
is also repeated three times.
The effects of later damage on these fields
will be corrected via the rules \( \Heal \) and \( \LocMaintain \);
when damage happens at the end of the work period, its effect is still limited to the location
of the damage rectangle.
The rest of the statement also follows from the error analysis of the \( \Compute \) rule.
 \end{proof}
 \end{sloppypar}

The lemma below follows easily from the above and from 
Condition~\ref{cond:time-mark} (Time Marking).

 \begin{lemma}\label{lem:dead-if-long}
  Assume that, for some number \( a \), \( \Damage^{*} \) does not
intersect \( \{x\}\times\rint{a-\Tls/2}{a+2\Tus} \).
  If \( \eta^{*}(x,\cdot) \) has no switching time during \( \rint{a}{a+2\Tus} \) then
\( \eta^{*}(x,a+) \) is vacant.
 \end{lemma}

  The following lemma infers about the present, not only about the past as
the Attribution Lemma.
For an 
island \( \lint{a_{0}}{a_{1}}\times\rint{u_{0}}{u_{1}} \), we call the rectangle
\begin{align}\label{eq:wake}
 \lint{a_{0}}{a_{1}}\times \rint{u_{0}}{u_{0}+4\tau_{2}}   
\end{align}
 its \df{healing wake}\index{damage!rectangle!healing wake}.

 \begin{lemma}[Present Attribution]\label{lem:present-attrib}
 \index{lemma@Lemma!Present Attribution}
 Assume that the live cell \( c_{0}=\pair{x_{0}}{t_{0}} \) in colony \( \cC(z_{0}) \)
 is not a germ. 
  Assume also that \( \cC \) at time \( t_{0} \) does not intersect the healing
wake of any island.
  Then one of the following cases holds.
 \begin{djenum}

  \item\label{i:present-attrib.center}
  \( c_{0} \) is a member cell, attributed to \( \cC \) which is covered by a domain
of member cells at time \( t_{0} \).
If \( \Q  <\Age(c_{0}) <\U -1-\Q  \) then every encoded package of this colony
has at most 4 errors.

  \item\label{i:present-attrib.new}
  \( c_{0} \) is a member cell from which a path of time projection at most
  \( \Q \tau_{2} \) leads back to a growth cell in \( \cC \).
  Assume, say, that it is a left growth cell; then it can be attributed to \( \cC+\Q \B  \).
  At time \( t_{0} \), if \( \cC+\Q \B  \) does not intersect the healing wake of an
island then \( \lint{x_{0}}{z_{0}+2\Q \B } \) is covered by a multi-domain.
If \( \Q  <\Age(c_{0}) <\U -1-\Q  \) then still every encoded package of this colony
has at most 4 errors.

  \item\label{i:present-attrib.ext}
    \( c_{0} \) is an outer cell, attributed to its originating colony, say \( \cC+ \Q  \B  \).
    If \( \cC+\Q \B  \) does not intersect the healing wake of an
    island then \( \lint{x_{0}}{z_{0}+2\Q \B } \) is covered at time \( t_{0} \)
    by a domain.

  \item\label{i:present-attrib.vanish}
  \( c_{0} \) is a member cell, and there is in \( \cC \) a backward path from
it, with time projection \( \le 2\splitT+(\Q +1)\tau_{2} \),
to a domain of doomed member cells covering \( \cC \).

 \end{djenum}
 \end{lemma}
  \begin{Proof}\
    \begin{step+}{present-attrib.old-member}
Suppose that \( c_{0} \) is a member cell.
  \end{step+}
  \begin{prooof}
    \begin{sloppypar}
      Recall the notation \( K,E_{0} \) from~\eqref{eq:J-K} and
Definition~\ref{def:attrib} (Attribution).
By Lemma~\ref{lem:attrib} (Attribution), \( c_{0} \) can be attributed to the originating
colony \( \cC \) which is either the colony of \( c_{0} \) or not.
Let time \( t_{1}\in K\setminus E_{0} \) be a time when \( \cC \) was covered by member cells.
Using Lemma~\ref{lem:colony-trace-back} (Colony Trace-back), 
we can go back from \( t_{1} \) to a time \( t_{2} \)
before age \( \computeStart \) in \( \cC \).
Then as in part~\ref {step:legality.non-free} of the proof of Lemma~\ref{lem:legality}
(Legality), we can follow the development of the colony forwards and see that it
forms a continuous domain together with its extension.
The locations containing encoded information have at most
4 errors---as seen in the same proof.      
    \end{sloppypar}

If the computation results in a nonvacant value for the represented big
cell then, if \( \cC=\cC(z_{0}) \) then case~\ref{i:present-attrib.center} holds.
The condition \( \Age>\Q \) guarantees that all cells of the colony belong to the same
work period.
Otherwise the represented field \( \Creating_{j} \) of the colony will be broadcast
into the field \( \Growing_{j} \) of its cells, as shown in 
Lemma~\ref{lem:update-loc-maint}.
The homogeneity of this latter field will be maintained by the healing rule and \( \LocMaintain \).
Thus, depending on the value of \( \Creating_{j} \) of the big cell, growth
will take place and the growth forms a continuous domain with the
originating colony until the age \( \U \) when growth cells turn into members,
and we have case~\ref{i:present-attrib.new}.
By the property of the growth rule
the colony encodes a latent big cell, so the statement about errors also holds.

Suppose that the computation results in a vacant value.
Then \( \Growing_{j} \) will be 0 everywhere but in the healable wake of
the damage.
Growth cannot start accidentally by a temporary wrong value
\( \Growing_{j}=1 \) in an endcell since there is enough time to correct
this value during the long waiting time of \( \rul{Grow} \). 
Also, all cells become doomed.
After \( \Age=1 \), any doomed cell dies, and
the whole colony decays within \( \Q \tau_{2} \) time units.
Before that, the colony is full.
After that, we have case~\ref{i:present-attrib.vanish}.
\end{prooof}

  \begin{step+}{present-attrib.outer}
  If \( c_{0} \) is an outer cell then we have case~\ref{i:present-attrib.ext}.
  \end{step+}
  \begin{pproof}
  Lemma~\ref{lem:attrib} (Attribution) attributes \( c_{0} \) to the originating colony which is
  covered by member cells during \( K\setminus E_{0} \).
  If \( c_{0} \) is a channel cell then this colony would not have time until \( t_{0} \)
  even to finish its computation, so it could not be destroyed.
  If \( c_{0} \) is a growth cell then the
  cell represented by this colony could not have become vacant, because
  then \( \Creating^{*}_{j} \) would have become 0, broadcast as \( \Growing=0 \) and so
  no growth would have occurred.
  
  It forms a continuous domain with its extension, until the age \( \U -1 \).  
  From that age on, they form a multi-domain.
  This could only change if the originating colony goes through a next work
  period and kills itself; however, there is not enough time for this:
  the definitions~\eqref{eq:destr-time}, \eqref{eq:compute-start} show that
  \( \destrT\Q \) is much smaller than \( \computeStart \).
 \end{pproof}

\end{Proof} 

\subsection{Robust media properties}

To prove Condition~\ref{cond:comput}\ref{i:comput.trans}
for big cells, assume that for some \( a \),
 \begin{equation}\label{eq:no-big-damage}
  \Damage^{*}\cap \lint{x-3\Q \B }{x+4\Q \B }\times\rint{a-\Tls/2}{a+3\Tus}=\eset.
\end{equation}
Recall the special switching times \( \sigma_{0},\sigma_{1},\sigma_{2} \)
of Definition~\ref{def:special-sw-times}.
By Lemma~\ref{lem:dead-if-long}, if there is no switching time
during \( \rint{a}{a+2\Tus} \) then \( \eta^{*}(x,a+) \) is vacant.
  
\begin{lemma}\label{lem:retr.trans}
  For a trajectory \( \eta \),
  Condition~\ref{cond:comput}\ref{i:comput.trans} holds for the decoded trajectory \( \eta^{*} \).
  \end{lemma}
  \begin{Proof}    
Let us take times \( v_{2}, v_{1} \) close from below \( \sigma_{2},\sigma_{1} \)
as in the proof of Lemma~\ref{lem:legality} (Legality).
From time \( v_{1} \) on, follow the development of colony \( \cC \).
The computation process can be treated similarly to the proof in 
that lemma; but in the communication with neighbor colonies, we must show that all retrieved
information can be attributed to a single time.
By event~\eqref{eq:no-big-damage}, the whole space-time area in
which this communication takes place contains at most one damage rectangle.
Repeated application of Lemma~\ref{lem:creation} (Creation) guarantees the success of
extending arms via the rule \( \Extend \), defined in Section~\ref{sec:growth},
over possible latent or germ cells or an opposing extension arm,
showing that in each direction, within a good time bound either the extension arm of the
colony \( \cC \)
reaches its neighbor \( \cC' \) (provided it exists) or the neighbor's  reaches it.
If there are extension arms from both colony \( \cC \) and a non-adjacent
neighbor colony \( \cC' \) then the strength ordering gives preference to one
side and therefore soon only (at most) one arm remains on each side.

\begin{sloppypar}
  In what follows, we call the space-time points \df{free} if they are outside
  the wake of the damage rectangle (see~\eqref{eq:wake}).
Cell \( x \) is called free at time \( t \) and time \( t \) is called free at cell \( x \)  if \( (x,t) \) is free.
We must find a time during
the communication period \( \rint{0}{\computeStart} \) of \( \cC \) introduced in
Definition~\ref{def:expansion-period} to which the retrieval from both neighbors
can be attributed.
Let \( a_{1} \) be the last free time when some cell of \( \cC \)
has age \( 0 \), and \( b_{1} \) the first free time when some such cell
has age \( \computeStart \), and let \(  J_{1}=\rint{a_{1}}{b_{1}} \).
By~\eqref{eq:safe-part} and the argument after it, we have 
\begin{align*}
 |J_{1}|>1.5\konst\rep\nP^{*}\Q\Tl,
\end{align*}
while the unsafe parts of the work period (the dwell period of the big cell simulated by \( \cC \))
have a total length \( < 1.5\rep\nP^{*}\Q\Tu \).
Retrieval by rule \( \Retrieve \) in Algorithm~\ref{alg:Retrieve-2} happens
in points of \( J_{1} \) (provided it did not happen earlier) if
a safe \( \Age \)  as defined in~\eqref{eq:Safe} has been seen in both neighbors.
We will find a sufficiently large subinterval of \( J_{1}\) when this condition holds.
Let \( m_{1} \) be the midpoint of \( J_{1} \).
\end{sloppypar}

\begin{step+}{retr.middle-live}
  Suppose that at time \( m_{1} \) the big cell \( x \) encoded by colony \( \cC \) has a left neighbor big cell,
  encoded by a colony \( \cC' \).
  Then there is a time interval \( J_{2}\subseteq J_{1} \)
  of size \(  \ge 0.7\konst\rep\Q\Tl \) in which the \( \Age_{-1} \) read by \( \cC \) in \( \Retrieve \)
  is safe. 
\end{step+}
\begin{pproof}
  It follows from Lemma~\ref{lem:present-attrib} that \( m_{1} \) is contained in a dwell period \( D \)
  of the big cell encoded by \( \cC' \), of length \( \ge\Tls \).
  By~\eqref{eq:safe-part}, \( D \) has a subinterval \( J' \)
  of size at least \(  1.5\konst\rep\nP^{*}\Q\Tl \) in which the \( \Age_{-1} \) read by \( \cC \) is safe.
  Without loss of generality assume that \( J' \) is above \( m_{1} \).
  In the worst case all the unsafe parts of \( D \) may come above \( m_{1} \): still,
  the  part \( J_{2} \) of \( J' \) falling into \( J_{1} \) has a length of at least
  \begin{align*}
   |J_{1}|/2-1.5\rep\nP^{*}\Q\Tu \ge 0.7\konst\rep\nP^{*}\Q\Tl
  \end{align*}
  if \( \konst \) is large enough.  
\end{pproof}

\begin{step+}{retr.middle-dead}
  Suppose that at time \( m_{1} \) the big cell \( x \) encoded by colony \( \cC \) has no left neighbor big cell.
  Then there is a time interval \( J_{2}\subseteq J_{1} \)  of size
  \(  \ge 0.4\konst\rep\nP^{*}\Q\Tl \) in which the \( \Age_{-1} \) read by \( \cC \) is safe. 
\end{step+}
\begin{pproof}
  Let \( J' \) be the largest subinterval of \( J_{1}\) containing \( m_{1} \) in which \( \cC \) has no left neighbor.
  Then either \( |J'|>|J_{1}|/3 \) or there is a subinterval of \( J_{1} \) of size \( \ge |J_{1}|/3 \) that
  either contains a whole work period of a left neighbor of \( \cC \) or is covered by it.
  In each case, similarly to the reasoning in part~\ref{retr.middle-live} above,
  a safe time subinterval \( J_{2} \) is found of size at least 
  \begin{align*}
   |J_{1}|/3-1.5\rep\nP^{*}\Q\Tu \ge 0.4\konst\rep\nP^{*}\Q\Tl
  \end{align*}
  if \( \konst  \) is large enough.
\end{pproof}

We found a time interval \( J_{2} \) of size \( 0.4\konst\rep\nP^{*}\Q\Tl \)
in which both \( \Age \) and \( \Age_{-1} \) are safe.
Repeating the argument we find a time interval \( J_{3}\subset J_{2} \)
of size \( 0.1\konst\rep\nP^{*}\Q\Tl \) in which also \( \Age_{1} \) is safe.
  This shows that the \( \Retrieve \) rule will succeed at some time \( t' \)
  during the work period of \( \cC \) and the simulation result can be attributed to this time point \( t' \) just
  as Condition~\ref{cond:comput}\ref{i:comput.trans} requires.
\end{Proof} 

\begin{lemma}[Growth]\label{lem:growth}\index{lemma@Lemma!Growth}
  For a trajectory \( \eta \), Condition~\ref{cond:comput}\ref{i:comput.creation}
  holds for the decoded trajectory \( \eta^{*} \).
  \end{lemma}
\begin{Proof}
This proof assumes that in case of conflict, growth to the right is preferred, but this choice 
is clearly arbitrary.    
As in Condition~\ref{cond:comput}\ref{i:comput.creation},
we assume that there is no non-vacant \( \eta^{*}(y,t) \) with
\( 0<|y-x|<\Q \B  \), \( t\in\clint{a}{a+3\Tus} \),
and that for every \( t \) in \( \clint{a}{a+3\Tus} \) at least one of
\( \eta^{*}(x-\Q \B ,t) \) and \( \eta^{*}(x+\Q \B ,t) \) is a potential
creator of \( x \).
\begin{step+}{growth.to-right}
  Suppose  \(   x-\Q \B \) is a potential creator in \( \eta^{*} \) at some time
  \( t\in\clint{a}{\Tus} \).
\end{step+}  
\begin{prooof}
  Tracing backward and forward the evolution of the colony of big cell \( x-\Q\B \),
  we find a work period \( \rint{t_{1}}{t_{2}} \) containing \( t \).
  If the computation does not kill the big cell \( x-\Q\B \) then it will create
  the big cell \( x \).
  Indeed, the only seeming obstacle to the growth to the right could be some
  non-germ cell \( z \) in its way.
  But Lemma~\ref{lem:present-attrib} (Present Attribution) shows that such a
  \( z \) must be a
  left extension cell of a live big cell \( y \) in \( \lint{x+\Q \B }{x+2\Q \B } \), and is
  therefore not stronger than the right growth it is preventing.

  Suppose that the computation kills the big cell \( x-\Q\B \).
  Then \( x-\Q\B \) could not become again a potential creator for time of length
  \( \ge\p_{0}\Tl^{*} \).
  By our assumption then by the time \( t_{2}<t+\Tu<2\Tu \), the big cell \( x+\Q\B \)
  must be a potential creator.
  As above, we would find that it has a whole work period \( \rint{t'_{1}}{t'_{2}} \)
  containing \( t_{2} \), and its computation cannot kill it because then
  neither \( x-\Q\B \) nor \( x+\Q\B \) would be a potential creator at time
  \( t'_{2} \).
  So \( x+\Q\B \) would attempt to create \( x \): let us show that it will succeed
 before time \( a+3\Tu \).
  As above, an obstacle to this could only be a non-germ cell \( z \) that is
  the right growth cell of a big cell \( y \) whose body overlaps that of
  \( x-\Q\B \).
  But this big cell \( y \) must be new and thus not become a creator for
  a time of length \( \ge\p_{0}\Tl^{*} \).
\end{prooof}

\begin{step+}{growth.to-left}
  Suppose  \(   x-\Q \B \) is not a potential creator in \( \eta^{*} \) at any time
  \( \clint{a}{\Tus} \).
\end{step+}  
\begin{prooof}
  Then \(   x+\Q \B \) is a potential creator in \( \eta^{*} \) at time \( a \).
  As above, we find a work period \( \rint{t_{1}}{t_{2}} \) of the colony
  \( \cC(x+\Q\B) \) containing \( a \).
  The computation does not kill the big cell \( x-\Q\B \) because then
  at time \( t_{2} \) neither \( x-\Q\B \) nor \( x+\Q\B \) would be a
  potential creator.
  So \( x+\Q\B \) will extend a growth arm to the left, which can only be stopped by some
  non-germ cell \( z \) in the way.
  As in the above reasoning, \( z \) must be a
  right growth cell of a live big cell \( y \) in \( \rint{x-2\Q \B }{x} \).
  Because right growth is preferred, this would succeed, creating a cell \( y \).
  If \( y=x \) then \( x \) has become non-vacant, so we are done.

  Let us show that other cases are not possible.
  We cannot have \( y=x-\Q\B \) because we assumed that
  \(   x-\Q \B \) is not a potential creator during this time.
  The body of \( y \) cannot overlap that of \( x \) (without being equal to it)
  because this was excluded by the original assumption.
  In the remaining case, \( y \) creates big cell \( y+\Q\B \) which overlaps
  the body of \( x \) while not equal to it, but this is again just what
  has been excluded.
  \end{prooof}
\end{Proof}

\subsection{The amplifier parameters}

Here we prove Lemma~\ref{lem:amp} (Amplifier), and with this,
as noted after the statement of that lemma,
we finish the proof of Theorems~\ref{thm:1dim.nonerg.var}
and~\ref{thm:1dim.stor.var} for the case of infinite space.
This lemma says that a uniform amplifier complex
can be built with the parameters defined in Section~\ref{sec:amp},
with large enough \( \rep \).
The main ingredient is the sequence of media \( M_{k} \)
as in~\eqref{eq:M_{k}} and the sequence of codes \( \fg_{k*} \), \( \Phi^{*}_{k} \),
defined in the course of the proof.
It remains to verify the properties of the amplifier complex listed in Lemma~\ref{lem:amp}.
Let us recall them here.
  \begin{cjenum}

  \item \( (\Payload^{k}) \) is an aggregated
field for \( (\fg_{k*}) \), as defined in Section~\ref{sec:amp}.

  \item The damage map of the simulation \( \Phi_{k} \) is defined as in
Section~\ref{sec:rob.damage}.

  \item \( \trans_{k} \) carries \( \PlTrans_{k} \) as in Definition~\ref{def:carries}.

  \item \( \Phi_{k} \) has \( \eps''_{k} \)-trickle-down.

 \end{cjenum}
  
  The first two properties follow immediately from the definition of
\( \trans_{k} \) and the code, given in the preceding sections.
  It has also been shown, by induction, that \( M_{k+1} \) is a robust medium
simulated by \( M_{k} \) via the code, with \( \eps_{k} \) as the error bound and
\( \Tl_{k} \), \( \Tu_{k} \) as the work period bounds.
In the construction and the proof we relied implicitly on the properties
proved in Lemma~\ref{lem:frame-conds}, mostly expressed as
inequalities.

To prove the two other properties
it needs to be shown yet that \( \eps'_{k} \) indeed serves as the bound in
the definition of \( M_{k}=\Rob(\cdots) \), and that the simulation has
\( \eps''_{k} \)-trickle-down.

  \( \eps'_{k+1} \) must bound the difference from 0.5 of the probability of
the new coin-toss of the simulated computation, and the simulation in the
work period must be shown to have \( \eps''_{k} \)-trickle-down.
  We must show that in a big cell transition, the probability that the
\( \Rand^{*} \) field is 1 is in \( \clint{0.5-\eps'_{k+1}}{0.5+\eps'_{k+1}} \).
  (We also have to show that the bounds on the probabilities are of the
form of sums and products as required in the definition of canonical
simulation, but this is automatic.)

  In case there is no island during the whole work period, the
field \( \Rand^{*} \) was computed with the help of the rule \( \rul{Randomize} \).
  This rule took the value \( X \) found in field \( \Rand \) of the base cell of
the colony and copied it into the location holding \( \Rand^{*} \).
  By the property of \( M_{k} \), the probability that \( X=1/2 \) is within
\( \eps'_{k} \) of 0.5.
  By its definition, \( \eps''_{k} \) upper-bounds the probability that any
island intersects the colony work period.
  Therefore the probability that \( \Rand^{*}\ne X \) can be bounded by
\( \eps_{k} \).
  Hence, the probability that the \( \Rand^{*} \) field is 1 is in
 \[
  \clint{0.5-\eps'_{k}-\eps''_{k}}{0.5+\eps'_{k}+\eps''_{k}}=
  \clint{0.5-\eps'_{k+1}}{0.5+\eps'_{k+1}}.
 \]
  As just noted, with probability \( \eps''_{k+1} \), no island
occurs during a colony work period.
Under such condition, the rule \( \RefreshPayload \) and \( \rul{Update-payloas} \)
work without a hitch and each
cell contains the information encoded by the code \( \fg_{*} \) from the state
of the big cell.
So the simulation has the \( \eps''_{k} \)-trickle-down property.

The number of steps in the work period
fits into \( \U'_{k} \) as long as the requirement~\eqref{eq:Ubd} of amplifier
parameters is satisfied.
The parts where this may not be obvious is the coding-decoding part of the program.
However the codes of Example~\ref{xmp:RS-code} we use, for a fixed number of errors,
can be computed in a linear number of algebraic operations (multiplication, division).
Indeed, the sets of equations to be solved involve only constant-size matrices,
and multiplications/divisions were explicitly allowed in our rule language
in Section~\ref{sec:rule-lang}.

 \section{Self-organization}\label{sec:sorg}

In the present paper, self-organization is not a goal in itself but a tool to
achieve reliability without a hierarchical initial configuration.
We achieve this goal via defining a kind of amplifier that while working is
creating more and more higher-level cells.

\subsection{Color control}\label{sec:color-control}

  Consider a robust medium
 \begin{equation}\label{eq:rob-med}
  M = \Rob(\trans, B,\Tl,\Tu,\eps,\eps',r).
\end{equation}
The field \( \Color \) will play an important role in self-organization.
Recall Definition~\ref{def:fit}, with two variants if the fitting relation between colors:
\ref{i:fit.sorg} and~\ref{i:fit.comp}.
From now on, we will focus on variant~\ref{i:fit.comp}, as
all reasoning can be transferred to the other variant, with simplifications.


\begin{notation}
For an interval \( I=\lint{a}{b} \) we define its \( d \)-neighborhood as%
 \glo{g.greek.cap@$\Gamma(I, d)$}%
 \begin{equation}\label{eq:Gg-interval}
  \Gamma(I, d)=\rint{a-d}{b+d}.
 \end{equation}
 We will also write $\Gamma(x,d)=\rint{x-d}{x+d}$.
\end{notation}

\begin{definition}[Color control]\label{def:color-controlling}\index{color control}
  Suppose we are given a medium \( M \), along with a history \( \eta \).
  A set \( E \) of space is \df{color-fitted} in \( \eta \) at time \( t \)
  if  \( \Color(x,t)\le\Color(y,t) \) for every pair of
  cells \( x<y \) whose body intersects \( E \).
  A space-time set is color-fitted if it is color-fitted at each time \( t \).
  Let 
\begin{align}\label{eq:J_t}
 J_{t}=\rint{t-(\p_{0}+2)\Tu}{t} .
\end{align}
The set \( E \) is \df{controlled} in \( \eta \) at time \( t \) if every subinterval \( K \) of
it of size \( \B \), is within distance \( <2\B \) of the body of a cell at that time,
and the rectangle \( K\times J_{t} \) intersects the (space-time) body of a cell.
It is \df{color-controlled} if in addition, the set \( E\times J_{t} \) is color-fitted.
It has color \( c \) if it is color-controlled, and each cell whose space-time
body intersects the set \( E\times J_{t} \) has color \( c \).
\end{definition}

The following lemma shows that if an interval of size \( 2\B \) contains a cell then
(in the absence of damage) it will never remain empty long.
Indeed, the cell can only be killed by other cells trying to create a new cell; this creation
fails only if others are in the way, and so on.
The \( d \)-neighborhood \( \Gamma(I, d) \)
of an interval \( I \) was defined in~\eqref{eq:Gg-interval}.

 \begin{lemma}\label{lem:lasting-control}
  Let \( t_{1}<t_{2} \) be times.
  Assume that interval \( K=\lint{x_{0}}{x_{0}+\B} \)
  is controlled at time \( t_{1} \), and
\begin{align*}
 \Gamma(K, 3\B)\times\rint{t_{1} - (\p_{0}+2)\Tu}{t_{2}}  
\end{align*}
is damage-free.
  Then \( K \) is controlled at time \( t_{2} \).
 \end{lemma}
\begin{Proof}

  Let \( I=\rint{x_{0}-\B}{x_{0}+\B} \), then a cell's body intersects \( K \)
  if and only if it is in \( I \).
  Let \( t\in\rint{t_{1}}{t_{2}} \).
  Suppose that \( K \) is controlled for all \( t'<t \),
  that is for all \( t'<t \),  each set \( I\times \rint{t'-(\p_{0}+2)\Tu}{t'} \) contains some cell.
 Let \( t_{0}=t-(\p_{0}+2)\Tu \).
  We want to show that
  \begin{cjenum}
  \item\label{i:color-control.1}
    \( I \) contains some cell \( z \) at  some time in \( \rint{t_{0}}{t} \),
  \item\label{i:color-control.2}
    at each such time the interval \( \Gamma(K,2\B) \) intersects a cell body.    
  \end{cjenum}
  The statement~\ref{i:color-control.1}
  may not hold only if a cell \( x_{1} \) in \( I \)
  disappears at time \( t_{0} \), so assume this happens.
  Then \( x_{1} \) must have been erased by 
  some cell \( y_{1} \) about to create an adjacent neighbor whose body
  overlaps with the body of \( x_{1} \).
  Without loss of generality assume
\begin{align}\label{eq:lasting-control-0}
  x_{1}\le x_{0},\quad  y_{1}+\B< x_{1}<y_{1}+2\B.
\end{align}
As \( x_{1}>x_{0}-\B \), we have 
\begin{align}\label{eq:lasting-control-1}
  x_{0}-3\B<y_{1}<x_{1}-\B\le x_{0}-\B .
\end{align}
If \( y_{1}\in  I \) then let \( z\gets y_{1} \), now assume it is not.
According to the rule \( \Adapt \), cell \( y_{1} \) must have had \( \Dying_{0}=0 \),
  \( \Creating_{1}=1 \) to be able to erase,
  and therefore survives until time \( t_{0}+\p_{1}\Tl \) as \( \Die(\p) \)
  always has \( \p\ge\p_{1} \) by Condition~\ref{cond:anim,kill}.
  Since according to Condition~\ref{cond:anim,kill} only the rule
\( \Create \) changes \( \Creating_{1} \), and only when a right adjacent neighbor has
disappeared, cell \( y_{1} \) keeps trying to create a right neighbor.

\begin{step+}{step:lasting-control.success}
Suppose \( y_{1} \) succeeds in creating \( y_{1}+\B \) before time \( t_{0}+2\Tu \).
\end{step+}
\begin{prooof}
  If \( y_{1}+\B>x_{0}-\B \) then set \( z\gets y_{1}+\B \).
  Suppose \( y_{1}+\B\le x_{0}-\B \).
  Then \( y_{1}+\B \) turns on \( \Creating_{1} \) within time \( \p_{0}\Tu \),
  and then tries to create \( y_{1}+2\B \).
  If it succeeds then, by~\eqref{eq:lasting-control-1} we can set \( z\gets y_{1}+2\B \).
  If it fails then a cell \( x_{2} \) with \( y_{1}+2\B<x_{2}<y_{1}+3\B \) is in the way,
  and we can set \( z\gets x_{2} \).
\end{prooof} 

\begin{step+}{step:lasting-control.failure}
  Suppose that \( y_{1} \) does not succeed creating \( y_{1}+\B \) before time \( t_{0}+2\Tu \).
\end{step+}
\begin{prooof}
  Then some cell \( x_{2} \) with  \( y_{1}+\B<x_{2}<y_{1}+2\B \) interferes.
  This cell was not there at time \( t_{0} \) as, by~\eqref{eq:lasting-control-0}, it
  would have overlapped with cell \( x_{1} \).
  Being so close to \( y_{1} \), cell \( x_{2} \) could not have
  arisen by spontaneous birth, hence it must have been created by \( x_{2}+\B \)
  which, as any creating cell,  must still be alive after the creating time,
  so we can set \( z\gets x_{2}+\B \).
\end{prooof} 
We have shown~\ref{i:color-control.1}.
The proof also shows that unless \( x_{1} \) survives \( y_{1} \) will, therefore for all
\( t' \) in \( \rint{t_{0}}{t} \), it is within \( 2\B \) of interval \( K \), which
also shows~\ref{i:color-control.2}.
\end{Proof}

\begin{lemma}[Lasting color control]\label{lem:lasting-color-control}\index{control!lasting}
Suppose that in a trajectory \( \eta \) of medium \( M=M_{k} \) of our amplifier,
for an intervals \( E \) and times $t_{1}<t_{2}$,  the interval
 \[
   E'=\Gamma(E, \B(t_{2}-t_{1})/\Tl)
 \]
  is color-controlled at time \( t_{1} \) and $I\times \rint{t_{1}-2\Tu}{t_{2}}$ is
  damage-free.
  Then \( E \) is color-controlled at time \( t_{2} \).
  \end{lemma}
 \begin{proof}
We have to prove that \( E \) is color-controlled at time \( t_{2} \).
Lemma~\ref{lem:lasting-control} implies that \( E \) is controlled at time \( t_{2} \); what
remains to show is that the rectangle $\Gamma(E,\B)\times\rint{t_{2}-3\Tu}{t}$ 
is color-fitted.
By the definition of control, in the space-time area considered, every interval
of size \( \B \) is within less than \( 2\B \) distance from some cell body,
and this prevents any germ cell from being born, due to
Condition~\ref{cond:comput}\ref{i:comput.legal-birth}.
Suppose that two non-fitting controlling cells was created: say, a cell with color 1 on
the left of a cell with color 0.
Then one could construct a path from both of these backwards, one with
color 1 and one with color 0.
These paths cannot cross, being of different color,
nor can reach outside \( E' \) during \( \rint{t_{1}}{t_{2}} \) since it needs
at least \( \Tl \) time units to move one cell width away from \( E \).
So we would find a cell with color \( 1 \) on the left of a cell
with color 0 in \( E' \) at time \( t_{1} \), contradicting the assumption.
\end{proof}

An interval will become controlled soon even if a cell is only nearby:

 \begin{lemma}\label{lem:arising-control}
   \( K=\lint{x_{0}}{x_{0}+\B} \) be a space interval and \( t_{1},t_{2} \)
   times with \( t_{2}=t_{0}+d\p_{0}\Tu \),  where
   \begin{align}\label{eq:d-bound}
 d < \p_{1}/\lambda\p_{0} .
\end{align}
   Assume that
\( \Gamma(K, (d+1)\B)\times\rint{t_{1} - 3\Tu}{t_{2}} \) is damage-free, and
the body of a cell \( x_{1} \) intersects \( \Gamma(K,d\B) \) at time \( t_{1} \).
  Then \( K \) is controlled at time \( t_{2} \).
\end{lemma}
\begin{proof}
  The proof is similar to that Lemma~\ref{lem:lasting-control}, so we only sketch it.
  Assuming \( x_{1} \) is the closest cell to \( K \), say from the left,
  it extends an arm of cells to the right, where each step of extension takes
  time \( \le\p_{0}\Tu \), for a total of \( \le d\p_{0}\Tu \).
  This can be stopped only in two ways.
  1: if a cell \( x_{2} \) closer to \( K \) emerges as an obstacle,
  in which case we continue the reasoning from \( x_{2} \),
  2:  If \( x_{1} \) gets killed by a cell \( y_{1} \) from the left, in which case
  we continue the reasoning from \( y_{1} \).
  We know that \( y_{1} \) cannot die for a time of at least \( \p_{1}\Tl \),
  and by the assumption~\eqref{eq:d-bound} we have \( d\p_{0}\Tu<\p_{1}\Tl \).
\end{proof}

\begin{definition}\label{def:germ-level}
  In a trajectory \( \eta \) of a medium \( M \), we will say that a space interval
  at time \( t \) is \df{germ-level}, if it has no colony cells.
\end{definition}

According to the lemma below, if only a small subinterval has possibly different color
in an otherwise color-controlled interval then the growth of this discrepancy in time
is limited.

\begin{lemma}\label{lem:robust-color-control}
  Let us be given a time \( t_{1} \), an interval \( J \), and
  \begin{align*}
    I\supset \Gamma(J,|J|), \quad  I'=\Gamma(I, \B^{*}),\quad
    t_{2}=t_{1}+\Tl^{*}.
\end{align*}
  In a trajectory \( \eta \), assume that the space-time rectangle
\(  I'\times\rint{t_{1}}{t_{2}} \) is free of damage.
Assume also that \( \eta(\cdot,t_{1}) \) is such that changing it in \( J \),
at time \( t_{1} \)  the whole interval \( I' \) will become color-controlled
and germ-level.
Then at time \( t_{2} \) there is an interval \( J'\supset J \), with
size \( \le 2|J| \) such that changing \( \eta \) in \( \Gamma(J,m\B) \),
the interval \( I \) becomes color-controlled.
If it contains any colony cells at this time, they encode cells of age at most \( 0 \).
\end{lemma}
\begin{proof}
  Let \( m=|J|/\B \).
  If there are no cells in \( J \) then repeated application of
  Lemma~\ref{lem:arising-control} implies
  that it will be populated with cells coming from outside it by
  time \( t_{1}+m\p_{1}\lambda\p_{0}\Tu \).
  Assume now that there are some cells in \( J \).
  Any colony cells  in \( J \) can only survive if they become
  part of a colony connected to either the left side or the right side,
  otherwise they will decay within time \( m\p_{2}\Tu \), as they are
  in a small interval at least one of whose ends is exposed.
  So assume we have a germ inside \( J \), of a color that is not fitting to
  its neighbors outside, say, on the left end.
  Then the left end of this germ cannot extend beyond the left end of \( J \).
  The distance to which the right end can extend is limited by the
  rule \( \GermGrow.\rActive \)
  of Algorithm~\ref{alg:Germ-grow.active} to \( m\B \).
  After the age \( \growEnd \), the germ starts decaying, and will disappear
  by time \( t_{1}+\Tl^{*} \).
  During this time, colony cells may intrude into the interval \( I' \), but
  in time \( \Tl^{*} \) they can enter to a distance at most \( \B^{*} \).

  Any colony cells in \( J \) that arose during the time indicated, belong to a brand new
  colony and as such, can only encode a big cell of age 0 that was not there at time \( t_{1} \).
  The time \( \Tl^{*} \) is not sufficient to increase this age even by 1.
\end{proof}
 
In some simulations, color ``trickles down''.

\begin{lemma}[Color trickle-down]\label{lem:color-trickle}\index{color!trickle}
  Let $\eta$ be a trajectory of $M$, define the times \( t_{0} \) and \( t_{1}=t_{0}+\Tus \).
  For interval $J$ let \(   J'=\Gamma(J, \lambda\B^{*}) \).
  If $J'$ is color-controlled at time \( t_{0} \) in $\Phi^{*}(\eta)$, and the area
$J'\times\rint{t_{1}-(\p_{0}+3)\Tus}{t_{1}}$ is damage-free in $\eta$, then
$J$ is color-controlled at time \( t_{1} \) in $\eta$,
and its colors can only be ones that were present in \( J \) at time \( t_{0} \).
  \end{lemma}
  \begin{proof}
Let \( K \) be any subinterval of \( J \) with \( |K|=\B \), and let
\( K'=\lint{a}{a+\Q\B} \)
be any subinterval of \( J' \) of size \( \Q\B \) containing \( K \).
As \( J' \) is color-controlled at time \( t_{0} \), the interval \( K' \) is intersected
by the body of a big cell \( x_{1} \) at some time \( t' \)
with \( t_{0}-(\p_{0}+2) \Tus< t'\le t_{0} \).
Without loss of generality assume that \( x_{1}\le a \).
By Lemma~\ref{lem:present-attrib} (Present Attribution) we can find a whole work
period \( \rint{u_{1}}{u_{2}} \) of \( x_{1} \) containing \( t' \) and
ending before time \( t_{1} \).
During this work period its colony will extend some arm to the right.
If it does not cover \( K' \) then again, by the Present Attribution lemma we would
find the left extension arm of another colony meeting this arm from the right.
This way, at some time \( t''\in\rint{u_{1}}{u_{2}} \) the whole interval \( K' \)
is covered by up to two intervals of adjacent small cells, meeting in a gap of size \( <\B \).
So then the interval \( K' \) is color-controlled in \( \eta \), and
by Lemma~\ref{lem:lasting-color-control}, this property will be conserved until time \( t_{1} \).
 \end{proof}

 \subsection{Colony birth} \label{sec:colony-birth}

  We will consider a trajectory \( \eta \).
  The final goal in this section is to find, in Lemma~\ref{lem:birth},
  with large probability, some germ will grow into a big cell.

 Recall the definition of \( \grFac \) in~\eqref{eq:germ-signif-larger}, and its use in
 the rule of germ growth.
 In what follows we will repeatedly rely on the following observation,
 for any \( \grFac>1 \).
 Let \( 0<L_{1}<\dots<L_{n} \) be such that
 \( L_{i+1}\ge\grFac L_{i} \) for all \( i \).
 Then 
\begin{align}\label{eq:geom-sum}
 L_{1}+\dots+ L_{n} \le \frac{\grFac}{\grFac-1}L_{n} .
\end{align}
  (We will use the fact that with the
  definition~\eqref{eq:germ-signif-larger} of \( \grFac \) 
  we have \( \frac{\grFac}{\grFac-1}=5 \).)

  \begin{lemma}\label{lem:germ-unexposed}
    \begin{sloppypar}
    In a trajectory \( \eta \), let \( c_{1}=\pair{x_{1}}{t_{1}} \) be a germ cell with
\( \GermSize = L \) and \( \Age>3\rep L \).
Assume that at time \( t_{1} \), with a germ cell body \( K=\lint{x_{1}}{x_{1}+2\B^{*}} \),
the space-time rectangle
  \begin{align}\label{eq:large-rect}
    \Gamma(K, 25\B^{*}) \times \Gamma(t_{1}, \Tus).
  \end{align}
is free of damage, and the domain containing \( c_{1} \) has an exposed edge.
Then there is an \( L' \) with \( \grFac L< L'\le 5\Q \) such that in 
\begin{align}\label{eq:germ-unexposed-area}
  \Gamma(x_{1}, 5 L' \B)
  \times\lint{t_{1}}{t_{1}+5\rep L'\tau_{2}}
\end{align}
there is a colony cell, or a germ cell with
\( \GermSize\ge L' \), or one with \( \Addr=e'_{j} \) for \( j=\pm 1 \)
(recall~\eqref{eq:e'_j}).
  \end{sloppypar}
\end{lemma}
When \(  \grFac L \) would be larger than \( 5\Q \) 
then \( \GermSize\ge \grFac L \) is not possible,
so only the case of a colony cell remains.
\begin{proof}
  Let us note that if there is a colony cell then by
  Lemma~\ref{lem:present-attrib} (Present Attribution) there is a whole extended
  colony containing it.
  
    Let us construct a path \( P_{1} \) backwards from \( c_{1} \), staying inside the germ,
    which we will call \( G_{1} \), to
    of an age that is at the beginning of the computation part which set \( \GermSize(c_{1}) \).
    Germ \( G_{1} \) must be without exposed edge during the first
    half of this computation part, because otherwise it
    would decay during the second half (which is long enough for this), contradicting
    the age of \( c_{1} \).
    Then at the end of the computation, at some time \( t_{0} \),
    the germ \( G_{1} \) has size \( \ge L \), and all its cells have the same \( \Dominant \)
    value, decided in the first step of the computation.

    Let us follow the development of \( G_{1} \) from time \( t_{0} \) forwards.
    At some time one of its edges gets exposed, so some
    protected edge \( c' \) was killed by an attack: without loss of generality,
    assume it is from the left, and consider the development of the
    attacker domain \( G_{2} \) after time \( t_{0} \).
    There are the following possibilities according to Algorithm~\ref{alg:Germ-grow.passive}:

    \begin{djenum}
    \item\label{i:unexposed.colony}
      \( G_{2} \) is an extended colony.      
    \item\label{i:unexposed.largest-germ}
      \( G_{2} \) has \( \GermEdge_{-1}=e'_{-1} \), and then
      naturally \( \GermSize \ge 5\Q/2 \).
    \item\label{i:unexposed.larger-germ} it has \( \GermSize > \grFac L \).
    \item\label{i:unexposed.attacking-germ}
      it has \( \Dominant=1 \) and \( \GermSize\ge L/\grFac \) while \( G_{1} \)
      has \( \Dominant=0 \).
    \end{djenum}

    In case~\ref{i:unexposed.colony}, this colony is still there at time \( t_{1} \), because
    the period \( \rint{t_{0}}{t_{1}} \) is too short for it to go through another work period
    and die.
    In the other cases, \( G_{2} \) must also be past the computation that led to its
    \( \GermSize \).
    Consider case~\ref{i:unexposed.largest-germ}.
    If \( G_{2} \) does not get exposed before time \( t_{1} \) then we are done.
    If it is exposed by an extended colony, we are back to case~\ref{i:unexposed.colony}.
    If it is exposed by a germ then we are back to cases~\ref{i:unexposed.larger-germ}
    or~\ref{i:unexposed.attacking-germ}, but with a germ size at least \( \grFac \) times
    larger than that of \( G_{1} \).
    Repeating this we may get to a distance \( L+\grFac L + \dots \grFac^{n}L \).
    With \( L'=\grFac^{n}L \), by~\eqref{eq:geom-sum} this distance is
    \( \le \frac{\rho}{\rho-1}L' \).

    Consider case~\ref{i:unexposed.attacking-germ}.
    If \( G_{2} \) gets exposed
    then we can repeat the whole above reasoning starting with
    \( G_{2} \); however, now it is attacking, so case~\ref{i:unexposed.attacking-germ}
    cannot occur again with the same \( \GermSize \).
    Now assume that \( G_{2} \) does not get exposed: call this the case X.    
    Then it can override \( G_{1} \), but
    it is possible that at the same time, case X occurs also on the right.
    Still, as a result, within time \( \tau_{2}\rep L \),
    at least half of the size of \( G_{1} \) will be added to
    one of the attacking germs, say \( G_{2} \),
    so after it finishes its growth (and new computation), it gets 
\begin{align*}
 \GermSize\ge L(1/2+1/\grFac)>L \grFac.
\end{align*}
  \end{proof}

  The following lemma strengthens the conclusion of Lemma~\ref{lem:germ-unexposed}.

  \begin{lemma}\label{lem:germ-unexposed-later}
Consider a trajectory \( \eta \)
with a germ cell \( c_{1}=(x_{1},t_{1}) \), \( K=\lint{x_{1}}{x_{1}+2\B^{*}} \)
and \( L=\GermSize(c_{1}) \),
and assume that the large space-time rectangle~\eqref{eq:large-rect} is free of damage.
Let \( 5\tau_{2}\rep L\le H\le \Tus \).
Then at some time \( t'\in\clint{t_{1}}{t_{1}+H} \) we can find a cell \( x' \)
with a domain of germ size \( L'\ge L \), where
  \begin{align}\label{eq:germ-unexposed-later}
    x' \in \Gamma(x_{1},5 L'\B).
\end{align}
If \( c'=(x',t') \) does not belong to an extended colony or to a  germ with maximal size
then \( t'=t_{1}+H \).
If \( L'\le\grFac L \) then \( c'=c_{1} \), and during \( \clint{t_{1}}{t_{1}+H} \)
the germ of \( c_{1} \) did not become exposed or have a neighbor that is a colony cell or
a germ cell with \( \GermSize\ge\grFac L \).
\end{lemma}
\begin{proof}
  Let \( t_{0}=t_{1}+H \).
  We will iterate the application of Lemma~\ref{lem:germ-unexposed}, with
  \( L_{1}= L \).
  If the germ of \( c_{1} \) does not get exposed or
  have a neighbor that is a colony cell or a germ cell with \( \GermSize\ge\grFac L \)
  during \( \clint{t_{1}}{t_{1}+H} \) then its germ has size \( \ge L \).
  As we will see this is the only possibility for \( L'<\grFac L \).
  
  Suppose now that one of these events does occur.
  Then Lemma~\ref{lem:germ-unexposed} gives a cell \( c_{2}=(x_{2},t_{2}) \)
  in the area~\eqref{eq:germ-unexposed-area}, without exposed edges,
  with \( \GermSize(c_{2})=L_{2}\ge \grFac L_{1} \),
  but with possibly \( t_{2}<t_{0} \).
  Follow the development of the germ of \( c_{2} \) forward.
  If it reaches its maximum size or 
  time \( t_{0} \) without getting an exposed edge then we are done.
  Otherwise, the germ of \( c_{2} \) gets an exposed edge at time \( t_{3} \).
  Applying Lemma~\ref{lem:germ-unexposed}, to \( c_{3}=(x_{2},t_{3}) \),
  we find a cell \( c_{4}=(x_{4},t_{4}) \) 
that is either a germ end-cell, or has \( \GermSize(c_{4})\ge\grFac L_{2} \).
Repeat the argument as long as necessary.
The relation~\eqref{eq:germ-unexposed-later} can be proved similarly to the corresponding
relation in Lemma~\ref{lem:germ-unexposed}.
\end{proof}

The following key lemma is the only one in this paper using a canonical
simulation with random rectangles as in Section~\ref{sec:canonical},
and then relying on the fact that \( \eta \) is a \emph{strong} trajectory.
Let
  \begin{align}\label{eq:no-birth-prob}
   \theta_{k}=e^{-\rep k^{2}}.
\end{align}

\begin{lemma}\label{lem:birth}
  Consider a trajectory \( \eta \) of \( M_{k} \),
  with a germ cell \( c_{1}=(x_{1},t_{1}) \), \( K=\lint{x_{1}}{x_{1}+2\B^{*}} \),
  and assume that the large space-time rectangle~\eqref{eq:large-rect}
  is free of damage.
  The probability that the rectangle contains no colony cell is bounded by \( \theta_{k} \).
  Moreover, there is a canonical simulation in the sense of Definition~\ref{def:canon-simul}
 providing this estimate.
\end{lemma}

\begin{Proof}
 \begin{step+}{step:birth.constr}
The following construction iterates that of Lemma~\ref{lem:germ-unexposed-later},
with cells \( c_{s}=(x_{s},t_{s}) \), \( L_{s}=\GermSize(c_{s}) \) and
time period of size
\begin{align}\label{eq:H_i}
  H_{s} = 10\lambda\tau_{2}\rep L_{s}
\end{align}
in step \( s \).
 \end{step+}
 \begin{prooof}
Starting with cell \( c_{1} \),
find a cell \( c' \) given by the conclusion of Lemma~\ref{lem:germ-unexposed-later},
call it \( c_{2}=(x_{2},t_{2}) \).
Repeat, getting to cells \( c_{3},c_{4},\dots \),
until cell \( c_{n} \) when
we reach a germ with maximal size, or a colony, or time \( t_{1}+2\Tus \).
Let us call step \( s \) a \df{growth step} if \( L_{s+1}\ge \grFac L_{s} \), and
\df{stalling step} otherwise.
As seen in Lemma~\ref{lem:germ-unexposed-later}, in a stalling step
a germ of \( c_{s} \) of size \( L_{s} \) does not get exposed, nor does it meet a
colony cell or a cell with \( \GermSize\ge\grFac L_{s} \).

Let us examine a stalling step  starting at time \( t_{s} \),
at the shorter end of the germ \( G_{s} \) of cell \( c_{s} \), that is
on the end in which direction its \( \GrowDir \) points; without loss
of generality, assume that it is the left end.
Let \( x_{s} \) be the cell of \( G_{s} \) with address \( \Q/2 \).
After \( t_{s} \), the germ \( G_{s} \) may grow, but since it will not reach \( \grFac L_{s} \),
it will stall soon, at some time \( u_{1} \).
One possibility is that the left end has reached its farthest possible address,
another that another germ is in the way: one with \( \GermSize \) between \( L_{s}/\grFac \)
and \( \grFac L_{s} \); we already excluded the other possibilities.
In case the left end has reached the farthest address,
we will turn attention to the right end, so
without loss of generality, assume that the left end does not reach the farthest address,
because another germ \( G' \) is in the way.

Let \( x' \) be the central cell of \( G' \), and
look at the next \( \lambda+2 \) times in which \( x_{s} \) or \( x' \) make
a random choice of \( \Dominant \).
Let \( \cE \) be the event that 
each choice is \( \Dominant(x_{s})=1 \) and \( \Dominant(x')=0 \).
We claim that if \( \cE \) occurs then 
\( G_{s} \) will advance, so the current step is not a stalling step.
To avoid some case distinctions, call \( x_{s}=y_{1} \), \( x'=y_{0} \).
Let \( j \) be such that the first choice of
\( y_{1-j} \), at some time \( \tau \),  falls between two choices of \( y_{j} \)
at times \( \tau' \) and \( \tau'' \).
If \( \tau \) happens at least \( 0.1\rep L \) age steps before \( \tau'' \)
them \( G_{s} \) will advance before \( \tau'' \).
Otherwise \( G_{s} \) will advance soon after \( \tau'' \).
Now, \( \tau \) occurs before \( \lambda+1 \) choices made by \( y_{j} \);
indeed, \( \Tu/\Tl\le\lambda \) implies that the stages of \( y_{1-j} \) are at most
\( \lambda \) times longer than those of \( y_{j} \).
So the current step is stalling only if \( \neg\cE \) holds, that is with probability
\begin{align*}
  \Prob(\neg\cE) \le 1-(1/2-\eps')^{\lambda+2}< 1-(1/2-\eps')^{5}<1-2^{-5},
\end{align*}
because \( \lambda\le 2 \) according to~\eqref{eq:lg-less-2}, and \( \eps' \) can
be chosen small.

by Condition~\ref{cond:comput}\ref{i:comput.rand} on robust media.
 \end{prooof} 

 \begin{step+}{step:birth.probab-estim}
   Let us estimate first the probability that there is no colony in~\eqref{eq:large-rect}
   informally (that is not in the sense of canonical simulation  yet).
 \end{step+}
 \begin{prooof}
   We assumed that the rectangle~\eqref{eq:large-rect} is damage-free; the straightforward
   estimation of the probability of damage will come later.
   Our construction gives in step \( s \) a cell \( c_{s}=(x_{s},t_{s}) \).
The number of growth steps is at most \( \log(5\Q)/\log\grFac \).
Recalling~\eqref{eq:germ-signif-larger}, this is \( < 4\log\Q \) when \( \Q \) is large.
As \( L_{s}<5\Q \), each  step  has time length \( \le 50\lambda\tau_{2}\Q \).
So the total number of steps is at least
\begin{align}\label{eq:n-lb-1}
   n\ge\Tus/25\lambda\p_{2}\Tu\Q.
\end{align}
We are using the parameters in~\eqref{eq:amp-fr}:
\begin{equation}\label{eq:frame-values}
  \begin{aligned}
  \Q &=\Q_{k}= \rep^{k+1},
\\ \bw &=1/\rep k^{2},
  \\ \U' &= \rep\Q/\bw,
  \\ \Tus &< 2\U'_{k}\Tu < 2\Tus,
  \\      \log\Q &= (k+1)\log\rep.
  \end{aligned}
\end{equation}
so by~\eqref{eq:n-lb-1} we have
\begin{align}
\nonumber  \Tus/\Q\Tu &>\U'/\Q=\rep/\bw=\rep^{2}k^{2},
  \\\label{eq:n-lb}
  n &\ge \rep^{2}k^{2}/25\lambda\p_{2}=:c\rep^{2}k^{2}.
\end{align}
 As we have seen the probability for any particular step \( s \) to be a stalling step is 
\( \le 1-2^{-5} \).
Using \( \ln (1+x)<1+x \) we have \( \ln (1-2^{-5})\le -2^{-5} \),
\begin{align}\label{eq:log-r}
   r := -\ln (1-2^{-5}) > 2^{-5}.
\end{align}
For any particular \( m \)-element subset of steps the probability to be stalling steps
is at most \( (1-2^{-5})^{m}=e^{-r m} \).
As seen, the number of stalling steps is at least \( n-4\log\Q \),
so the probability of no colony is 
\begin{equation}\label{eq:prob-bd}
  \begin{aligned}
  &\le\sum_{i>n-4\log\Q}\binom{n}{n-i}e^{-r(n-i)}
    = \sum_{i<4\log\Q}\binom{n}{i}e^{-r(n-i)}
  \\ &<   4\log\Q \cdot n^{4\log\Q}e^{-r(n-4\log\Q)}
       = 4\log\Q e^{-(rn-4r\log\Q+4\ln n \log \Q)},    
  \end{aligned}
\end{equation}
because the term inside the sum is increasing with \( i \).
Given~\eqref{eq:n-lb} and~\eqref{eq:frame-values},
for large \( \rep \) the right-hand side is \( <4\log\Q e^{-rn/2} \),
and the sum of~\eqref{eq:prob-bd} for all \( n>c\rep^{2}k^{2} \) is
\begin{align*}
    \le 4\log\Q \sum_{n>c\rep^{2}k^{2}}e^{- r n/2}
\le 4(k+1)\log\rep\cdot e^{-r c\rep^{2}k^{2}/2}\sum_{n> 0}e^{-rn/2}.
\end{align*}
The multipliers of \( e^{-r c \rep^{2}k^{2}/2} \) don't change the rate of growth much,
so we can bound the probability by \( \theta_{k} \) as defined in~\eqref{eq:no-birth-prob}.
\end{prooof} 
\begin{step+}{step:birth.canon}
  Now we repeat the construction above, turning it into a canonical simulation.
\end{step+}
\begin{prooof}
  As we constructed the sequence of cells \( c_{s}=(x_{s},t_{s}) \), \( s=1,2,\dots \),
  the times \( t_{s} \) were chosen on the basis of the trajectory \( \eta \) before \( t_{s} \),
  so these times \( t_{s} \) are stopping times according to Definition~\ref{def:stopping-time}.
  It is not important to upper-bound the number of steps, but \( n^{*}=2\Tus/\Tu \)
  will serve.
  Following the reasoning in part~\ref{step:birth.probab-estim},
  choose in every possible way a subset \( S \) of size \( > n^{*}-4\log\Q \) of the set
  \( \{1,\dots,n^{*}\} \) of indices.
  Number these subsets as \( S_{1},S_{2},\dots,S_{N} \).
  Subset \( S_{i} \) represents the assumption that the steps \( s\in S_{i} \) are
  the stalling steps.
  In order to define the canonical simulation, we repeat the reasoning of
  part~\ref {step:birth.constr} above.

  Let us fix one of the sets \( S=S_{i} \).
  For each element \( s\in S \) in increasing order,
  let us examine a \emph{presumed} stalling step  starting at time \( t_{s} \),
  at the shorter end of the germ \( G_{s} \) of cell \( c_{s} \).
  Without loss of generality, assume it is the left end.
Let \( u_{1} \) be the start of a growth part of a new work stage
of \( G_{s} \) after \( t_{s} \) (started by the center cell of \( G_{i} \)).
After \( u_{1} \), the germ \( G_{s} \) may grow.
If it reaches \( \grFac L_{s} \), or meets a colony cell or a germ cell
with \( \GermSize\ge\grFac L_{s} \), then \( s \) is not a stalling step: we set
\begin{align}\label{eq:bail-out}
\beta_{S,s}=\omega,\quad W_{i,r}=\emptyset,   
\end{align}
and proceed to \( j+1 \).
We will call this action \df{bailing out}.
Another possibility is that the left end has reached its farthest possible address,
another that another germ is in the way: one with \( \GermSize \) between \( L_{s}/\grFac \)
and \( \grFac L_{s} \); we already excluded the other possibilities.
In case the left end has reached the farthest address, we will turn attention to the right end, so
without loss of generality, assume that the left end does not reach the farthest address,
because another germ \( G' \) is in the way.
Let \( x_{s} \) be the central cell of \( G_{s} \) and \( x' \) be that of \( C' \).
Now we look at the first \( m=\lambda+2 \) times
\( \tau^{*}_{1},\tau^{*}_{2},\dots \) when \( x_{2} \) or \( x' \)
make random choices, let \( y_{j}\in \{x_{s},x'\} \) be the element making its choice.
For every possible binary sequence \( \bb = b_{1},\dots,b_{m} \) we will bail out
unless \( y_{j}=x_{s} \) if and only if \( b_{i}=1 \).
We will also bail out if the choice is \( \Dominant(y_{j})=j \) for all \( j \).
Both of these decisions can be made at the last time \( \tau_{m} \), if not earlier.
For the remaining cases,
\begin{align*}
\beta_{S,s,\bb,j}:=\alpha(\var{rand},f_{j}),\quad W_{S,s}:=W_{0}(y_{j},\tau^{*}_{j}),
\end{align*}
where \( \alpha(\var{rand},d) \) and \( W_{0}(x,t) \) were
defined in Condition~\ref{cond:comput}\ref{i:comput.rand}.
This construction results in the required canonical simulation.
Summing up the non-bailed out bounds results in the same estimate as the more
intuitive one given in parts~\ref{step:birth.constr} and~\ref{step:birth.probab-estim} above.
\end{prooof} 

\end{Proof}

We now strengthen the conclusion of Lemma~\ref{lem:birth} for the cases when 
more is known about interval we start with.

\begin{lemma}\label{lem:color-lift}
  Let us be given an interval \( K=\lint{x_{1}}{x_{1}+2\B^{*}} \), and let
\begin{align*}
  t_{2}=t_{1}+26\p_{0}\Tus ,\quad
  K'=\Gamma(K,2\cdot (26\p_{0}+2)\lambda\B^{*}).
\end{align*}
In a trajectory \( \eta=\eta^{k} \) of medium \( M=M_{k} \),
assume that the space-time rectangle
\(  K'\times\rint{t_{1}}{t_{2}} \) is free of damage.
Assume also that there is an interval \( J \) of size \( <|\B^{*}|/3 \)
in which when changing \( \eta \) at time \( t_{1} \) 
the whole interval \( K' \) can be made color-controlled in \( M \).
Then at time \( t_{2} \), except for an event of probability \( \theta \)
defined in~\eqref{eq:no-birth-prob},
the interval \( K \) is color-controlled in trajectory \( \eta^{*} \) of \( M^{*} \).
If \( K' \) is germ-level for \( \eta \) at time \( t_{1} \) then \( K \)
is germ-level for \( \eta^{*} \) at time \( t_{1} \).
\end{lemma}
\begin{proof}
  Assume \( \Gamma(J,\B^{*})\subseteq K' \);
  if this is not true then we can simply ignore \( J \).
  Let \( K'' \) be the larger of the two intervals or \( K'\setminus J \).
  By Lemma~\ref{lem:birth}, except for an event of probability \( \theta \),
  before time \( t_{1}+\Tus \) a colony arises in \( K'' \), not
  farther than \( 25\B^{*} \) from \( K \).
  Then applying Lemma~\ref{lem:arising-control}, to \( \eta^{*} \),
  until time \( t_{2} \) the interval \( K \) becomes controlled in \( M^{*} \).
  Indeed, given that by Lemma~\ref{lem:robust-color-control},
  any non-fitting germ that started in \( J \) 
  can grow only by a factor of 2 in each time interval \( \Tl^{*} \), 
  the size of such a germ remains much smaller than \( \B^{*} \), and will
  be overridden by the colonies.

  Cells that are not color-fitted can only intrude over germ cells if they are colony cells.
  In time \( \Tl^{*} \) they can cover at most one colony length, so in time
  \( 26\p_{0}\Tus \), this amounts to \( \le 26\p_{0}(\Tus/\Tl^{*})=26\p_{0}\lambda \)
  colony lengths, still away from \( K \).

  Let us show that \( K \) will be germ-level in \( \eta^{*} \).
  At time \( t_{1} \), there were no colony cells yet in \( K' \),
  so even if a new colony occurs right after \( t_{1} \), it starts simulating a
  germ cell of \( M^{*} \), with age 0.
  Every dwell period of such a cell takes time at least \( \Tl^{*} \), so
  until time \( t_{2} \) it can have at most \( 26\p_{0}\lambda \) dwell periods, a constant
  less than the age \( U^{*} \) needed to become a colony cell of \( M^{*} \).  
\end{proof}

\subsection{Lifting the simulation level}
   
\begin{definition}[Probabilistic color control]\label{df:qualif-control}
Let us define the constants
\begin{align}\label{eq:C_1,C_2}
   C_{1}=(26\p_{0}+2)\lambda, \quad C_{2}=26\p_{0}+2.
\end{align}
  Let ``blue'' be some color.
  For some parameter $\sigma>0$, let \( \cE \) be an event derived
  from the behavior of trajectory \( \eta \) of medium \( M \) before time \( t \).
  We will say that trajectory $\eta$ is \df{germ-level blue}
  \df{modulo} \( (\sigma,\cE) \) at time \( t \) on the finite or infinite interval \( H \) if
for every set of space intervals $I_{1},\dots,I_{n}\subset H$,
of size \( C_{1}\B \), at distance \( \ge 2 C_{1}\B \)
from each other, the probability that \( \neg\cE \) holds and
none of the \( I_{j} \) is germ-level and blue is at most $\sigma^{n}$.
When \( H \) is the whole space then we will have \( \cE=\emptyset \), so
the ``modulo \( \cE \)'' qualification can be omitted.

For an amplifier, the parameter \( \theta_{k} \) was
defined in~\eqref{eq:no-birth-prob}.
Let 
\begin{equation}\label{eq:sg}
    \sigma_{k} =
                    \begin{cases}
                      \theta_{1} &\text{ if } k=1,
                      \\ 2\theta_{k-1}&\text{ if }k>1.
                    \end{cases}
 \end{equation}
As usual, for level \( k \) we will simply write \( \sigma=\sigma_{k} \), 
\( \sigma^{*}=\sigma_{k+1} \), and so on.
\end{definition}

A special case of the lemma below says that if \( H \) is the whole space then
for every trajectory $\eta$ of $M$, if \( H \) is
germ-level blue in \( \eta \) at time $t_{1}$ modulo \( \sigma \) then it is
germ-level blue for the trajectory \( \eta^{*} \) of \( M^{*} \) modulo \( \sigma^{*} \)
at time $t_{1}+C_{2}\Tus$.
The statement is a little more complex if \( H \) is not the whole space,
because of intervals \( I_{j} \) that may fall close to its boundary.

\begin{lemma}[Self-organization]\label{lem:sorg}
  For a trajectory $\eta$ of one of the media $M=M_{k}$, of our amplifier, and
  \( \cE \) an event depending on the behavior of \( \eta \) before \( t_{1} \).
For a finite or infinite space interval \( H \), assume that 
  \( H':=\Gamma(H,C_{1}\B^{*}) \) is germ-level blue modulo \( (\sigma,\cE) \)
  in \( \eta \) at time $t_{1}$, and let $t_{2}=t_{1}+ C_{2}\Tus$.
Let \( \cE^{*} \) be the event that \( \cE \) holds, or not all (0,1 or 2)
intervals of size \( C_{1}\B^{*} \) at the ends of \( H' \) are blue in \( \eta \)
at time \( t_{1} \).
Then \( H \) is germ-level blue modulo \( (\sigma^{*},\cE^{*}) \) at time \( t_{2} \)
for the trajectory \( \eta^{*}  \) in \( M^{*} \).
  \end{lemma}

\begin{Proof}
Let \( I_{1},\dots,I_{n}\subset H \) be any system of intervals of
size \( C_{1}\B^{*} \), at distance \( \ge 2 C_{1}\B^{*} \) from each other.
  We need to bound the probability by \( (\sigma^{*})^{n} \) that \( \neg\cE^{*} \) holds
  and none of the intervals \( I_{j} \) are good in \( \eta^{*} \).
  It is sufficient to look at one of these intervals---call it \( I \), with
  \begin{align*}
   I'=\Gamma(I,C_{1}\B^{*}).
\end{align*}
  Indeed, as the space-time rectangles of the argument
  for all of them are disjoint, and the probability estimates come from a canonical
  simulation, the bounds \( \sigma^{*} \) on the probability will multiply.
  Here is an outline of the argument.
  
  We will work within the space-time rectangle
  \begin{align}\label{eq:sorg-rectangle}
   I' \times \rint{t_{1}}{t_{1}+ C_{2}\Tus}.
\end{align}
\begin{itemize}
\item 
  Estimate the probability of the event \( \cF_{1} \)
  that damage occurs in the rectangle~\eqref{eq:sorg-rectangle}.
\item Estimate the probability of the event \( \cF_{2} \) that at time \( t_{1} \) it is not true
  that the interval \( I' \) is germ-level blue with the possible exclusion of a 
  subinterval \( J \) of size \( C_{1}\B \).
\item Show that assuming \( \neg(\cF_{1}\cup\cF_{2}) \),
  except for an event \( \cF_{3} \) of probability \( \le\theta \)
  the interval \( I \) turns germ-level blue for \( M^{*} \) by time \( t_{2} \).
\end{itemize}
  
Condition~\ref{cond:restor} (Restoration Property) gives the upper bound 
\begin{equation}\label{eq:sorg.q_{1}}
  \Prob(\cF_{1})
  \le \frac{3 C_{1}\B^{*} \cdot  C_{2}\Tus)}{\B\Tu}\eps\le 4 C_{1}C_{2}\Q\U'\eps.
\end{equation}
Assuming \( \neg\cE^{*} \) implies \( \neg\cE \) as well.
Let \( \cF_{2} \) be the event that \( \neg\cE \) holds and there is no interval
\( J \) of size \( C_{1}\B \) such that, setting \( J'=\Gamma(J,2 C_{1}\B) \),
the set \( I\setminus J' \) is blue at time \( t_{1} \).
\begin{step+}{step:sorg.first}
  We have \( \Prob(\cF_{2})\le (C_{1}\Q\sigma)^{2} \).
  \end{step+}
  \begin{pproof}
 Let \(  J_{i} = i C_{1}\B+\lint{0}{C_{1}\B} \).
    Under the assumption of \( \neg\cE \), which we will omit to write
    in the estimates below,
    for some \( i,h \) with \( |i-h|\ge 4 \), let \( \cF_{i,h} \) be the event that neither
\( J_{i} \) nor \( J_{h} \) is germ-level and blue, and both intersect \( I' \)..
The germ-level blueness
of \( \eta \) at time \( t_{1} \) modulo \( \sigma \)
implies \( \Prob(\cF_{i,h})\le\sigma^{2} \), and hence
\( \Pbof{\exists i,h: \cF_{i,h}}\le  (3 C_{1}\Q)^{2}\sigma^{2} \).
If there is no such pair \( i,h \) then it is easy to see that there is an \( i \)
such that \( I\setminus J'_{i} \) is blue.
\end{pproof} 

\begin{step+}{step:sorg.lift}
  Let \( \cF_{3} \) be the event that \( \neg(\cE\cup\cF_{1}\cup\cF_{2}) \) holds
  and \( I \) does not become germ-level blue in \( \eta^{*} \).
  Then \( \Prob(\cF_{3})\le\theta\), where \( \theta \) was
  defined in~\eqref{eq:no-birth-prob}.
\end{step+}
\begin{pproof}
  By \( \neg\cF_{2} \) there is a \( J \) of size \( C_{1}\B \) such that 
  \( I'\setminus J \) is germ-level blue at time \( t_{1} \), so
  we invoke Lemma~\ref{lem:color-lift}.
\end{pproof} 

Summarizing, the probability that under condition \( \neg\cE \) the interval \( I \) does
not become blue in \( \eta^{*} \) is bounded by
\begin{align}\label{eq:sorg-sum}
  \Prob(\cF_{1}\cup\cF_{2}\cup\cF_{3})
  \le 4C_{1}C_{2}\Q\U'\eps + (C_{1}\Q\sigma)^{2}+\theta.
\end{align}
Let us show that the latter is \( \le\sigma^{*} \),
using the definitions in~\eqref{eq:sg}, those
in~\eqref{eq:frame-values}, and the inequality~\eqref{eq:exp-error}.
We have, for large enough \( \rep \) and small enough \( \eps \):
\begin{align*}
   4 C_{1}C_{2}\Q_{k}\U'_{k}\eps_{k}&=4C_{1}C_{2}\rep^{2k+4}k^{2}\eps_{k}
 <\theta_{k}/2,
  \\  (C_{1}\Q_{k}\sigma_{k})^{2}&= C_{1}^{2}\rep^{2k+2}e^{-2\rep k^{2}}
                                  < \theta_{k}/2,
\end{align*}
so~\eqref{eq:sorg-sum} is bounded by \( 2\theta=\sigma^{*} \).
\end{Proof}

\subsection{Computing supported by self-organization}

Here we will prove the remaining main theorems, which rely on self-organization.

\begin{Proof}[Proof of Theorem \protect\ref{thm:1dim.nonerg.var}]
  The field that will be remembered will be \( \Color \).
  For the size \( N \) of our space, let us define the finite or infinite level \( K \):
 \[
  K=\bigvee\setof{ k : 5\B_{k} \le N}.
 \]
  For a certain value of the field $\Color$ that we call blue,
we create an initial configuration $\eta^{1}(\cdot,0)$ with
blue latent cells covering the whole space.
Let the trajectories $\eta^{1}$, $\eta^{k}$ be defined as before.
The proof will show that by time \( C_{2}\Tu_{k} \) the space will be
germ-level blue modulo \( \sigma_{k} \); from this via trickle-down
we will be able to control the probability of any particular space-time point
$\pair{x_{1}}{t_{1}}$ not being blue.

Let \( \cF_{1} \) be the event \( \var{Blue}\ne\eta^{1}(x_{1},t_{1}).\Color^{1} \).   
Eventually, we want to upper-bound \( \Prob(\cF_{1}) \).
Let us consider first the case of infinite space, that is \( K=\infty \).
  Let 
 \begin{align*}
   v_{1}   &= t_{1},       & v_{k+1} &= v_{k}-\Tu_{k+1},
   \\  I_{1} &= \{x_{1}\},&  I_{k+1} &= \Gamma(I_{k}, \lambda\B_{k+1}),
   \\   u_{1}  &= 0, & u_{k+1} &= u_{k}+C_{2}\Tu_{k+1}.
 \end{align*}
 By Lemma~\ref{lem:sorg}, for each \( k\le K \),
 in the trajectory \( \eta^{k} \), the whole space
 is germ-level blue at time \( u_{k} \) modulo \( \sigma_{k} \).
Let \( \cF_{k} \) be the event that \( I_{k} \) is blue at time \( v_{k} \) in \( \eta^{k} \).

  \begin{step+}{step:sorg.stor.F-k}
    We have \(   \Prob(\cF_{k+1}\setminus\cF_{k})\lem \U'_{k}\Q_{k}\eps_{k} \),
    where we used the notation introduced in~\eqref{eq:lem}.
 \end{step+}
  \begin{pproof}
    If damage does not occur in \( \eta^{k} \) in
    \( I_{k+1}\times\rint{v_{k}-(\p_{0}+2)\Tu_{k+1}}{v_{k}} \)
    and \( \cF_{k+1} \) holds then Lemma~\ref{lem:color-trickle} (Color trickle-down)
    implies \( \cF_{k} \).
    The right-hand side above upper-bounds the probability of damage there.
    Now
    \begin{align*}
\U'_{k}\Q_{k}\eps_{k} < \eps^{2^{n-2} + 2^{(n-3)/2}}
       <\eps^{2^{n-2}},   
\end{align*}
\end{pproof}

  Let
 \begin{align*}
      n = \bigwedge\setof{k: v_{k+1} \le u_{k+1} \text{ or } k = K},
 \quad t_{2} = u_{n}.
 \end{align*}
If \( K=\infty \) then we always have \( v_{n+1}\le u_{n+1} \), hence
\begin{align*}
 v_{n+1} &=v_{n}-\Tu_{n+1} \le u_{n+1},\quad u_{n}<v_{n}, 
\\  v_{n} &\le u_{n+1}+\Tu_{n+1},\quad v_{n}-u_{n}\le(C_{2}+1)\Tu_{n+1}.
\end{align*}
   Let \(  I' = \Gamma(I_{n}, (C_{2}+1)\B_{n} \Tu _{n+1}/\Tl_{n}) \).
  Let \( \cB \) be the event that \( I' \) is blue in \( \eta^{n} \) at time \( u_{n} \).
  Since \( \eta^{n} \) is blue at time \( u_{n} \) modulo \( \sigma_{n} \) for each
\( n \), we have 
\begin{align*}
 \Prob(\neg\cB)\lem \U'_{n}\sigma_{n}.
\end{align*}
We claim \(  \Prob(\cB\setminus\cF_{n})\lem (\U'_{n})^{2}\eps_{n} \).
 Indeed, if damage does not occur in \( \eta^{k} \)  in \( I'\times\rint{u_{n}}{v_{n}} \)
 then \( \cB \) and Lemma~\ref{lem:lasting-color-control} imply \( \cF_{n} \).
 And the right-hand side upper-bounds the probability of damage occurring there.  
Thus, 
\begin{align*}
  \Prob(\neg\cF_{1})&\le \sum_{k=1}^{n-1}\Prob(\cF_{k+1}\setminus\cF_{k})
  + \Prob(\neg\cB)+\Prob(\cB\setminus\cF_{n}) 
\\ &\lem \sum_{k=1}^{n-1} \U'_{k}\Q_{k}\eps_{k}
     +\U'_{n}\sigma_{n} + (\U'_{n})^{2}\eps_{n}.
\end{align*}
With the parameters of~\eqref{eq:amp-fr}, large enough \( \rep \)
and small enough \( \eps \) we have
\begin{align*}
     &\sum_{k=1}^{n-1}\Prob(\cF_{k+1}\setminus\cF_{i}) = O(\eps),
\\  &\U'_{n}\sigma_{n}=\rep^{n+3}n^{2}e^{-\rep n^{2}}<e^{-\rep n^{2}/2},
  \\ &(\U'_{n})^{2}\eps_{n} = \rep^{2n+6}n^{4}\eps^{2^{n-2} + 2^{(n-3)/2}}
       <\eps^{2^{n-2}},
\end{align*}
for any somewhat large \( n \) and small \( \eps \),
giving \( \Prob(\neg\cF_{1})\lem e^{-\rep n^{2}/2}+\eps \).
Let us upper-bound \( e^{-\rep n^{2}/2} \) by lower-bounding \( n \).
We have
\begin{align*}
   t_{1}\le (C_{2}+1)\sum_{k=1}^{n+1}\Tu_{k}\lem \Tu_{n+1}.
\end{align*}
From~\eqref{eq:B_k-expl} for large \( k \):
\begin{align*}
  t_{1} &\le \Tu_{n+1}= \rep^{(n+2)(n+3)/2-2}(n!)^{2},
  \\  \ln t_{1}&< 0.6 n^{2}\ln\rep,
  \\  n^{2} &>(1/0.6\ln\rep) \ln t_{1},
  \\ e^{-\rep n^{2}}&\le t_{1}^{-\rep/0.6\ln\rep},
\end{align*}
which leads to the estimate in Theorem~\ref{thm:1dim.nonerg}.

 \begin{step+}{step:sorg.stor.fin}
  Suppose \( u_{n+1}<v_{n+1} \), hence \( n= K \) and the space is finite.
 \end{step+}
 \begin{prooof}
Given \( \B_{n}\le N\le 5 \B_{n} \) and the explicit expression for \( \B_{n} \)
and  (within constant factor of) \( \Tu_{n} \)
in~\eqref{eq:B_k-expl}, the estimate for \( \eps_{n} \)
in~\eqref{eq:exp-error} and the definition~\eqref{eq:h_0}, we have for large \( n,N \),
where for a later application we have been a little more generous:
\begin{align}
\nonumber  N &<\B_{k+1}= \rep^{(n+1)(n+2)/2-1},
\\\nonumber\log N &< (n+1)^{2}\log\rep/2,
\\\label{eq:k-explicit}  n-2 &> (\log N)^{1/2}(2/\log\rep)^{1/2}-4,
\\\nonumber 2^{n-2} &> c^{(\log N)^{1/2}}=h_{0}(N,c) \text{ for an appropriate } c>1,
\\\nonumber \eps_{n} &< \eps^{2^{n-2}} \le \eps^{h_{0}(N,c)}.
\end{align}
  Let \( \cB \) be the event that the whole space \( \bbZ_{N} \) is blue in
  \( \eta^{n} \) at time \( u_{n} \).
Since \( \eta^{n} \) is blue at time \( u_{n} \) modulo \( \sigma_{n} \),
and \( 5\B_{n+1}>N \), we have
\begin{align*}
 \Prob(\neg\cB)\lem \Q_{n}\sigma_{n}.
\end{align*}
The event \( \cB\setminus \cF_{n} \) occurs only if there is any
damage in \( \eta^{n} \) between times \( u_{n} \) and \( v_{n} \).
We can upper-bound the probability of this within a constant factor by
\begin{align*}
  \eps_{n}N t_{1}/\B_{n}\Tu_{n} \le \eps_{n} t_{1}/\Tu_{n}\le \eps_{n}t_{1}
  \le\eps^{h_{0}(N,c)}t_{1},
\end{align*}
which is what Theorem~\ref{thm:1dim.nonerg} states.
\end{prooof} 
\end{Proof}

\begin{Proof}[Proof of Theorem~\ref{thm:1dim.comp.var}]
  This also proves Theorem~\ref{thm:1dim.comp.finite}, and the remark
  at the last part of the proof also proves Theorem~\ref{thm:1dim.stor.var} for the
  finite space.
  Let \( \trans \) be a given standard computing transition function over an alphabet \( \Sigma \).
  Construct an amplifier as in Lemma~\ref{lem:amp} with
  the  transition function \( \PlTrans_{1}(\cdot)=\trans(\cdot) \).
  This amplifier defines cell body sizes \( \B_{k} \) for all \( k \).
  Given a string \( \rho\in\Sigma_{0}^{n} \) in the input-output alphabet, let
  \begin{align*}
    l=\bigvee\setOf{k}{\B_{k}\le n}.
  \end{align*}
  Create a germ of cells of \( M_{l} \) covering the space interval \( \lint{0}{n} \),
  with address \( \Q_{l}/2 \) somewhere inside it.
  Fill its \( \Payload \) track first with \( * \), then write over the \( \Payload.\Input \)
  track the string \( \rho \).
  Now apply consecutively the encodings \( \phi_{*l} \), \( \phi_{*(l-1)} \), \( \dots \),
  to get a string \( \rho_{l} \) of \( \bbS_{1} \).
  The code \( \psi_{*} \) whose existence is stated in the theorem is defined as
  \( \psi_{*}(\rho)=\rho_{l} \),
  which gives \( n_{l}:=|\rho_{l}|<2n \) as required by the theorem.
  From here, with  \( \Init \) as in Definition~\ref{def:Init},
 we obtain the initial configuration of \( M_{1} \)
 as \( \eta(\cdot,0)=\Init_{\psi_{*}}(\rho) \).
 For readability, we will call the color \( -1 \) ``blue'', color 0 ``white'' and
 color 1 ``red''.
 First, for simplicity, assume that the space is infinite.

By the construction, the interval \( \lint{0}{n_{l}} \) is already an \( l \)-level white germ.
First we will show how the space outside it will also promptly self-organize to level \( l \):
of color blue on the left and red on the right.
From then on, the self-organization continues on the whole space, and
at the same time, the computation in the white germ proceeds reliably,
lifting to a higher level in every step of self-organization.
Finally, the trickle-down argument shows that at any time
the result can be read out at the lowest level.
 
 For \( k=1,2,\ldots \), let
 \begin{align*}
   a_{k}&=\sum_{i=2}^{k}C_{1}\B_{i},\quad
          v_{k} = \sum_{i=2}^{k}C_{2}\Tu_{i},
\\ J_{k} &= \lint{-a_{k}}{n_{l}+a_{k}},          
\quad  H^{-}_{k}=\opint{-\infty}{-a_{k}},\quad H^{+}_{k}=\lint{n_{l}+a_{k}}{\infty}, 
 \end{align*}
 then \( J_{l}=\lint{0}{n_{l}} \),  \( H^{-}_{1}=\opint{-\infty}{0} \).
For \( j=1,2,3 \) let \( \cA_{k,j} \) be the event that  in \( \eta^{k} \) at time \( v_{k} \),
the interval \( \lint{- j a_{k+1}}{-a_{k}} \) is not germ-level blue or
\( \lint{n_{l}- a_{k}}{n_{l}+ j a_{k+1}} \) is not germ-level red.
 Let
\begin{align*}
 \cD_{k}=\bigcup_{i= 1}^{k}\cA_{i,1}.
\end{align*}
 Note that \( \cA_{1,1}=\cD_{1}=\emptyset \).
\begin{step+}{step:1dim.comp.outside}
 For each \( k \), in \( \eta^{k} \)  at time \( v_{k} \),
  the set \( H^{-}_{k}\setminus H^{-}_{k+1} \)
  is germ-level blue and \( H^{+}_{l}\setminus H^{+}_{l+1} \) is germ-level red
 modulo \( (\sigma_{k},\cD_{k-1}) \).
\end{step+}
   \begin{sloppypar}
 \begin{pproof}
   Interval \( H^{-}_{1} \) is germ-level blue at time \( 0 \) in \( \eta^{1} \) by definition.
Lemma~\ref{lem:sorg} implies that \( H^{-}_{2} \) is germ-level blue in
 \( \eta^{2} \)  at time \( v_{2} \) modulo \( (\sigma_{2},\cD_{1}) \), further
\( H^{-}_{3} \) is germ-level blue in
\( \eta^{3} \)  at time \( v_{3} \) modulo \( (\sigma_{3},\cD_{2}) \), and so on.
The same argument works for \( H^{+}_{k} \).     
\end{pproof} 
   \end{sloppypar}

For \( k\ge l \), let \( \cF_{k} \)
be the event that there is damage in \( \eta^{k} \) in the space-time rectangle 
\begin{align*}
 \Gamma(J_{l},3 a_{k+1})\times \rint{0}{v_{k+1}}.
\end{align*}
\begin{step+}{step:1dim.comp.middle}
  Assuming \( \neg \cF_{l} \), the interval \( J_{l} \) is part of a white germ
 at time \( v_{l} \) in \( \eta^{l} \).
\end{step+}
\begin{pproof}
 We started from the interval \( J_{l} \) of level \( l \) cells surrounded by level 1 latent cells.
In the absence of  damage in \( \eta^{l} \), these cells perform according to the program,
hence up to time \( v_{l}<2 \Tu_{l+1} \), none of them would die.
\end{pproof} 

\begin{step+}{step:1dim.comp.boundary}
From \( \neg(\cA_{l,3}\cup\cF_{l}) \) follows that in \( \eta^{l} \)
at time \( v_{l}+2\Tu_{l+1} \),  interval \( J_{l} \) is part of a white germ,
and \( \Gamma(J_{l},2 a_{l+1}) \) is germ-level color-controlled.
\end{step+}
\begin{pproof}
At time \( v_{l} \)  in \( \eta^{l} \), as part~\ref{step:1dim.comp.middle} shows,
  \( \neg\cF_{l} \) implies that  \( J_{l} \) is part of a white germ.
And \( \neg\cA_{l,3} \), implies by definition that
\( \lint{-3 a_{l+1}}{-a_{l}} \) is germ-level blue and
\( \lint{n_{l}-a_{l}}{n_{l}+3 a_{l+1}} \) is germ-level red.
We can now apply the argument of Lemma~\ref{lem:robust-color-control}
simultaneously with 
with \( I=\lint{-3 a_{l+1}}{n_{l}} \),  \( J=\lint{-a_{l}}{0} \),
and with \( I=\lint{0}{n_{l}+3 a_{l+1}} \), \( J=\lint{n_{l}}{n_{l}+a_{l}} \),
to imply that in \( \eta^{l} \) at time \( v_{l}+2\Tu_{l+1} \),
the whole interval \( \Gamma(J_{l},2 a_{l+1}) \) becomes germ-level color-controlled.
\end{pproof} 

\begin{step+}{step:1dim.comp.next}
From \( \neg(\cA_{l,3}\cup\cF_{l}) \) follows that in \( \eta^{l+1} \)
at time \( v_{l+1} \), the interval \( J_{l} \) is part of a white germ,
and interval \( \Gamma(J_{l},a_{l+1}) \) is germ-level color-controlled,
 except for an event \( \cE_{l} \) of probability
\begin{align*}
   \le (2 C_{1}+7)\theta_{l},
\end{align*}
where \( \theta_{l} \)  was defined in~\eqref{eq:no-birth-prob}.
\end{step+}
\begin{pproof}
  From \( \cF_{l} \) follows that \( J_{l} \) is part of a white germ also in
  \( \eta^{l+1} \).
By part~\ref{step:1dim.comp.boundary},
from \( \neg(\cA_{l,3}\cup\cF_{l}) \) follows that in \( \eta^{l} \) at time
\( v_{l}+2\Tu_{l+1} \),
the interval \( \Gamma(J_{l},2 a_{l+1}) \) is germ-level color-controlled.
Note that \( v_{l+1}=v_{l}+2\Tu_{l+1}+26\p_{0}\Tu_{l+1} \).
Applying Lemma~\ref{lem:color-lift}, for case \( \eta=\eta^{l} \),
\( t_{2}-t_{1}=26\p_{0}\Tu_{l+1} \),
to any subinterval \( K \) of length \( 2\B_{l+1} \) of \( \Gamma(J_{l},a_{l+1}) \),
we find that \( K \) becomes color-controlled in \( \eta^{l+1} \) except for
an event of probability \( \theta_{l} \).
If each of the \( r \)  intervals of the form
\( \lint{i\B_{l+1}}{(i+2)\B_{l+1}} \) 
intersecting \( \Gamma(J_{l},a_{l+1}) \) becomes color-controlled then so does
\( \lint{-a_{l+1}}{n_{l}-a_{l+1}} \).
Let \( \cE_{l} \) be the event that any one of them does not become color-controlled,
then \( \Prob(\cE_{l})\le r\theta_{l} \).
The construction gave \( |J_{l}|<5 \B_{l+1} \),
we know \( a_{l+1}<2C_{1}\B_{l+1} \), hence
\( |\Gamma(J_{l},a_{l+1})|< (4 C_{1}+5)\B_{l+1} \), \( r\le 4 C_{1}+9 \).
\end{pproof} 

Let us handle now first the case of infinite space \( \bbZ \).

\begin{step+}{step:1dim.comp.above}
  Let \( k>l \).
  Assume that in \( \eta^{k} \) at time \( v_{k} \)
the interval \( J_{l} \) is part of a white germ, and
the interval \( \Gamma(J_{l},a_{k}) \) is germ-level color-controlled.
Then from \( \neg(\cA_{k,2}\cup\cF_{k}) \) follows that in \( \eta^{k+1} \)
at time \( v_{k+1} \), the interval \( J_{l} \) is part of a white germ,
and interval \( \Gamma(J_{l},a_{k+1}) \) is germ-level color-controlled,
 except for an event \( \cE_{k} \) of probability \( \le (2 C_{1}+3)\theta_{k} \).
\end{step+}
\begin{pproof}
  (This proof is similar to the one above.)
    From \( \cF_{k} \) follows that \( J_{l} \) is part of a white germ also in
  \( \eta^{k+1} \).
  From \( \neg\cA_{k,2} \) follows that in \( \eta^{l} \) at time \( v_{l} \)
also the interval \( \Gamma(J_{l},2 a_{k+1}) \) is germ-level color-controlled.
 Applying Lemma~\ref{lem:color-lift}, for case \( \eta=\eta^{k} \),
\( t_{2}-t_{1}=C_{2}\Tu_{k+1} \),
to any subinterval \( K \) of length \( 2\B_{k+1} \) of \( \Gamma(J_{l},a_{k+1}) \),
we find that \( K \) becomes color-controlled in \( \eta^{k+1} \) except for
an event of probability \( \theta_{k} \).
If each of the \( r \)  intervals of the form \( \lint{i\B_{k+1}}{(i+2)\B_{k+1}} \) 
intersecting \( \Gamma(J_{l},a_{l+1}) \) becomes color-controlled then so does
\( \Gamma(J_{l},a_{l+1}) \).
Let \( \cE_{k} \) be the event that any one of them does not become color-controlled;
then \( \Prob(\cE_{k})\le r\theta_{k} \).
The construction gave \( |J_{l}|<\B_{k+1} \),
we know \( a_{k+1}<2C_{1}\B_{k+1} \), hence
\( |\Gamma(J_{l},a_{l+1})|< (4 C_{1}+1)\B_{l+1} \), \( r\le 4 C_{1}+5 \).
\end{pproof} 

For \( k\ge l \) let
\begin{align*}
   \cD'_{k} = \bigcup_{i=l}^{k}\cA_{l,2}, \quad
  \cG_{k}=\cD_{l-1}\cup\cA_{l,3}\cup\cD'_{k-1}\cup\bigcup_{i=l}^{k-1}\cF_{i}
           \cup\bigcup_{i=l}^{k-1}\cE_{i},
\end{align*}
where the events \( \cE_{i} \) were defined in the proofs of part~\ref{step:1dim.comp.next}
and part~\ref{step:1dim.comp.above}.

\begin{step+}{step:1dim.comp.combine}
  Let \( k> l \).
  Assuming \( \neg\cG_{k} \), in \( \eta^{l+1} \) at time \( v_{l+1} \)
the interval \( \Gamma(J_{l},a_{k}) \) is germ-level color-controlled,
and the interval \( J_{l} \) is part of a white germ.
Also,
\begin{equation}\label{eq:G-bound}
\begin{aligned}
  \Prob(\cG_{k})&\le 2\sum_{i=1}^{l-1}\Q_{i}\sigma_{i}+6\Q_{l}\sigma_{l}
                  + 4\sum_{i=l}^{k-1}\Q_{i}\sigma_{i}
  \\ &+ 33\sum_{i=l}^{k-1}\Q_{i}\U'_{i}\eps_{i}
       +  (2 C_{1}+7)\theta_{l} + \sum_{i=l+1}^{k-1}(4 C_{1}+5)\theta_{i}=O(\sigma_{1}).
\end{aligned}
\end{equation}
\end{step+}
\begin{pproof}
  The statement before the probability estimate,
  for \( k=l+1 \), follows from part~\ref{step:1dim.comp.next}.
  For \( k>l+1 \), it follows by induction from part~\ref{step:1dim.comp.above}.
  It remains to prove \eqref{eq:G-bound}.
From the definition, \( \cD_{k}=\bigcup_{i=2}^{k}(\cA_{i,1}\setminus\cD_{i-1}) \).  
From
\begin{align*}
  \Prob(\cA_{i+1,j}\setminus\cA_{i,j})&\le 2 j \Q_{k}\sigma_{k}
\end{align*}
we have
\(  \Prob(\cD_{k})\le 2\sum_{i=1}^{k-1}\Q_{i}\sigma_{i} \),
\(  \Prob(\cD'_{k})\le 4\sum_{i=l}^{k-1}\Q_{i}\sigma_{i} \).                     
Given \( n_{l}<5\B_{l+1} \),
\( a_{k+1}< 2\B_{k+1} \),  \( v_{k+1}\le 2\Tu_{k+1} \), we have
\begin{align*}
  \Prob(\cF_{k})\le 11 (\B_{k+1}/\B_{k})\cdot 2(\Tu_{k+1}/\Tu_{k})
  \le 33 \Q_{k}\U'_{k}\eps_{k}.
\end{align*}
The estimates for \( \cE_{k} \) were given in parts~\ref{step:1dim.comp.next} 
and~\ref{step:1dim.comp.above}.
\end{pproof} 

In the reasoning below, we assume \( \neg\cG_{k} \).
Between times \( v_{k} \) and \( v_{k+1} \), we have a white germ
containing \( J_{l} \), around the germ cell with address \( \Q_{k}/2 \).
Besides lifting from level \( k \) to level \( k+1 \),
the trajectory \( \eta^{k} \) performs its payload computation;
for each level \( i\le k \), this proceeds without faults
up to time \( v_{i+1} \), in the area \( \Gamma(J_{l},a_{i+1}) \)
which is all the area that matters, as the germ does not grow beyond it.
But we are interested in the output of level \( 1 \), so trickle-down will be invoked.
According to the rules \( \Compute \) of Algorithm~\ref{alg:Compute} and
\( \UpdatePayload \) of Algorithm~\ref{alg:Update-payload},
payload computation happens on level \( k \)
representing germ cells of level \( k+1 \) with addresses in \( \lint{0}{\Q_{k+1}} \).
We will estimate the number of payload computation steps performed, but first assume that
time \( t' \) is the end of a dwell period of a germ cell \( x' \) of \( \eta^{k} \),
in which \( x=x'+a \) for some \( a\in\lint{0}{\B_{k}} \),
and \( \eta^{k}(x',t').\Payload(a) \) is supposed to have the value \( \zeta(x,t) \).
Under the assumption \( \neg\cG_{k} \) indeed
\( \zeta(x,t)=\eta^{k}(x', t').\Payload(a) \).
A trickle-down reasoning just like in Section~\ref{sec:trickle-down}, in the proof of
Theorem~\ref{thm:remember-inf}, with \( D_{k}=\Tu_{k} \), as done for
information storage in the discussion after Lemma~\ref{lem:amp} (Amplifier),
will bring the equality down to level 1:
\begin{align*}
  \Pbof{\zeta(x,t).\Output\ne \eta(x, t').\Output}&\le\Prob(\cG_{k})+
                                                    \eps_{k}+\sum_{i=1}^{k-1}\eps''_{i}
\\ &= O(\sigma_{1})=O(e^{-\rep}).
\end{align*}

Let us now estimate \( t' \) in terms of \( t \).
Let \( t'_{k} \) be the time by which the first white cell of \( M_{k} \)
will be created, then we have \( v_{k-1}< t'_{k} \le v_{k} \).
Between times \( t'_{k} \) and \( t'_{k+1} \), 
there are enough computation steps to simulate in \( \trans \)
the \( \B_{k+1}-\B_{k} \) steps of the expansion of an area of
size \( \B_{k} \) to one of size \( \B_{k+1} \).
This shows that if \( \B_{k}< t \le \B_{k+1} \) then
\( t_{k}\le t'<t_{k+1} \).
Recall the definition of bandwidth \( \bw_{k} \) in~\eqref{eq:amp-fr}.
During the time interval \( \rint{t'_{k}}{t'_{k+1}} \), starting at time \( t'_{k} \)
until the payload of a germ of size \( \B_{k+1} \) is filled,
each work period of a germ cell simulates at least a constant times
\( \bw_{k}\B_{k} \) steps of \( \trans \).
So if \( t>\B_{k} \) then to simulate \( t \) steps of \( \trans \) we need,
within constant factor,
\( t/\B_{k}\bw_{k} \) work periods, each of length \( \le\Tu_{k} \), so 
\begin{align*}
   t \Tu_{k}/\B_{k}\bw_{k}.
\end{align*}
For \( t\le\B_{k} \) we need less, so we are conservative when going with this
estimate for all \( t \).
Noting that \( \prod_{i<k}\U'_{i} \) is within constant factor of \( \Tu_k \),
and using~\eqref{eq:B_k-expl}, we have within constant factor
\begin{align*}
   \Tu_{k}/\B_{k}\bw_{k} = \rep^{k}(k!)^{2}< 2^{2.1  k\log k}.
\end{align*}
Using~\eqref{eq:B_k-expl} we have for large \( t \),
\begin{align}
  \nonumber
  t&>B_{k}=\rep^{k(k+1)/2-1},
\\\nonumber \log t &> \log\rep\cdot k^{2}/2,
  \\\label{eq:k-bound}  k &< (\log t)^{1/2}(2/\log\rep)^{1/2},
\\\nonumber    2^{2.1 k\log k} &< c^{ (\log t)^{1/2}\log\log t},
\end{align}
for a constant \( c \) depending on \( \rep \).
So we can estimate \( t'\le t\cdot c^{(\log t)^{1/2}\log\log t} \)
as Theorem~\ref{thm:1dim.comp} states.

In case of the finite space \( \bbZ_{N} \) we can only hope to simulate the
computation \( \zeta(\cdot,t) \) while it fits into a space of size \( N \),
so let
\begin{align*}
   K = \bigvee\setOf{k}{B_{k}< N}.
\end{align*}
\begin{step+}{step:1dim.comp.finite}
  In case of a space \( \bbZ_{N} \) for finite \( N \), adding \( \eps^{h_{0}(N,c)}t \)
  to the estimate with a constant \( c>1 \) will do.
  The same argument works also to finish the proof of Theorem~\ref{thm:1dim.stor.var}
  for a finite space.
\end{step+}
\begin{pproof}
  After the level has been raised to \( K \), it will not be raised any longer.
  In every time period of size  \( \Tu_{K} \), the computation may fail
  if a fault of level \( K \) occurs, that is with probability \( \eps_{K} \).
  In time \( t \) the probability of this happening at least once can be (generously)
  be upper-bounded by \( \eps_{K}t \).
  We estimate, similarly to~\eqref{eq:k-explicit}, for \( k=K \):
  \begin{align*}
 \eps_{K} &< \eps^{h_{0}(N,c)}
\end{align*}
for an appropriate \( c>1 \).
\end{pproof} 

\end{Proof} 

\section{Some applications and open problems}\label{sec:concl}

\subsection{Non-periodic Gibbs states}

  Consider spin systems\index{spin system} in the usual sense
(generalizations of the Ising model\index{Ising model}).
  All 2-dimensional spin systems hitherto known were known to have only a
finite number of extremal Gibbs states\index{state!Gibbs} (see
for example \cite{AizenIsing80}): thus, theoretically, the amount of storable
information on an $n\times n$ square lattice did not grow with the size of
the lattice.
  In 3 dimensions this is not true anymore, since we can stack independent
2-dimensional planes: thus, in a cube $C_{n}$\glo{c.cap@$C_{n}$} of size
$n$, we can store $n$ bits of information.
  More precisely, the information content of $C_{n}$ can be measured by the
dimension of the set of vectors\glo{s.greek@$\sigma$}
 \[
    \tupof{\mu\setof{\sigma(x) = 1} : x\in C_{n}}
 \]
  where $\mu$ runs through the set of Gibbs states.
  This dimension can be at most $O(n^{d-1})$ in a $d$-dimensional lattice,
since the Gibbs state on a cube is determined by the distribution on its
boundary.
  The stacking construction shows that storing
$\Og(n^{d-2})$\glo{o.greek.cap@$\Og$} bits of
information is easy.
  We can show that $\Og(n^{d-1})$ is achievable: in particular, it is
possible to store an infinite sequence in a 2-dimensional spin system in
such a way that $n$ bits of it are recoverable from any $n$ sites with
different $x$ coordinates.

  For this, we apply a transformation from~\cite{GKLMEquil86} (see 
also~\cite{BennGriIrrev85}): a probabilistic cellular automaton $M$ in $d$
dimensions gives rise to an equilibrium system $M'$ in $(d+1)$ dimensions.
  Essentially, the logarithms of the local transition probabilities define
the function $J$ and histories of $M$ become the
space-configurations of $M'$.
  Non-ergodicity of $M$ corresponds to phase transition in $M'$.
  In the cellular automaton of Theorem~\ref{thm:1dim.stor.var}, each infinite
sequence $\rg$ gives rise to a history storing the bits of
$\rg$ in consecutive cells.
  Now, each of these histories gives rise to a separate
Gibbs state belonging to one and the same potential.

  Let us note that though the Gibbs system is defined in terms of an energy
function \( H(\sigma) \)\glo{h.cap@\( H(\sigma) \)}, it is not helpful to represent this
energy function in terms of a temperature\index{temperature}
\( T \)\glo{thm:cap@\( T \)} as \( H(\sigma)/T \).
  The reason is that the individual terms of \( H(\sigma) \) do not depend
linearly on the error probability \( \eps \) (or any function of it).
  In fact, we believe it can be proved that if an artificial \( T \) is
introduced (with \( T=1 \) for a certain sufficiently small value of \( \eps \))
then a slight decrease of \( T \) destroys the phase transition.

 \subsection{Some open problems}

 \paragraph{Turing machines}\index{Turing machine}
  It is an interesting question (asked by Manuel Blum) whether a reliable
Turing machine can be built if the tape is left undisturbed, only the
internal state is subject to faults.
By now this question has been answered affirmatively, in~\cite{Turing21}.
The hierarchical construction shares many similarities with the one
of this paper, but there is quite a number of additional details.

\paragraph{Non-hierarchical construction in continuous time}

The simple three-dimensional reliable cellular automaton in~\cite{GacsReif3dim88}
(with its simpler and stronger proofs in~\cite{BermSim88} and~\cite{GacsToom95})
relies strongly on being in discrete time.
Each plane supposed to store copies of a single symbol needs to switch instantly
in each step.
In fact, currently no (proved) construction is known in any dimension of a simple reliable
cellular automaton in continuous time whose work is non-hierarchical.
There is a simple way to simulate a discrete-time cellular automaton by a continuous-time one,
described in Section~\ref{sec:sync} and called the ``marching soldiers'' scheme,
and also used in this paper.
However, without the hierarchy, its issues caused by randomness and faults are not yet
solved.
Besides the possibly large local delays during which the self-correction stalls while the
faults continue to occur, in dimensions higher than one the faults can also create deadlocks
in the synchronization scheme.
An important step was taken by Cook and Winfree with a rule that, under plausible conditions,
eliminates these deadlocks: see the report~\cite{3Dasync09}.
Further progress depends on the ability to bound the probability of large local delays.

\paragraph{Does a 1-dimensional solution have to be ``hierarchical''?}

Over the years, there have been several attempts to define a non-ergodic one-dimensional
cellular automaton that is ``simple''.
All these attempts ended up similar to~\cite{GacsKurdLevIslands78}, which is
ergodic (at least under some small noise, as shown in~\cite{ParkThes96}).
By ``complexity'' I don't just refer to the number of states, rather to the fact that
the history requires inhomogenities (to exist or to arise) on a larger-and-larger scale.
Ideally, we will see a result confirming the intuition that this sort of
complexity is necessary, but as a first step, we need a
mathematically formulated conjecture.

\paragraph{Relaxation time as a function of space size}
Recall relaxation time in Definition~\ref{def:relaxation-time}.
  Consider now Toom's medium\index{Toom rule} as a typical example of a
medium non-ergodic for \( m=\infty \).
  It follows from Toom's proof that, for small enough fault
probability \( \eps \), we have
 \( \lim_{m\to\infty}r_{m}(0,1/3)=\infty \), that is the increase of space
increases the relaxation time\index{relaxation time} (the length of time
for which the rule keeps information) unboundedly.
  The speed of this increase is interesting since it shows the
durability of information as a function of the size of the cellular
automaton in which it is stored.
  Toom's original proof gives only \( r_{m}(0,1/3)>cm \) for some constant
\( c \).
  The proof in~\cite{BermSim88}, improving on~\cite{GacsReif3dim88},
gives \( r_{m}(0,1/3)>e^{cm} \) for some constant \( c \), and this is essentially
the meaning of saying that Toom's rule helps remember a bit of information
for exponential time.

  So far, the relaxation times of all known nontrivial non-ergodic
media (besides Toom's, the ones in~\cite{Gacs1dim86} and \cite{Gacs2dim89})
depend exponentially on the size of the space.
  It is an interesting question whether this is necessary.
  In a later work, we hope to show that this is not the case and that
there are non-ergodic media such that \( r_{m}(n,1/3)<m^{c} \) holds for all
\( n \), for some constant \( c \).
  The main idea is that since the medium will be able to perform an
arbitrary reliable computation this computation may involve recognizing the
finiteness of the space rather early (in time \( m^{c} \)) and then erasing all
information.

 \paragraph{Relaxation time as a function of observed area}
  Lemma~\ref{lem:unif} seems to suggest that the issue of information
loss in a cellular automaton is solved by the question of ergodicity
for \( m=\infty \) where \( m \) is the space size.
  This is not so, however: as noted together with Larry Gray, Leonid
Levin and Kati Marton, we must also take the dependence of
\( r_{m}(n, \dg) \) on \( n \) into account.
  As time increases we may be willing to use more and more cells to
retrieve the original information.
  Even if \( r_{\infty}(n,\dg)<\infty \) for each \( m \), we may be satisfied
with the information-keeping capability of the medium if, say,
\( r_{m}(n,1.9)>e^{cn} \) for some constant \( c \).

  In some mixing systems, the dependence of \( r_{m}(n, \dg) \) on \( n \) is
known.

\begin{example}[The contact process.]\index{contact process}
  The contact process is a one-dimensional \( \CCA \) as defined in 
Section~\ref{sec:Markov}.
  Let us note that this this process is not noisy: not all local
transition rates are positive.
  The process has states 0,1.
  In trajectory \( \eta(x,t) \), state \( \eta(x,t)=1 \) turns into 0 with
rate 1.  
  State \( \eta(x,t)=0 \) turns into 1 with rate
\( \lambda(\eta(x-1,t)+\eta(x+1,t)) \).
  It is known that this process has a critical rate\index{critical rate}
\( \lambda_{c}\in\opint{0}{\infty} \) with the following properties.
  
  If \( \lambda\le \lambda_{c} \) then the process is mixing with the invariant
measure concentrated on the configuration \( \xi \) with \( \xi(x)=0 \) for
all \( x \).
  If \( \lambda>\lambda_{c} \) then the process is non-ergodic.

  If \( \lambda<\lambda_{c} \) then it is known (see Theorem 3.4 in Chapter VI 
of~\cite{Liggett85}) that the order of magnitude of \( r_{\infty}(n,\dg) \) 
is \( \log n \).

  If \( \lambda=\lambda_{c} \) then the convergence is much slower, with an order
of magnitude that is a power of \( n \) (see Theorem 3.10 in Chapter VI
of~\cite{Liggett85} and \cite{BezuidenhoutGrimmett90}).
\end{example}

  We believe it possible to construct a medium that is mixing for
\( m=\infty \) but loses information arbitrary slowly: for any computable
function \( f(n) \), and a constant \( c \) there is a medium with
 \[
  r_{m}(n,1.9) > f(n)\wedge e^{cm}
 \]
  for finite or infinite \( m \).
  Notice that this includes functions \( f(n) \) like \( e^{e^{e^{n}}} \).
  Such a result could be viewed as an argument against the relevance of the
non-ergodicity of the infinite medium for practical information
conservation in a finite medium.
  The construction could be based on the ability to perform arbitrary
computation reliably and therefore also to destroy locally identifiable
information arbitrarily slowly.

\subsubsection*{Acknowledgement}
  I am very thankful to Robert Solovay for reading parts of the paper and
finding important errors.
  Larry Gray revealed his identity as a referee and gave generously of his
time to detailed discussions: the paper became much more readable (yes!) as
a result.

\bibliographystyle{amsplain}
 \bibliography{reli,gacs-publ}

 \renewcommand{\memjustarg}[1]{\hyperpage{#1}}
\renewcommand{\indexname}{Notation}
\renewcommand{\indexname}{Index}

\end{document}

%% file: tex-figures/running-gap-fig.tex
 \setlength{\unitlength}{1mm}
 \begin{picture}(100, 120)(-10, -5)
    \put(0,0){\circle{1}}
    \put(0,0){\line(0,1){10}}
    \put(0,10){\circle{1}}
    \put(0,10){\line(1,2){5}}
    \put(5,20){\circle{1}}
    \put(5,20){\line(0,1){10}}
    \put(5,30){\circle{1}}
    \put(5,30){\line(0,1){10}}
    \put(5,40){\circle{1}}
    \put(5,40){\line(1,2){5}}
    \put(10,50){\circle{1}}
    \put(10,50){\line(1,2){5}}
    \put(15,60){\circle{1}}
    \put(15,60){\line(0,1){10}}
    \put(15,70){\circle{1}}
    \put(15,70){\line(1,2){5}}
    \put(20,80){\circle{1}}
    \put(20,80){\line(1,2){5}}

    \put(20,0){\circle{1}}
    \put(20,0){\line(1,2){5}}
    \put(25,10){\circle{1}}
    \put(25,10){\line(0,1){10}}
    \put(25,20){\circle{1}}
    \put(25,20){\line(1,1){5}}
    \put(30,25){\circle{1}}
    \put(30,25){\line(1,1){5}}
    \put(35,30){\circle{1}}
    \put(35,30){\line(1,0){5}}
    \put(40,30){\circle{1}}
    \put(40,30){\line(1,2){5}}
    \put(45,40){\circle{1}}
    \put(45,40){\line(1,1){5}}
    \put(50,45){\circle{1}}
    \put(50,45){\line(1,0){5}}
    \put(55,45){\circle{1}}
    \put(55,45){\line(1,0){5}}
    \put(60,45){\circle{1}}
    \put(60,45){\line(1,1){5}}
    \put(65,50){\circle{1}}
    \put(65,50){\line(1,1){5}}
    \put(70,55){\circle{1}}
    \put(70,55){\line(1,0){5}}
    \put(75,55){\circle{1}}
    \put(75,55){\line(1,0){5}}
    \put(80,55){\circle{1}}
    \put(80,55){\line(1,1){5}}
    \put(85,60){\circle{1}}
    \put(85,60){\line(0,1){10}}
    \put(85,70){\circle{1}}
    \put(85,70){\line(1,1){5}}
    \put(90,75){\circle{1}}
    \put(90,75){\line(1,2){5}}
    \put(95,85){\circle{1}}
    \put(95,85){\line(1,1){5}}

    \put(20,0){\line(1,0){55}}
    \multiput(25,0)(5,0){10}{\circle*{1}}

    \put(75,0){\circle*{1}}
    \put(75,0){\line(-1,2){5}}
    \put(70,10){\circle*{1}}
    \put(70,10){\line(0,1){10}}
    \put(70,20){\circle*{1}}
    \put(70,20){\line(1,2){5}}
    \put(75,30){\circle*{1}}
    \put(75,30){\line(0,1){10}}
    \put(75,40){\circle*{1}}
    \put(75,40){\line(0,1){10}}
    \put(75,50){\circle*{1}}
    \put(75,50){\line(0,1){10}}
    \put(75,60){\circle*{1}}
    \put(75,60){\line(1,2){5}}
    \put(80,70){\circle*{1}}
    \put(80,70){\line(0,1){10}}
    \put(80,80){\circle*{1}}
    \put(80,80){\line(0,1){10}}

    \put(105, 30){\vector(0,1){10}}
    \put(108, 35){\makebox(0,0)[l]{time}}

    \put(8, 35){\makebox(0,0)[l]{$l(t)$}}
    \put(45, 35){\makebox(0,0)[l]{$r(t)$}}
    \put(78, 35){\makebox(0,0)[l]{$P_{1}(t)$}}

    \put(75, -3){\makebox(0,0)[t]{$\pair{x_{0}}{v_{0}}$}}
    \put(20, -3){\makebox(0,0)[t]{$\pair{y_{0}}{v_{0}}$}}
    \put(82, 100){\makebox(0,0)[r]{$\pair{x_{n}}{y_{n}}$}}

 \end{picture}
